\documentclass[a4paper, 10pt]{article}

\usepackage{multicol}
\usepackage{euscript}
\usepackage{amssymb}
\usepackage{amsmath}
\usepackage{amsthm}
\usepackage[dvips]{graphicx}
\usepackage{psfrag}
\newcommand{\unpo}{\mathcal{O}(1)}
\newcommand{\ee}{\varepsilon}
\newcommand{\ue}{u^{\varepsilon}}
\newcommand{\erreN}{\mathbb{R}^N}

\newcommand{\ve}{v^{\varepsilon}}
\newcommand{\un}{u^{\nu}}
\newcommand{\tr}{\tilde{r}_i}

\newcommand{\B}{\text{\rm \ss}}

\newtheorem{teo}{Theorem}[section]
\newtheorem{pro}{Proposition}[section]
\newtheorem{lem}{Lemma}[section]
\newtheorem{hyp}{Hypothesis}

\theoremstyle{definition}
\newtheorem{ex}{Example}[section]
\newtheorem{say}{Definition}[section]
\newtheorem{rem}{Remark}[section]

\begin{document}
\numberwithin{equation}{section}
\bibliographystyle{plain}

\title{The boundary Riemann solver coming from
the real vanishing viscosity approximation}
\author{\it Stefano Bianchini \footnote{SISSA-ISAS, via Beirut 2-4 34014 Trieste, Italy } and Laura V. Spinolo \footnote{Northwestern University,
2033 Sheridan Road
Evanston, IL 60208-2370, USA} \\
        \footnotesize \it Email:
         bianchin@sissa.it,  spinolo@math.northwestern.edu}
\date{}
\maketitle \vspace{0.3cm}
\begin{abstract}
We study the limit of the hyperbolic-parabolic approximation
\begin{equation*}
\left\{
\begin{array}{lll}
      \ve_t + \tilde{A} \big( \ve, \, \ee \ve_x \big) \ve_x =
      \ee \tilde{B}(\ve ) \ve_{xx} \qquad \ve \in \mathbb{R}^N\\
      \tilde \B(\ve (t, \, 0)) \equiv \bar g \\
      \ve (0, \, x) \equiv \bar{v}_0. \\
\end{array}
\right.
\end{equation*}
The function $\tilde \B$ is defined in such a way to guarantee that the initial boundary value problem is well posed even if $\tilde B$ is not invertible.
The data $\bar g$ and $\bar v_0$ are constant.

When $\tilde B$ is invertible, the previous problem takes the simpler form 
\begin{equation*}
\left\{
\begin{array}{lll}
      \ve_t + \tilde{A} \big( \ve, \, \ee \ve_x \big) \ve_x =
      \ee \tilde{B}(\ve ) \ve_{xx} \qquad \ve \in \mathbb{R}^N\\
      \ve (t, \, 0) \equiv \bar v_b \\
      \ve (0, \, x) \equiv \bar{v}_0. \\
\end{array}
\right.
\end{equation*}
Again, the data $\bar v_b$ and $\bar v_0$ are constant. The conservative case is included in the previous formulations.

It is assumed convergence of the $\ve$, smallness of the total variation and other technical hypotheses and it is provided a complete characterization of the limit. 

The most interesting points are the following two.  

First, the boundary characteristic case is considered, i.e. one eigenvalue of $\tilde A$ can be $0$. 

Second, as pointed out before we take into account the possibility that $\tilde B$ is not invertible. To deal with this case, we take as hypotheses conditions that were introduced by Kawashima and Shizuta relying on physically meaningful examples. We also introduce a new condition of block linear degeneracy. We prove that, if it is not satisfied, then pathological behaviours may occur. 
\vspace{0.3cm}

\vspace{0.3cm}


\end{abstract}
\section{Introduction}
\label{sec_introduction} The aim of this work is to describe the limit of the parabolic approximation 
\begin{equation}
\label{eq_introduction_parabolic_approx}
\left\{
\begin{array}{lll}
      \ve_t + \tilde{A} \big( \ve, \, \ee \ve_x \big) \ve_x  =
      \ee \tilde{B}(\ve ) \ve_{xx} \qquad \ve \in \mathbb{R}^N\\
      \tilde \B(\ve (t, \, 0)) \equiv \bar g \\
      \ve (0, \, x) \equiv \bar{v}_0 \\
\end{array}
\right.
\end{equation}
for $\ee \to 0^+$.
In the previous expression, the function $\tilde \B$ is needed to have well posedness in the case the matrix $\tilde{B}$ is not invertible. The function $\tilde \B$ is defined in Section \ref{sec_hypotheses}. 
When the matrix $\tilde{B}$ is indeed invertible, system \eqref{eq_introduction_parabolic_approx} takes the simpler form 
\begin{equation}
\label{e:intro:parapp_inv}
\left\{
\begin{array}{lll}
      \ve_t + \tilde{A} \big( \ve, \, \ee \ve_x \big) \ve_x  =
      \ee \tilde{B}(\ve ) \ve_{xx} \qquad \ve \in \mathbb{R}^N\\
      \ve (t, \, 0)  \equiv \bar{v}_b \\
      \ve (0, \, x) \equiv \bar{v}_0. \\
\end{array}
\right.
\end{equation}
Even if the equations in \eqref{eq_introduction_parabolic_approx} and \eqref{e:intro:parapp_inv} are not necessarily in conservation form, nevertheless the conservative case 
$$
   \ve_t + f (\ve )_x = \ee \Big( \tilde{B}(\ve ) \ve_x \Big)_x  
$$
is included in the previous formulation. Indeed, one can define 
$$
    \tilde{A} (\ve, \, \ee \ve_x ) : = D f ( \ve ) - \ee \Big( \tilde{B} (\ve ) \Big)_x
$$
and obtain an equation like the one in \eqref{eq_introduction_parabolic_approx} or \eqref{e:intro:parapp_inv}.

In the present paper we assume that when $\ee \to 0^+$ the solutions $\ve$ converge, in the sense that will be specified in the following, to a unique limit $v$. 
Since in both \eqref{eq_introduction_parabolic_approx} and \eqref{e:intro:parapp_inv} the initial and boundary data are constant, then the limit $v$ solves a so called boundary Riemann problem, i.e. an hyperbolic initial boundary value problem with constant data. Results in \cite{AmaCol} show that the study of boundary Riemann problems is a key point to determine the semigroup of solutions for an hyperbolic initial boundary value problem.  We will came back to this point at the end of the introduction.

The goal of the work is to determine the value of $v(t, \, x)$ for a.e. point $(t, \, x)$.  In particular, we determine the value of the trace $\bar{v}$ 
of the limit on the axis $x=0$. The reason why this is interesting is the following. Let us focus for simplicity on the case in which 
$\tilde{B}$ is invertible, i.e. on \eqref{e:intro:parapp_inv}. It is known that, in general, the trace of the limit, which we denote by $\bar{v}$, is different from $\bar{v}_b$, the boundary datum imposed in \eqref{e:intro:parapp_inv}.  The relation between $\bar{v}_b$ and $\bar{v}$ is relevant for the study of hyperbolic initial boundary value problems and was first investigated (as far as we know) in \cite{DubLeFl}. Also, in \cite{Gisclon:etudes}  it was proved that the value of $\bar{v}$ in general depends on the choice of the matrx $\tilde{B}$. In other words, if in \eqref{e:intro:parapp_inv} one keeps fixed $\bar{v}_0$, $\bar{v}_b$
and the function $\tilde{A}$ and changes only the matrix $\tilde{B}$, then in general the value of $\bar{v}$ will change, even if the system is in conservation form.

The most interesting points in the work are the following two. First, we cover the characteristic case, which occurs when an eigenvalue of the
matrix $\tilde{A}$ can attain the value $0$. The non characteristic case occurs when none of the eigenvalues of $\tilde A$ can attain the value $0$. The characteristic case is more complicated to handle than the non characteristic one. Loosely speaking, the reason is the following. Suppose that the $k$-th eigenvalue can assume the value $0$. Then in the non linear case we do not know a priori if the waves of the $k$-th family are entering or are leaving the domain.

The second point is that we cover the case of a non invertible viscosity matrix $\tilde{B}$. To tackle this case we assume that the so called Kawashima Shizuta condition holds. We also introduce a new condition of block linear degeneracy. We provide a counterexample which shows that, if this condition is violated, then there may be pathological behaviours, in the sense that will be specified in the following. 

The exposition is organized as follows. In Section \ref{s:intro:over} we give an overview of the paper, introducing the main ideas involved in the analysis. In Section \ref{sec_hypotheses} we discuss the hypotheses assumed in the work. In Section \ref{sec_B_invertible} we give a characterization of the limit of the parabolic approximation \eqref{e:intro:parapp_inv}, i.e. when the viscosity matrix $\tilde{B}$ is invertible. Finally, in Section \ref{sec_B_non_invertible} we discuss the limit of \eqref{eq_introduction_parabolic_approx} when the matrix $\tilde{B}$ is singular.


\subsection{Overview of the paper}
\label{s:intro:over}

\subsubsection{Section \ref{sec_hypotheses}: hypotheses}
\label{subs:intro:hyp}
Section \ref{sec_hypotheses} describes the hypotheses exploited in the work and it is divided into three parts.

Section \ref{subsec_hypotheses_invertible} describes the hypotheses assumed in the case the matrix $\tilde{B}$ is invertible. These hypotheses were already considered in several previous works and they are automatically satisfied when the system admits a dissipative entropy.

In Section \ref{subsec_hypotheses_non_inv} we discuss the hypotheses assumed in the case the matrix $\tilde{B}$ is singular. These hypotheses can be divided into two groups.

The first group is composed by conditions that were already exploited in several previous works
(e.g. in \cite{Kaw:phD, Kaw:Large, KawShi:systems,
KawShi:normal,Serre-Zum, Rousset:residual, Rousset:char}). In particular one assumes that there exists a regular and invertible change of variables $u = u(\ve)$ such that the following holds. If $v$ satisfies 
$$
    \ve_t + \tilde{A} \big( \ve, \, \ee \ve_x \big) \ve_x  =
      \ee \tilde{B}(\ve ) \ve_{xx}, 
$$
then $u$ satisfies  
\begin{equation}
\label{e:inttro:rescaled}
   E(u) u_t + A (u, \, u_x) u_x = B (u) u_{xx}.
\end{equation}
The matrix $B$ has constant rank $r$ and admits the block decomposition
\begin{equation}
\label{e:intro:b}
  B(u) =
  \left(
  \begin{array}{ll}
        0 & 0 \\
	0 & b(u) \\
  \end{array}
  \right)
\end{equation}  
for a suitable $b(u) \in \mathbb{M}^{r \times r}$. Also, $B$ and $A$ satisfy suitable hypotheses that are reasonable from the physical point of view since they were introduced in \cite{Kaw:phD, Kaw:Large, KawShi:systems,
KawShi:normal} relying on examples with a physical meaning. In particular, we assume that the so called Kawashima Shizuta condition is satisfied.

Apart from these hypotheses, we introduce a new condition of block linear degeneracy, which is the following. Let 
\begin{equation}
\label{e:intro:block}   
   A(u, \, u_x)=
  \left(
  \begin{array}{ll}
        A_{11} (u) & A_{12} (u) \\
	A_{21}(u) & A_{22}(u, \, u_x) \\
  \end{array}
  \right)
  \qquad 
  E(u) =
  \left(
  \begin{array}{ll}
        E_{11} (u) & E_{12} (u) \\
	E_{21}(u) & E_{22}(u) \\
  \end{array}
  \right)
\end{equation}
be the block decomposition of $A$ corresponding to \eqref{e:intro:b}, namely $A_{11}$ and $E_{11}$ belong to $ \mathbb{M}^{(N-r) \times (N-r)}$ and $A_{22}$ and $E_{22}$ belong to $\mathbb{M}^{r \times r}$. The condition of block linear degeneracy says that, for every given real number $\sigma$, the dimension of the kernel of $[A_{11} (u) - \sigma E_{11} (u)]$ is constant with respect to $u$. In other words, the dimension of the kernel may vary as $\sigma$ varies, but it cannot change when $u$ varies. 

Block linear degeneracy is not just a technical condition. Indeed, in Section \ref{subsubsec_hypotheses_two_examples} we discuss counterxamples which show how, when the block linear degeneracy is violated, one can have pathological behaviors. More precisely, we exhibit examples in which block linear degeneracy does not hold and there is a solution of \eqref{e:inttro:rescaled} which is not $\mathcal{C}^1$. These can be considered a pathological behaviour since one usually expects that the parabolic approximation has a regularizing effect.

On the other side, block linear degeneracy is not an optimal condition, in the following sense. It is possible to show that block linear degeneracy is satisfied by the Navier Stokes equation written using Lagrangian coordinates, but it is not satisfied by the Navier Stokes equation written in Eulerian coordinates. On the other side the two formulations of the Navier Stokes equation are equivalent, provided that the density of the fluid is strictly positive. This remark was first proposed by Frederic Rousset in \cite{Rousset:ns} and it suggests that it is interesting to look for a condition strong enought to prevent pathological behaviours but at the same time sufficiently weak to be satisfied by the Navier Stokes equation written in Eulerian coordinates. This problem is tackled in the forthcoming paper \cite{BiaSpi:per}. 

In Section \ref{subsub_hyp_datum} it is defined the function $\B$ which is used to define the boundary condition in \eqref{eq_introduction_parabolic_approx}. The point here is the following. 

Consider an hyperbolic initial boundary value problem, and for simplicity let us focus on the conservative case: 
\begin{equation}
\label{e:intro:hyp}
\left\{
\begin{array}{lll}
       u_t + f(u)_x = 0 \\
       u (0, \, x) = \bar{u}_0 \\
       u (t, \, 0) = \bar{u}
\end{array}
\right.
\end{equation}
It is known that if one assigns a datum $\bar{u} \in \mathbb{R}^N$ (i.e. if one assigns $N$ boundary conditions), then the initial boundary value problem  
\eqref{e:intro:hyp} may be ill posed, in the following sense. In general, there is no function $u$ which is a solution of 
$$
  u_t + f(u)_x = 0 
$$
in the sense of distributions, which 
assumes the initial datum for $t =0$ and satisfies 
$$
   \lim_{x \to 0^+ } u (t, \, x) = \bar{u}
$$
for almost every $t$. Also, it is know that a necessary condition to obtain a well posed problem is to assign a number of conditions on the boundary which is in general smaller then $N$. Assume the boundary is non characteristic, i.e. assume that none of the eigenvalues of the jacobian $D f (u)$ can attain the value $0$. In this case, one can impose on the boundary datum a number of conditions equal to the number of positive eigenvalues of $Df(u)$.

We can now came back to the parabolic equation
\begin{equation}
\label{e:intro:pareq}
  u_t + A(u, \, u_x ) u_x = B(u) u_{xx}.
\end{equation}
Let us write $u = ( u_1, \, u_2)^T$, where $u_1 \in \mathbb{R}^{N-r}$, $u_2 \in \mathbb{R}^r$. Here $r$ is the rank of $B$, as in \eqref{e:intro:b}. With this notations equation \eqref{e:intro:pareq} can be rewritten as 
\begin{equation}
\label{e:intro:pareqsplit}
\left\{
\begin{array}{ll}
       u_{1 t} + A_{11} u_{1 x} + A_{2 1} u_{2 x} = 0 \\
        u_{2 t } + A_{1 2 } u_{1 x} + A_{2 2} u_{2 x}  =  b u_{2 xx} \\ 
\end{array}
\right.
\end{equation}
Roughly speaking, the reason why one has to introduce the function $\B$ is the following. Let $n_{11}$ be the number of strictly negative eigenvalues of the block $A_{11}$, which itself is a $(N-r) \times (N-r)$ matrix. Also, let $q$ denote the dimension of the kernel of $A_{11}$. The second line in 
\eqref{e:intro:pareqsplit} contains a second order derivative of $u_2$ and hence $u_2$ can be seen as a parabolic component. On the other side, only first derivatives of $u_1$ appear and hence $u_1$ can be seen as an hyperbolic component. Actually, there is an interaction between the two components (this is ensured by the Kawashima-Shizuta condition). On the other side, because of the hyperbolic component one is not completely free to assign the boundary condition in \eqref{eq_introduction_parabolic_approx}. As pointed out in previous works (e.g in \cite{Rousset:residual}) the number of conditions one can impose on the boundary is $N-n_{11}-q$. Indeed, one can impose $r$ conditions on $u_2$. On the other hand, one can impose on $u_1$ a number of conditions equal to the number of positive eigenvalues of $A_{11}$, i.e. to $N-r-n_{11} - q$. Summing up one obtains exactly $N-n_{11} - q$. 

Thus, the function $\tilde \B$ in \eqref{eq_introduction_parabolic_approx} takes values in $\mathbb{R}^{N-n_{11} -q}$ and $\bar{g}$ is a fixed vector in $\mathbb{R}^{N-n_{11}-q}$. The precise definition of $\tilde \B$ is given in section \ref{subsub_hyp_datum} and it is such that the initial boundary value problem \eqref{eq_introduction_parabolic_approx} is well posed. In Section \ref{sub_Bsing_boundary_datum} it is given a more general definition for the function $\B$. 

Section \ref{subsubsec_hypotheses_two_examples} discusses three examples. The first two show that, if the condition of block linear degeneracy is violated, then there may be solution of 
$$
   u_t + A(u, \, u_x ) u_x = B(u) u_{xx}
$$
exhibiting pathological behaviors, in the sense explained before. More precisely, the first example deal with steady solutions   
\begin{equation}
\label{e:intro:bl}
  A(u, \, u_x ) u_x = B(u) u_{xx},
\end{equation}
while the second one deals with travelling waves, 
$$
   [ A (u, \, u' ) - \sigma E (u) ] u' = B u''.
$$
In the previous expression, $\sigma$ represents the speed of the wave and it is a real parameter. Finally, the third example in Section \ref{subsubsec_hypotheses_two_examples} shows that if the rank of the matrix $B$ is not constant, then there may be solutions of \eqref{e:intro:bl} exhibiting pathological behaviours of the same kind discussed before. 

Section \ref{subsec_hypotheses_introduction} discusses the hypotheses that are assumed in both cases, when the matrix $\tilde{B}$ in \eqref{eq_introduction_parabolic_approx} is invertible and when it is not. It is assumed that the system is strictly hyperbolic (see Section \ref{subsec_hypotheses_introduction}  for a definition of strict hyperbolicty). Also, it is assumed that when $\ee \to 0^+$ the solutions of $\ve$ of \eqref{eq_introduction_parabolic_approx} converge to a unique limit. Also, it is assumed that the approximation is stable with respect to the initial and the boundary data and that the limit has finite propagation speed. 

We refer to Section \ref{subsec_hypotheses_introduction} for the exact statement of the hypotheses, here instead we underline another point. The proof of the convergence of $\ve$ in the case of a generic matrix $\tilde{B}$ is still an open problem. However, there are results that provide a justification of our hypotheses. In particular, in \cite{Gisclon:etudes} it is proved the local in time convergence in the case $\tilde{B}$ is invertible, but in general different from the identity. Moreover, in \cite{AnBia} the authors proved the global in time convergence in the case of an artificial viscosity ($\tilde{B}(\ve)$ is identically equal to $I_N$).   
The analysis in \cite{AnBia} exploits techniques that were introduced in \cite{BiaBre:BV, BiaBre:case, BiaBre:center, BiaBrevv}
to deal with the Cauchy problem. In \cite{AnBia} it is proved the same kind of convergence we assume in the present properties. Also, other properties we assume here (stability of the approximation, finite propagation speed of the limit) are as well proved in \cite{AnBia}.  Analogous results were proved in \cite{Spi:2bvp} for a special class of problems with $2$ boundaries.   

Also, we point out that there are several works that study the stability of the approximation in the case of a very general viscosity matrix $\tilde{B}$. Actually, the literature concerning this topic is very wide and hence we will quote only works that concern specifically initial boundary value problems: \cite{Serre-Zum, Rousset:inviscid, Rousset:residual, Rousset:char}.

\subsubsection{Section \ref{sec_B_invertible}: the characterization of the hyperbolic limit in the case of an invertible viscosity matrix}
\label{sus:intro:Binv} Section \ref{sec_B_invertible} discusses the characterization of the limit of \eqref{e:intro:parapp_inv} when the matrix $\tilde{B}$ is invertible. Actually, because of the hypotheses we assume in Section \ref{subsec_hypotheses_invertible}, we study the equivalent (in the sense specified therein) problem 
\begin{equation}
\label{e:intro:parapp_inv2}
\left\{
\begin{array}{lll}
      E(\ue) \ue_t + {A} \big( \ue, \, \ee \ue_x \big) \ue_x  =
      \ee {B}(\ue ) \ue_{xx} \qquad \ue \in \mathbb{R}^N\\
      \ue (t, \, 0)  \equiv \bar{u}_b \\
      \ue (0, \, x) \equiv \bar{u}_0. \\
\end{array}
\right.
\end{equation}
Also, Section \ref{sec_B_invertible} is divided into fours parts.

Section \ref{subsec_Binv_preliminary} collects preliminary results that are needed in the following.

Section \ref{subsec_Binv_riemannsolver} gives a quick review of some results concerning the characterization of the limit in the Riemann problem. These results were introduced in \cite{Bia:riemann}. 

The Riemann problem is a Cauchy problem with a piecewise constant 
inial datum with a single jump. Let us focus for simplicity on the conservative case:  
\begin{equation}
\label{e:iintro:hypc}
\left\{
\begin{array}{ll}
       u_t + f(u)_x =0 \\
       u (0, \, x) = 
       \left\{
       \begin{array}{ll}
              u^- & x < 0 \\
	       u^+ & x \ge 0
       \end{array}
       \right. 
\end{array}
\right.
\end{equation}
A solution of \eqref{e:iintro:hypc} was first described in \cite{Lax} assuming  some technical  hypotheses (i.e. that all the fields are either genuinely non 
linear or linearly degenerate).  Since we will need in the following, here we briefly review the ideas exploited in \cite{Lax} to obtain a solution of  \eqref{e:iintro:hypc}.

Denote by 
$$
   \lambda_1 (u) < \dots < \lambda_N (u) 
$$
the eigenvalues of the jacobian $D f (u)$ and by $r_1 (u) \dots r_N (u)$ the corresponding eigenvectors. For simplicity, we assume that all the fields are genuinely non linear. In this case, one can chose the orientation of $r_i$ in such a way that $ \nabla  \lambda_i \cdot r_i > 0$. 
For every $i = 1 \dots N$ we denote by $\mathcal{S}^i_{s_i} (u^+)$ a curve in $\mathbb{R}^N$ which is parameterized by $s_i$ and which enjoys the following property.  For every value  $\mathcal{S}^i_{s_i} (u^+)$ there exists a speed 
$\lambda$ close to $\lambda_i (u^+)$ such that the Rankine Hugoniot condition is satisfied, i.e. 
\begin{equation}
\label{e:intro:rh}
    f (u^+ ) - f \Big( \mathcal{S}^i_{s_i} (u^+) \Big) = \lambda \Big[ u^+ - \mathcal{S}^i_{s_i} (u^+)  \Big].
\end{equation}  
Also, let $\mathcal{R}^i_{s_i} (u^+)$ be the integral curve of $r_i (u)$ starting at $u^+$, in other words 
$ \mathcal{R}^i_{s_i} (u^+)$  is the solution of the Cauchy problem 
$$
\left\{
\begin{array}{lll}
          \displaystyle{ \frac{d}{ds}  \mathcal{R}^i = r_i ( \mathcal{R}^i_{s} (u^+) ) }\\
          \\
            \mathcal{R}^i_{0} (u^+) = u^+ \\           
\end{array}
\right.
$$
Finally, let  $\mathcal{T}^i_{s_i} (u^+)$ be defined as follows:
$$
     {T}^i_{s_i} (u^+)  = 
     \left\{
     \begin{array}{ll}
                 \mathcal{R}^i_{s_i} (u^+) & s_i \ge 0 \\
                  \mathcal{S}^i_{s_i} (u^+) & s_i < 0 \\
     \end{array}
     \right.
$$
In \cite{Lax} it is proved that  ${T}^i_{s_i} (u^+)$ 
is a $\mathcal{C}^2$ curve such that 
$$
   \frac{d  {T}^i_{s_i} (u^+)}{ d s_i} \Bigg|_{ s_i = 0} = r_i ( u^+ ).
$$
In the following, we will say that  ${T}^i_{s_i} (u^+)$  is the i-th curve of admissible states. Indeed, every state  ${T}^i_{s_i} (u^+)$  can be connected to $u^+$ by either a rarefaction wave or a shock which is admissible in the sense of Liu.
 In other words, if $s_i \ge 0$ then the solution of the Riemann problem 
\begin{equation}
\label{e:intro:riei}
\left\{
\begin{array}{ll}
     u_t + f(u)_x =0 \\
       u (0, \, x) = 
       \left\{
       \begin{array}{ll}
              u^+ & x \ge 0 \\
	     \mathcal{T}^i_{s_i} (u^+)  &  x < 0
       \end{array}
       \right. 
\end{array}
\right.
\end{equation}
is 
\begin{equation}
\label{e:Intro:rare}
    u(t, \, x ) = 
    \left\{
    \begin{array}{lllll}
               u^+   &   \quad  x   \ge t \lambda_i (u^+)  \\
               \\
              \displaystyle{ T^i_{s_i}     (u^+) }            & \quad  \displaystyle{ x  \leq   t   \lambda_i  \Big( T^i_{s_i} (u^+)      \Big) }   \\
              \\
             \displaystyle{ T^i_{s} (u^+)}  &  \quad   \displaystyle{  t \lambda_i \Big( T^i_{s_i} (u^+) \Big)  < x <  t \lambda_i (u^+), \quad x =  t   \lambda_i \Big(T^i_{s} (u^+) \Big)} \\
    \end{array}
    \right.
\end{equation}
The meaning of the third line is that $u(t, \, x )$ is equal to ${T}^i_{s} (u^+)$ (the value assumed by the curve at the point $s$) when $x$ is exactly equal to $\lambda_i$ evaluated at the point 
${T}^i_{s} (u^+)$. The value of $u$ is well defined because of the condition $\nabla \lambda_i \cdot r_i > 0$.  

On the other side, if $s_i < 0$ then the solution of the Riemann problem \eqref{e:intro:riei} is 
\begin{equation}
\label{e:Intro:shock}
    u(t, \, x ) = 
    \left\{
    \begin{array}{lll}
               u^+   &   \quad  x   \ge t \lambda   \\
               \\
              \displaystyle{ T^i_{s_i}     (u^+) }            & \quad  x  \leq   t   \lambda.  \\
    \end{array}
    \right.
\end{equation}
The speed $\lambda$ satisfies the Rankine Hugoniot condition and it is close to $\lambda_i (u^+)$.

In this way one obtains $N$ curve of admissible states $\mathcal{T}^1_{s_1} (u^+)$, $\dots$ ${T}^N_{s_N} (u^+)$. To define the solution of the Riemann problem \eqref{e:iintro:hypc} one can proceed as follows. Consider the function
\begin{equation}
\label{e:intro:rsolver}
    \psi (s_1 \dots s_N, \, u^+ ) = {T}^1_{s_1} \circ \dots T^{N-1}_{s_{N-1}} \circ T^N_{s_N} (u^+).
\end{equation}
With the notation $T^{N-1}_{s_{N-1}} \circ T^N_{s_N} (u^+)$ we mean that the starting point of $T^{N-1}_{s_{N-1}} $ is $ T^N_{s_N} (u^+)$, i.e. that
$$
    T^{N-1}_{0} \circ T^N_{s_N} (u^+) = T^N_{s_N} (u^+).
$$
It is proved in \cite{Lax} that the map $\psi$ is locally invertible with respect to $s_1 \dots s_N$. In  other words, the values of $s_1 \dots s_N$ are uniquely determined if one imposes
$$
    u^- = \psi (s_1 \dots s_N, \, u^+ ), 
$$ 
 at least if $u^-$ is a close enough to $u^+$. One takes $u^-$ as in \eqref{e:iintro:hypc} and obtains the values $s_1 \dots s_N$. Indeed, we assume, here and in the following, 
 $|u^+ - u^-| <<1$. Once $s_1 \dots s_N$ are known, one can obtain the limit gluing together pieces like \eqref{e:Intro:rare} and \eqref{e:Intro:shock}. 
 
 The construction in \cite{Lax} was extended in \cite{Liu:riemann} to more general systems. Also, in \cite{Bia:riemann} it was given a characterization of the limit of 
 \begin{equation}
 \label{e:intro:limitrie}
\left\{
\begin{array}{lll}
      \ue_t + {A} \big( \ue, \, \ee \ue_x \big) \ue_x  =
      \ee {B}(\ue ) \ue_{xx} \qquad \ue \in \mathbb{R}^N\\
      \ue (0, \, x) =   
       \left\{
       \begin{array}{ll}
              u^+ & x \ge 0 \\
	       u^- &  x < 0
       \end{array}
       \right. 
\end{array}
\right.
\end{equation} 
 when $|u^+ - u^-|$ is sufficiently small (under hypotheses slightly different from the ones we consider here).

The construction works as follows. Consider a travelling profile 
for 
$$
   u_t + {A} \big( u, \, u_x \big) u_x  =
       {B}(u ) u_{xx},
$$
i.e. $u (x - \sigma t)$ such that
$$
    B (u ) u'' = \Big( A(u, \, u' )  - \sigma E  ( u )     \Big) u'.
$$
In the previous, expression, the speed of the wave $\sigma$ is a real parameter.
Then $u$ solves 
\begin{equation}
\label{eq_intro_system_trav} \left\{
\begin{array}{lll}
      u' = v \\
      B(u )v' = \Big( A(u, \, v) - \sigma E(u) \Big)v \\
      \sigma' = 0.
\end{array}
\right.
\end{equation}
The point $( u^+, \, \vec{0}, \, \lambda_i (u^+)  )$ is an equilibrium for \eqref{eq_intro_system_trav}. Also, one can prove that a center manifold around $( u^+, \, \vec{0}, \, \lambda_i (u^+)  )$  has dimension $N+2$.   

We recall here that every center manifold is invariant with respect to \eqref{eq_intro_system_trav} and moreover satisfies the following property: let $(u^0, \, p^0, \, \sigma^0)$ belong to $\mathcal{M}^{c}$ 
and denote by $(u(x), \, p(x), \, \sigma (x)$ the orbit starting at $(u^0, \, p^0, \, \sigma^0)$. Then 
$$
    \lim_{x \to + \infty} \Big(u(x), \, p(x), \, \sigma (x)  \Big) e^{-  c x  / 2} 
    = (\vec 0, \, \vec 0, \, 0)  \qquad    \lim_{x \to -  \infty} \Big(u(x), \, p(x), \, \sigma (x)  \Big)  e^{-  c x  / 2} 
    = (\vec 0, \, \vec 0, \, 0) . 
$$
The constant $c$ in the previous expression is strictly positive, depends on the matrix $E^{-1}A$ and it is defined in Section \ref{sec_hypotheses} (Hypothesis \ref{hyp_hyperbolic_I}).  

Fix a center manifold $\mathcal{M}^c$. 
If $(u, \, v, \, \sigma)$ is a solution to \eqref{eq_intro_system_trav} laying on $\mathcal M^c$, then 
\begin{equation}
\label{eq_intro_reduction} \left\{
       \begin{array}{lll}
             u' = v_i \tilde{r}_i (u, \, v_i, \, \sigma_i) \\
             v_i' = \phi_i(u, \, v_i, \, \sigma_i) v_i \\
             \sigma_i' =0. \\
\end{array}
\right.
\end{equation}
The functions $\tilde{r}_i$ and $\phi_i$ are defined in Section \ref{subsub_Binv_speed+}. 

The construction of $T^i_{s_i} u^+$ works as follows. Fix $s_i>0$ and consider  
 the following fixed point problem, defined on a interval $[0, \, s_i]$: 
\begin{equation}
\label{eq_intro_fixedpt} \left\{
\begin{array}{lll}
      u(\tau) = u^+ + {\displaystyle \int_0^{\tau} \tr(u(\xi), \, v_i(\xi), \,
      \sigma_i(\xi))d \xi }  \\
      v_i (\tau) = f_i (\tau, \, u, \, v_i, \, \sigma_i ) -
      \mathrm{conc} f_i (\tau, \, u, \, v_i, \, \sigma_i )\\
      \sigma_i(\tau)=  {  \frac{1}{c_E  (\bar{u}_0) }  \displaystyle \frac{d}{d \tau}
      \mathrm{conc} f_i (\tau, \, u, \, v_i, \, \sigma_i )}. \\
\end{array}
\right.
\end{equation}
We have used the following notations:
$$
    f_i (\tau) = \int_0^{\tau} \tilde{\lambda}_i [u_i, \, v_i, \, \sigma_i] (\xi) d \xi,
$$
where 
$$
    \tilde{\lambda}_i [u_i, \, v_i, \, \sigma_i] (\xi) = \phi_i \Big( u_i (\xi), \, v_i (\xi), \, \sigma_i (\xi) \Big) + c_E (\bar{u}_0) \sigma.
$$
Also,  $\mathrm{conc} f_i$ denotes the concave
envelope of the function $f_i$: 
$$
    \mathrm{conc} f_i  (\tau)=  \inf  \{ h (s): \; h \mathrm{ \; is \; concave}, \; h (y) \ge f_i (y) \; \forall \, y \in [0, \, s_i] \} .
$$
One can show that the fixed point problem \eqref{eq_intro_fixedpt} admits a unique solution. 
 
The link between \eqref {eq_intro_fixedpt} and \eqref{eq_intro_system_trav} is the following: let $(u_i, \, v_i, \, \sigma_i)$ satisfy \eqref{eq_Binv_fixedpt}. Assume that $v_i < 0$ on $]a, \, b[$ and that 
$v_k (a) = v_k (b) = 0$. Define $ \alpha_i (\tau) $ as the solution of the Cauchy problem 
\begin{equation*}
\left\{ 
\begin{array}{ll}
            \displaystyle{ \frac{d \alpha_i}{ d \tau}   =  -  \frac{1}{v_i (\tau)}}\\
           \alpha_i ( a+ b / 2) = 0 \\
 \end{array}
\right.
\end{equation*}
 then $(u_i \circ \alpha_i, \, v_i \circ \alpha_i, \, \sigma_i \circ \alpha_i)$ is a solution to \eqref{eq_Binv_reduction} satisfying 
 $$
     \lim_{x \to - \infty} u_i \circ \alpha_i (x) = u_i (a) \qquad \lim_{x \to + \infty} u_i \circ \alpha_i (x) = u_i (b).
 $$
Thus, $ u_i (a)$ and $ u_i (b)$ are connected by a travelling wave profile. 

On the other side, if $v_i \equiv 0$ on the interval $[c, \, d]$, then the following holds. Consider $\mathcal{R}^i_{s_i} u (c)$, the integral curve of $r_i (u)$ such that $\mathcal{R}^i_0 u(c) = u (c)$. Then $u(d)$ lays on $\mathcal{R}^i_{s_i} u (c)$, thus $u(c)$ and $u(d)$ are connected by a rarefaction or by a contact discontinuity. 

If $s_i<0$, one considers a fixed problem like \eqref{eq_intro_fixedpt}, but instead of the concave envelope of $f_i$ one takes the convex envelope:
$$
    \mathrm{conv} f_i  (\tau)=  \sup \{ h (s): \; h \mathrm{ \; is \; convex}, \; h (y) \leq f_i (y) \; \forall \, y  \in [0, \, s_i]\} .
$$
Again, one can prove the existence of a unique fixed point $(u_i, \, v_i, \, \sigma_i)$. 

The curve $T^i_{s_i} u^+$ is defined setting 
$$
  T^i_{s_i} u^+ : = u (s_i). 
$$
This curve contains states that are connected to $u^+$ by rarefaction waves and shocks with speed close to $\lambda_i(u^+)$.

If $u^- = T^i_{s_i} u^+$, then the limit of the approximation \eqref{e:intro:limitrie} is 
\begin{equation}
\label{e:intro:s+:limit}
     u  (t, \, x) = 
     \left\{
     \begin{array}{ll}
                 u^+   &  x \leq \sigma_i (0) t \\
                 u_i (s)      &  x = \sigma_i (s) t \\
                 u_i (s_i)=T^i_{s_i} u^+    &  x \ge \sigma_i (s_i) t \\
     \end{array}
    \right.
\end{equation}
In the previous expression, $\sigma_i$ is given by \eqref{eq_intro_fixedpt} and it is a monotone non increasing function. 

It can be shown that in the case of conservative systems with only genuinely non linear or linearly degenerate fields the i-th curve of admissible states $T^i_{s_i} u^+$ defined in \cite{Bia:riemann} coincides with the one described in \cite{Lax}.  
Once $T^i_{s_i} u^+$ is known, then one defines $\psi$ as in \eqref{e:intro:rsolver} and find the limit gluing together pieces like \eqref{e:intro:s+:limit}. 
\vspace{1cm}

In Section \ref{subsubsec_Binv_nonchar} we give a characterization of the limit of 
\begin{equation}
\label{e:intro:bdatinv}
\left\{
\begin{array}{lll}
      E(\ue) \ue_t + {A} \big( \ue, \, \ee \ue_x \big) \ue_x  =
      \ee {B}(\ue ) \ue_{xx} \qquad \ue \in \mathbb{R}^N\\
      \ue (t, \, 0)  \equiv \bar{u}_b \\
      \ue (0, \, x) \equiv \bar{u}_0. \\
\end{array}
\right.
\end{equation}
when the boundary is not characteristic, i.e. when none of the eigenvalues of $E^{-1}A$ can attain the value $0$. The idea is to costruct a locally invertible map $\phi $ which describes all the states that can be connected to $\bar{u}_0$, in the sense that is specified in the following. Loosely speaking, the map $\phi$ represents for the initial boundary value problem what the map $\psi$ defined in \eqref{e:intro:rsolver} represents  for the Cauchy problem. Once $\phi$ is defined, one takes $\bar{u}_b$ as in \eqref{e:intro:bdatinv} and imposes 
$$
   \bar{u}_b = \phi (\bar{u}_0, \, s_1 \dots s_N).
$$
If $\bar{u}_0$ and $\bar{u}_b$ are sufficiently close, then this relation uniquely determines the values of $s_1 \dots s_N$. Once these values are known, then the limit $u(t, \, x)$ is uniquely determined. More precisely, one can define the value of $u(t, \, x)$ for a.e. $(t, \, x)$. We will come to this point later. 

The construction of the map $\phi$ works as follows. Denote by $\lambda_1 (u) \dots \lambda_N(u)$ the eigenvalues of $E^{-1}(u) A(u, \, 0)$ and by 
$r_1 (u) \dots r_N(u) $ the associated eigenvectors. Also, assume that for every $u$, 
$$
  \lambda_1 (u) < \dots \lambda_n (u) < - c < 0 < c < \lambda_{n +1 } < \dots  \lambda_N(u)
$$  
In other words, $n$ is the number of negative eigenvalues of $E^{-1}(u)A(u, \, 0)$. These eigenvalues are 
real because this is one of the hypotheses listed in Section \ref{sec_hypotheses} (Hypothesis \ref{hyp_hyperbolic_I}). 

For $i=n +1  \dots N$ consider the i-th curve of admissible states. Fix $N-n$ parameters $s_{n+1} \dots s_N$ and define
$$
  \bar{u} = T^{n+1}_{s_{n+1}} \circ \dots T^{N-1}_{s_{N-1}}\circ T^N_{s_N} \bar{u}_0. 
$$ 
As before, the notation $T^{N-1}_{s_{N-1}}\circ T^N_{s_N} \bar{u}_0$ means that 
the starting point of   $T^{N-1}_{s_{N-1}}$ is  \\ $T^N_{s_N} \bar{u}_0$, $T^{N-1}_{0}= T^N_{s_N} \bar{u}_0$.
Thanks to the results in \cite{Bia:riemann} we quoted before, $\bar{u}$ is connected to $\bar{u}_0$ by a sequence of rarefaction waves and shocks with stricly positive speed.

To complete the construction, one considers steady solutions of 
$$
  E(u) u_t + A (u, \, u_x ) u_x = B (u) u_{xx}
$$
i.e. couples $(U, \, p)$ such that 
\begin{equation}
\label{e:intro:steady}
\left\{
\begin{array}{ll}
       U' = p \\
       B (u) p' = A(u, \, p) p \\
\end{array}
\right.
\end{equation}
The point $(\bar{u}, \, 0)$ is an equilibrium for \eqref{e:intro:steady}. As shown in section \ref{subsubsec_Binv_nonchar}, the stable manifold around $(\bar{u}, \, 0)$ has dimension $n$, i.e. has dimension equal to the number of strictly negative eigenvalues of $E^{-1} A$. Also, the following holds. Let $ \psi $ be a map that parameterizes the stable manifold, then $\psi$  takes values into $\mathbb{R}^N \times \mathbb{R}^N$ and it is defined on a space of dimension $n$. To underline the dependence on $\bar{u}$ we will write $\phi(\bar{u}, \, s_1 \dots s_n)$. Denote by $\pi_u$ the projection
\begin{equation*}
  \begin{split}
        \pi_u : 
  &     \mathbb{R}^N \times \mathbb{R}^N \to \mathbb{R}^N \\
  &     (u, \, p) \mapsto u  \\
  \end{split}
\end{equation*}
Fix $n$ parameters $s_1 \dots s_n$ and consider $\pi_u\circ \psi (\bar{u}, \, s_1 \dots s_n)$. Consider the problem
$$
\left\{
\begin{array}{lll}
        A (u, \, u_x ) u_x = B (u) u_{xx} \\
       u (0) = \pi_u \circ \psi (\bar{u}, \, s_1 \dots s_n),
\end{array}
\right.
$$
then there exists a unique solution of this problem such that 
$$
  \lim_{x \to + \infty} u (x) = \bar{u}. 
$$
Setting $\ue (x) : = u ( x / \ee)$, one finds a solution of 
$$
   A (\ue, \, \ee \ue_x ) \ue_x = B (\ue) \ue_{xx}
$$
such that $\ue (0) = \pi_u \circ \psi (\bar{u}, \, s_1 \dots s_n)$ and for every $x >0$, 
$$
   \lim_{\ee  \to 0^+} \ue (x) = \bar{u}.
$$
In this sense, we say that there is a loss of boundary condition when passing to the hyperbolic limit, because the boundary condition $\pi_u \circ \psi (\bar{u}, \, s_1 \dots s_n)$ disappears in the limit. We point out that the idea of studying steady solutions to take into account the loss of boundary condition was already exploited in many previous work, e.g in \cite{Gisclon:etudes}. 

To complete the characterization of the limit, we define $\phi$ as 
$$
  \phi (\bar{u}_0, \, s_1 \dots s_N) = \pi_u \circ \psi 
  \Big(T^{n+1}_{s_{n+1}} \circ \dots T^{N-1}_{s_{N-1}}\circ T^N_{s_N} \bar{u}_0, \, s_1 \dots s_n  \Big).
$$
In section \ref{subsubsec_Binv_nonchar} we prove that the map $\phi$ is locally invertible, i.e. the values of $s_1 \dots s_N$ are completely determined is one sets  
$$
  \bar{u}_b = \phi (\bar{u}_0, \, s_1 \dots s_N).
$$
We take the same $\bar{u}_b $ as in \eqref{e:intro:bdatinv}. As pointed out at the beginning of the paragraph, once $s_1 \dots s_N$ are known, then the value of the self similar solution $u(t, \, x)$ is determined for a.e. $(t, \, x)$. One can indeed glue together pieces like \eqref{e:intro:s+:limit}. In particular, it turns out that the trace of $u (t, \, x)$ on the axis $x=0$ is 
$$
  \bar{u} = T^{n+1}_{s_{n+1}} \circ \dots T^{N-1}_{s_{N-1}}\circ T^N_{s_N} \bar{u}_0.
$$
As we have already underlined before, the relation between $\bar{u}_b$ (the boundary datum in \eqref{e:intro:bdatinv}) and $\bar{u}$ is interesting for the study of hyperbolic initial boundary value problems. 
\vspace{1cm}

In Section \ref{subsubsec_Binv_char} we give a characterization of the limit of the parabolic approximation 
\begin{equation}
\label{e:intro:bdatinv2}
\left\{
\begin{array}{lll}
      E(\ue) \ue_t + {A} \big( \ue, \, \ee \ue_x \big) \ue_x  =
      \ee {B}(\ue ) \ue_{xx} \qquad \ue \in \mathbb{R}^N\\
      \ue (t, \, 0)  \equiv \bar{u}_b \\
      \ue (0, \, x) \equiv \bar{u}_0. \\
\end{array}
\right.
\end{equation}
when the boundary is characteristic, i.e. one eigenvalue of $E^{-1} A$ can attain the value $0$. The characterization of the limit works as follows. We costruct a locally invertible map $\phi $ which describes all the states that can be connected to $\bar{u}_0$. Once $\phi$ is defined, one takes $\bar{u}_b$ as in \eqref{e:intro:bdatinv2} and imposes 
$$
   \bar{u}_b = \phi (\bar{u}_0, \, s_1 \dots s_N).
$$
If $\bar{u}_0$ and $\bar{u}_b$ are sufficiently close, then this relation uniquely determines the values of $s_1 \dots s_N$. Once these values are known, then the limit $u(t, \, x)$ is uniquely determined.

Formally, the idea is the same as in Section \ref{subsubsec_Binv_nonchar}. However, the construction of the map $\phi$ is definitely more complicated in the boundary chararacteristic case. 

Roughly speaking, the reason is the following. Let $\lambda_1 (u) \dots \lambda_N (u)$ be the eigenvalues of $E^{-1}(u)A(u, \, 0)$. They are real by Hypothesis \ref{hyp_hyperbolic_I}. Assume 
$$
  \lambda_1 (u) < \dots \lambda_{k-1} (u) < -c < \lambda_k (u) < c < \lambda_{k +1 } < \dots \lambda_N(u), 
$$
where $c$ is a suitable positive constant. In other words, there are at least $k-1$ strictly negative eigenvalues, $N-k$ strictly positive eigenvalues and one eigenvalue close to $0$.

Define 
$$
  \bar{u}_k = T^{k+1}_{s_{k+1}} \circ \dots T^{N-1}_{s_{N-1}}\circ T^N_{s_N} \bar{u}_0,
$$
then $\bar{u}_k$ is connected to $\bar{u}_0$ by rarefaction waves and shocks with stricly positive speed. We now want to define the k-th curve of admissible states. To define $T^k_{s_k} \bar{u}_k$, we might try to consider the fixed point problem 
\begin{equation}
\label{e:intro:kappa}
 \left\{
\begin{array}{lll}
      u(\tau) = \bar{u}_k + {\displaystyle \int_0^{\tau} \tilde{r}_k(u(\xi), \, v_k(\xi), \,
      \sigma_k(\xi))d \xi }  \\
      v_k (\tau) = f_k (\tau, \, u, \, v_k, \, \sigma_k ) -
      \mathrm{conc} f_k (\tau, \, u, \, v_k, \, \sigma_k )\\
      \sigma_k(\tau)=  {  \frac{1}{c_E  (\bar{u}_0) }  \displaystyle \frac{d}{d \tau}
      \mathrm{conc} f_k (\tau, \, u, \, v_k, \, \sigma_k )}, \\
\end{array}
\right.
\end{equation}
where $\tilde{r}_k$, $f_k$ and $c_E$ are the same as in \eqref{eq_intro_fixedpt}. However, if we consider \eqref{e:intro:kappa} we are not doing the right thing. Indeed, we might have that the speed $\sigma_k$ is negative at a certain point $\tau$. Since eventually we want to define the limit $u(t, \, x)$ as in \eqref{e:intro:s+:limit}, we want $\sigma_k$ to be greater then $0$. 

Another problem is the following. Consider the system satisfied by steady solutions of 
$$
  E(u) u_t + A (u, \, u_x ) u_x = B (u) u_{xx},
$$
i.e. consider 
\begin{equation}
\label{e:intro:steady2}
  \left\{
  \begin{array}{ll}
       U' = p \\
       B (u) p' =  A(u, \, p) p \\
  \end{array}
  \right.
\end{equation}
Also, consider the equilibrium point $(\bar{u}_k, \, 0)$. Let 
$\mathcal{M}^s$ be the stable manifold of \eqref{e:intro:steady2} around $(\bar{u}_k, \, 0)$. For simplicity, assume $\lambda_k (\bar{u}_k)  = 0$. Then, there might be a solution $(U, \, p)$ such that 
$$
  \lim_{x \to + \infty} \Big( U(x), \, p(x) \Big) = \Big( \bar{u}_k, \, 0\Big)
$$ 
but $(U, \, p)$ does {\sl not} belong to the stable manifold. However, this kind of solution should be taken into account when we study the loss of boundary condition. 

To tackle these difficulties one can proceed as follows. Instead of the fixed point problem  \eqref{e:intro:kappa}, one considers  
\begin{equation}
\label{e:intro:kappaok}
 \left\{
\begin{array}{lll}
      u(\tau) = \bar{u}_k + {\displaystyle \int_0^{\tau} \tilde{r}_k(u(\xi), \, v_k(\xi), \,
      \sigma_k(\xi))d \xi }  \\
      v_k (\tau) = f_k (\tau, \, u, \, v_k, \, \sigma_k ) -
      \mathrm{mon} f_k (\tau, \, u, \, v_k, \, \sigma_k )\\
      \sigma_k(\tau)=  {  \frac{1}{c_E  (\bar{u}_0) }  \displaystyle \frac{d}{d \tau}
      \mathrm{mon} f_k (\tau, \, u, \, v_k, \, \sigma_k )}. \\
\end{array}
\right.
\end{equation}
In the previous expression, $ \mathrm{mon} f_k $ denotes the monotone concave envelope of the function $f_k$, 
$$
   \mathrm{mon} f_k (\tau) = \inf \Big\{  g (\tau): \; g (s ) \ge f_k( s ) \; \forall s, \;  g \; \mathrm{concave \; monotone \; non \; decreasing \; in } \; [0, \, s_i] \;    \Big\}.
$$
Some properties of the monotone envelope are discussed in Section \ref{subsubsec:Binv_monotone}, here we stress that $\mathrm{mon} f_k$ is a concave and non decreasing function, thus the solution $\sigma_k$ of \eqref{e:intro:kappaok}  is always non negative. 

Also, the following holds.  Denote by $(u_k, \,v_k, \, \sigma_k)$ the solution of \eqref{e:intro:kappaok} (existence and uniqueness are proved in Section \ref{subsubsec_Binv_char}).  Define 
\begin{equation}
\label{e:intro:bars}
     \bar{s} = \min \{  s: \; \sigma_k (s) = 0 \} 
\end{equation}
and
$$
     \underline{s} = \max \{  s: \; \sigma_k (s) = 0, \; v_k (s) =0  \}.  
$$
 Assume $0 < \bar{s} \leq \underline{s} < s_k$. Then $u_k (\bar{s})$ is connected to $\bar{u}_k$ by a sequence of rarefaction and shocks with positive speed. Also, one can show that there exists a steady solution $U$,
 $$
     A(U, \, U_x) U_x = B(U ) U_{xx}
 $$
 such that $U(0) = u_k (s_k)$ and 
 $$
     \lim_{x \to +  \infty} U(x) = u_k (\underline{s}). 
 $$
 However, in general this solution {\sl does not } belong to the stable manifold of system \eqref{e:intro:steady2}. This means that considering system \eqref{e:intro:kappaok} we also manage to take into account the converging steady solutions we were missing considering just the stable manifold of \eqref{e:intro:steady2}.

Heuristically, to complete the construction one should consider the stable manifold of  \eqref{e:intro:steady2} and hence take into account the steady solution that, for $x$ that goes to $+ \infty$, converge to $u_k (\underline{s})$ with fast exponential decay, in the sense specified in Section \ref{subsubsec_Binv_char}. Actually, the situation is more complex. The reason, loosely speaking, is the following . There may be a solution $U$ that converges to $u_k (\underline{s})$ and such that some of its components converge with fast exponential decay, but other components converge  more slowly. This possibility is not covered if we consider only the solutions laying on the stable manifold and those given by \eqref{e:intro:kappaok}. To take into account this possibility some technical tools are introduced. More precisely, one considers suitable manifolds:  center stable manifold and uniformly stable manifold. The existence of these manifolds is a consequence of results in \cite{KHass}, but some of the most important properties are recalled in Section \ref{subsubsec_Binv_char}. 

Eventually, one manages to define a locally invertible function $\phi$. One then takes $\bar{u}_b$ as in \eqref{e:intro:bdatinv2} and imposes 
$$
   \bar{u}_b = \phi (\bar{u}_0, \, s_1 \dots s_N).
$$
If $\bar{u}_0$ and $\bar{u}_b$ are sufficiently close, then this relation uniquely determines the values of $s_1 \dots s_N$. Once these values are known, then the limit $u(t, \, x)$ is uniquely determined and can be obtained gluing together pieces like \eqref{e:intro:s+:limit}. In particular, it turns out that the trace of $u$ on the $x$ axis is $u_k (\bar{s})$, where $u_k$ solves \eqref{e:intro:kappaok} and $\bar{s}$ is given by \eqref{e:intro:bars}.

\subsubsection{Section \ref{sec_B_non_invertible}: the characterization of the hyperbolic limit in the case of a singular viscosity matrix}
\label{sus:intro:Bsing}
In Section \ref{sec_B_non_invertible} we discuss the characterization of the limit of \eqref{eq_introduction_parabolic_approx} when the matrix $\tilde{B}$ is not invertible. Actually, because of the hypotheses introduced in Section \ref{sec_hypotheses}, one studies the limit of 
\begin{equation}
\label{e:intro:bsing:pa}
\left\{
\begin{array}{lll}
      \ue_t + {A} \big( \ue, \, \ee \ue_x \big) \ue_x  =
      \ee {B}(\ue ) \ue_{xx} \qquad \ue \in \mathbb{R}^N\\
       \B(\ue (t, \, 0)) \equiv \bar g \\
      \ue (0, \, x) \equiv \bar{u}_0 \\
\end{array}
\right.
\end{equation}
This system is equivalent to \eqref{eq_introduction_parabolic_approx} in the sense specified in Section \ref{subsec_hypotheses_non_inv}.  Also, $\B = \tilde \B \Big( u (\ve )\Big) $: a precise definition is given Section \ref{sec_hypotheses}. In particular, $\B$ ensures that the initial boundary value problem \eqref{e:intro:bsing:pa} is well posed. Section \ref{sec_B_non_invertible} is divided into several parts. 

In Section \ref{sub_Bsing_preliminary} we 
introduce some preliminary results. The point here is the following. In Section \ref{sec_B_invertible} we give a characterization of the hyperbolic limit  when the viscosity matrix is invertible. A key point in the analysis is the study of travelling waves 
\begin{equation}
\label{e:intro:tw}  
  [ A(U, \, U' ) - \sigma E(U)] U' = B(U) U'' 
\end{equation}
and of steady solutions
\begin{equation}
\label{e:intro:steady3}   
   A(U, \, U' ) U' = B(U) U'' .
\end{equation}
To give a characterization of the hyperbolic limit when the viscosity matrix is not invertible, we have to study again systems \eqref{e:intro:tw} and \eqref{e:intro:steady3}. However, being the viscosity matrix $B$ singular, a technical difficulty arises. Let us focus, for simplicity, on the case of steady solutions. If $B$ is invertible, we can write 
\begin{equation}
\label{e:intro:steadys}  
  \left\{ 
  \begin{array}{ll}
         U' = p \\
	 p' = B(u)^{-1}A(u, \, p) p \\
  \end{array}
  \right.
\end{equation}
In this way, we write system \eqref{e:intro:steady3} in an explicit form. On the other side, if the matrix $B$ is singular, additional work is required to reduce \eqref{e:intro:steady3} in a form like \eqref{e:intro:steadys}. This is indeed done in Section \ref{sub_Bsing_preliminary}. What we actually obtain is not a $2N$-dimensional first order system like \eqref{e:intro:steadys}, but a lower dimensional first order system. The exact dimension depends on the structure of the matrix $A$, in the sense specified in Section \ref{sub_Bsing_preliminary}. 

In Section \ref{subsub_Bsing_travelling_waves} we review the characterization of the hyperbolic limit in the case of a Riemann problem, i.e. the limit of    
\begin{equation*}
\left\{
\begin{array}{lll}
      \ue_t + {A} \big( \ue, \, \ee \ue_x \big) \ue_x  =
      \ee {B}(\ue ) \ue_{xx} \qquad \ue \in \mathbb{R}^N  \\
      \ue (0, \, x) =
      \left\{ 
      \begin{array}{ll}
            u^+ & x \ge 0 \\
	    u^- & x < 0 \\
      \end{array}
      \right.
\end{array}
\right.
\end{equation*}
As for the case of an invertible $B$ (Section \ref{subsub_Binv_speed+}), the key point in the analysis is the description of the i-th curve of admissible states $T^i_{s_i} u^+$. However, there are technical difficulties due to the fact that $B$ is not invertible. Actually, in Section \ref{subsub_Bsing_travelling_waves} we only give a sketch of the construction, and we refer to \cite{Bia:riemann} for the complete analysis. 

In Section \ref{subsub_Bsing_dimension} we introduce a technical lemma. The problem is the following. Consider a steady solution \eqref{e:intro:steady3} and assume that it is written in an explicit form like \eqref{e:intro:steadys}. This is possible thanks to the considerations carried on in Section \ref{sub_Bsing_preliminary}. Given an equilibrium point for this new system, consider the stable manifold around that equilibrium point. For reasons explained in Section \ref{sub_Bsing_Riemann_solvernonchar}, we need to know the dimension of this stable manifold. Lemma \ref{lem_Bsing_crucial} ensures that the dimension of the stable manifold is equal to $n - n_{11} - q$, where $n$ is the number of strictly negative eigenvalues of $A$, $n_{11}$ is the number of strictly negative eigenvalues of the block $A_{11}$ and $q$ is the dimension of the kernel of the block $A_{11}$. The block $A_{11}$ is defined by \eqref{e:intro:block}. Lemma \ref{lem_Bsing_crucial} gives an answer to a question left open in \cite{Rousset:char}.

In Section \ref{sub_Bsing_Riemann_solvernonchar} we discuss the characterization of the hyperbolic limit of 
\begin{equation}
\label{e:intro:bsing:pa2}
\left\{
\begin{array}{lll}
      \ue_t + {A} \big( \ue, \, \ee \ue_x \big) \ue_x  =
      \ee {B}(\ue ) \ue_{xx} \qquad \ue \in \mathbb{R}^N\\
       \B(\ue (t, \, 0)) \equiv \bar g \\
      \ue (0, \, x) \equiv \bar{u}_0 \\
\end{array}
\right.
\end{equation}
when the matrix $B$ is singular, but the boundary is not characteristic, i.e. none of the eigenvalues of $E^{-1}(u) A(u, \, 0)$ can attain the value $0$. 

To provide a characterization, we construct a map $\phi (\bar{u}_0, \, s_{n_{11} + q+ 1} \dots s_N)$ which enjoys the following properties. It takes values in $\mathbb{R}^N$ and when $ s_{n_{11} + q+ 1} \dots s_N$ vary it describes all that states that can be connected to $\bar{u}_0$, in the sense specified in the following. Note that similar functions are constructed in Section \ref{subsubsec_Binv_nonchar} and Section \ref{subsubsec_Binv_char} and  they are used to describe the limit of the parabolic approximation when the viscosity is {\sl invertible}. However, in those cases the functions depend on $N$ variables, while in the case of a singular matrix $B$ the function $\phi$ depends    
only on $(N-n_{11} - q)$ variables.

Also, in the case of a singular viscosity we have to compose $\phi$ with the function $\B$, which is used to assign the boundary datum in \eqref{e:intro:bsing:pa2}. It is possible to show that the map $\B \circ \phi$ is locally invertible, in other words the values of  $ s_{n_{11} + q+ 1} \dots s_N$ are uniquely determined if one imposes 
\begin{equation}
\label{e:intro:bphi}
   \B \circ \phi ( s_{n_{11} + q+ 1} \dots s_N) = \bar{g},
\end{equation}
provided that $| \B (\bar{u}_0) - \bar{g} |$ is small enough. We will plug the same $\bar{g}$ as in \eqref{e:intro:bsing:pa2}. Once the values of $s_{n_{11} + q+ 1} \dots s_N$ are known, the limit $u(t, \, x)$ can be determined a.e. $(t, \, x)$ in the same way as in Sections \ref{subsubsec_Binv_nonchar}  and \ref{subsubsec_Binv_char}. In particular, one can determine exactly the value of the trace of the limit $u$ on the axis $x=0$. As pointed out before, this is important for the study of hyperbolic initial boundary value problems.

The construction of the map $\phi$ works as follows.  Denote as before by $n$ the number of the eigenvalues of $E^{-1}(u) A(u, \, 0)$ that are strictly negative, since the boundary is not characteristic  the number of strictly positive eigenvalues is $N-n$.  For $ i = n+1 \dots N$, let $T^i_{s_i}$ be the i-th curve of admissible states, whose construction is reviewed in Section \ref{subsub_Bsing_travelling_waves}. Fix $s_{n +1} \dots s_N$ and define
\begin{equation}
\label{e:intro:traceBsing}
         \bar{u} = T^{n +1}_{s_{n +1}} \circ \dots T^{N-1}_{s_{N-1}} \circ T^N_{s_N} \bar{u}_0.
\end{equation}
As before, with the notation $T^{N-1}_{s_{N-1}} \circ T^N_{s_N} \bar{u}_0$ we mean that the starting point of $T^{N-1}_{s_{N-1}}$ is $T^N_{s_N} \bar{u}_0$, i.e.
$$
    T^{N-1}_{{0}} \circ T^N_{s_N} \bar{u}_0 =  T^N_{s_N} \bar{u}_0. 
$$   
The value $\bar{u}$ is connected to $\bar{u}_0$ by a sequence of rarefactions and shocks with strictly positive speed. 

To complete the construction, we consider steady solutions of 
$$
     E(u) u_t + A(u, \, u_x) u_x = B(u) u_{xx},
$$ 
i.e. we consider solutions of 
\begin{equation}
\label{e:intro:steadybsing}
    A(u, \, u_x) u_x = B(u) u_{xx}
\end{equation}
In Section \ref{subsub_Bsing_prel_normal} we discuss this system and we explain how to write it as first order O.D.E..  Also, in Section \ref{subsub_Bsing_dimension} we study 
 the stable manifold of \eqref{e:intro:steadybsing}  around an equilibrium point such that $u = \bar{u}$ and $u_x = 0$. In particular, we prove that this manifold has dimension $n-n_{11} - q$. 
 Let $\psi (\bar{u}, s_{n_{11} + q} \dots s_n)$ a map that parameterizes the stable manifold (we are also putting in evidence the dependence on $\bar{u}$). For every  $ s_{n_{11} + q} \dots s_n$ there exists a solution $u$ of \eqref{e:intro:steadybsing} such that $u (0) = \psi (\bar{u}, s_{n_{11} + q} \dots s_n)$ and 
 $$
     \lim_{x \to + \infty } u(x) = \bar{u}.
 $$
setting $\ue (x ) : = u ( x / \ee)$ one obtains a steady solution of 
$$
     E(\ue) \ue_t + A(\ue, \, \ee \ue_x) \ue_x = B(\ue) \ue_{xx},
$$
such that for every $x >0$
$$
   \lim_{\ee \to 0^+} \ue (x) = \bar{u}.
$$
In other words, we experience again a possibile loss of boundary condition from the parabolic approximation to the hyperbolic limit.

The map $\phi$ is then defined as follows:
$$
     \phi ( \bar{u}_0, \, s_{n_{11} + q+ 1} \dots s_N) : = \psi \Big(T^{N-1}_{s_{0}} \circ T^N_{s_N} \bar{u}_0 =  T^N_{s_N} \bar{u}_0, \,      s_{n_{11} + q +1 } \dots s_n     \Big). 
$$
In Section \ref{sub_Bsing_boundary_datum} we prove that $\B \circ \phi$ is locally invertible with respect to $ s_{n_{11} + q+ 1} \dots s_N$.  

Here, instead, we point out the following.  The fact that the dimension of the stable manifold of \ref{e:intro:steadybsing} has dimension $n - n_{11} - q$ is proved in 
Section \ref{subsub_Bsing_dimension} and it is not, a priori, obvious. However, $n - n_{11} - q$  is exactly the number that makes things works if one wants the function $\B \circ \phi$ to be locally invertible, in the sense that $\B \circ \phi$ takes values in $\mathbb{R}^{N - n_{11} - q}$ and hence to be locally invertible it should depend on $N - n_{11} - q$ variables. Since we have to take into account $N-n$ curve of admissible states, we are left with  $n - n_{11} - q$ variables for the stable manifold.
\vspace{1cm}

In Section \ref{subsubsec_Binv_nonchar} we provide a characterization of the limit of 
\begin{equation}
\label{e:intro:pappb}
\left\{
\begin{array}{lll}
      \ue_t + {A} \big( \ue, \, \ee \ue_x \big) \ue_x  =
      \ee {B}(\ue ) \ue_{xx} \qquad \ue \in \mathbb{R}^N\\
       \B(\ue (t, \, 0)) \equiv \bar g \\
      \ue (0, \, x) \equiv \bar{u}_0 \\
\end{array}
\right.
\end{equation}
when the matrix $B$ is singular and the boundary is characteristic, i.e one of the eigenvalues of $A$ can attain the value $0$. Actually this case requires no new ideas, in the sense that
one can combine the techniques that are described in Section \ref{subsubsec_Binv_char} (to deal with the fact that  the boundary is characteristic) and in Section \ref{subsub_Bsing_nonchar} (to deal with the singularity of 
the matrix $B$).
\vspace{1cm}

Finally, in Section \ref{sub_Bsing_boundary_datum}  we introduce a technical lemma which guarantees that the function $\B \circ \phi$ in \eqref{e:intro:bphi}
is locally invertible. We also provide a more general definition for the function $\B$. Such a definition still ensures that the initial boundary value problem 
\eqref{e:intro:pappb} is well posed and that the map $\B \circ \phi$ is locally invertible.

 \subsection{Unquoted references}
 \label{s:intro:ur}
Existence results for hyperbolic initial value problems were obtained 
 in \cite{Good} and
\cite{SabTou:mixte} relying on an adaptation of the Glimm scheme of approximation introduced
in \cite{Gli}.

These results were later improved by mean of wave front
tracking techniques. These techniques were used  in a
series of papers
(\cite{Breft,Bre:Gli,BreCol,BreCraPi,BreG,BreLF,BreLew,BreLY}) to
establish the well posedness of the Cauchy problem.  A comprehensive
account of the stability and uniqueness results for the Cauchy
problem for a system of conservation laws can be found in
\cite{Bre}. There are many references for a general introduction to system of conservation laws, for example \cite{Daf:book} and
to \cite{Serre:book}. As concerns initial boundary value problems,
in \cite{Ama} the existence results in \cite{Good} and
\cite{SabTou:mixte} were substantially improved, while well posedness results were obtained in \cite{DonMar}.  

In \cite{AmaCol} it is studied the limit of the wave front tracking approximation. In particular, the authors extended to initial boundary value problems the definition of Standard Riemann Semigroup. Such a notion was introduced in 
\cite{Bre:Gli} for the Cauchy problem.  Roughly
speaking,  the analysis in  \cite{AmaCol}  guarantees that  to identify the semigroup of solutions it is enough to consider the behaviour of the semigroup in the case the initial and the boundary data are constant.
It hence one of the most important motivations for the study of the parabolic approximation \eqref{eq_introduction_parabolic_approx}. 
\vspace{1cm}

\noindent {\bf Acknowledgments:} 
the authors express their gratitude to Denis Serre for having carefully read a 
previous version of the work and for his useful remarks. Also,  
the authors wish to thank Jeffrey
Rauch for useful observations concerning the definition of boundary
condition for the hyperbolic-parabolic system. Finally, we would like to thank Fabio Ancona for his useful remarks concerning the concave envelope of a function. 
\section{Hypotheses}
\label{sec_hypotheses} This section introduces the hypotheses
exploited in the work.

The exposition is organized as follows. In Section
\ref{subsec_hypotheses_invertible}  we introduce the hypotheses we impose on system \eqref{e:intro:parapp_inv} when the viscosity matrix   
$\tilde{B}$ is invertible. In Section \ref{subsec_hypotheses_non_inv} we discuss the hypotheses exploited when the matrix $\tilde B$ in \eqref{eq_introduction_parabolic_approx} 
is singular. Finally, in Section \ref{subsec_hypotheses_introduction} we introduce the hypotheses needed in both the cases, when the matrix $\tilde B$ in \eqref{eq_introduction_parabolic_approx} 
is invertible and when it is not. 

In Section \ref{subsub_hyp_datum} we also discuss three examples which show that, if some of the conditions introduced are not satisfied, then there may be pathological behaviors, in the sense specified there.  
Also, we define the function $\tilde \B$ used to assign the boundary condition in \eqref{eq_introduction_parabolic_approx}. 

In \cite{BiaSpi:per} it will be discussed a way to extend to a more general setting the condition of block linear degeneracy, which is the third condition in Hypothesis \ref{hyp_Bsing}.

\subsection{The hypotheses assumed when the viscosity matrix is invertible}
\label{subsec_hypotheses_invertible} In this section it is
considered the system
\begin{equation}
\label{eq_hypotheses_system}
       v_t + \tilde{A}(v, \, v_x) v_x = \tilde{B}(v) v_{xx}
\end{equation}
in the case the viscosity matrix $\tilde{B}(v)$ is invertible. We assume the following.
\begin{hyp}
\label{hyp_invertible}
       There exists an (invertible) smooth change of variables 
       $v = v(u)$ such that \eqref{eq_hypotheses_system}
       is equivalent to system
       \begin{equation}
       \label{eq_hypotheses_equivalent}
             E(u) u_t + A(u, \, u_x) u_x =
             B(u) u_{xx},
       \end{equation}
       where
       \begin{enumerate}
       \item for any $u$, the matrix
       $E(u)$
       is real, symmetric and positive definite: there exists a constant $c_E(u)$
       such that
       \begin{equation*}
              \forall \, \xi \in \mathbb{R}^N, \;
              \langle E(u) \xi, \, \xi \rangle
              \ge c_E(u) |\xi|^2.
       \end{equation*}
       \item for any $u$,
       the matrix $A(u, \, \, 0)$ is symmetric.
       \item for any $u$, the
       viscosity matrix
       ${B}(u)$
       is real and there
       exists a constant $c_B(u)>0$
       such that
       \begin{equation*}
              \forall \, \xi \in \mathbb{R}^N, \;
              \langle B(u) \xi, \, \xi \rangle
              \ge c_B(u) |\xi|^2.
       \end{equation*}
       \end{enumerate}
\end{hyp}
In particular, the matrix $E(u)$ may be the differential $Dv(u)$
of the change of coordinates.

The initial boundary value problem
\eqref{e:intro:parapp_inv}  is equivalent to
\begin{equation}
\label{eq_hypotheses_parabolic_approximation} \left\{
\begin{array}{lll}
      \ue_t + A \big( \ue, \, \ee \ue_x \big)\ue_x =
      \ee B(\ue ) \ue_{xx} \qquad \ue \in \mathbb{R}^N \\
      \ue (t, \, 0) \equiv \bar{u}_b \\
      \ue (0, \, x) \equiv \bar{u}_0, \\
\end{array}
\right.
\end{equation}
where $v(\bar{u}_0) = \bar{v}_0$ and $v(\bar{u}_b) = \bar{v}_b$.
Indeed, one can verify that, thanks to the invertibility of the
viscosity matrix, it is possible to assign a full boundary
condition $\bar{u}_b \in \mathbb{R}^N$.

      In the case system \eqref{eq_hypotheses_system} is in conservation form
      \begin{equation}
      \label{eq_hypotheses_conservation_form}
             v_t + f(v)_v = \Big( \tilde{B}(v)v_x \Big)_x,
      \end{equation}
      Hypothesis
      \ref{hyp_invertible} is guaranteed by the presence of a
      dissipative and convex entropy.

      For completeness we recall 
      that the entropy $\eta$ is dissipative for
      \eqref{eq_hypotheses_conservation_form}
      if for any $v$ there exists a constant
      $c_D(v)$ such that
      \begin{equation}
      \label{eq_hypotheses_dissipative}
            \forall \, \xi \in \mathbb{R}^N, \;
            \langle D^2 \eta(v) \xi, \, \tilde{B}(v)\xi \rangle
            \ge c_D(v) |\xi|^2
      \end{equation}
      It is known (see for example \cite{GodRav}) that if system
      \eqref{eq_hypotheses_system}
      admits a convex entropy $\eta$, then there exists a
      symmetrizing and
      invertible change of variable $v= v(u)$. More precisely,
      if the inverse function $u(v)$ is defined by $u = \nabla \eta
      (v)$, then $v$ satisfies system
      \eqref{eq_hypotheses_equivalent}, with $A(u, \, 0)$ symmetric and given by
      \begin{equation*}
            A(u, \, 0) = Df \big( v(u) \big)
            \Big( D^2 \eta \big( v(u) \big) \Big)^{-1}-
            \frac{\partial}{\partial x} \bigg( B \big( v(u) \big)
            \Big( D^2 \eta \big( v(u) \big) \Big)^{-1}
            \bigg) \bigg|_{u_x =0} = Df \big( v(u) \big)
            \Big( D^2 \eta \big( v(u) \big) \Big)^{-1}
      \end{equation*}
      Moreover,
      \begin{equation*}
             E(u) = \Big( D^2 \eta \big( v(u) \big) \Big)^{-1}
      \end{equation*}
      is symmetric (being the inverse of an hessian) and
      positive definite by convexity, while the dissipation
      condition \eqref{eq_hypotheses_dissipative} guarantees that
      \begin{equation*}
            B(u) = \tilde{B} \big( v(u) \big)
            \Big( D^2 \eta \big( v(u) \big) \Big)^{-1}
      \end{equation*}
      is positive definite. Hence
      Hypothesis \ref{hyp_invertible} is satisfied.

\subsection{The hypotheses
assumed when the viscosity matrix is not invertible}
\label{subsec_hypotheses_non_inv} The aim of this section is to
introduce the hypotheses that will be needed in Section
\ref{sec_B_non_invertible} to study system
\begin{equation}
\label{eq_hypotheses_system_Bsing}
       v_t + \tilde{A}(v, \, v_x) v_x = \tilde{B}(v) v_{xx}
\end{equation}
in the case the viscosity matrix is singular.
\begin{hyp}
\label{hyp_Bsing}
       There exists an (invertible) smooth change of coordinates
       $v = v(u)$ such that \eqref{eq_hypotheses_system_Bsing}
       is equivalent to system
       \begin{equation}
       \label{eq_hypotheses_true_system_Bsing}
             E(u) u_t +  A(u, \, u_x) u_x =
             B (u) u_{xx},
       \end{equation}
       where
       \begin{enumerate}
       \item the matrix $B(u)$ has constant rank $r < N$ and
       admits the block decomposition
       \begin{equation}
       \label{eq_hypotheses_Bsing_B}
             B(u) =
             \left(
             \begin{array}{cc}
                   0 & 0 \\
                   0 & b(u) \\
             \end{array}
             \right)
       \end{equation}
       with $b(u) \in \mathbb{M}^{r \times r}$.
       Moreover, for every $u$ there exists a constant $c_b
       (u)>0$ such that
       \begin{equation*}
              \forall \, \xi \in \mathbb{R}^r, \;
              \langle b(u) \xi, \, \xi \rangle
              \ge c_b(u) |\xi|^2.
       \end{equation*}
       \item for any $u$
       the matrix ${A}(u, \, \, 0)$ is symmetric. Moreover,
       the block decomposition of $A(u, \, u_x)$ corresponding to
       \eqref{eq_hypotheses_Bsing_B} takes the form
       \begin{equation}
       \label{eq_hypotheses_Bsing_A}
              \left(
             \begin{array}{cc}
                   A_{11}(u) & A_{12}(u) \\
                   A_{21}(u) & A_{22}(u, \, u_x) \\
             \end{array}
             \right),
       \end{equation}
       with $A_{11} \in \mathbb{M}^{(N-r) \times (N-r)}$,
       $A_{12} \in \mathbb{M}^{(N-r) \times r}$,
       $A_{21} \in \mathbb{M}^{r \times (N-r) }$ and
       $A_{22} \in \mathbb{M}^{r \times r}$. Namely, only in the block
       $A_{22}$ depends on $u_x$.
       \item a condition of block linear degeneracy holds.
       More precisely, for any $\sigma \in \mathbb{R}$ the dimension of the kernel 
       $\text{ker} [A_{11}(u) - \sigma E_{11}(u )]$ does not depend on $u$,
       but only on $\sigma$.

       This holds in particular for the limit cases
       in which the dimension of $\text{ker} (A_{11} - \sigma E_{11})$ is equal to
       $N-r$ or to $0$.
       \item the so called Kawashima condition holds: for any $u$,
       \begin{equation*}
              \text{ker} \big( B(u) \big) \cap \{\mathrm{eigenvectors \; of \; }
              E^{-1}(u) A(u, \,
              0)\}= \emptyset
       \end{equation*}
       \item for any $u$, the matrix
       $E(u)$
       is real, symmetric and positive definite: there exists a constant $c_E(u)$
       such that
       \begin{equation*}
              \forall \, \xi \in \mathbb{R}^N, \;
              \langle E(u) \xi, \, \xi \rangle
              \ge c_E(u) |\xi|^2.
       \end{equation*}
       In the following, we will denote by
       \begin{equation}
       \label{e:hyp:e}
              E(u) =
              \left(
             \begin{array}{cc}
                   E_{11}(u) & E_{12}(u) \\
                   E_{21}(u) & E_{22}(u) \\
             \end{array}
             \right)
       \end{equation}
       the block decomposition of $E(u)$ corresponding to
       \eqref{eq_hypotheses_Bsing_B}.
       \end{enumerate}
\end{hyp}
The change of variable $v=v(u)$ guarantees that the
initial-boundary value problem
\eqref{eq_introduction_parabolic_approx} is equivalent to
\begin{equation}
\label{eq_hypotheses_parabolic_approximation_Bsing} \left\{
\begin{array}{lll}
      E(\ue) \ue_t + A \big( \ue, \, \ee \ue_x \big) \ue_x  =
      \ee B(\ue ) \ue_{xx} \qquad \ue \in \mathbb{R}^N\\
      \B(\ue (t, \, 0)) \equiv \bar g \\
      \ue (0, \, x) \equiv \bar{u}_0, \\
\end{array}
\right.
\end{equation}
where $v(\bar{u}_0) = \bar{v}_0$, $v(\bar{u}_b) = \bar{v}_b$,
$\B(\ue) = \tilde \B(v(\ue))$.

Relying on the block decomposition of the viscosity matrix
described in Hypothesis \ref{hyp_Bsing}, in Section
\ref{subsub_hyp_datum} the explicit definition of the function
$\B$ is introduced. Moreover, in Section
\ref{sub_Bsing_boundary_datum} such a definition is extended to a
more general formulation.

In order to justify the assumptions summarized in Hypothesis
\ref{hyp_Bsing} it will be made reference to the works of
Kawashima and Shizuta, in particular to \cite{Kaw:phD} and to
\cite{KawShi:normal}.

In particular, in \cite{KawShi:normal} it is assumed that the
system in conservation form
\begin{equation}
\label{eq_hypotheses_conservation}
       v_t + f(v)_x = \Big( \tilde{B}(v) v_x \Big)_x
\end{equation}
admits a convex and dissipative entropy $\eta$ which moreover
satisfies
\begin{equation*}
       \Big( D^2 \eta (v) \Big)^{-1} \tilde{B}(v)^T=
       \tilde{B}(v) \Big( D^2 \eta (v) \Big)^{-1}.
\end{equation*}
If one performs the change of variables defined by $w = \nabla
\eta (v)$, finds that system \eqref{eq_hypotheses_system_Bsing} is
equivalent to
\begin{equation*}
        \hat{E}(w) w_t +  \hat{A}(w) w_x =
             \Big( \hat{B} (w) w_{x} \Big)_x,
\end{equation*}
with $\hat{A}$ and $\hat{B}$ symmetric.

It is then introduced the assumption \vspace{0.2cm}

\hspace*{2cm} {\it Condition N:} the kernel of $\hat{B}(w)$ does
not depend on $w$. \vspace{0.2cm}

\noindent and it is proved that {\it Condition N} holds if and
only if there exists a change of variable $w= w(u)$ which ensures
that system \eqref{eq_hypotheses_conservation} is equivalent to
\begin{equation*}
      E(u) u_t +  A(u, \, u_x) u_x =
      B (u) u_{xx}
\end{equation*}
for some $E(u)$ that satisfies condition 5 in Hypothesis
\ref{hyp_Bsing} and some $A(u, \, u_x)$ and $B(u)$ as in condition
2 and 1 respectively.

Moreover, it is shown that {\it Condition N} is verified in the
case of several systems with physical meaning. \vspace{0.2cm}

On the other side, the fourth assumption in Hypothesis
\ref{hyp_Bsing} is the so called Kawashima condition and was
introduced in \cite{Kaw:phD}. Roughly speaking, its meaning is to
ensure that there exists an interaction between the parabolic and
the hyperbolic component of system
\eqref{eq_hypotheses_system_Bsing} and hence to guarantee the
regularity of a  solution \eqref{eq_hypotheses_system_Bsing}.

Examples
\ref{ex_kernel} and \ref{ex_travelling} in Section
\ref{subsubsec_hypotheses_two_examples}  show that, if the condition of block linear degeneracy is violated, then one encounters pathological behaviours, in the following sense. One may find a solution of  
\eqref{eq_hypotheses_system_Bsing} that are not continuously differentiable. This is a pathological behaviour in the sense that, when one introduces a parabolic approximation, one expects a regularizing effect. 
In Section \ref{subs:intro:hyp} we explain why it is interesting to look for an extension of the condition of block linear degeneracy to a more general setting. This problem will be tackled in \cite{BiaSpi:per}.

Finally, Example \ref{ex_rank}  show that if the first condition in Hypothesis \ref{hyp_Bsing}
is violated, then one can have pathological behaviours like the one described before. In other words, if the rank of $B$ can vary, then one may find a solution of \eqref{eq_hypotheses_true_system_Bsing} 
which is not continuously differentiable. 

\subsubsection{An explicit definition of boundary
datum for the parabolic problem} \label{subsub_hyp_datum} Thanks
to Hypothesis \ref{hyp_Bsing} system \eqref{e:intro:parapp_inv}
is equivalent to 
\begin{equation}
\label{e:hyp:datum:sys}
\left\{
\begin{array}{lll}
      E(u) u_t + A \big( u, \,  u_x \big) u_x  =
      B(u ) u_{xx} \qquad u \in \mathbb{R}^N\\
      \B (u) (t, \, 0)) \equiv \bar g \\
      u (0, \, x) \equiv \bar{u}_0, \\
\end{array}
\right.
\end{equation}
where $\B (\ue ) = \tilde \B [ \ve ( u ) ]$. In this section we define the function $\B$. Once $\B (u)$ is known, one can obtain $\tilde \B (v)$  
exploiting the fact that the map $v(u)$ is invertible. In Section \ref{sub_Bsing_boundary_datum} we extend this definition to a more 
general setting.

Let $r$ be, as in the statement of Hypothesis \ref{hyp_Bsing}, the rank of the matrix $B$. Decompose $u$ as $ u = (u_1, \, u_2 )^T$, where $u_1 \in \mathbb{R}^{N-r}$ and $u_2 \in \mathbb{R}^r$. then the equation 
$$
          E(u) u_t + A \big( u, \,  u_x \big) u_x  =
      B(u ) u_{xx}     
$$
can be written as 
\begin{equation}
\label{e:hyp:parhyp}
\left\{
\begin{array}{ll}
         E_{11} u_{1t} +  E_{12} u_{2t} + A_{11} u_{1x} + A_{12} u_{2x} =0 \\
               E_{21} u_{1t} +  E_{22} u_{2t} + A_{21} u_{1x} + A_{22} u_{2x} = b u_{2xx} \\
\end{array}
\right.
\end{equation}
The blocks $E_{11}$, $A_{11}$ and so on are as in \eqref{eq_hypotheses_Bsing_A} and \eqref{e:hyp:e}.

In the first line of \eqref{e:hyp:parhyp} only first order derivatives appear, while in the second line there is a second order derivative $u_{2 xx}$. In this sense, $u_1$ can be regarded as the {\sl hyperbolic} component of  
\eqref{e:hyp:parhyp}, while $u_2$ is the {\sl parabolic} component.  As explained in section \ref{subs:intro:hyp}, one can impose $r$ boundary conditions on $u_2$, while one can impose on $u_1$ a number of boundary conditions equal to the number of eigenvalues of $E_{11}^{-1}A_{11}$ with strictly positive real part. 

We recall that we denote by $n_{11}$ the number of strictly negative eigenvalues of
$A_{11}$, by $q$ the dimension of $\text{ker} A_{11}$ and by $n$
the number of strictly negative eigenvalues of $A$. One can prove that the number of eigenvalues  of $E_{11}^{-1}A_{11}$  with strictly negative real part is equal to $n_{11}$ (see Lemma \eqref{lem_Binv_dimension} in Section  
\ref{sssec:tran}, which was actually introduced in \cite{BenSerreZum:Evans}). Also, the dimension of the kernel of $E_{11}^{-1}A_{11}$   is $q$. 

Let
$\vec{\zeta}_i(u, \, 0) \in \mathbb{R}^{N-r}$ be an eigenvector of
$E_{11}^{-1} A_{11}(u)$ associated to an eigenvalue
$\eta_i $ with non positive real part. Let $Z_i (u, \, 0) \in \mathbb{R}^N$ be defined
by
 \begin{equation}
 \label{eq_bsing_bd_zeta_1}
       Z_i:=
       \left(
       \begin{array}{cc}
              \vec{\zeta}_i \\
              0
       \end{array}
       \right)
\end{equation}
and finally let
\begin{equation*}
             \mathcal{Z}(u):= \text{span}  \langle Z_1(u, \, 0),
             \dots, Z_{n_{11}+q}(u, \, 0)
             \rangle.
\end{equation*}
Finally, let 
\begin{equation*}
      \mathcal{W}(u):= \text{span}  \Big\langle
      \left(
      \begin{array}{cc}
             0 \\
             \vec{e}_1 \\
      \end{array}
      \right),
      \dots,
      \left(
      \begin{array}{cc}
             0 \\
             \vec{e}_r \\
      \end{array}
      \right),
      \left(
      \begin{array}{cc}
             \vec{w}_1 \\
             0 \\
      \end{array}
      \right),
      \dots,
      \left(
      \begin{array}{cc}
             \vec{w}_{N-r-n_{11}-q} \\
             0 \\
      \end{array}
      \right)
      \Big\rangle,
\end{equation*}
where $\vec{e}_i \in \mathbb{R}^r$ are the vectors of a basis in
$\mathbb{R}^r$ and $\vec{w}_j \in \mathbb{R}^{N-r}$ are the
eigenvectors of $E_{11}^{-1} A_{11}$ associated to  eigenvalues with strictly positive real part.

Since $\mathcal{W}(u) \oplus \mathcal{V}(u)= \mathbb{R}^N$, every
$u \in \mathbb{R}^N$ can be written as
\begin{equation}
\label{E:boundecom}
      u = u_w + u_z, \quad u_w \in \mathcal{W}(u),
      u_z \in \mathcal{Z}(u).
\end{equation}
in a unique way.  

We define $\B$ as follows.
\begin{say}
\label{def_bc}
      The function $\B$ which gives the boundary condition in  \eqref{eq_hypotheses_parabolic_approximation_Bsing}
      is
      \begin{equation}
      \label{eq_hyp_beta}
      \begin{split}
           \B: \, & \mathbb{R}^N \to \mathbb{R}^{N - n_{11} -q} \\
      &     u  \; \; \; \; \mapsto u_{z}, \\
      \end{split}
      \end{equation}
      where $u_z$ is the component of $u$ along $\mathcal{W}(u)$,
      according to the decomposition \eqref{E:boundecom}.
\end{say}
We refer to Section \ref{sub_Bsing_boundary_datum} for an extension of this definition.

\subsubsection{Examples}
\label{subsubsec_hypotheses_two_examples}
\begin{ex}
\label{ex_kernel} This example deals with the equation
\begin{equation}
\label{eq_hypotheses_sigma_system}
      A(u, \, u_x) u_x = B(u) u_{xx},
\end{equation}
in the case the condition of block linear degeneracy (the third in
Hypothesis \ref{hyp_Bsing}) does not hold.

System \eqref{eq_hypotheses_sigma_system} is satisfied by the
steady solutions of \eqref{eq_hypotheses_true_system_Bsing}.
One therefore expects the solution of
\eqref{eq_hypotheses_sigma_system} to have good regularity
properties: the example which is going to be discussed show that
if the condition of block linear degeneracy does not hold a
function $u$ satisfying \eqref{eq_hypotheses_sigma_system} may
have the graph illustrated by Figure \ref{fig_example_kernel} and
hence be not $\mathcal{C}^1$. Moreover, the figure suggests that a
pathological behavior typical of the solution of the porous-media
equation may occur, namely it may happen that a solution is
different from zero on a interval and then vanishes identically.
\begin{figure}
\begin{center}
\caption{the graph of the function $u_1(x)$ in Example
\ref{ex_kernel}} \psfrag{u}{$u_1$} \psfrag{x}{$x$}
\psfrag{u0}{$u_1(0)$} \psfrag{v}{$u_1(0)$}
\label{fig_example_kernel}
\includegraphics[scale=0.6]{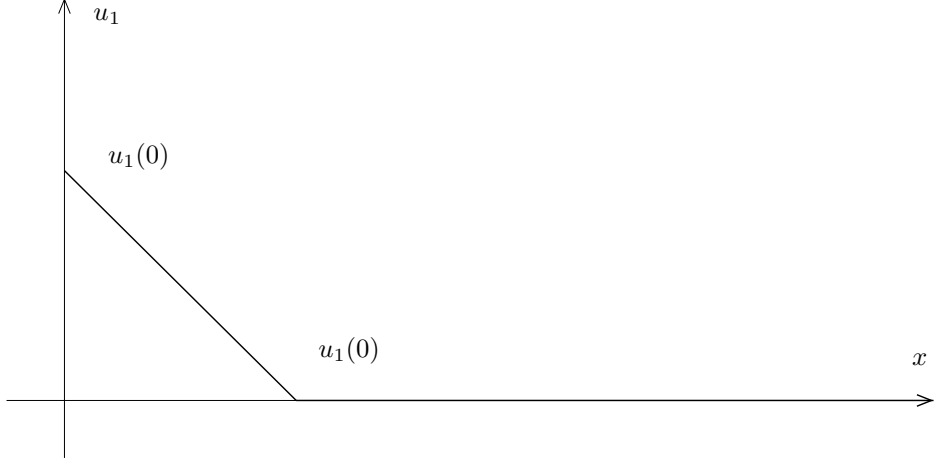}
\end{center}
\end{figure}
Let $u= (u_1, \, u_2 )^T$ and let 
\begin{equation*}
       B(u) =
       \left(
      \begin{array}{ll}
             0 & 0 \\
             0 & 1 \\
      \end{array}
      \right).
\end{equation*}
and
\begin{equation*}
      A(u, \, u_x) =
      \left(
      \begin{array}{ll}
             u_1 & 1 \\
             1 & 0, \\
      \end{array}
      \right).
\end{equation*}
Then the assumption of block linear degeneracy is
not satisfied in a neighborhood of $u_1=0$: we therefore impose the limit
conditions
\begin{equation*}
      \lim_{x \to + \infty} u_1 (x) =0 \qquad
      \lim_{x \to + \infty} u_2 (x) =0.
\end{equation*}
In this case equation \eqref{eq_hypotheses_sigma_system} writes
\begin{equation}
\label{eq_hypotheses_example_kernel}
       \left\{
      \begin{array}{ll}
            u_1 \big( u_{1 x} + 1 \big) = 0 \\
            u_{1} = u_{2x}
      \end{array}
      \right.
\end{equation}
Since the equation satisfied by $u_1$ admits more solutions in a
neighborhood of $u_1 =0$, to introduce a selection principle we
proceed as follows.

We consider a solution  $u^{\nu}= (u_1^{\nu}, \, u^{\nu}_2 )^T$ of \eqref{eq_hypotheses_sigma_system} 
such that 
$$
    \lim_{x \to + \infty }u^{\nu}_1 (x) = \nu.
$$
The parameter $\nu$ is positive and we study the limit $\nu \to 0^+$. Fix the initial datum $u_{1 0} > 0$. The component $u^{\nu}_1$ satisfies the Cauchy problem
\begin{equation}
\label{e:ex:cauchypb}
\left\{
\begin{array}{ll}
          u^{\nu}_{1x} = \displaystyle{ \frac{\nu - u^{\nu}_1 }{ u^{\nu}_1} } \\
          u^{\nu}_1 (0) = u^{\nu}_{10} \\
\end{array}
\right.
\end{equation}
If $\nu$ is sufficiently small then $u_{10} > \nu$ and hence $u_{1x}$ is always negative. In other words, $u^{\nu}_1$ is a monotone non increasing function which satisfies $\nu < u^{\nu}_1 (x) \leq u_{10}$ for every $x$.  Also, note that if $\nu_1 < \nu_2$ then, for every $x$,
\begin{equation}
\label{e:hyp:ex:mono}
    u_{1}^{\nu} (x) \leq u_2^{\nu} (x). 
\end{equation}
This can be deduced by a comparison argument applied to the Cauchy problem \eqref{e:ex:cauchypb}.  Indeed, if $u^{\nu}_1> 0$ then
$$
    \displaystyle{ \frac{\nu - u^{\nu_1}_1 }{ u^{\nu_1}_1} }  < \displaystyle{ \frac{\nu - u^{\nu_2}_1 }{ u^{\nu_2}_1} }. 
$$ 
In particular, from \eqref{e:hyp:ex:mono} we deduce that for every $x$ $u^{\nu}_1(x)$ is monotone decreasing with respect to $\nu$ and hence it admits limit $\nu \to 0^+$.

Denote by $u_1$ the pointwise limit of $u^{\nu}_1$ for $\nu \to 0^+$:  we claim that $u_1$ satisfies  
$$
               u_1 \big( u_{1 x} + 1 \big) = 0.
$$
Indeed, let 
$$
    x_0 : = \min \{ x: \; u^1 (x) =0  \}. 
$$
If $x < x_0$, then by monotoniticity for every $\nu$ and for every $y<x$ $u_{10} \ge u^{\nu}_1 (y )  \ge  u_1 (x) >0$. Also, if $\nu$ is sufficiently small then 
$u^{\nu}_1 (x) >0$ and by monotonicity $u^{\nu}_1(y) \ge u^{\nu}_1 (x) >0$.
Consider the relation
$$
    u^{\nu}_1 (x) = u_{10} + \int_{0}^x \frac{\nu - u^{\nu}_1 (y)}{ u^{\nu}_1 (y)} dy.
$$
We take the limit $\nu \to 0^+$ and, applying Lebesgue's dominated convergence theorem, we get
$$
    u_1 (x) = u_{10} - x.
$$
On the other side, if $x \ge x_0$ then by monotonicity $u_1 (x) \leq 0$. On the other side, $ u_1 (x) \ge 0 $ for every $x$. We conclude that
$u_1 (x) =0$ if $x \ge x_0$. In other words, $x_0 = u_{10}$
and 
$$
    u_1 (x ) = 
    \left\{
    \begin{array}{ll}
                u_{10} - x & x \leq u_{10} \\
                0 & x \ge u_{10}.
    \end{array}
    \right.
$$
This function ha the graph illustrated in Figure \ref{fig_example_kernel} and it is not continuously differentiable. In this sense, we encounter a pathological behaviour. 

\begin{rem}
\label{rem_dim_man}
      An alternative interpretation of the previous
      considerations is the following.

      Consider the family of systems
      \begin{equation}
      \label{eq_hyp_ex_ker1_fam}
      \left\{
      \begin{array}{lll}
            u_1 u_{1x} + u_{2x} =0 \\
            u_{1x} = u_{2xx} \\
            u_1(0) = \upsilon \qquad
            \lim_{x \to + \infty} u_1(x) = y \qquad
            \lim_{x \to + \infty} u_2 (x) =0
      \end{array}
      \right.
      \end{equation}
      parametrized by the limit value $y \in \mathbb{R}$.
      Let $\mathcal{U}_1 (y)$ be the set of the values $\upsilon$
      such that \eqref{eq_hyp_ex_ker1_fam} admits a solution.
      The ODE satisfied by $u_1$ is
      \begin{equation*}
            u_{1x} = \frac{y- u_1}{u_1}
      \end{equation*}
      when $u_1 \neq 0$. From a standard analysis it turns out
      that when $y>0$, then $\mathcal{U}_1 (y) = ]0, \, +
      \infty[$. When $y=0$, $\mathcal{U}_1 (0)= [0, \, + \infty[$,
      while when $y < 0$, $\mathcal{U}_1 (y) = \{ y \}$.

      In other words, $\mathcal{U}_1$ is a manifold of dimension
      $1$ when $y > 0$, a manifold with boundary and with
      dimension $1$ when $y =0$, while when $y < 0$ the dimension
      of the manifold drops to zero, i.e. the manifold reduces to
      a point.
\end{rem}
\end{ex}
\begin{ex}
\label{ex_travelling} Example \ref{ex_kernel}  shows 
that, if the condition of blcok linear degeneracy  (the third in Hypothesis
\ref{hyp_Bsing}) is not satisfied, then one can find a steady solution of
\begin{equation}
\label{e:ex2:sys}
         E (u) u_t + A(u, \, u_x) u_x = B(u) u_{xx} 
\end{equation} 
which is not continuously differentiable. 

As we see in the following sections, steady solutions are important when studying the limit of the parabolic approximation \eqref{e:hyp:datum:sys}.  Other solutions that play an important role 
are travelling wave solutions, i.e $U$ such that 
\begin{equation}
\label{eq_hyp_ex_ker2_trav}
    \Big( A (U, \, U' ) - \sigma E (U) \Big) U' = B(U) U''.
\end{equation}
In the previous expression, $\sigma$ is a real parameter and it is the speed of the travelling wave. More precisely, in the following sections we study travelling waves such that $\sigma$ is close to an eigenvalue of the matrix 
$E^{-1} A$. 

Let
\begin{equation}
\label{eq_hyp_ex_ker2_aeb} A(u, \, u_x) := \left(
\begin{array}{ccc}
       u_1 & 1       & 0 \\
       1      & 1  & 0 \\
       0      & 0       & 0 \\
\end{array}
\right) \qquad B(u) := \left(
\begin{array}{ccc}
       0      & 0       & 0 \\
       0      & 1       & 0 \\
       0      & 0       & 1 \\
\end{array}
\right),
\end{equation}
which are obtained from a system in conservation form
\begin{equation*}
      u_t + f(u)_x = \Big( B(u) u_x \Big)_x
\end{equation*}
where
\begin{equation*}
      f(u) = \Big( u^2_1 / 2 + u_2, \, u_1 + u_2, \, 0 \Big)^T
\end{equation*}
and $B(u)$ is the constant matrix defined by
\eqref{eq_hyp_ex_ker2_aeb}.
\begin{figure}
\begin{center}
\caption{the graph of the function $u_1(x)$ in Example
\ref{ex_travelling}} \psfrag{2}{$2$} \psfrag{u}{$u_1(x)$}
\psfrag{x}{$x$} \psfrag{z}{$u_1$} \label{fig_2}
\includegraphics[scale=0.6]{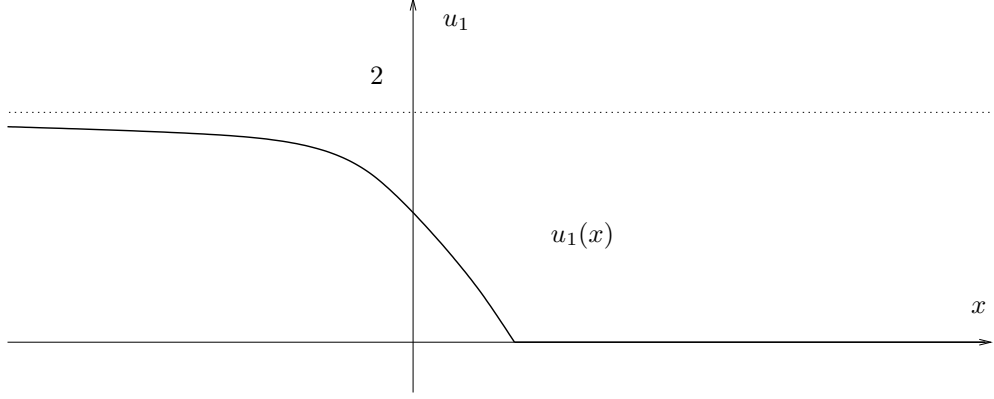}
\end{center}
\end{figure}

The matrix $A$ defined by
\eqref{eq_hyp_ex_ker2_aeb} has an eigenvalue identically equal to
zero: we will focus on the travelling waves with speed $\sigma =0$.

System \eqref{eq_hyp_ex_ker2_trav} can then be rewritten as
\begin{equation}
\label{eq_hyp_ex_ker2_syst} \left\{
\begin{array}{lll}
       u_1 u_{1x} + u_{2x} =0 \\
       u_{1x} + u_{2x} = u_{2xx}\\
        0= u_{3 xx}.
\end{array}
\right.
\end{equation}
As in the case considered in Example \ref{ex_kernel}, the
condition of block linear degeneracy is not satisfied in a
neighborhood of $u_1 =0$ and therefore we study travelling wave
solutions such that
\begin{equation*}
      \lim_{x \to + \infty} u_1 (x) = 0 \qquad
      \lim_{x \to + \infty} u_2 (x) = 0 \qquad
      \lim_{x \to + \infty} u_3 (x) = 0.
\end{equation*}

From system  \eqref{eq_hyp_ex_ker2_syst} and from the previous condition one obtains
\begin{equation*}
\left\{
\begin{array}{lll}
       u_1^2 /2  + u_{2} =0 \\
       u_{1} +  u_{2} = u_{2x}\\
       u_3 \equiv 0.
\end{array}
\right.
\end{equation*}
Hence in the following we ignore the third component of the
solution and we study only the first two lines of the system,
which can be rewritten as
\begin{equation}
\label{eq_hyp_ex_ker2_plus} \left\{
\begin{array}{ll}
       u_2 = - u_1^{2} / 2 \\
       u_1 u_{1x} =  u^2_1 /2  -u_1.
\end{array} \right.
\end{equation}
If $u_1 \neq 0$, the second line is equivalent to
\begin{equation*}
      u_{1x} = u_1 /2   -1.
\end{equation*}
Fix $u_{10}$ such that $0 < u_{10} < 2$. We impose $u_1 (0) = u_{10}$ and we obtain
\begin{equation*}
       u_1 (x) = 2 +( u_{10} - 2 ) e^{ x /2 },
\end{equation*} 
This solution can be extended as far as
$u_1
>0$: when $u_1$ reaches the value zero,
\eqref{eq_hyp_ex_ker2_plus} admits more that one solution.

In order introduce a selection principle, we
proceed as in Example \ref{ex_kernel} and we introduce a family of solutions $u^{\nu} = (u^{\nu}_1, \, u^{\nu}_2, \, u^{\nu}_3)$ of 
 \eqref{eq_hyp_ex_ker2_trav} such that 
\begin{equation}
\label{eq_hyp_ex_ker2_lim}
      \lim_{x \to + \infty} \un_1 (x) = \nu \qquad
       \lim_{x \to + \infty} \un_2 (x) = 0 \qquad
       \lim_{x \to + \infty} \un_3 (x) = 0.
\end{equation}
We also impose $u^{\nu}_1 (0) = u_{10}$. 
In the previous expression,  $\nu$ is a small and positive parameter and we study limit $\nu \to 0^+$. One can repeat the same considerations as in the previous example and conlude that when $\nu \to 0^+$ the solution $u^{\nu}_1$ converges pointwise to 
\begin{equation}
\label{eq_hyp_ex_ker2_sol}
      u_1 (x) =
      \left\{
      \begin{array}{ll}
             u_1 (x) = 2 +( u_1(0) - 2 ) e^{x /2 } \qquad x \leq
             2 \log \Big( 2/ \big(2 - u_{10}  \big) \Big) \\
             0 \qquad \qquad \qquad \qquad
             \qquad \qquad \qquad x >   2 \log
             \Big( 2/ \big(2 - u_{10} \big) \Big)
      \end{array}
      \right.
\end{equation}
This function  has the graph illustrated in Figure \ref{fig_2} and it is not continuously differentiable. 
\end{ex}

\begin{ex}
\label{ex_rank} The aim of this example is to show that if in a
system of the form \eqref{eq_hypotheses_true_system_Bsing} the
rank of the matrix $B$ is not constant, thus contradicting the
first assumption of Hypothesis \ref{hyp_Bsing}, then pathological
behaviors of the same kind described before may appear.
More precisely, we find a steady solution which is not continuously differentiable.

We consider the system in conservation form
\begin{equation*}
     u_t + f(u)_x = \Big( B(u) u_x \Big)_x
\end{equation*}
with
\begin{equation*}
       f(u_1, \, u_2)=
       \left(
       \begin{array}{ll}
             u_2 \\
             u_1 \\
       \end{array}
       \right)
       \qquad
       B(u_1, \, u_2)=
       \left(
       \begin{array}{ll}
             \gamma u_1^{\gamma - 1 } & 0  \\
             0 & 1 \\
       \end{array}
       \right)
\end{equation*}
for some $\gamma >3$. In this case, the rank of $B(u)$ drops from
$2$ to $1$ when $u_1$ reaches $0$: to find out pathological
behaviors it seems therefore natural to study the solution in a
neighborhood of $u_1 =0$. More precisely, the attention is
focused on steady solutions
\begin{equation}
\label{eq_hypotheses_rank_boundary_layers}
       f(u)_x = \Big( B(u) u_x \Big)_x
\end{equation}
such that
\begin{equation}
\label{eq_hypotheses_rank_limit_cond}
       \lim_{x \to + \infty} u_1 (x) =0 \qquad u_1 (x) \ge 0
       \qquad
       \lim_{x \to + \infty} u_2 (x)=0
\end{equation}
and it will be shown that the graph of the first component $u_1$
has the shape illustrated in Figure \ref{fig_example_rank}. Hence
$u_1$ is not continuously differentiable and presents a behaviour
like the one typical of the solutions of the porous-media
equation: it is different from zero on a interval and then
vanishes identically.

In this case equation \eqref{eq_hypotheses_rank_boundary_layers}
writes
\begin{equation}
\label{eq_hypotheses_boundary_layers_limit}
      \left\{
      \begin{array}{ll}
             u_{2x}= \Big( u_1^{\gamma} \Big)_{xx} \\
             u_1 = u_{2x} \\
      \end{array}
      \right.
\end{equation}
and hence after some computations one obtains
\begin{equation}
\label{eq_hypotheses_rank_solution}
      \left\{
      \begin{array}{lll}
             u_{1x}= - \displaystyle{ \sqrt{
                      \frac{2}{\gamma(\gamma +1)} u_1^{(3 - \gamma) }}} \\
             \\
             u_2 = \displaystyle{\int_{+ \infty }^{x} u_1 (y) dy}. \\
      \end{array}
      \right.
\end{equation}
The equation satisfied by $u_1$ admits more than one solution in a
neighborhood of $u_1 = 0$.

To introduce a selection principle we consider the matrix 
\begin{equation*}
      B^{\nu}(u^{\nu})=
      \left(
       \begin{array}{cc}
             \nu + \gamma (\un_1)^{\gamma - 1 } & 0  \\
             0 & 1 \\
       \end{array}
       \right),
\end{equation*}
and the equation
\begin{equation}
\label{e:ex3:sys}
       \un_{xx} = \Big( B^{\nu}(\un) \Big)^{-1}
       f(\un)_x.
\end{equation}
In the previous expression, $\nu$ is a positive parameter and we study the limit $\nu \to 0^+$.  We keep the limit conditions 
\eqref{eq_hypotheses_rank_limit_cond} fixed. In other words, this time we do not perturb the limit condition (as in the previous examples), but the equation itself. Note that now $B^{\nu}$ is invertible. The solution of \eqref{e:ex3:sys} with limit conditions \eqref{eq_hypotheses_rank_limit_cond} satisfies 
\begin{equation}
\label{eq_hypotheses_rank_orbit}
       ( \un_2 )^2 = \nu (\un_1)^2 + \frac{2 \gamma }{\gamma +1}
       (\un_1)^{\gamma + 1}.
\end{equation}
Indeed, from
\begin{equation*}
      \left\{
      \begin{array}{ll}
            \un_{2 x} = \Big( \nu \un_{1x} +
            \gamma (\un_1 )^{\gamma -1  } \un_{1x}
            \Big)_{x} \\
            \un_{1 x} = \un_{2 xx}
      \end{array}
      \right.
\end{equation*}
one obtains integrating
\begin{equation*}
      \left\{
      \begin{array}{ll}
            \un_{2 } =  \nu \un_{1x} +
            \gamma (\un_1 )^{\gamma -1  } \un_{1x}  \\
            \un_{1 } = \un_{2 x}.
      \end{array}
      \right.
\end{equation*}
Multiplying the first line by $u_1 = u_{2x}$ and then integrating
again one obtains \eqref{eq_hypotheses_rank_orbit}.  From 
\eqref{eq_hypotheses_rank_orbit} we get
$$    
    u^{\nu}_{2x} = - \frac{         ( 2 \nu + \gamma (u^{\nu}_1)^{\gamma}) u^{\nu}_{1x}       }{ \displaystyle{\sqrt{  
                                                                                                                                      \nu u^{\nu}_1 + \frac{  2 \gamma }{   \gamma + 1 }  
                                                                                                                                                                                                          (u^{\nu}_1)^{\gamma + 1 }   }}}.
$$
Taking $u^{\nu}_{2x} = u^{\nu}_1$ one eventually gets 
$$
    u^{\nu}_{1x}= -   \frac{  \sqrt{  \nu u^{\nu}_1 + \frac{2 \gamma }{   \gamma + 1 } ( u^{\nu}_1)^{\gamma + 1 }       } }{   2 \nu + \gamma (u^{\nu}_1)^{\gamma}   }.
$$
Fix $u_{10} >0$ and consider the Cauchy problem
$$
   \left\{
   \begin{array}{lll}
                            u^{\nu}_{1x}= -   \displaystyle{ \frac{  \sqrt{  \nu u^{\nu}_1 + \frac{2 \gamma }{   \gamma + 1 } ( u^{\nu}_1)^{\gamma + 1 }       } }{   2 \nu + \gamma (u^{\nu}_1)^{\gamma} } }\\
                            u^{\nu}_1 (0) = u^{\nu}_{10}  \\
   \end{array}
   \right.
$$
One can then exploit the same considerations  
as in the previous examples: in particular, the fact that for every $x$ $u^{\nu}_1(x)$ is monotone decreasing with respect to $\nu$ follows from 
$$
     \frac{\partial}{\partial \nu } \Bigg( -   \displaystyle{ \frac{  \sqrt{  \nu u^{\nu}_1 + \frac{2 \gamma }{   \gamma + 1 } ( u^{\nu}_1)^{\gamma + 1 }       } }{   2 \nu + \gamma (u^{\nu}_1)^{\gamma} } } \Bigg) 
     \Bigg|_{\nu=0, \; u^{\nu}_1 > 0} < 0
$$
if $\gamma < 3$. One can the prove (proceeding as in the previous examples) that 
when $\nu \to 0^+$,  $u_1^{\nu}$ converges pointwise to a function $u_1$ satisfying 
\begin{equation*}
      u_1 (x) =
       \left\{
      \begin{array}{lll}
             \displaystyle{ \sqrt[\gamma -1]{
                  \frac{(\gamma -1)^2}{2 \gamma(\gamma +1)} (x- x_0)^2}}
                  \qquad x \leq x_0 \\
             \\
             0 \qquad \qquad  \qquad  \qquad \qquad \qquad x > x_0, \\
      \end{array}
      \right.
\end{equation*}
where
\begin{equation*}
      x_0 = \sqrt{\frac{2 \gamma (\gamma +1 )}{(\gamma - 1)^2 }
       u_{10}^{\gamma -1}}.
\end{equation*}
This function has the graph illustrated in Figure \ref{fig_example_rank} and it is not continuously differentiable. 

\begin{figure}
\begin{center}
\caption{the graph of the function $u_1 (x)$ in Example
\ref{ex_rank}} \psfrag{u}{$u_1 (x)$} \psfrag{v}{$u_1(0)$}
\psfrag{o}{$x_0$} \psfrag{x}{$x$} \label{fig_example_rank}
\includegraphics[scale=0.6]{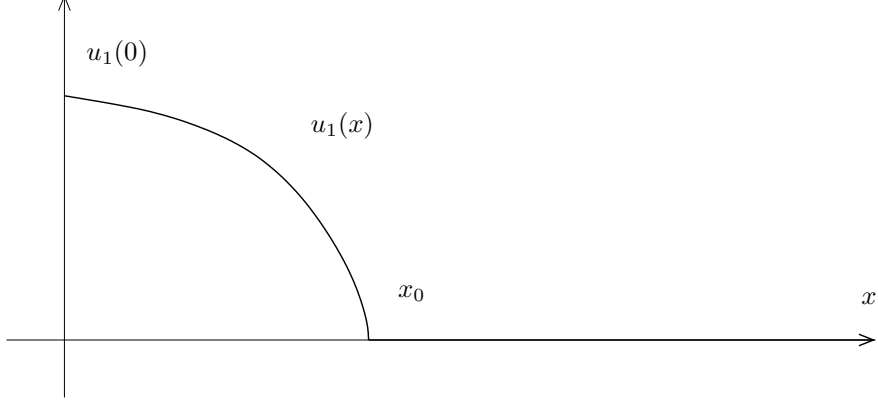}
\end{center}
\end{figure}

\end{ex}
\subsection{General hypotheses}
\label{subsec_hypotheses_introduction} This section introduces the
general hypotheses required in both cases, i.e. when the viscosity matrix $\tilde B$ in \eqref{eq_introduction_parabolic_approx}
is invertible and when it is singular. In the statements of the
hypotheses we actually make  reference to the formulation \eqref{eq_hypotheses_system} and 
\eqref{e:hyp:datum:sys}. Thus, in particular, we consider the same value $\bar u_0 $ as in 
 \eqref{eq_hypotheses_system} and 
\eqref{e:hyp:datum:sys}. 

First of all, we assume strict hyperbolicity:
\begin{hyp}
\label{hyp_hyperbolic_I}
            There exists $\delta >0$ such that, if $u$ belongs to a neighbourhood of $\bar{u}_0$
            of size $\delta$ then
      all the eigenvalues of the matrix $E^{-1}(u)A(u, \, 0)$ are real.
      Moreover, there exists a constant $c>0$ such that
      \begin{equation*}
            \inf_{u}
            \big\{ \big| \lambda_i (u, \, 0) -
            \lambda_j(u, \, 0) \big| \big\} \ge c > 0 \quad
            \forall \, i \neq j.
      \end{equation*}
\end{hyp}

We also introduce an hypothesis of convergence: 
\begin{hyp}
\label{hyp_convergence}
      Let
      \begin{equation*}
      \left\{
      \begin{array}{lll}
             E(\ue ) \ue_t + A(\ue, \, \ee \ue_x ) \ue_x = \ee B(\ue )
             \ue_{xx} \\
             \B(\ue(t, \, 0)) = g(t)  \\
              \ue (0, \, x) =
             \bar{u}_0(x)
      \end{array}
      \right.
      \end{equation*}
      be a parabolic initial boundary value problem such that 
      $g$, $\bar{u}_0 \in L^1_{loc}$ and
      \begin{equation}
      \label{eq_hypotheses_small_bv}
            \big| g(0) - \B \big( \bar{u}_0(0) \big) \big| , \;
            \mathrm{Tot Var} \big\{
            g \big\}, \; \mathrm{Tot Var} \big\{
            \bar{u}_0 \big\} \leq \delta
      \end{equation}
      For a suitable constant $\delta << 1$. 
      Then we assume that there exists a time $T_0$, which depends only on the
      bound $\delta$ and on the matrices $E$, $A$ and $B$, such
      that
      \begin{equation*}
            \mathrm{Tot Var} \big\{
            \ue (t, \, \cdot) \big\} \leq C \delta \quad \forall \, \ee,
            \; t \in [ 0, \, T_0 ].
      \end{equation*}
      We also assume
      directly a result of convergence and of uniqueness of the
      limit:
     \begin{equation*}
            \forall \, t \in [ 0, \, T_0 ] \quad  \ue (t, \, \cdot) \to u(t) \; L^1_{loc}
            \;  \mathrm{when} \; \ee \to
            0^+ \qquad \mathrm{Tot Var} \big\{
            u(t) \big\} \leq C \delta.
     \end{equation*}
\end{hyp}

We point out that the uniqueness of the limit is actually implied by
the next hypotheses and the uniqueness results for the Standard
Riemann Semigroup with boundary: it is made reference to
\cite{AmaCol} for the extension of the definition of SRS to
initial boundary value problems, while an application of this
notion to prove the uniqueness of the limit of the vanishing
viscosity solutions was introduced in \cite{BiaBrevv} in the case
of the Cauchy problem.

\begin{figure}
\caption{the finite propagation speed in the hyperbolic limit}
\begin{multicols}{2}
\raggedright \psfrag{t}{$t_0$} \psfrag{a}{$\beta(t-t_0)$}
\psfrag{b}{$b - \beta t $} \psfrag{c}{$b$} \psfrag{s}{$t$}
\psfrag{u}{$t$} \psfrag{f}{$a$} \psfrag{e}{$a + \beta t $}
\label{fig_finite_prop}
\includegraphics[scale=0.3]{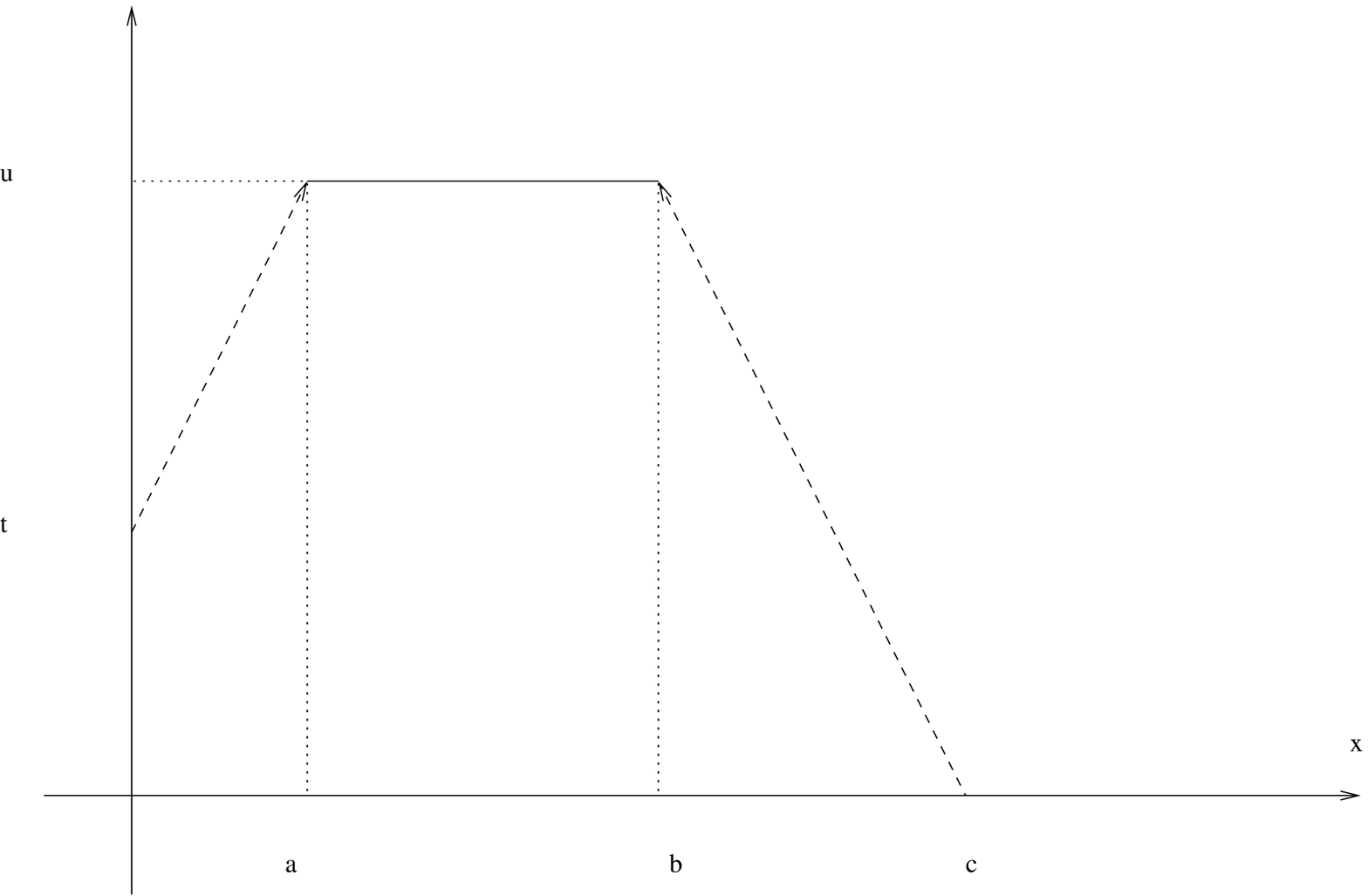}
\includegraphics[scale=0.3]{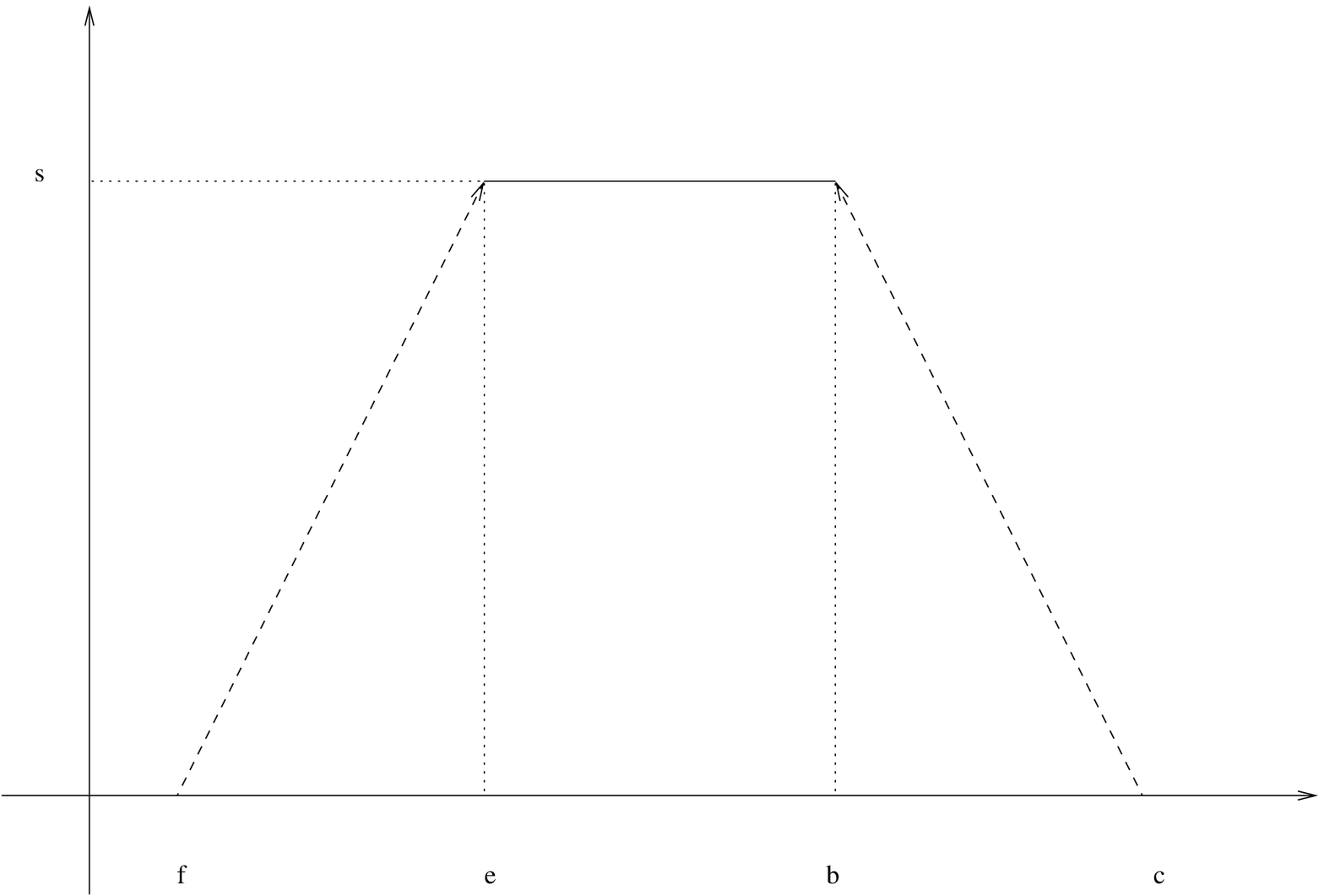}
\end{multicols}
\end{figure}

It is also assumed:
\begin{hyp}
\label{hyp_small_bv}
       It holds
       \begin{equation*}
              |\B(\bar{u}_0) - \bar g| < \delta
       \end{equation*}
       for the same constant $\delta < < 1$ that appears in
       Hypothesis \ref{hyp_convergence}.
\end{hyp}
An hypothesis of stability with respect to $L^1$ perturbations in
the initial and boundary data is also introduced:
\begin{hyp}
\label{hyp_stability}
       There exists a constant $L >0$ such that
       the following holds.

       Two families of parabolic initial boundary value
       problems are fixed:
       \begin{equation}
       \label{eq_hypotheses_2_systems}
             \left\{
       \begin{array}{lllll}
              u^{1 \, \ee}_t + A(u^{1 \, \ee}, \, \ee u^{1 \, \ee}_x)
              u^{1 \, \ee}_x =
       \ee B(u^{1 \, \ee}) u^{1 \, \ee}_{xx} \\
       \B(u^{1 \, \ee} (t, \, 0)) = \bar g^1(t). \\
       u^{1 \, \ee} (0, \, x) = \bar{u}^1_0 (x)  \\
       \end{array}
       \right. \qquad
             \left\{
       \begin{array}{lllll}
              u^{2 \, \ee}_t + A(u^{2 \, \ee}, \, \ee u^{2 \, \ee}_x)
              u^{2 \, \ee}_x =
             \ee B(u^{2 \, \ee}) u^{2 \, \ee}_{xx} \\
       \B(u^{2 \, \ee} (t, \, 0)) = \bar g^2(t). \\
       u^{2 \, \ee} (0, \, x) = \bar{u}^2_0 (x)  \\
       \end{array}
       \right.
       \end{equation}
       with $\bar{u}^1_0, \; \bar g^1$ and $\bar{u}^2_0, \; \bar g^2$ in $L^1_{loc}$
       and satisfying the assumption
       \eqref{eq_hypotheses_small_bv}.

       Then for all $t \in [0, \, T_0]$ and for all $\ee >0$ it
       holds
       \begin{equation*}
             \|u^{1 \, \ee}(t)  - {u}^{2 \, \ee}(t) \|_{L^1} \leq
             L \bigg( \|\bar{u}^1_0 - \bar{u}^2_0\|_{L^1} +
             \|\bar g^1 - \bar g^2 \|_{L^1} \bigg).
       \end{equation*}
\end{hyp}
From the stability of the approximating solutions and from the
$L^1_{loc}$ convergence one can deduce the stability of the limit.
More precisely, let $u^1$ and $u^2$ the limits of the two
approximations defined above, then
\begin{equation*}
      \|u^{1}(t)  - {u}^{2}(t) \|_{L^1} \leq
             L \bigg( \|\bar{u}^1_0 - \bar{u}^2_0\|_{L^1} +
             \|\bar g^1 - \bar g^2 \|_{L^1} \bigg).
\end{equation*}
Finally, it is assumed that in the hyperbolic limit there is a
{\sl finite propagation speed} of the disturbances:
\begin{hyp}
\label{hyp_finite_pro}
       There exist constants $\beta, c > 0$ such that the following
       holds.

       Let $\bar{u}^1_0, \; \bar g^1$ and $\bar{u}^2_0, \;
       \bar g^2$ in $L^1_{loc}$ and bounded, let $u^{\ee \, 1}$ and $u^{\ee \, 2}$
       be the solutions of
       \eqref{eq_hypotheses_2_systems} and $u^1$ and $u^2$ the
       corresponding limits. If
       \begin{equation*}
              \bar g^1 (t) = \bar g^2 (t) \quad \forall \, t
              \leq t_0 \qquad
               \bar{u}_0^1 (x) = \bar{u}_0^2 (x) \quad \forall \, x \leq b
       \end{equation*}
       then
       \begin{equation*}
              |u^{\ee 1}(x, \, t) - u^{\ee 2} (x, \, t)| \leq \mathcal{O}
              \Big( e^{- c \min\big\{|x-\max\{0, \,
              \beta( t- t_0)\}|,|x-b + \beta t|\big\}/\epsilon} \Big) \quad
              \forall x \in [\max\{0, \, \beta( t- t_0)\}, \, b - \beta t
              ].
       \end{equation*}
       Analogously, if
       \begin{equation*}
               \bar{u}_0^1 (x) = \bar{u}_0^2 (x) \quad \forall \, x \in [a, \,
               b]
       \end{equation*}
       then
       \begin{equation*}
             |u^{\ee 1}(x, \, t) - u^{\ee 2} (x, \, t)| \leq \mathcal{O}
             \Big( e^{- c \min\big\{|x-a - \beta t|,|x-b + \beta t|\big\}/\epsilon} \Big)
             \quad \forall \, x \in [a + \beta t, \, b - \beta
             t].
       \end{equation*}
\end{hyp}
\begin{rem}
Hypothesis \ref{hyp_finite_pro} may appear a bit technical. It implies the finite propagation speed of disturbances in the hyperbolic limit: an heuristic representation of this phenomenon is illustrated in  Figure \ref{fig_finite_prop}. Loosely speaking, the reason why we need finite propagation speed is the following. We have to be sure that the limit  of \eqref{e:hyp:datum:sys} in the case of a generic couple of data  $(\bar{u}_0, \, \bar{g})$ can be obtained gluing  together the limit one obtains in the case of {\sl cooked up} data, namely data connected by travelling wave profiles.  
\end{rem}
\section{Characterization of the hyperbolic limit in the case of an invertible viscosity matrix}
\label{sec_B_invertible} 
The aim of this section is to provide a characterization of the limit of the parabolic approximation 
 \eqref{eq_introduction_parabolic_approx} when the viscosity matrix $\tilde{B}$ is invertible.
The precise
hypotheses that are assumed are listed in Section 
\ref{subsec_hypotheses_invertible} and 
\ref{subsec_hypotheses_introduction}: in particular, these
hypotheses guarantee that it is sufficient to study system
\begin{equation}
\label{eq_Binv_the_system}
      \left\{
      \begin{array}{ll}
             E(\ue) \ue_t + A \big( \ue, \, \ee \ue_x \big)  \ue_x=
             \ee B(\ue ) \ue_{xx} \qquad \ue \in \mathbb{R}^N\\
             \ue (t, \, 0) \equiv \bar{u}_b \\
             \ue (0, \, x) \equiv \bar{u}_0, \\
      \end{array}
      \right.
\end{equation}
where the matrices $E$, $A$ and $B$ satisfy Hypotheses
\ref{hyp_invertible} and \ref{hyp_hyperbolic_I}.

The exposition is organized as follows. In Section
\ref{subsec_Binv_preliminary} we discuss some preliminary results. More precisely, in Section \ref{sssec:tran} 
we recall some transversality results discussed in  \cite{BenSerreZum:Evans}. In Section \ref{subsubsec:Binv_monotone} we 
recall the definition of monotone envelope of a function and we introduce some related results.
In Section
\ref{subsec_Binv_riemannsolver} we review some results in 
\cite{Bia:riemann}. Namely, we give a characterization of the limit of parabolic approximation 
\begin{equation*}
\left\{
\begin{array}{ll}
       E(\ue) \ue_t + A \big( \ue, \, \ee \ue_x \big) \ue_x =
             \ee B(\ue ) \ue_{xx} \\
      \ue (0, \, x) = 
      \left\{
      \begin{array}{ll}
                  u^-  &  x \leq 0 \\
                  \bar{u}_0 & x >0. \\      
      \end{array}
      \right.
\end{array}
\right.             
\end{equation*}

In Section
\ref{subsub_Bsing_nonchar} it is given a 
characterization of the limit of the parabolic approximation
\eqref{eq_Binv_the_system} in the case of a non characteristic
boundary, i.e. when none of the eigenvalues of $E^{-1} A$ can attain the value $0$. 
The case of a characteristic boundary occurs when one of the eigenvalues of $E^{-1} A$ can attain the value $0$ and it is definitely more complicated than the previous one.
The characterization of the limit \eqref{eq_Binv_the_system} in the case of a characteristic boundary is discussed in Section
\ref{subsubsec_Binv_char}. 
\subsection{Preliminary results}
\label{subsec_Binv_preliminary} 
\subsubsection{Transversality results}
\label{sssec:tran}
The following lemma is discussed in \cite{Gisclon:etudes}. However, for completeness we will repeat the proof.
\begin{lem}
\label{lem_Binv_dimension}
      Assume that Hypotheses \ref{hyp_invertible} and \ref{hyp_hyperbolic_I}
      hold.
      Then for any given $u$ and $u_x$
      \begin{enumerate}
      \item $B^{-1}(u) A(u, \, u_x) \vec{\xi}=0 \; \Longleftrightarrow \;
      A(u, \, u_x) \vec{\xi}=0 \; \Longleftrightarrow \;
      E^{-1}(u) A(u, \, u_x) \vec{\xi}=0 \qquad \vec{\xi} \in \mathbb{R}^N$
      \item the number of
      eigenvalues of the
      matrix $B^{-1}(u) A(u, \, 0)$ with negative (respectively positive)
      real part is equal to the number of
      eigenvalues of
      $E^{-1}(u) A(u, \, 0)$ with negative (respectively
      positive) real part.
      \end{enumerate}
\end{lem}
\begin{proof}
The first point is actually an immediate observation.

To simplify the notation, in the proof of the second point, we
will write $A$, $B$ and $E$ instead of $A(u, \, 0)$, $B(u)$ and
$E(u)$ respectively. One can define the continuous path
\begin{equation*}
\begin{split}
&     F: [0, \, 1] \to \mathbb{M}^{N \times N} \\
&     \qquad \; \, \quad s \mapsto (1- s) B + s E,
\end{split}
\end{equation*}
which satisfies the following condition: for every $s \in [0, \,
1]$, $F(s)$ is positively definite and hence invertible. Indeed,
\begin{equation}
\label{eq_Binv_def+}
      \forall \, \vec{\xi} \in \erreN \quad
      \langle F(s) \vec{\xi}, \, \vec{\xi} \rangle
      \ge \min \{ c_B(u), \, c_E(u) \} |\vec{\xi} |^2.
\end{equation}
Moreover, by classical results (see for example \cite{Knopp:book}) from
the continuity of the path $F(s)^{-1}A$ it follows the continuity
of the corresponding paths of eigenvalues $\lambda_1 (s) \dots
\lambda_N (s)$.

Let $k-1$ be the number of eigenvalues of $F(0)^{-1}A= B^{-1} A$
with negative real part:
\begin{equation*}
      Re (\lambda_1(0) ) < \dots Re (\lambda_{k-1}(0)) < 0 =
      Re (\lambda_k(0)) < Re (\lambda_{k+1}(0)) < \dots Re (\lambda_N (0)).
\end{equation*}
Because of the continuity of $\lambda_1 (s) \dots \lambda_i (s)$,
to prove that the number of negative eigenvalues of $F(1)^{-1}A=
E^{-1}A$ is $(k-1)$ it is sufficient to prove that for any $i=1
\dots (k-1)$, $\lambda_i(s)$ cannot cross the imaginary axis.
Moreover, the case $\lambda_i(s)=0$ is excluded because otherwise
the corresponding eigenvector $\vec{\xi}$ should satisfy
$F(s)^{-1}A \vec{\xi}(s)\equiv 0$ for any $s$ and hence $\lambda_i
(s) \equiv 0$. It remains therefore to consider the possibility
that $\lambda_i(s)$ is purely imaginary.

By contradiction, assume that there exist $\bar{s} \in [0, \, 1]$,
$\lambda \in \mathbb{R} \setminus \big\{ 0 \big\}$,  $r \in
\mathbb{C}^N$ such that
\begin{equation}
\label{eq_Binv_pi_eigenvalue}
       F(\bar{s})^{-1} A r = i \lambda r:
\end{equation}
just to fix the ideas, it will be supposed that $\lambda > 0$.
Moreover, let
      \begin{equation*}
             r= r_1 + i r_2, \quad r_1, \; r_2 \in \mathbb{R}^N.
      \end{equation*}
      Then from \eqref{eq_Binv_pi_eigenvalue} one gets
      \begin{equation*}
      \begin{split}
      &       A r_1 = - \lambda F(\bar{s}) r_2 \\
      &       A r_2 =  \lambda F(\bar{s}) r_1 \\
      \end{split}
      \end{equation*}
      and hence from \eqref{eq_Binv_def+} and from the symmetry of
      $A$ it follows
      \begin{equation*}
      \begin{split}
      &       \langle A r_1, \ r_2  \rangle \leq - \lambda
              \min \{c_{B}(v), \, c_E(v) \} |r_2|^2 \\
      &       \langle A r_1, \ r_2  \rangle =
              \langle A r_2, \ r_1  \rangle \ge  \lambda
              \min \{c_{B}(v), \, c_E(v) \} |r_1|^2 , \\
      \end{split}
      \end{equation*}
      which is a contradiction since $\lambda > 0$. This
      ends the proof of Lemma \ref{lem_Binv_dimension}.
\end{proof}
In the following, given a matrix $C$ we denote by $V^s(C)$ the direct sum of all eigenspaces associated with eigenvectors with strictly negative real part,
by $V^c (C)$  the direct sum of the eigenspace associated with the eigenvector with zero real part, and $V^u (C)$ the  the direct sum of all eigenspaces associated with eigenvectors with strictly positive real part. Lemma \ref{lem_Binv_dimension} ensures that, if
Hypothesis \ref{hyp_invertible} holds, then for all $u, \, u_x$
\begin{equation*}
\begin{split}
&      
       V^c \big( E^{-1}(u) A(u, \, 0) \big) =
       V^c \big( A(u, \, 0) \big) =
       V^c \big( B^{-1}(u) A(u, \, 0) \big) \\
&      dim V^{u}\big( E^{-1}(u) A(u, \, 0) \big)=
       dim V^{u}\big( B^{-1}(u) A(u, \, 0) \big)
       \\
&      dim V^{s}\big( E^{-1}(u) A(u, \, u_x) \big)=
       dim V^{s}\big( B^{-1}(u) A(u, \, u_x) \big) \\
\end{split}
\end{equation*}
\begin{lem}
\label{lem_Binv_zero intersection}
      Assume that Hypotheses
      \ref{hyp_invertible} and \ref{hyp_hyperbolic_I} hold, then
      \begin{equation*}
             V^{u}\Big( E^{-1}(u) A(v, \, 0) \Big) \cap
             V^{s}\Big( B^{-1}(u) A(v, \, 0) \Big) = \{ 0 \}.
      \end{equation*}
\end{lem}
Thanks to the previous lemma, which is a direct application of
Lemma 7.1 in \cite{BenSerreZum:Evans}, one can conclude that
actually
\begin{equation}
\label{eq_Binv_direct_sum}
      V^{u} \Big( E^{-1}(u) A(u, \, 0) \Big) \oplus
      V^c \Big( E^{-1}(u) A(u, \, 0) \Big) \oplus
      V^{s}\Big( B^{-1}(u) A(u, \, 0) \Big) = \erreN.
\end{equation}

\subsubsection{Some results about the monotone envelope of a function}
\label{subsubsec:Binv_monotone}
The aim of this section is to collect some results that will be needed in Section \ref{subsubsec_Binv_char}. Proposition \ref{p:reg_conc} is very similar to results that were introduced, in a much more general form, for example in \cite{BallKiKr} and \cite{GreRav}. However, since the situation discussed here is slightly different, for completeness we give a proof. Also, results analogous to Propositions \ref{p:es_mon} and \ref{p:concint} are discussed in \cite{AnMar:front_track}, but a  proof is given here for completeness. 

In the following, $\mathrm{conc}_{[0, \, s] } f $ will denote the concave envelope of the function $f$ on the interval $[0, \, s]$, namely
\begin{equation}
\label{e:conc}
         \Big(          \mathrm{conc}_{[0, \, s] } f    \Big) (\tau) :  = \inf{\Big\{  h( \tau): \; h (t) \ge f(t) \quad \forall \, t \in [0, \, s], \; h \; \mathrm{ is \; concave  }      \Big\}}.
\end{equation}
We will consider only the case of a function $f \in \mathcal{C}^{1, \, 1}_k ([0, \, s])$, i.e.  $f \in \mathcal{C}^{1, \, 1}$ and $f'$ is Lipschitz continuous with Lipschitz constant smaller or equal to $k$. 

The symbol $\mathrm{mon}_{[0, \, s] } f $ will denote the monotone envelope of the function $f$ on the interval $[0, \, s]$, i.e.
\begin{equation}
\label{e:mono}
         \Big(          \mathrm{mon}_{[0, \, s] } f \Big) (\tau) :  = \sup{\Big\{  h( \tau): \; h (t) \ge f(t) \quad \forall \, t \in [0, \, s], \; h \; \mathrm{ is \; concave \; and \; nondecreasing  }      \Big\}}.
\end{equation}.

\begin{lem}
\label{l:reg_concave}
          Let $f \in \mathcal{C}^{1, \, 1}_k ([0, \, s])$, then 
          \begin{enumerate}
          \item $   \mathrm{conc}_{[0, \, s] } f  $ is continuous and satisfies 
          $\|      \mathrm{conc}_{[0, \, s] } f   \|_{\mathcal{C}^0}      \leq     \| f   \|_{\mathcal{C}^0}  $.
          \item If $\tau \in ]0, \, s[$ and
          $$
              \Big(          \mathrm{conc}_{[0, \, s] } f \Big) (\tau) = f(\tau),
           $$
           then  $\mathrm{conc}_{[0, \, s] } f $ is differentiable at $\tau$ and 
           $$
                  \Big(          \mathrm{conc}_{[0, \, s] } f \Big) ' (\tau) = f' (\tau) .
           $$
           \item If $    \tau \in ]0, \, s[ $ and 
             $$
                \Big(          \mathrm{conc}_{[0, \, s] } f \Big) (\tau) > f(\tau),
           $$
           then there exists $a, \, b$ such that $0 \leq a < \tau < b \leq s $ such that 
           $$
             \forall \;   t \in [a, \, b], \quad \Big(     \mathrm{conc}_{[0, \, s] } f     \Big)(t) = \frac{(t- a ) f(b) + (b - t) f (a) }{b  - a}.
           $$
           \item $  \mathrm{conc}_{[0, \, s] } f (0) = f(0)$, $ \mathrm{conc}_{[0, \, s] }f (s) = f (s)$.
    \end{enumerate}
\end{lem}
\begin{proof}
By definition, 
\begin{equation}
\label{e:Binv:concint:bi}
    f(\tau) \leq  \mathrm{conc}_{[0, \, s] } f (\tau) \quad \forall \, \tau \in [0, \, s]
\end{equation}
and hence $ - \| f  \|_{\mathcal{C}^0}  \leq  \mathrm{conc}_{[0, \, s] } f (\tau)$. Moreover, the constant function $ \| f \|_{   \mathcal{C}^0  } $ is concave and greater or equal to $f$ and hence  
$\mathrm{conc}_{[0, \, s] } f (\tau) \leq \| f \|_{   \mathcal{C}^0  }$. Combining the two inequalities,  one gets 
$\| \mathrm{conc}_{[0, \, s] } f \|_{L^{\infty}}  \leq \|   f  \|_{\mathcal{C}^0}$. Moreover, $\mathrm{conc}_{[0, \, s] } f $ is a concave function and hence it is continuous in all the 
inner points of  its domain, thus in $]0, \, s[$. Later on we prove that it is also continuous at $t=0$ and $t=s$. 

To prove the second point, assume that $ \tau \in ]0, \, s[ $ and that  $ \mathrm{conc}_{[0, \, s] } f (\tau) = f(\tau)$. If $h >0$, then thanks to \eqref{e:Binv:concint:bi} 
$$
     \frac{ f(\tau + h ) - f (\tau)     }{    h}   \leq \frac{ \mathrm{conc}_{[0, \, s] } f (\tau   +   h) - \mathrm{conc}_{[0, \, s] } f (\tau) }{h}
$$
and hence 
$$
    f ' ( \tau ) =  \lim_{h \to 0^+}     \frac{ f(\tau + h ) - f (\tau)     }{    h}    \leq  \lim_{h \to 0^+}    \frac{ \mathrm{conc}_{[0, \, s] } f (\tau   +   h) - \mathrm{conc}_{[0, \, s] } f (\tau) }{h}
$$
Note that $ f'(\tau)$ exists because $f \in \mathcal{C}^{1 \, 1}_k$, while
$$
    \lim_{h \to 0^+}    \frac{ \mathrm{conc}_{[0, \, s] } f (\tau   +   h) - \mathrm{conc}_{[0, \, s] } f (\tau) }{h}
$$
exists because the difference quotient is a monotone non increasing function. By the same reason,  
$$
     \lim_{h \to 0^-}    \frac{ \mathrm{conc}_{[0, \, s] } f (\tau   +   h) - \mathrm{conc}_{[0, \, s] } f (\tau) }{h} \leq  \lim_{h \to 0^-}    \frac{ \mathrm{conc}_{[0, \, s] } f (\tau   +   h) - \mathrm{conc}_{[0, \, s] } f (\tau) }{h}.
$$
Moreover, taking $ h <0 $ one can proceed as before and get
$$
    f ' ( \tau ) =  \lim_{h \to 0^-}     \frac{ f(\tau + h ) - f (\tau)     }{    h}    \ge  \lim_{h \to 0^-}    \frac{ \mathrm{conc}_{[0, \, s] } f (\tau   +   h) - \mathrm{conc}_{[0, \, s] } f (\tau) }{h}
$$
and hence putting all the previous considerations together one gets that $ \big(   \mathrm{conc}_{[0, \, s] } f  \big)'  (   \tau)  $ exists and it is equal to $f ' (  \tau  ) $.  

To prove the third point, take $\tau$ such that $ \mathrm{conc}_{[0, \, s] } f (\tau) >  f(\tau)$ and take $a$ and $b$ as follows:
$$
    a : = \inf \{ t < \tau: \;    f (t)   <  \mathrm{conc}_{[0, \, s] } f ( t ) \} 
$$
and 
$$
    b :=  \sup \{ t >  \tau: \;    f (t)   <  \mathrm{conc}_{[0, \, s] } f ( t ) \}.  <  s
$$
To simplify the notations in the following we will also assume 
$$
    \mathrm{conc}_{[0, \, s] } f ( a) =  \mathrm{conc}_{[0, \, s] } f ( b ) = 0,
$$
the general case being analogous by subtracting a linear function. 

First of all, note that $ \mathrm{conc}_{[0, \, s] } f ( t )  \ge 0$ for every $t \in [a, \, b]$ because $ \mathrm{conc}_{[0, \, s] } f $ is a concave function and hence it is greater or 
equal to the segment joining two of its values.  We want to prove that also $ \mathrm{conc}_{[0, \, s] } f ( t ) \leq 0$. More precisely, we will consider the function
\begin{equation*}
h ( t ) : =
\left\{
\begin{array}{lll}
             \mathrm{conc}_{[0, \, s] } f ( t)  &  t \leq a \\
             0 &  a \leq t \leq b \\
               \mathrm{conc}_{[0, \, s] } f ( t) & t \ge b. \\                
\end{array}
\right.
\end{equation*}
and we will prove that $\mathrm{conc}_{[0, \, s] } f ( t )  \leq h (t)$.
Since the function $h$ is concave, it is enough to show that $h (t) \ge f ( t ) $ for every $t \in [a, \, b]$. By contradiction, 
assume that there exists $\tau \in [a, \, b]$ such that $f (\tau) >0$. Assume that $\tau$ is exactly a point at which $f$ assumes it maximum on $[a, \, b]$. If $\tau=a$ or $\tau=b$ then there is nothing to prove because $f(t) \leq 0$ for every $t$ in $[a, \, b]$.
Moreover, let $\tau_m$ denote a point at which the maximum of $     \mathrm{conc}_{[0, \, s] } f $ on $[a, \, b] $ is assumed: since $     \mathrm{conc}_{[0, \, s] } f $ is concave, then it is non decreasing for $t \leq \tau_m$ and non increasing for $t \ge \tau_m$. Moreover, by \eqref{e:Binv:concint:bi}  $   f (\tau) \leq  \mathrm{conc}_{[0, \, s] } f (\tau_m)$. Consider the function $z$ defined as follows:
\begin{equation*}
z ( t ) : =
\left\{
\begin{array}{lll}
             \mathrm{conc}_{[0, \, s] } f ( t)  &  t \leq t_a \\
             f(\tau) &  a \leq t \leq b \\
               \mathrm{conc}_{[0, \, s] } f ( t) & t \ge t_b, \\                
\end{array}
\right.
\end{equation*}
where 
$$
    t_a = \min \{ t \in [a, \, b]: \,   \mathrm{conc}_{[0, \, s] } f ( t) = f (\tau)   \}
$$
and 
$$
    t_b = \max \{ t \in [a, \, b]: \,   \mathrm{conc}_{[0, \, s] } f ( t) = f (\tau)   \}.
$$
Putting all the previous considerations together it turns out that $ t_a \leq \tau_m \leq t_b$ and hence $z$ is a concave function which is greater or equal to $f$. Thus, 
$ \mathrm{conc}_{[0, \, s] } f \leq z$ and hence in particular 
$$
      \mathrm{conc}_{[0, \, s] } f  (\tau) \leq z (\tau) = f (\tau) .
$$
This contradicts the assumption $ \mathrm{conc}_{[0, \, s] } f (t) > f (t)$ on $]a, \, b[$.

Note that if $a >0$, then automatically $f (a) =   \mathrm{conc}_{[0, \, s] } f  (  a  ) $. Now we want to prove that, even if $a = 0$,
$f (0) =   \mathrm{conc}_{[0, \, s] } f  (  0  ) $. By contradiction,  assume that $f (0) <  \mathrm{conc}_{[0, \, s] } f  (  0 ) =0$, then 
consider the function 
\begin{equation*}
w ( t ) : =
\left\{
\begin{array}{lll}
            f (0)  + \|   f ' \|_{\mathcal{C}^0} t  & t \leq \bar{t} \\
             0 &  a \leq t \leq b \\
               \mathrm{conc}_{[0, \, s] } f ( t) & t \ge \bar{t}, \\                
\end{array}
\right.
\end{equation*}
where $\bar{t}$ is the point at which $   f (0)  + \|   f ' \|_{\mathcal{C}^0} \bar{t} = 0$. Then $w$ is a concave function greater or equal to $f$ and hence 
$w \ge  \mathrm{conc}_{[0, \, s] } f $, which contradicts the assumption $f (0) <  \mathrm{conc}_{[0, \, s] } f  (  0 ) $. 

In an entirely similar way one proves that $f (b) =  \mathrm{conc}_{[0, \, s] } f  (  b ) $  even if $b = s$.  This completes the proof of the third point. 

To prove the fourth point we proceed as follows. From the points $2$ and $3$ we know that $ \mathrm{conc}_{[0, \, s] } f $ is a bounded function that satisfies 
$$
    \Big|   \Big(     \mathrm{conc}_{[0, \, s] } f   \Big)' (\tau)  \Big| \leq \|   f ' \|_{\mathcal{C}^0 } \quad \forall \, \tau \in \, ]0, \, s[
$$
Thus, it has bounded total variation and hence the following limit exists:
$$
    \lim_{x \to 0^+ }  \mathrm{conc}_{[0, \, s] } f  (x) : = L.
$$
We want to prove $L = f (0)$. Let
$$
     S : = \Big\{ \tau  \in ]0, \, s[ \, : \; f(\tau) < \mathrm{conc}_{[0, \, s] } f  (\tau)  \Big\}.
$$
The set $S$ is open and hence 
$$
     S = \bigcup_{n > 1 } ]a_n , \, b_n[
$$
for suitable sequences $\{ a_n \}$ and $\{ b_n \}$.

If 
$$
      \inf \{ a_n \} =0,
$$
there are two possibilities. If the infimum is actually a minimum, then 
$$
    f(\tau) = \frac{f (b) - f (0) }{b}  \quad  \forall \; \tau \in ]0, \, b[
$$
for some $b>0$. One can then exploit the same function $w$ as before to prove that $f(0) =L$. If 
there exists  a subsequence $a_{n_k} \to 0^+$, then 
$$
     L = \lim_{k \to + \infty} \mathrm{conc}_{[0, \, s] } f  (a_{n_k} ) =  \lim_{k \to + \infty} f (a_{n_k}) = f (0).
$$
Finally, if
$$
      \inf \{ a_n \} >0,
$$
then
$$
     f (\tau) = \mathrm{conc}_{[0, \, s] } f  (\tau) \quad \forall \; \tau \in ]0, \, \inf \{ a_n \}[,
$$
 then by the continuity of $f$  $L= f (0)$. 
 
 Analogous considerations guarantee that $ \mathrm{conc}_{[0, \, s] } f $ is continuous at $t=s$ and that $ \mathrm{conc}_{[0, \, s] } f (s) = f(s)$.
\end{proof}

Because of points 2 and 3 of Lemma \ref{l:reg_concave}, the function $\Big( \mathrm{conc}_{[0, \, s] } f \Big)'$  is continuous 
at every point in $]0, \, s[$. Also, relying on considerations analogous to those performed in the last part of the proof of 
Lemma \ref{l:reg_concave}, one can prove that $\Big( \mathrm{conc}_{[0, \, s] } f \Big)'$ is continuous at $t=0$ and $t=s$. Note that it is actually 
enough to prove that the following limits exist:
$$
    \lim_{\tau \to 0^+} \Big( \mathrm{conc}_{[0, \, s] } f \Big)' (\tau) \qquad   \lim_{\tau \to s^-} \Big( \mathrm{conc}_{[0, \, s] } f \Big)' (\tau). 
$$
Indeed, one can then apply the theorem of the limit of the derivative. We conclude that the following holds true:
\begin{pro}
\label{p:reg_conc}
          If $f \in \mathcal{C}^{1, \, 1}_k ([0, \, s])$ , then $   \mathrm{conc}_{[0, \, s] } f  \in  \mathcal{C}^{1, \, 1}_k ([0, \, s])$.
\end{pro}
The following result describes the relation between the concave and the monotone envelope. The proof exploits considerations similar to those used to prove Lemma  \ref{l:reg_concave} and it will be therefore omitted.  
\begin{lem}
\label{l;mon}
          Let $f \in \mathcal{C}^{1, \, 1}_k ([0, \, s])$, then 
          \begin{equation*}
          \mathrm{mon}_{[0, \, s] } f  (\tau) =
          \left\{
          \begin{array}{ll}
                     \mathrm{conc}_{[0, \, s] } f (\tau)               & \mathrm{if } \;     \tau \leq \tau_0 \\
                     \mathrm{conc}_{[0, \, s] } f (\tau_0)           & \mathrm{if }   \;   \tau >   \tau_0. \\
          \end{array}
          \right.
          \end{equation*} 
          The value $\tau_0$ is given by
          \begin{equation}
          \label{e:tau0}
                    \tau_0 : = max \Big\{ t \in [0, \, s]:   \;  \big(    \mathrm{conc}_{[0, \, s] } f  \big) ' (t) \ge 0    \Big\}.
          \end{equation}
          If $ \big(    \mathrm{conc}_{[0, \, s] } f  \big) ' (t)$ is negative for every $t$, then we set $\tau_0 = 0$.
\end{lem}
Combining Lemma \ref{l;mon} and Proposition \ref{p:reg_conc} one deduces that if $f$ belongs to $ \mathcal{C}^{1, \, 1}_k ([0, \, s])$ then 
$$
     \mathrm{mon}_{[0, \, s] } f  \in  \mathcal{C}^{1, \, 1}_k ([0, \, s_1]).
$$

The following propostion collects some estimates that will be exploited in Section \ref{subsubsec_Binv_char}:
\begin{pro}
\label{p:es_mon}
          Let $ f, \, g   \in \mathcal{C}^{1, \, 1}_k ([0, \, s])$. Then 
          \begin{enumerate}
          \item 
          $$
               \|    \mathrm{conc}_{[0, \, s] } f  \|_{\mathcal{C}^0 } \leq \| f  \|_{\mathcal{C}^0 }   \qquad 
               \|   \big(   \mathrm{conc}_{[0, \, s] } f   \big) '  \|_{\mathcal{C}^0 }  \leq \|   f '   \|_{\mathcal{C}^0 }    .
          $$
          \item 
          \begin{equation}
         \label{e:monfg}
                     \| \mathrm{conc}_{[0, \, s] } f  -  \mathrm{conc}_{[0, \, s] } g    \|_{\mathcal{C}^0 }   \leq \| f -  g  \|_{\mathcal{C}^0 }        
                     \qquad 
                           \|  \big( \mathrm{conc}_{[0, \, s] } f    -  \mathrm{conc}_{[0, \, s] } g \big)'    \|_{\mathcal{C}^0 }   \leq \| f'  -  g '  \|_{\mathcal{C}^0 }      
          \end{equation}
          \end{enumerate}
\end{pro}
\begin{proof}
To prove the first point, note that 
$$
     \mathrm{conc}_{[0, \, s] } f  (t )  \leq  \mathrm{conc}_{[0, \, s] }  ( f - g ) (t)+  \mathrm{conc}_{[0, \, s] } g (t)   \quad \forall \, t  
$$
because $\mathrm{conc}_{[0, \, s] }  ( f - g ) (t)+  \mathrm{conc}_{[0, \, s] } g (t) $ is a concave function greater or equal to $f - g +g = f$.  Thus, by the first point in Lemma \ref{l:reg_concave}, 
$$
        \mathrm{conc}_{[0, \, s] } f  (t) -    \mathrm{conc}_{[0, \, s] } g (t) \leq    \| \mathrm{conc}_{[0, \, s] } (f   - g )    \|_{\mathcal{C}^0   }  \leq
        \| f - g \|_{  \mathcal{C}^0    }    \quad \forall \, t \in [0, \, s]
$$
and by analogous considerations 
$$
        \mathrm{conc}_{[0, \, s] } g  (t) -    \mathrm{conc}_{[0, \, s] } f (t)  \leq
        \| f - g \|_{  \mathcal{C}^0    }    \quad \forall \, t \in [0, \, s]
$$
This concludes the proof of the first point. 

To prove the second point, let us first make a preliminary consideration:
\begin{equation}
\label{e:Binv:concint:pc}
         \big(    \mathrm{conc}_{[0, \, s] }  g \big) ' (s) \neq g ' (s) \implies \; \exists \; \alpha < s:  \; \big(    \mathrm{conc}_{[0, \, s] } g   \big) ' (t) = \mathrm{constant} \; \forall \, 
         t \in  [ \alpha, \, s]  
\end{equation}
Indeed, if there is $t_n \to s^-$ such that $g(t_n) =    \mathrm{conc}_{[0, \, s] }  g (t_n)$ then $g'(t_n) = \big(    \mathrm{conc}_{[0, \, s] } g \big)' (t_n)$ and hence 
$\big(    \mathrm{conc}_{[0, \, s] } g \big)' (s)  = g ' (s) $. In particular, if $  \big(    \mathrm{conc}_{[0, \, s] }  g  \big) ' (s) \neq g ' (s) $ then $    \mathrm{conc}_{[0, \, s] }  g> g $ on an interval 
$[\alpha, \, s] $ and thus exploiting Lemma \ref{l:reg_concave}, point 2, one conclude with \eqref{e:Binv:concint:pc}. 

Let $\tau$ be a point at which the maximum of $\big|   \big(     \mathrm{conc}_{[0, \, s] }  f -    \mathrm{conc}_{[0, \, s] }  g \big) '   \big|$ is 
assumed.  Just to fix the ideas, suppose
$$
      \big(        \mathrm{conc}_{[0, \, s] }  f -    \mathrm{conc}_{[0, \, s] }  g \big)'  (\tau)  \ge 0.
$$ 
Because of the second point in Lemma \ref{l:reg_concave}, if $ \big(   \mathrm{conc}_{[0, \, s] }  f   \big) ' (\tau)  \neq f ' (\tau)$, then there are $a, \; b $ such that 
$\tau \in [a, \, b]$ and 
$$
    \big(    \mathrm{conc}_{[0, \, s] } f \big)' (  t ) = \frac{ f (a) - f  (b) }{b - a} \quad \forall \, t \in [a, \,b].
$$
Moreover, since $ \big(    \mathrm{conc}_{[0, \, s] }  g \big)' $ is non increasing, it is not restrictive to suppose $\tau = b$. Also, because of the first point in the statement of Lemma \ref{l:reg_concave}, if $b< s$ then $ \big(   \mathrm{conc}_{[0, \, s] }  f   \big) ' (b)  =  f' (b)$. Thus, if $ \big(   \mathrm{conc}_{[0, \, s] }  f   \big) ' (b)  \neq f ' (b)$ then $b =s$. 
If $ \big(   \mathrm{conc}_{[0, \, s] }  g \big) ' (s) = g' (s) $ then 
$$
   \big(    \mathrm{conc}_{[0, \, s] }  f  -    \mathrm{conc}_{[0, \, s] }  g \big)' ( s ) \leq  \big( f - g \big)' (s). 
$$
Indeed,
$$
    f (s)  =    \mathrm{conc}_{[0, \, s] }  f (s)  \qquad    \mathrm{conc}_{[0, \, s] }  f (t) \ge f (t) \quad  \forall \, t \leq b
$$
and hence $\big(    \mathrm{conc}_{[0, \, s] }  f  \big)' ( s )  \leq f' (s) $.

On the other side, if $ \big(   \mathrm{conc}_{[0, \, s] }  g \big) ' (s) \neq g' (s) $ then combining \eqref{e:Binv:concint:pc} and the second part of Lemma \ref{l:reg_concave} one gets
that there exists $\alpha < s$ such that 
$$
     \big(       \mathrm{conc}_{[0, \, s] }  g \big)' (s)  = \frac{ g (  s) - g (  \alpha) }{s-\alpha}.  
$$
Moreover, $\alpha \ge a$. Indeed, if $\alpha < a$ then $ \big( \mathrm{conc}_{[0, \, s] }  f  \big)' (\alpha)  >  \big(   \mathrm{conc}_{[0, \, s] }  f \big) ' (b)$ and hence $b$ would be no more a point at which the maximum of $\big(   \mathrm{conc}_{[0, \, s] }  f -   \mathrm{conc}_{[0, \, s] }  g\big)' $ is assumed. By Lagrange's theorem, there exists $x \in [\alpha, \, s]$ such that
\begin{equation*}
\begin{split}
    \big( f - g \big)' (x) 
&    =  \frac{       [  f (s) - g (s) ]  - [f (\alpha)  - g (\alpha)]    }{ s - \alpha} \\  
&    \ge \frac{ [  f (s)  - g(s)] -  [  \mathrm{conc}_{[0, \, s] }      f  (\alpha)  - g (\alpha)] }{ s - \alpha}  \\
&     = \frac{ f (s)  - \mathrm{conc}_{[0, \, s] } f(\alpha)} {s - \alpha}  - \big( \mathrm{conc}_{[0, \, s] } g \big)' (s) \\
&    =  \big(  \mathrm{conc}_{[0, \, s] }  f \big)' (x) - \big( \mathrm{conc}_{[0, \, s] }  \big)' (x) \\
\end{split}
\end{equation*}
This concludes the proof of the second point in the case $\big( \mathrm{conc}_{[0, \, s] }  f \big)' ( \tau)  \ge  \big( \mathrm{conc}_{[0, \, s] }  f \big)' ( \tau)  $. The opposite case  $\big( \mathrm{conc}_{[0, \, s] }  f \big)' ( \tau)  <  \big( \mathrm{conc}_{[0, \, s] }  f \big)' ( \tau)  $
is entirely analogous.
\end{proof}

The following result concerns the dependence of the concave envelope from the interval:
\begin{pro}
\label{p:concint}
          Let $s_1 \leq s_2$, $f \in \mathcal{C}^{1, \, 1}_k ([0, \, s_2])  $ and assume that $\| f'    \|_{  \mathcal{C}^0 ([0, \, s_2])  }   \leq C_1 \delta_0$. Then there are constants $C_2$ and $C_3$ such that 
          \begin{equation}
          \label{e:Binv:concint:sta1}
                \|  \mathrm{conc}_{[0, \, s_1] } f - \mathrm{conc}_{[0, \, s_2] } f   \|_{  \mathcal{C}^0 ([0, \, s_1])  }
                \leq C_2 \delta_0 (s_2 -  s_1 )          
           \end{equation}
           and 
           \begin{equation}
           \label{e:Binv:concint:sta2}
               \|  ( \mathrm{conc}_{[0, \, s_1] } f  -  \mathrm{conc}_{[0, \, s_2] } f  )'  \|_{  \mathcal{C}^0 ([0, \, s_1])  }   
                       \leq C_3 \delta_0 (s_2 -  s_1 ) .
         \end{equation}        
\end{pro}
\begin{proof}
We will proceed in several steps: 
\begin{enumerate}
\item For every $\tau \in [0, \, s_1]$,
\begin{equation}
\label{e:Binv:concint:o}
         f(\tau) \leq \mathrm{conc}_{[0, \, s_1]}  f  (\tau)  \leq \mathrm{conc}_{[0, \, s_2]}  f  (\tau)   
\end{equation}
Indeed, the restriction of $\mathrm{conc}_{[0, \, s_2]}  f  (\tau) $ to the interval $[0, \, s_1]$ is a concave and non decreasing function which is also greater then $f$. From \eqref{e:Binv:concint:o} one deduces that, for every $\tau \in [0, \, s_1]$, the following implication holds: 
\begin{equation}
\label{e:Binv:monint:i}
              f(\tau)   =   \mathrm{conc}_{[0, \, s_2]}  f  (\tau)   \implies      \mathrm{conc}_{[0, \, s_1]}  f  (\tau)  =  \mathrm{conc}_{[0, \, s_2]}  f  (\tau).        
\end{equation}
\item It is not restrictive to assume 
\begin{equation}
\label{e:Binv:concint:a} 
          \exists   \; \bar{s} \in ]  s_1, \, s_2]: \;  f (  \bar{s})  >  \frac{d}{d \tau }\Big(     \mathrm{conc}_{[0, \, s_1]}  f    \Big) \Big|_{\tau = s_1}  (\bar{s}  - s_1)+
            \mathrm{conc}_{[0, \, s_1]}  f (s_1)    
\end{equation}
Indeed, if \eqref{e:Binv:concint:a} is not satisfied then 
$$
    h (\tau) = 
    \left\{
    \begin{array}{ll}
                 \mathrm{conc}_{[0, \, s_1]}  f  (\tau)  &  \tau \in [0, \, s_1] \\
                 \displaystyle{ \frac{d}{d \tau }\Big(     \mathrm{conc}_{[0, \, s_1]}  f    \Big) \Big|_{\tau = s_1}  [\bar{s}  - s_1] +
            \mathrm{conc}_{[0, \, s_1]}  f (s_1)  }   & \tau \in ]s_1, \, s_2] \\  
    \end{array}
    \right.
$$
is a concave function which is  greater or equal to $f$ on $[0, \, s_2]$. Thus, 
$$
    h(\tau) \ge                  \mathrm{conc}_{[0, \, s_2]}  f  (\tau)  \quad \forall \; \tau \in [0, \, s_2]. 
$$  
In particular, on $[0, \, s_1]$ 
$$
                     \mathrm{conc}_{[0, \, s_1]}  f  (\tau)  \ge                  \mathrm{conc}_{[0, \, s_2]}  f  (\tau)  
$$
and thus by \eqref{e:Binv:concint:o} the two coincide.
\item Before proving \eqref{e:Binv:concint:sta1} we introduce other preliminary observations. 

Given two arbitrary points  $a$ and $b$ in $[0, \, s_2]$, let $\ell$ the line
$$
    \ell_{a\, b} (s) = \frac{1}{b - a} \bigg(   f(a) \Big(b -s \Big) +  f(b)   \Big(   s - a \Big)   \bigg). 
$$
Then 
\begin{equation}
\label{e:concint:l}
            \mathrm{conc}_{[0, \, s_1]}  \Big(  f - \ell_{a \, b}  \Big) =   \Big( \mathrm{conc}_{[0, \, s_1]}  f   \Big)  - \ell_{a \, b} 
            \qquad 
             \mathrm{conc}_{[0, \, s_2]}  \Big(  f - \ell_{a \, b}  \Big) =   \Big( \mathrm{conc}_{[0, \, s_2]}  f   \Big)  - \ell_{a \, b} . 
\end{equation}
Indeed, 
$$
        \mathrm{conc}_{[0, \, s_1]}  \Big(  f - \ell_{a \, b}  \Big) \leq   \Big( \mathrm{conc}_{[0, \, s_1]}  f   \Big)  - \ell_{a \, b}  
$$
because $  \Big( \mathrm{conc}_{[0, \, s_1]}  f   \Big)  - \ell_{a \, b}  $ is a concave function bigger than $f - \ell_{a \, b} $. On the other side, 
$$
      \Big( \mathrm{conc}_{[0, \, s_1]}  f   \Big)  \leq \ell_{a \, b}  +  \mathrm{conc}_{[0, \, s_1]}  \Big(  f - \ell_{a \, b}  \Big) 
$$
since $\ell_{a \, b}  +  \mathrm{conc}_{[0, \, s_1]}  \Big(  f - \ell_{a \, b}  \Big) $ is a concave function bigger then $f - \ell_{a \, b}  + \ell_{a \, b}  = f$. 

Equalities \eqref{e:concint:l} imply, in particular, that in proving  \eqref{e:Binv:concint:sta1} and  \eqref{e:Binv:concint:sta2} one can fix two points $a$ and $b$ in $[0, \, s_2]$ and assume that 
\begin{equation}
\label{e:Binv:concint:ab}
    f(a) =0 
    \qquad f(b) =0.
\end{equation}
Indeed, instead of $f$ one can consider $f - \ell_{a \, b}$, which satisfies \eqref{e:Binv:concint:ab}. Thanks to  \eqref{e:concint:l}, $f$ satisfies \eqref{e:Binv:concint:sta1} and  \eqref{e:Binv:concint:sta2} if and only if $f - \ell_{a \, b}$ does. 
\item Define 
$$
    {s}^{\ast}: = \max\{   s\in [0, \, s_1]: \;    \mathrm{conc}_{[0, \, s_1]}   f  (s) = \mathrm{conc}_{[0, \, s_2]}   f (s)     \}. 
$$ 
If $s \leq s ^{\ast}$, then $ \mathrm{conc}_{[0, \, s_1]}   f  (s) = \mathrm{conc}_{[0, \, s_2]}   f (s)   $. Indeed, by contradiction assume 
that there exists a $\tau \in [0, \, s^{\ast}]$ such that 
\begin{equation}
\label{e:binv:concint:c}
     \mathrm{conc}_{[0, \, s_1]}   f  (\tau) <  \mathrm{conc}_{[0, \, s_2]}   f (\tau).   
\end{equation}
By \eqref{e:Binv:concint:o} this implies that $f (\tau) < \mathrm{conc}_{[0, \, s_2]}   f (\tau)$ and hence, thanks to the third part of Lemma \ref{l:reg_concave}, there are $\tau_a$ and 
$\tau_b$ such that $\tau_a < \tau< \tau_b$ and $ \mathrm{conc}_{[0, \, s_2]}   f $ on $[\tau_a, \, \tau_b]$ is a linear. Moreover, 
$
    f(\tau_a) = \mathrm{conc}_{[0, \, s_2]}   f (\tau_a)
$
and 
$
    f(\tau_b) = \mathrm{conc}_{[0, \, s_2]}   f (\tau_b).
$
Because of \eqref{e:Binv:monint:i}, this implies 
$
    f(\tau_a) = \mathrm{conc}_{[0, \, s_1]}   f (\tau_a).
    $
Also, because of the previous observation, it is not restrictive to assume 
$f (\tau_a )= f (\tau_b)$. Since $\mathrm{conc}_{[0, \, s_1]}   f $ is a concave function, if $\mathrm{conc}_{[0, \, s_1]}   f (\tau_a) =0$ and 
$\mathrm{conc}_{[0, \, s_1]}   f (\tau) <0$, then $\mathrm{conc}_{[0, \, s_1]}   f (t) <0$ for every $t \ge \tau$. 
If $\tau_b > s_1$, then $\mathrm{conc}_{[0, \, s_1]}   f < \mathrm{conc}_{[0, \, s_2]}   f $ 
on $] \tau_a, \, s_1$ and hence $s^{\ast} \leq \tau_a < \tau$, which contradicts our hypothesis. On the other side, if $\tau_b \leq s_1$ then from $  f(\tau_b) = \mathrm{conc}_{[0, \, s_2]}   f (\tau_b)$ and \eqref{e:Binv:monint:i} we deduce      
$
    f(\tau_b) = \mathrm{conc}_{[0, \, s_1]}   f (\tau_b).
$
Being a concave function, $\mathrm{conc}_{[0, \, s_1]}   f $ on $[\tau_a, \, \tau_b]$ is above the line joining 
$\mathrm{conc}_{[0, \, s_1]}   f (\tau_a)$ and $\mathrm{conc}_{[0, \, s_1]}   f (\tau_b)$. 
Thus, in particular, $\mathrm{conc}_{[0, \, s_1]}   f (\tau)  \ge  \mathrm{conc}_{[0, \, s_2]}   f (\tau)  $, which contradicts \eqref{e:binv:concint:c}.
\item We now prove \eqref{e:Binv:concint:sta1}. 
Let $s^{\ast}$ as in the previous step. If $s \in ] s^{\ast}, \, s_1 ]$, then $\mathrm{conc}_{[0, \, s_2]}   f > f $ and hence there are $s_a$, $s_b$ such that $[s^{\ast}, \, s_1]  \subseteq 
[ s_a, \, s_b] \subseteq [ 0, \, s_2]$ and on $[s_a, \, s_b]$ $\mathrm{conc}_{[0, \, s_2]}   f $ is the linear interpolation between $f(s_a)$ and $f(s_b)$. Relying on \eqref{e:concint:l} , one can assume $f(s_a) = f(s_b) = 0$, getting $\mathrm{conc}_{[0, \, s_2]}   f \equiv 0$ on $[s_a, \, s_b]$. Thus, $\mathrm{conc}_{[0, \, s_1]}   f (s^{\ast}) = 0$ and $\mathrm{conc}_{[0, \, s_1]}   f <0$ on $]s^{\ast}, \, s_1]$. This implies 
$$
    \big( \mathrm{conc}_{[0, \, s_1]}   f  \big)' (s^{\ast}) \leq 0. 
$$
Since  $\big( \mathrm{conc}_{[0, \, s_2]}   f \big) '$ is the derivative of a concave function, then it is non increasing and 
$$
    \big( \mathrm{conc}_{[0, \, s_1]}   f  \big)' (s) \leq 0  \quad \forall \, s \in [s^{\ast}, \, s_1]. 
$$
In particular,  $\mathrm{conc}_{[0, \, s_1]}   f $ is non increasing on $[s^{\ast}, \, s_1]$ and hence its minimum on $[s^{\ast}, \, s_1]$ is assumed in $s_1$. Thanks to the previous step, 
$$
    \max_{ s \in [0, \, s_1]  } \{  | \mathrm{conc}_{[0, \, s_1]}   f   -  \mathrm{conc}_{[0, \, s_1]}   f     |    \} =  
      \max_{ s \in [s^{\ast}, \, s_1]  } \{  | \mathrm{conc}_{[0, \, s_1]}   f   -  \mathrm{conc}_{[0, \, s_1]}   f     |    \} , 
$$
which is equal to $-  \mathrm{conc}_{[0, \, s_1]}   f   (s_1)$  by the previous observations. By the forth part of Lemma \ref{p:reg_conc},  $ \mathrm{conc}_{[0, \, s_1]}   f   (s_1) = f(s_1)$. 
To conclude, one can observe that 
$$
    -  f(s_1) = f(s_b) - f (s_1) \leq C_1 \delta (s_b - s_1 ) \leq C_1 \delta (s_2 - s_1).
$$
\item Before proving \eqref{e:Binv:concint:sta2} we introduce the following observations. Let $g \in \mathcal{C}^{1, \, 1}_k[a, \, b]$ such that $g(a)= g(b)=0 $. Then 
\begin{equation}
\label{e:Binv:concint:g}
         g (s) \ge \frac{k}{2}  (s - a) (s- b) \quad \forall \, s \in [a, \, b].
\end{equation}
Indeed, consider for simplicity the case $a=0$, $b =1$. Let $\bar{s}$ denote a point at which $g$ assumes its minimum and suppose, just to fix the ideas, that $\bar{s} \leq 1/2$. Denote by $h$ the first derivative of  $g$: $h$ is Lipschitz continuous function with Lipschitz constant smaller than $k$ and such that $h(\bar s) =0$.Thus, 
 $h(s) \ge k (s - \bar{s})$ if $s \leq \bar{s}$. Moreover,  Let us consider the Cauchy problem for $g$ at $s=0$: 
\begin{equation*}
\left\{ 
\begin{array}{ll}
           g'(y) = h(y)  \ge k (s - \bar{s}) \\
           g(0) =0.
\end{array}
\right.
\end{equation*}
By a comparison argument one gets
\begin{equation}
\label{e:Binv:concint:c1} 
        g(s) \ge k s^2 /2 - k \bar{s} s  \quad \forall \, s \in [0, \, \bar{s}].
\end{equation}        
Moreover, 
\begin{equation}
\label{e:Binv:concint:c2}
           \forall \, s \in [0, \, 1], \; g(s) \ge g (\bar{s}) \ge - \frac{k}{2} \bar{s}^2.  
\end{equation}
To complete the proof of \eqref{e:Binv:concint:g}, let us consider to cases separately. If $g ' (1 /2) \leq 0$, consider the Cauchy problem 
\begin{equation*}
\left\{ 
\begin{array}{ll}
           h'(s) = g''(s) \leq k \\
           h(1 /2)\leq 0.
\end{array}
\right.
\end{equation*}
Then $h (s) \leq k ( s - 1 /2 )$ and thus the solution of 
\begin{equation*}
\left\{ 
\begin{array}{ll}
           g'(s) = h(s)  \leq k (s - 1 /2) \\
           g(1) =0.
\end{array}
\right.
\end{equation*}
satisfies 
\begin{equation}
\label{e:Binv:concint:c3}
          g (s) \ge k s (s -1) /2  \quad \forall s \in [1/2, \, 1] .
\end{equation}
Define 
\begin{equation*}
           \psi (s) =
           \left\{ 
           \begin{array}{lll}
                      k s^2 /2 - k \bar{s} s & s \leq \bar{s} \\
                      - k \bar{s}^2 /2          & \bar{s} \leq s \leq 1 /2 \\
                      k s (s -1) /2               &  s > 1 /2.
           \end{array}
           \right.
\end{equation*}
By direct check one gets $\psi (s) \ge      k s (s -1) /2   $ for every $s \in [0, \, 1]$. Putting together \eqref{e:Binv:concint:c1}, \eqref{e:Binv:concint:c2} and \eqref{e:Binv:concint:c3}
one gets therefore 
$$
     g(s) \ge \psi (s) \ge   k s (s -1) /2   \quad \forall \, s \in [0, \, 1]. 
$$

If $g'(  1/2) >  0$ then 
\begin{equation}
\label{e:Binv:conincint:c4}
            g(s) \ge   k s (s -1) /2   \quad \forall \, s \in [1/2, \, 1].     
\end{equation}
Indeed, thanks to \eqref{e:Binv:concint:c2} $g (1 /2 ) \ge  - k /8$. Assume by contradiction 
that  
\begin{equation}
\label{e:Binv:concint:c5}
          g (s) <  k s (s -1) /2   
\end{equation}
for some $s \in [1/2, \, 1]$. Define 
$$
     s_n : = \max \{   s: \; g(s) \ge k s (s -1) /2  \}. 
$$          
Since $g (1 /2 ) \ge  - k /8$ and $g(s_n) = k s_n (s_n -1) /2 $, then there exits $\tau  \in [1/2, \, s_n]$ such that 
$$
    g '(\tau) \leq k (\tau - 1/2).
$$
Consider the Cauchy problem 
\begin{equation*}
\left\{ 
\begin{array}{ll}
           h'(s) = g''(s)  \leq k  \\
           h(\tau) = g'(\tau) \leq  k (\tau - 1/2),
\end{array}
\right.
\end{equation*}
then $h (s ) \leq  k (s - 1/2)$ for $s \ge \tau$. Then the solution of the backward Cauchy problem 
\begin{equation*}
\left\{ 
\begin{array}{ll}
           g'(s) = h(s)    \leq k (s - 1/2)  \\
           g(1) = 0
\end{array}
\right.
\end{equation*} 
satisfies $g (s) \ge  k s (s -1) /2 $ for every $s \in [\tau, \, 1]$. This contradicts \eqref{e:Binv:concint:c5} 
and hence \eqref{e:Binv:conincint:c4} holds. One can then exploit the same function $\psi$ considered in the case $g'(1/2) \leq 0$ and conclude.

\item We now prove \eqref{e:Binv:concint:sta2}. Let ${s}^{\ast}$ be as in the previous steps, then $  \mathrm{conc}_{[0, \, s_2]}   f  =   \mathrm{conc}_{[0, \, s_1]}   f $ 
on $[0, \,  s^{\ast}]$ and hence 
$$
     \big(   \mathrm{conc}_{[0, \, s_1]}   f \big) ' (s) =  \big(  \mathrm{conc}_{[0, \, s_2]}   f  \big) ' (s) \quad s \in [0, \, s^{\ast}[ .  
$$   
Thus, 
\begin{equation}
\label{e:Binv:concint:m}
    \|   \big(   \mathrm{conc}_{[0, \, s_1]}   f \big) '   -   \big(  \mathrm{conc}_{[0, \, s_2]}   f  \big) ' \|_{  \mathcal{C}^0  [0, \, s_1]} = \max_{s \in [s^{\ast}, \, s^1]} 
    |  \big(   \mathrm{conc}_{[0, \, s_1]}   f \big) ' (s)  -  \big(  \mathrm{conc}_{[0, \, s_2]}   f  \big) ' (s)|.
\end{equation}
As in step 4, denote by $s_a$ and $s_b$ the values such that 
$$
 \big(  \mathrm{conc}_{[0, \, s_2]}   f  \big) ' (s) = \frac{f(s_a) - f(s_b)} {s_a - s_b } \quad s \in ]s_a, \, s_b[. 
$$
As before it is not restrictive to assume $f(s_a) = f(s_b) =0$ and in this case 
$$
     \big(  \mathrm{conc}_{[0, \, s_1]}   f  \big) ' (s) \leq 0 \quad s \in [s^{\ast}, \, s_1]
$$
Since it is a non increasing function, 
the maximum in \eqref{e:Binv:concint:m} is given by 
$-  \big(  \mathrm{conc}_{[0, \, s_2]}   f  \big) ' (s_1)$.
Define 
$$
     s^0: = \sup  \{   s \in [    0, \, s_1 [: \;      \mathrm{conc}_{[0, \, s_1]}   f (s)  = f (s) \}.
$$

If $s^0 =0$, then for every $s \in [0, \, s_1[$
$$
    f(s) <   \mathrm{conc}_{[0, \, s_1]}   f (s)  \leq   \mathrm{conc}_{[0, \, s_2]}   f (s) .
$$
In particular, one has $s_a =0$ and hence $f(0) =0$. By the third part of Lemma \ref{l:reg_concave}
$$
    \big(   \mathrm{conc}_{[0, \, s_1]}   f \big) '(s_1) = \frac{f(s_1)  - f(0)}{ s_1} = \frac{  f(s_1)}{s_1}. 
$$
Now we can apply the previous step to $f$ defined on $[0, \, s_b]$, thus obtaining
$$
    - \frac{  f(s_1)}{s_1} \leq  \frac{k}{2 s_1}  s_1 (s_b - s_1)   \leq \frac{k}{2} (s_b - s_1) \leq \frac{k}{2} (s_2  - s_1).  
$$

If $s^0 = s_1$, then $\big(   \mathrm{conc}_{[0, \, s_1]}   f  \big)' (s_1) = f ' (s_1) $. Consider two subcases separately. If $s_b < s_2$, then 
also  $\big(   \mathrm{conc}_{[0, \, s_2]}   f  \big)' (s_2) = f ' (s_2) $ and hence  $f'  (s_2) =0$ because we are assuming $f(s_a) = f(s_b) = 0$.
Then 
$$
    - f'(s_1) = f'(s^b) - f ' (s_1) \leq k (s^b - s_1) \leq k (s_2 -s_1).
$$
On the other side, if $s_b = s_2$ one can proceed as follows. Since $f'(s_1)\leq 0$, $f (s_1) \leq 0$ and $f (s_2) =0$, then there exists $\tilde{s} \in [s_1, \, s_2]$ such that $f'(\tilde{s} )= 0$. Then 
$$
    - f'(s_1) = f'(\tilde{s}) - f ' (s_1) \leq k (\tilde{s} - s_1) \leq k (s_2 -s_1).
$$

We are now left to deal with the case $0 < s^0 < s_1$. One then has 
\begin{equation}
\label{e:Binv:concint:u}
      - \big(  \mathrm{conc}_{[0, \, s_1]}   f  \big)'  (s_1) = \frac{     f (s^0)   - f (s_1)}{  s_1  - s^0}
\end{equation}
Because of the implication \eqref{e:Binv:monint:i}, $s_a \leq s^0$ and hence $f(s^0) \leq 0$. Moreover, $f(s^0) =  \mathrm{conc}_{[0, \, s_1]}   f (s^0)$, $f(s^1) =  \mathrm{conc}_{[0, \, s_1]}   f (s_1)$ and $\big(  \mathrm{conc}_{[0, \, s_1]}   f \big) '  (s) \leq 0$ on $[  s^0, \, s_1]$. Thus, 
$$
   f (s^1) \leq f (s^0) \leq 0.
$$
Since $f (s_b) =0$, there exists $\check{s} \in [s_1, \, s_b]$ such that  $f(\check{s}) = f (s^0) $. Applying the previous step to the function $f - f (s^0)$ defined on the interval $[s^0, \, \check{s}]$, one gets $ f(s^1)  - f(s^0) \ge k ( s_1 - s^0) (s_1 - \check{s}) / 2$. Inserting in  \eqref{e:Binv:concint:u}, one gets 
$$
     - \big(  \mathrm{conc}_{[0, \, s_1]}   f  \big)'  (s_1) \leq \frac{k}{2} (\check{s}  - s_1)  \leq \frac{k}{2} (s_2  - s_1).
$$ 
This concludes the proof of \eqref{e:Binv:concint:sta2} and hence of the lemma. 
\end{enumerate}
\end{proof}

Combining Lemma \ref{l;mon} and Proposition \ref{p:monint} one finally gets 
\begin{pro}
\label{p:monint}
          Let $s_1 \leq s_2$, $f \in \mathcal{C}^{1, \, 1}_k ([0, \, s_1])  $ and assume that $\| f'    \|_{  \mathcal{C}^0 ([0, \, s_2])  }   \leq C_1 \delta_0$. Then there are constants $C_2$ and            $C_3$ such that 
          $$
                \|  \mathrm{mon}_{[0, \, s_1] } f - \mathrm{mon}_{[0, \, s_2] } f   \|_{  \mathcal{C}^0 ([0, \, s_1])  }
                \leq C_2 \delta_0 (s_2 -  s_1 )          
           $$
           and 
           $$
               \|  ( \mathrm{mon}_{[0, \, s_1] } f  -  \mathrm{mon}_{[0, \, s_2] } f  )'  \|_{  \mathcal{C}^0 ([0, \, s_1])  }   
                       \leq C_3 \delta_0 (s_2 -  s_1 ) .
          $$        
\end{pro}

\subsection{The hyperbolic limit in the case of an 
invertible viscosity matrix} \label{subsec_Binv_riemannsolver}
\subsubsection{The hyperbolic limit in the case of a Cauchy problem}
\label{subsub_Binv_speed+} 
In this section for completeness we will give a quick review of the construction of the Riemann solver for the Cauchy problem. We refer to 
\cite{Bia:riemann} for the complete analysis. The goal is to characterize the limit of 
\begin{equation}
\label{e:BInv:s+:app}
\left\{
\begin{array}{ll}
          E( \ue) \ue_t + A (\ue, \, \ue_x ) \ue_x = \ee B(  \ue ) \ue_{xx} \\
          \ue (0, \,,  x) = 
          \left\{
          \begin{array}{ll}
                     u^-    &  x < 0\\
                     \bar{u}_0     & x  \ge 0  \\
          \end{array}
          \right.
\end{array}
\right.
\end{equation}
under Hypotheses \ref{hyp_invertible}, \ref{hyp_hyperbolic_I}, \ref{hyp_small_bv} and \ref{hyp_stability}.

The construction works as follows. Consider a travelling profile 
$$
    B (U ) U'' = \Big( A(U, \, U' )  - \sigma E  ( U )     \Big) U',
$$
then $U$ solves 
\begin{equation}
\label{eq_Binv_system_trav} \left\{
\begin{array}{lll}
      u' = v \\
      B(u )v' = \Big( A(u, \, v) - \sigma E(u) \Big)v \\
      \sigma' = 0.
\end{array}
\right.
\end{equation}
Let $\lambda_i (\bar{u}_0)$ be the $i$-th eigenvalue of $E^{-1} (\bar{u}_0) A(\bar{u}_0, \, 0)$ and let $r_i (\bar{u}_0)$ be the corresponding eigenvalue. If one linearizes the previous system around the equilibrium point $( \bar{u}_0, \, \vec{0}, \, \lambda_i (\bar{u}_0)  )$ one finds 
\begin{equation*}
     \left(
     \begin{array}{ccc}
                0 & I_N & 0 \\
                0 &  A(\bar{u}_0, \, 0) - \lambda_i ( \bar u_0 )E(  \bar{u}_0) & 0 \\
		0 & 0 & 0 \\
     \end{array}
     \right)      
\end{equation*}
The generalized eigenspace corresponding to the $0$ eigenvector is 
\begin{equation*}
         V^c = \Big\{
            (u, \, x r_i(\bar{u}_0), \,  \sigma
             ): u \in \erreN, \, x, \, \sigma \in \mathbb{R}
         \Big\}.
\end{equation*}
It is then possible to define a  {\sl center} manifold which is parameterized by $ V^c$: we refer to \cite{KHass} for an extensive analysis, here we will just recall some of the fundamental properties of a center manifold.  

Every center manifold is invariant with respect to \eqref{eq_Binv_system_trav} and moreover satisfies the following property: let $(u^0, \, p^0, \, \sigma^0)$ belong to $\mathcal{M}^{c}$ 
and denote by $(u(x), \, p(x), \, \sigma (x))$ the orbit starting at $(u^0, \, p^0, \, \sigma^0)$. Then 
$$
    \lim_{x \to + \infty} \Big(u(x), \, p(x), \, \sigma (x) \Big)  e^{-  c x  / 2} 
    = (0, \, 0, \, 0)  \qquad    \lim_{x \to -  \infty} \Big(u(x), \, p(x), \, \sigma (x) \Big) e^{-  c x  / 2}  
    = (0, \, 0, \, 0) . 
$$
The constant $c$ is the separation speed defined in Hypothesis \ref{hyp_hyperbolic_I}. A center manifold $\mathcal{M}^{c}$ is defined in a neighbourhood of the equilibrium point $(\bar{u}_i, \, \vec{0}, \, \lambda(\bar{u}_i)$ and it is tangent to $V^c$ at $(\bar{u}_i, \, 0, \, \lambda(\bar{u}_i))$.  

Also, the following holds. Define 
$$
    V^{su} : = \Big\{
            (\bar{u}_0, \,  \sum_{ j \neq i} x_j r_j(\bar{u}_0), \,  \sigma
             ): u \in \erreN, \, x_j, \, \sigma \in \mathbb{R}
         \Big\},
$$
where $r_1 (\bar{u}_0) \dots r_{i-1} (\bar{u}_0), \, r_{i +1} (\bar{u}_0) \dots r_N (\bar{u}_0)$ are eigenvectors of   $E^{-1} (\bar{u}_0) A(\bar{u}_0, \, 0)$ different than $r_i (\bar{u}_0)$. 
If we write $\mathbb{R}^{2N + 1} = V^c  \oplus V^{su}$, then the map 
$$
    \Phi_{c}:  V^c \to  V^c \oplus V^{us}
$$
that parameterizes $\mathcal{M}^{c}$ can be chosen in such a way that if $\pi_{c}$ is the projection from $ V^c \oplus V^{us}$ onto $ V^c $, then $\pi_{c} \circ \Phi_{c}$ is the identity.

Fix a center manifold $\mathcal{M}^c$.
Putting all the previous considerations together, one gets that a point $(u, \, v, \, \sigma)$ belongs $\mathcal{M}^c$ if and only if 
\begin{equation}
\label{eq_app_center}
      v= v_i r_i(\bar{u}_0) + \sum_{j \neq i}  \psi_j
      (u, \, v_i, \, \sigma_i) r_j (\bar{u}_0).
\end{equation}
Since
all the equilibrium points $(u, \, 0, \, \sigma_i)$ lay on the
center manifold, it turns out that when $v_i=0$, then
\begin{equation*}
      \psi_j (u, \, 0, \, \sigma_i) =0 \quad \forall \, j, \; u, \;
      \sigma_i.
\end{equation*}
Hence for all $j$, $\psi_j (u, \, v_i, \, \sigma_i) = v_i
\phi_j(u, \, v_i, \, \sigma_i)$ for a suitable regular function
$\phi_j$. The relation \eqref{eq_app_center} can be therefore
rewritten as
\begin{equation*}
      v= v_i \Big( r_i(\bar{u}_0) + \sum_{j \neq i}
      \phi_j (u, \, v_i, \, \sigma_i) r_j(\bar{u}_0) \Big) : =
      v_i \tilde{r}_i (u, \, v_i, \, \sigma_i).
\end{equation*}
Because of the tangency condition of the center manifold to the
center space, it holds
\begin{equation*}
      \tilde{r}_i \Big(\bar{u}_0, \, 0, \, \lambda_i(\bar{u}_0 \Big)= r_i ( \bar{u}_0).
\end{equation*}
Inserting the expression found for $v$ into
\eqref{eq_Binv_system_trav}, one gets
\begin{equation*}
      v_{i x }  B(u) \tilde{r}_i + v_i^2 B(u) D \tilde{r}_i
      \tilde{r}_i + v_i v_{i x } B(u) \tilde{r}_{i v } =
      v_i \Big( A(u, \, v_i \tilde{r}_i)  - \sigma_i E(u) \Big)
      \tilde{r}_i.
\end{equation*}
Considering the scalar product of the previous expression with
$\tilde{r}_i$ one obtains
\begin{equation*}
      v_{ix } \Big( \langle \tilde{r}_i, \, B(u) \tilde{r}_i \rangle +
      v_i \langle \tilde{r}_i, \, B(u) \tilde{r}_{i v } \rangle \Big) =
      v_i \Big( \langle
      \tilde{r}_i, \, \Big(A(u, \, v_i \tilde{r}_i)- \sigma_i E(u)\Big)\tilde{r}_i \rangle -
      v_i \langle \tilde{r}_i, \, B(u) D \tilde{r}_i \tilde{r}_i \rangle
      \Big).
\end{equation*}
Hence setting
\begin{equation*}
\begin{split}
&      c_i (u, \, v_i, \, \sigma_i) := \langle \tilde{r}_i, \, B(u) \tilde{r}_i \rangle +
      v_i \langle \tilde{r}_i, \, B(u) \tilde{r}_{i v } \rangle  \\
&      a_i (u, \, v_i, \, \sigma_i):= \langle
      \tilde{r}_i, \, \Big(A(u, \, v_i \tilde{r}_i)- \sigma_i E(u)\Big)\tilde{r}_i \rangle -
      v_i \langle \tilde{r}_i, \, B(u) D \tilde{r}_i \tilde{r}_i \rangle \\
\end{split}
\end{equation*}
one can define
\begin{equation*}
      \phi_i(u, \, v_i, \, \sigma_i) : = \frac{ a_i (u, \, v_i, \, \sigma_i)}{ c_i (u, \, v_i, \,
      \sigma_i)}.
\end{equation*}
The fraction is well defined since $c_i(\bar{u}_0, \, 0, \,
\lambda_i(\bar{u}_0)) \ge c_B(\bar{u}_0)>0$ and hence $c_i$ is strictly
positive in a small neighborhood. The constant $c_B$  is as in the third condition in Hypothesis \ref{hyp_invertible}.

Thus, 
$$
   \frac{\partial \phi_i}{ \partial \sigma_i } \bigg|_{\big(\bar u_0, \, 0, \, \lambda_i(\bar u_0)\big)}=
      - \frac{a_i }{c_i^2} \,
      \frac{\partial c_i}{ \partial \sigma_i}+
      \frac{1}{c_i } \,
      \frac{\partial a_i}{ \partial \sigma_i}. 
$$
Since 
$$
   a_i(\bar{u}_0, \, 0, \, \lambda_i(\bar{u}_0)) \ge c_B(\bar{u}_0) = 0 
$$
then 
\begin{equation}
\label{eq_Binv_negative}
\begin{split}
      \frac{\partial \phi_i}{ \partial \sigma_i } \bigg|_{\big(\bar u_0, \, 0, \, \lambda_i(\bar u_0)\big)}
&      =      \frac{1}{c_i  } \,
      \bigg(
            \langle \tilde{r}_{i \sigma}, \, \Big(A(\bar{u}_0, \, 0)- \lambda_i E(\bar{u}_0) \Big)r_i (\bar u_0)\rangle
            - \langle r_i(\bar{u}_0), \, E(\bar{u}_0) r_i(\bar{u}_0)
            \rangle
            \bigg) \\
&     \quad + \frac{1}{c_i  } \,
      \bigg(
            \langle {r}_i(\bar{u}_0), \, \Big(A(\bar{u}_0, \, 0)- \lambda_i E(\bar{u}_0) \Big)\tilde{r}_{i \sigma} \rangle
      \bigg) \\
&      = - c_E(\bar{u}_0)< 0. \\
\end{split}
\end{equation}
In the previous computations, we have exploited the symmetry of
$A(\bar{u}_0, \, 0)$ and $E(\bar{u}_0)$ and hence the fact that
\begin{equation*}
       \left \langle {r}_i(\bar{u}_0, \, 0), \, \Big(A(\bar{u}_0, \, 0)- \lambda_i E(\bar{u}_0) \Big)\tilde{r}_{i \sigma}
       \right \rangle =  \left \langle \Big(A(\bar{u}_0, \, 0)- \lambda_i E(\bar{u}_0) \Big) r_i(\bar{u}_0), \, \tilde{r}_{i \sigma}
       \right \rangle =0.
\end{equation*}
Also, the constant $c_E$ in \eqref{eq_Binv_negative} is the same as in Hypothesis \ref{hyp_invertible}.

In conclusion, system \eqref{eq_Binv_system_trav} restricted to $\mathcal{M}^c$ can be rewritten as 
\begin{equation}
\label{eq_Binv_reduction} \left\{
       \begin{array}{lll}
             u' = v_i \tilde{r}_i (u, \, v_i, \, \sigma_i) \\
             v_i' = \phi_i(u, \, v_i, \, \sigma_i) v_i \\
             \sigma_i' =0. \\
\end{array}
\right.
\end{equation}
One actually studies the following fixed point problem, defined on a interval $[0, \, s_i]$: 
\begin{equation}
\label{eq_Binv_fixedpt} \left\{
\begin{array}{lll}
      u(\tau) = \bar{u}_0 + {\displaystyle \int_0^{\tau} \tr(u(\xi), \, v_i(\xi), \,
      \sigma_i(\xi))d \xi }  \\
      v_i (\tau) = f_i (\tau, \, u, \, v_i, \, \sigma_i ) -
      \mathrm{conc} f_i (\tau, \, u, \, v_i, \, \sigma_i )\\
      \sigma_i(\tau)=  {  \frac{1}{c_E  (\bar{u}_0) }  \displaystyle \frac{d}{d \tau}
      \mathrm{conc} f_i (\tau, \, u, \, v_i, \, \sigma_i )}. \\
\end{array}
\right.
\end{equation}
We have used the following notations:
$$
    f_i (\tau) = \int_0^{\tau} \tilde{\lambda}_i [u_i, \, v_i, \, \sigma_i] (\xi) d \xi,
$$
where 
$$
    \tilde{\lambda}_i [u_i, \, v_i, \, \sigma_i] (\xi) = \phi_i \Big( u_i (\xi), \, v_i (\xi), \, \sigma_i (\xi) \Big) + c_E (\bar{u}_0) \sigma.
$$
Also,  $\mathrm{conc} f_i$ denotes the concave
envelope of the function $f_i$: 
$$
    \mathrm{conc} f_i  (\tau)=  \inf  \{ h (s): \; \mathrm{h \; is \; concave}, \; h (y) \ge f_i (y) \; \forall \, y  \in [0, \, s_i ]\} .
$$
The link between \eqref {eq_Binv_fixedpt} and \eqref{eq_Binv_system_trav} is the following: let $(u_i, \, v_i, \, \sigma_i$ satisfy \eqref{eq_Binv_fixedpt}. Assume that $v_i < 0$ on $]a, \, b[$ and that 
$v_k (a) = v_k (b) = 0$. Define $ \alpha_i (\tau) $ as the solution of the Cauchy problem 
\begin{equation*}
\left\{ 
\begin{array}{ll}
            \displaystyle{ \frac{d \alpha_i}{ d \tau}   =  -  \frac{1}{v_i (\tau)}} \\
           \alpha ( a+ b / 2) = 0 \\
 \end{array}
\right.
\end{equation*}
 then $(u_i \circ \alpha_i, \, v_i \circ \alpha_i, \, \sigma_i \circ \alpha_i)$ is a solution to \eqref{eq_Binv_reduction} satisfying 
 $$
     \lim_{x \to - \infty} u_i \circ \alpha_i (x) = u_i (a) \qquad \lim_{x \to + \infty} u_i \circ \alpha_i (x) = u_i (b).
 $$
Thus, $ u_i (a)$ and $ u_i (b)$ are connected by a travelling wave profile.

As shown in \cite{Bia:riemann},  \eqref{eq_Binv_fixedpt} admits a unique continuous solution $(u_i, \, v_i, \, \sigma_i)$. Also,  one can show that $u_k  (0)$ and $u_k (s_i)$ are connected by a sequence of rarefaction and travelling waves with speed close to $\lambda_i (\bar{u}_0)$. If $u (t, \, x )$ is the limit of 
\begin{equation*}
\left\{
\begin{array}{ll}
          E( \ue) \ue_t + A (\ue, \, \ue_x ) \ue_x = \ee B(  \ue ) \ue_{xx} \\
          \ue (0, \,,  x) = 
          \left\{
          \begin{array}{ll}
                     \bar{u}_0  & x < 0  \\
                       u_i (s_i) & x \ge 0  \\
          \end{array}
          \right.
\end{array}
\right.
\end{equation*}
then as $\ee \to 0^+$ 
\begin{equation}
\label{e:Binv:s+:limit}
     u  (t, \, x) = 
     \left\{
     \begin{array}{ll}
                 \bar{u}_0    &  x \leq \sigma_i (0) t \\
                 u_i (s)      &  x = \sigma_i (s) t \\
                 u_i (s_i)   &  x \ge \sigma_i (s_i) t \\
     \end{array}
    \right.
\end{equation}
The $i$-th {\sl curve of admissible states} is defined setting 
$$
     T^i_{s_i} \bar{u}_0:= u_i (s_i). 
$$
 
If $s_i<0$, one considers a fixed problem like \eqref{eq_Binv_fixedpt}, but instead of the concave envelope of $f_i$ one takes the convex envelope:
$$
    \mathrm{conv} f_i  (\tau)=  \sup \{ h (s): \; \mathrm{h \; is \; convex}, \; h (y) \leq f_i (y) \; \forall \, y  \} .
$$
Again, one can prove the existence of a unique fixed point $(u_i, \, v_i, \, \sigma_i)$. 

Also, in \cite{Bia:riemann} it is proved that the curve $T^i_{s_i}$ is Lipschitz continuos with respect to $s_i$ and it is differentiable at $s_i =0$ with derivative 
given by $\vec{r}_i ( \bar{u}_0)$.  Moreover, the function $T^i_{s_i}$ is Lipschitz continuos with respect to $\bar{u}_0$.

Consider the composite function 
$$
    \psi (\bar{u}_0, \, s_1 \dots s_N) = T^1_{s_1} \circ \dots T^N_{s_N} \bar{u}_0
$$  
With the previous expression we mean that the starting point for $T^{N-1}_{s}$ is not $\bar{u}_0$ but  $T^N_{s_N} \bar{u}_0$.
Thanks to the previous steps, the map $\psi$
is Lipschitz continous with respect to $s_1 \dots s_N$ and differentiable at $s_1 = 0$, $\dots s_N = 0$. 
The column of the jacobian are $r_1 (\bar{u}_0)  \dots r_N (\bar{u}_0)$.  Thus the jacobian is invertible and hence, exploiting the extension of the implicit function theorem discussed in \cite{Cl} (page 253), the map $\phi ( \bar{u}_0, \cdot) $ is invertible in a neighbourhood of $(s_1 \dots s_N) = ( 0 \dots 0 )$. 
In other words, if $u^-$ is fixed and is sufficiently close to $\bar{u}_0$, then the values of $s_1 \dots s_N$ are uniquely determined by the equation 
\begin{equation*}
    u^- = \psi (\bar{u}_0, \, s_1 \dots s_N).
\end{equation*}
Taking the same $u^-$ as in \eqref{e:BInv:s+:app}, one obtains the parameters $(s_1 \dots s_N)$ which can be used to reconstruct the hyperbolic limit $u$ of \eqref{e:BInv:s+:app}. Indeed, once $(s_1 \dots s_N)$  are known then $u$ can be obtained  gluing together solutions like \eqref{e:Binv:s+:limit}.

\subsubsection{The hyperbolic limit in the non characteristic case}
\label{subsubsec_Binv_nonchar} 
The goal of this section is to provide a characterization of the limit of the parabolic approximation 
 \eqref{eq_Binv_the_system} in the case of a non characteristic boundary, i.e. when none of the eigenvalues of $E^{-1} (u) A( u, \, u_x)$
 can attain the value $0$. More precisely, we assume the following.
\begin{hyp}
\label{hyp_nonchar}
         Let $\lambda_1 (u) \dots \lambda_N(u)$ be the eigenvalues of the matrix $E^{-1} (u)A(u, \, 0)$. Then there exists a constant $c$ such that 
         for every $u$
           \begin{equation}
          \label{eq_Binv_nonchar_separation}
                  \lambda_1 (u) < \dots < \lambda_{n}
                 ( u) \leq  -  \frac{c}{2}  < 0 <  \frac{c}{2} \leq
                 \lambda_{n +1 }(u) < \dots <
                  \lambda_N (u).     
         \end{equation}
\end{hyp}
Thus, in the following $n$ will denote the number of eigenvalues with strictly negative
negative real part and $N-p$ the number of eigenvalues with strictly positive real part. 

\begin{figure}
\begin{center}
\psfrag{ub}{$\bar{u}_b$} \psfrag{x}{$x$} \psfrag{t}{$t$}
\psfrag{u0}{$\bar{u}_0$} \psfrag{bu}{$\bar{u}$}
\psfrag{T}{$T^N_{s_N} \bar{u}_0$} \label{fig_states_Binv_nonchar}
\caption{the solution of the boundary Riemann problem when the
viscosity is invertible and the boundary is not characteristic}
\includegraphics[scale=0.6]{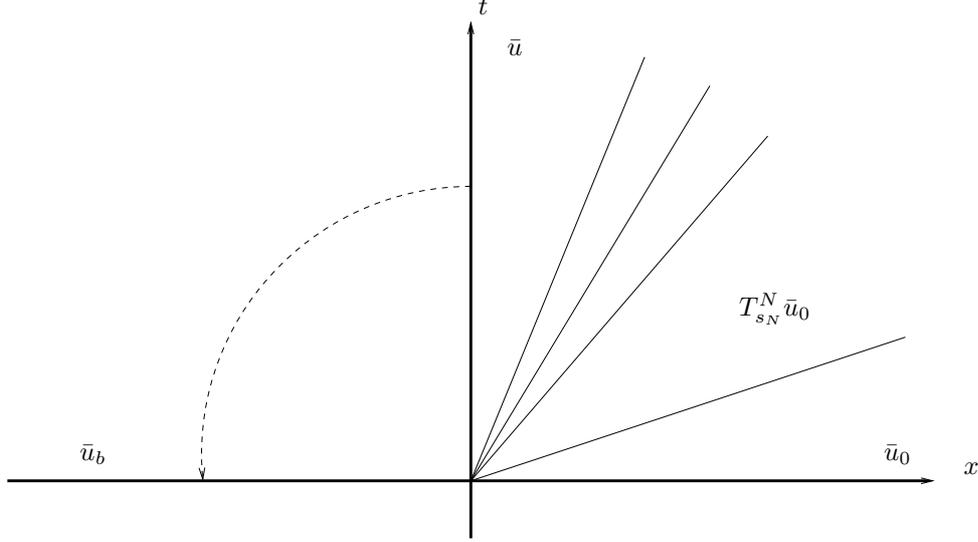}
\end{center}
\end{figure}

To give a characterization of the limit of  \eqref{eq_Binv_the_system} we will proceed as follows. We wil construct a map 
$ \phi( \bar{u}_0, s_1 \dots s_N)$ which is a {\sl boundary Riemann solver} in the sense that as $(s_1 \dots s_N)$ vary, describes states that can be connected to $\bar{u}_0$.
We will the show that the map $\phi$ is locally invertible. Hence, given $\bar{u}_0$ and $\bar{u}_b$ sufficiently close, the values of $(s_1 \dots s_N)$ are uniquely determined by the equation
$$
    \bar{u}_b = \phi ( \bar{u}_0, \, s_1 \dots s_N).
$$
Once $(s_1 \dots s_N)$ are known the limit of \eqref{eq_Binv_the_system}  is completely characterized. The construction of the map $\phi$ is divided in some steps:  
\begin{enumerate}
 \item {\sl Waves with positive speed }

Consider the Cauchy datum $\bar{u}_0$, fix $(N-k)$ parameters $(s_{N-k} \dots s_N)$
and consider the value
$$
    \bar{u} = T^{N-k}_{s_{N-k}} \circ \dots T^N_{s_N} \bar{u}_0.
$$
The curves $T^{N-k}_{s_{N-k}} \dots T^N_{s_N}$ are, as in Section  \ref{subsub_Binv_speed+}, the {\sl curves of admissible states} introduced in \cite{Bia:riemann}. 
The state $\bar{u}_0$ is then connected to $\bar{u}$ by a sequence of rarefaction and travelling waves with positive speed.

\item {\sl Boundary layers }

We have now to characterize the set of values $u$ such that the following problem admits a solution:
\begin{equation*}
\left\{
\begin{array}{lllll}
      A(U, U_x) U_x = B(U) U_{xx} \\
      U(0) = u  \\
      \lim_{x \to + \infty} U(x) = \bar{u}.
\end{array}
\right.
\end{equation*}
We have thus to study system 
\begin{equation}
\label{eq_Binv_stableman} \left\{
\begin{array}{ll}
       U_x = p \\
       p_x =B(U)^{-1} A(U, \, p) p \\
\end{array}
\right.
\end{equation}
Consider the equilibrium point $(\bar{u}, \, 0)$, linearize at that point and denote by $V^s$ the stable space, i.e. the eigenspace associated to the eigenvalues with strictly negative real part.  Thanks to Lemma
\ref{lem_Binv_dimension}, the dimension of $V^s$ is equal  to the number of negative eigenvalues of
$E^{-1}  (\bar{u}) A(\bar{u}, \, 0)$, i.e. to $n$.  Also, $V^s$ is given by
$$
    V^s = \Big\{    \big( \bar{u}  + \sum_{i = 1}^n  \frac{x_i}{  \mu_i  (\bar{u}) } \vec{\chi}_i   (\bar{u}), \,   \sum_{i = 1}^n  x_i \vec{\chi}_i  (\bar{u})  \big), \; x_1 \dots x_n \in \mathbb{R}  \Big\},
$$
 where $\mu_1  (\bar{u}) \dots \mu_n   (\bar{u})$ are the eigenvalues of $B^{-1} (\bar{u})   A(  (\bar{u}), \, 0) $ with negative real part  and $\vec{\chi}_1  (\bar{u})  \dots \vec{\chi}_n( \bar u)$ are the  
 corresponding eigenvectors. 
 
 Denote by $\mathcal{M}^s$ the stable manifold, which is parameterized by $V^s$. Also, denote by $\phi_s$ a parameterization of $\mathcal{M}^s$:
 $$
     \phi_s :   V^s \to \mathbb{R}^N. 
 $$
 Let $\pi_u$ be the projection 
\begin{equation*}
\begin{split}
     \pi_u: \,
&    \erreN \times \erreN
      \to \erreN \\
&    (u, \, p) \mapsto u 
\end{split}
\end{equation*}
If $u \in \pi_u \big( \phi_s  (s_1 \dots s_n)  \big)$ for some $s_1 \dots s_n$, then system  \eqref{eq_Binv_stableman} admits a solution. Note that thanks to classical results about the stable manifold (see eg \cite{KHass}) the map $\phi_s$ is differentiable and hence also $\pi_u \circ \phi_s$ is differentiable. In particular, the stable manifold is tangent at $s_1 =0 \dots s_n =0 $ to $V^s$ and hence  the columns of the jacobian of $\pi_u \circ \phi_s$ computed at  $s_1 =0 \dots s_n = 0$ are $\vec{\chi}_1 \dots \vec{\chi}_n(\bar u) $. 

Note that the map $\pi_u \circ \phi_s $ actually depends also on the point $\bar{u}$ and it  does in a Lipschitz continuos way:
$$
    |  \pi_u \circ \phi_s  ( \bar{u}_1, \, s_1 \dots s_n ) - \pi_u \circ \phi_s ( \bar{u}_2, \, s_1 \dots s_n ) | \leq L | \bar{u}_1 - \bar{u}_2  |
$$

\item {\sl Conclusion}

Define the map $\phi$ as follows:
\begin{equation}
\label{e:Binv:nonchra:solver}
    \phi ( \bar{u}_0, \, s_1 \dots s_n ) = \pi_u \circ \phi_s \Big(   T^{N-k}_{s_{N-k}} \circ \dots T^N_{s_N} \bar{u}_0, \, s_1 \dots s_n  \Big)
\end{equation}
From the previous steps it follows that $\phi$ is Lipschitz continuos  and that it is differentiable at $s_1 =0 \dots s_N =0$. Also, the columns of the jacobian   
are $\vec{\chi}_1(\bar{u}_0) \dots \vec{\chi}_n  (\bar{u}_0)  , \, \vec{r}_{n + 1} (\bar{u}_0)  \dots \vec{r}_N $. Thus, thanks to Lemma \ref{lem_Binv_dimension}, the jacobian is invertible. 
One can thus exploit the extension of the implicit function theorem discussed in \cite{Cl} (page 253) and conclude that the map $\phi ( \bar{u}_0, \cdot) $ is invertible in a neighbourhood of $(s_1 \dots s_N) = ( 0 \dots 0 )$. 
In particular, if one takes $\bar{u}_b$ as in \eqref{eq_Binv_the_system} and assumes that $| \bar{u}_0 - \bar{u}_b |$ is sufficiently small, then the values of $s_1 \dots s_N$ are uniquely determined by the equation 
\begin{equation}
\label{eq_Binv_noncharsolver}
      \bar{u}_b =   \phi ( \bar{u}_0, \, s_1 \dots s_n )
\end{equation}
Once the values of $s_1 \dots s_N$ are known, then the limit $u (t, \, x )$ can be reconstructed. In particular, the trace of $u$ on the axis $x =0$ is given by
\begin{equation}
\label{e:Binv:nonchar:trace}
    \bar{u} : =  T^{n+1}_{s_{n+1}} \circ \dots T^N_{s_n} \bar{u}_0.
\end{equation}
The self similar function $u$ is represented in Figure \ref{fig_states_Binv_nonchar} and can be obtained gluing together pieces like \eqref{e:Binv:s+:limit} . 
 \end{enumerate}

\begin{teo}
\label{pro_Binv_nonchar}
           Let Hypotheses \ref{hyp_invertible}, \ref{hyp_hyperbolic_I}, \ref{hyp_convergence}, \ref{hyp_small_bv}, \ref{hyp_stability} , \ref{hyp_finite_pro} and \ref{hyp_nonchar} hold. Then 
             there  exists $\delta >0$ small enough such that the following holds. If $|  \bar{u}_0 - \bar{u}_b| << \delta$, then the limit of the parabolic approximation 
             \eqref{eq_Binv_the_system} 
            satisfies 
$$
    \bar{u}_b = \phi (\bar{u}_0, \, s_1 \dots s_N)
$$
         for a suitable vector $(s_1 \dots s_N)$. The map $\phi$ is defined by \eqref{e:Binv:nonchra:solver}. Given $\bar{u}_0$ and $\bar{u}_b$, one can invert $\phi$ and determine uniquely 
         $     (s_1 \dots s_N)$. Once $(s_1 \dots s_N)$ are known one can determine a.e. $(t, \, x)$ the value $u(t, \, x)$ assumed by the limit function. In particular, the trace $\bar{u}$ of the hyperbolic limit in the axis $x =0$ is given by \eqref{e:Binv:nonchar:trace}.     
\end{teo}
\subsubsection{The hyperbolic limit in the boundary characteristic case}
\label{subsubsec_Binv_char} 
The aim of this section is to provide a characterization of the limit of the parabolic approximation 
\eqref{eq_Binv_the_system}  when the matrix $B$ is invertible, but the boundary is characteristic, i.e. one eigenvalue of 
$A (u, \, u_x)$ can attain the value $0$. 

Let $\delta$ be the bound on the total variation introduced in Hypothesis \ref{hyp_convergence}. Moreover, as in the previous section we will denote by $\lambda_1 (u) \dots \lambda_N (u)$ the eigenvalues of the matrix $E^{-1}(u) A(u, \, 0)$. The eigenvalue $\lambda_k(u)$ is Lipschitz continuous with respect to $u$ and hence there exists a suitable constant $M$ such that, if $|\lambda_k (\bar{u}_0) |> M \delta$, then $|\lambda_k (u) | > 0$ for every $u$ in a neighbourhood of $\bar{u}_0$ of size $\delta$,  $| u - \bar{u}_0 | \leq \delta$. Thus, in particular, if $\lambda_k (\bar{u}_0) > M \delta$ then Hypothesis \ref{hyp_nonchar} is satisfied for some $c$ and hence the boundary is non characteristic. In this section we will therefore make the following assumption: 
\begin{hyp}
\label{hyp_char} 
         Let $\delta$ the same constant as in Hypothesis \ref{hyp_convergence} and let $M$ the constant introduced before. Then 
         \begin{equation}
          \label{e:Binv:char:lambda:piccino}
                     |   \lambda_k  (   \bar{u}_k  )  | \leq M \delta.
          \end{equation}
\end{hyp}
Note that, because of strict hyperbolicity (Hypothesis
\ref{hyp_hyperbolic_I}) the other eigenvalues of $E^{-1}(u)A(u, \, 0)$ are well separated from zero, in the following sense: there exist a constant $c
>0$ such that for all $u$ satisfying $|\bar u_0 - u | \leq \delta$,
it holds
\begin{equation}
\label{eq_Binv_char_separation}
      \lambda_1 (u) < \dots < \lambda_{k -1} (u)
      \leq
      -c < 0 < c \leq \lambda_{k+1} (u) < \dots < \lambda_N (u).
\end{equation}
The notation used in this section is therefore the following: $(k-1)$ denotes the number of
strictly negative eigenvalues of $E^{-1}A(u, \, 0)$, while $(N-k)$ is the number of strictly positive eigenvalues. 
\vspace{1cm}

The study of the limit is more complicated in the boundary characteristic case than in the case considered in the previous section. The main ideas exploited in the analysis 
are described in the introduction, here instead we will give the technical details of the construction.  As in the previous section, the characterization of the limit will be as follows: we will construct a map $\phi( \bar{u}_0, \, s_1 \dots s_N)$ which describes all the states that can be connected to $\bar{u}_0$. We will then show that the map is  locally invertible: 
thus, given $\bar{u}_0$ and $\bar{u}_b$ sufficiently close, the values of $(s_1 \dots s_N)$ are uniquely determined by the equation
$$
    \bar{u}_b = \phi ( \bar{u}_0, \, s_1 \dots s_N).
$$
Once $(s_1 \dots s_N)$ are known one can proceed as in the previous section and determine a.e. $(t, \, x)$ the value $u(t, \, x)$ assumed by the limit function. 
The construction of the map $\phi$ will be given in several steps: 
\begin{enumerate}

\item {\sl Study of the waves with uniformly positive speed}

Consider the Cauchy datum $\bar{u}_0$, fix $(N-k)$ parameters $(s_{N-k} \dots s_N)$
and consider the value
$$
    \bar{u}_k = T^{N-k}_{s_{N-k}} \circ \dots T^N_{s_N} \bar{u}_0.
$$
The curves $T^{N-k}_{s_{N-k}} \dots T^N_{s_N}$ are, as in Section  \ref{subsub_Binv_speed+}, the {\sl curves of admissible states} introduced in \cite{Bia:riemann}. 
The state $\bar{u}_0$ is then connected to $\bar{u}_k$ by a sequence of rarefaction and travelling waves with uniformly positive speed. The speed is uniformly positive in the sense of 
\eqref{eq_Binv_char_separation}. 

\item {\sl Reduction on the center stable manifold}

Consider system 
\begin{equation}
\label{e:Binv:char:sta}
\left\{
\begin{array}{lll}
            u_x = p \\
            p_x = B(u)^{-1} \Big(    A(u, \, p)  - \sigma E(u) \Big)  p \\
           \sigma_x = 0 \\
\end{array}
\right.
\end{equation}
and the equilibrium point $\big(  \bar{u}_k, \, 0, \, \lambda (\bar{u}_k) \big)$.  Linearizing the system around $\big(\bar{u}_k, \, 0, \, \lambda(\bar{u}_k)\big)$ one obtains 
\begin{equation}
\label{e:Binv:char:jacobian}
      \left(
      \begin{array}{ccc}
            0 & I_N ^ 0 \\
            0 & B^{-1}(\bar{u}_k) \Big(  A(\bar{u}_k, \, 0) - \lambda_k (\bar{u}_k E (\bar{u}_k)   \Big) \\
            0 & 0 &`0
      \end{array}
      \right)
\end{equation}
Thanks to \eqref{eq_Binv_char_separation}, the matrix $  E( \bar{u}_k) - \lambda_k (\bar{u}_k) A( \bar{u}_k, \, 0)  $
has 1 null eigenvalue and $(N-k )$ strictly negative eigenvalues.  Because of Lemma \ref{lem_Binv_dimension} one can then conclude that the matrix 
\eqref{e:Binv:char:jacobian} has $N+1$ null eigenvalues and $(N-k)$ eigenvalues s with strictly negative real part..

Let $V^c$ be the kernel of \eqref{e:Binv:char:jacobian}, $V^s$ the eigenvspace associated eigenvalues s with strictly negative real part. and $V^u$ the eigenspace associated to eigenvalues s with strictly positive real part. There exists a so called {\sl center stable} manifold $\mathcal{M}^{cs}$ which is parameterized by $V^s \oplus V^c$: we refer to \cite{KHass} for an extensive analysis, here we will just recall some of the fundamental properties of a center stable manifold.  If we write $\mathbb{R}^{2N + 1} = V^c \oplus V^s \oplus V^u$, then the map 
$$
    \Phi_{cs}:  V^s \oplus V^c \to  V^c \oplus V^s \oplus V^u
$$
that parameterizes $\mathcal{M}^{cs}$ can be chosen in such a way that if $\pi_{cs}$ is the projection from \\ $ V^c \oplus V^s \oplus V^u$ onto $ V^c \oplus V^s$, then $\pi_{cs} \circ \Phi_{cs}$ is the identity.

Every center stable manifold is invariant with respect to \eqref{e:Binv:char:sta} and moreover satisfies the following property: let $(u^0, \, p^0, \, \sigma^0)$ belong to $\mathcal{M}^{cs}$ 
and denote by $(u(x), \, p(x), \, \sigma (x)$ the orbit starting at $(u^0, \, p^0, \, \sigma^0)$. Then 
$$
    \lim_{x \to + \infty} \Big(u(x), \, p(x), \, \sigma (x) e^{-  c x  / 2} \Big) 
    = (0, \, 0, \, 0). 
$$
The constant $c$ is the separation speed defined in \eqref{eq_Binv_char_separation}. We are also assuming that $\lambda_k (\bar{u}_k)$ is bigger than $- c / 4$: this is a consequence of Hypothesis 
\ref{hyp_char} provided that $\delta$ is sufficiently small. A center stable manifold $\mathcal{M}^{cs}$ is defined in a neighbourhood of the equilibrium point $(\bar{u}_k, \, 0, \, \lambda(\bar{u}_k)$ and it is tangent to $V^c \oplus V^s$ at $(\bar{u}_k, \, 0, \, \lambda(\bar{u}_k)$.  

In the following, we will fix a center stable manifold $\mathcal{M}^{cs}$ and we will focus on the solutions of \eqref{e:Binv:char:sta} 
that lay on $\mathcal{M}^{cs}$. The center and the stable space of \eqref{e:Binv:char:jacobian} are given respectively by
$$
    V^c = \Big\{   (u, v^k_{cs}  \vec{r}_k (\bar{u}_k), \, \sigma   ):   \;  u \in \mathbb{R}^N, \;   v^k_{cs}, \; \sigma \, \in  \mathbb{R}  \Big\}  
$$
and 
$$
    V^s = \Big\{   (\bar{u}^k   +  \sum_{ i < k} \frac{v^i_{cs}}{\mu_i}  \vec{\chi}_i (\bar{u}_k) ,   \sum_{ i < k} v^i_{cs}  \vec{\chi}_i (\bar{u}_k), \, \lambda_k (\bar{u}_k)   ) \Big\} . 
$$
In the previous expression, $r_k (\bar{u}_k)$ is a unit vector in the kernel of $\Big( A(\bar{u}_k, \, 0) - \lambda_k (\bar{u}_k) \Big)$, while $\vec{\chi}_1 \dots \vec{\chi}_{k-1}$ are the eigenvectors of $B^{-1} (\bar{u}_k) \Big(   A(\bar{u}_k, \, 0) - \lambda_k (\bar{u}_k)    \Big)$ associated to the eigenvalues $\mu_1 \dots \mu_{k-1}$ s with strictly negative real part. Since $\Phi_{cs} \circ \pi_{cs}$ is the identity, then a point $(u, \, p, \, \sigma)$ belongs to the  manifold $\mathcal{M}^{cs}$ is and only if 
$$
    p =  \sum_{ i < k} v^i_{cs} \vec{\chi}_i (\bar{u}_k) + v^k_{cs} \vec{r}_k + \sum_{i > k}  \phi^i_{cs} (u, \, v^1_{cs}  \dots v^k_{cs}, \, \sigma ) \vec{\chi}_i (\bar{u}_k),
$$
where $\vec{\chi}_{k +1} (\bar{u}_k) \dots \vec{\chi}_N (\bar{u}_k)$ are the eigenvectors of $B^{-1} (\bar{u}_k) \Big(   A(\bar{u}_k, \, 0) - \lambda_k (\bar{u}_k)    \Big)$ associated to eigenvalues with strictly positive real part.. Since all the equilibrium points $(u, \, 0, \, \sigma)$ belong to the center stable manifold, then we must have $p=0$ when $v^1_{cs} \dots v^k_{cs}$ are all $0$. Thus, since $\vec{\chi}_{k +1} (\bar{u}_k) \dots \vec{\chi}_N (\bar{u}_k)$  are linearly independent,     
 $\phi^{i}_{cs}(u, \, \vec{0}, \, \sigma ) =0$ for all
$i= k + 1\dots N$ and for every $u$ and $\sigma$. This means, in particular, that for every $i=k+1 \dots N$ there exists a $k$-dimensional vector
$\tilde{\phi}^i_{cs}(u, \, v_{cs}^1 \dots p_{cs}^k, \, \sigma)$ such that
\begin{equation*}
      \phi^{i}_{cs}(u, \, v_{cs}^1 \dots v_{cs}^k) = \langle
      \tilde{\phi}_{cs}^{i}(u, \, v_{cs}^1 \dots v_{cs}^k, \, \sigma ), \,
       v_{cs}^1 \dots v_{cs}^k  \rangle.
\end{equation*}
Define 
\begin{equation*}
   R_{cs} (u, \,  v_{cs}^1 \dots v_{cs}^k, \, \sigma ) : = 
  \Big(
           \vec{\chi}_1 |
           \dots |
           \vec{\chi}_{k -1} | \vec{r}_k
      \Big)
       +
      \Big(
           \vec{\chi}_{k+1} |
           \dots |
           \vec{\chi}_N
      \Big)
      \left(
      \begin{array}{ccc}
             \tilde{\phi}_{cs}^{k+1}(u, \, v_{cs}^1 \dots v_{cs}^k, \, \sigma ) \\
             \hline
             \dots \\
             \hline
                  \tilde{\phi}_{cs}^{N}(u, \, v_{cs}^1 \dots v_{cs}^k, \, \sigma )
      \end{array}
      \right),
\end{equation*}
then a point $(u, \, p, \, \sigma)$ belongs to $\mathcal{M}^{cs}$ if and only if
\begin{equation}
\label{e:Binv:char:pcs}
    p = R_{cs}  (u, \, V_{cs}, \, \sigma ) V_{cs},
\end{equation}
where $V_{cs}$ denotes  $ ( v_{cs}^1 \dots v_{cs}^k )^T$. Since $\mathcal{M}^{cs}$ is tangent to $V^c \oplus V^s$ at $(\bar{u}_k, \, \vec{0}, \, \lambda_k(\bar{u}_k))$, then 
$$
    R_{cs} (\bar{u}_k, \, \vec{0}, \, \lambda_k(\bar{u}_k)) =   \Big(
           \vec{\chi}_1 |
           \dots |
           \vec{\chi}_{k -1} | \vec{r}_k
      \Big). 
$$
Plugging  $ p = R_{cs}  (u, \, V_{cs}, \, \sigma ) V_{cs} $ into system \eqref{e:Binv:char:stable} one gets 
\begin{equation}
\label{eq_app_derivative}
      \bigg( (D_u R_{cs}) R_{cs} V_{cs} +
             (D_{v}R_{cs} ) V_{cs x}
      \bigg) V_{cs} + R_{cs} V^{cs}_x = B^{-1} \Big(  A  (u, \, R_{cs} V_{cs}) - \sigma E(u)  \Big)R_{cs}
      V_{cs}.
\end{equation}
Moreover, the matrix
\begin{equation*}
      D(u, \, V_{cs}) = (d_{i h}) \qquad d_{ih} = \sum_{j=1}^{k}
      \frac{\partial R_{cs ij}}{\partial V_h} V_j
\end{equation*}
satisfies
\begin{equation*}
      D(u, \, 0)=0 \qquad \Big( (D_{v}R_{cs} ) V_{csx}
      \Big) V_{cs} = D(u, \, V^{cs}) V^{cs}_x.
\end{equation*}
Hence \eqref{eq_app_derivative} can be rewritten  as
\begin{equation}
\label{eq_app_derivative2}
       \bigg( R_{cs} + D(u, \, V_{cs}) \bigg) V_{csx} =
       \bigg( B^{-1}  \Big(  A   -  \sigma E \Big)    - (D_u R_{cs}) R_{cs} V_{cs}
       \bigg)  R_{cs} V_{cs}.
\end{equation}
Since 
$$
     \bigg( R_{cs} + D(u, \, V_{cs}) \bigg) (\bar{u}_k, \, \vec{0}, \, \lambda_k(\bar{u}_k)) =    \Big(
           \vec{\chi}_1 |
           \dots |
           \vec{\chi}_{k -1} | \vec{r}_k
      \Big),
$$
then the columns of $R_{cs} + D(u, \, V_{cs}) $ are linearly independent in a small enough neighbourhood. Thus one can find a matrix $L(u, \, V_{cs}, \, \sigma)$ such that  
$$
    L(u, \, V_{cs}, \, \sigma)   \bigg( R_{cs} + D(u, \, V_{cs}) \bigg) = I_k. 
$$
Multiplying
\eqref{eq_app_derivative2} by $L(u, \, V_{cs}, \, \sigma )$ one finds
\begin{equation}
\label{e:Binv:char:vcsx}
       V^{cs}_x = \bigg(L B^{-1}   ( A - \sigma E) R_{cs} - L
       (D_u R_{cs} )R_{cs} V_{cs}
       \bigg)V_{cs}.
\end{equation}
Define 
\begin{equation}
\label{e:Binv:char:lambdacs}
       \Lambda_{cs}(u, \, V_{cs}, \, \sigma ):=  L B^{-1}   ( A - \sigma E) R_{cs} - L
       (D_u R_{cs}) R_{cs} V_{cs},
\end{equation}
then by construction $  \Lambda_{cs}(\bar{u}_k, \, \vec{0}, \, \lambda_k (\bar{u}_k) )$ is a $(k \times k)$ diagonal matrix with one eigenvalue equal to $0$ and $N-k$ eigenvalues
with strictly negative real part.
Also,  if $\vec{\ell}_k = (0 \dots 0, \, 1)$ and $\vec{e}_k = \vec{\ell}_k^T$, then 
\begin{equation}
\label{e:Binv:char:Lambdacssigma}
         \vec{ \ell}_k    \Bigg( \frac{\partial \Lambda_{cs} }{\partial \sigma } \Bigg|_{u = \bar{u}_k, \, V_{cs} = \vec{0}, \, \sigma = \lambda_k (\bar{u}_k ) }  \Bigg) \vec{e}_k = 
         - \bar{k}  <0
\end{equation} 
for a suitable strictly positive constant $\bar k$. This is a consequence of conditions 1 and 3 in Hypothesis \ref{hyp_invertible}.

In conclusion, the solutions of system \eqref{e:Binv:char:sta} laying on $\mathcal{M}^{cs}$ satisfy
\begin{equation}
\label{e:Binv:char:reduced}
\left\{
\begin{array}{lll}
           u_x =  R_{cs} (u, \, V_{cs}, \, \sigma ) V_{cs}  \\
           V_{cs x} =      \Lambda_{cs}(u, \, V_{cs}, \, \sigma ) V_{cs} \\
           \sigma_x =0
\end{array}
\right.
\end{equation}
where $    \Lambda_{cs}(u, \, V_{cs}, \, \sigma )$ is defined by \eqref{e:Binv:char:lambdacs}.
\item {\sl Analysis of the uniformly stable component of \eqref{e:Binv:char:reduced}}

Linearizing system \eqref{e:Binv:char:reduced} around the equilibrium point $(\bar{u}_k, \, \vec{0}, \, \lambda_{k} (\bar{u}_k) )$ one obtains the matrix
\begin{equation}
\label{e:Binv:char:jaccs} 
\left(
\begin{array}{lll}
            0 &      R_{cs} (\bar{u}_k, \, \vec 0, \, \lambda_k (\bar{u}_k ) ) & 0 \\
            0  &    \Lambda_{cs}(\bar{u}_k, \, \vec{0}, \, \lambda_k (\bar{u}_k) )  & 0\\
            0 & 0 & 0 \\
\end{array}
\right),
\end{equation}
which has $k-1$ distinct eigenvalues with strictly negative real part and the eigenvalue $0$ with multiplicity $N+1$. Also, the manifold 
$$
    E : = \Big\{   (u, \, \vec{0}, \, \sigma ), \; u \in \mathbb{R}^N, \; \sigma \in \mathbb{R} \Big\}
$$ 
is entirely constituted by equilibria.

Then there exists a {\sl uniformly stable manifold} $\mathcal{M}^{us}_{E}$ which is invariant for \eqref{e:Binv:char:reduced} and is characterized by the following property: given $(u^0, \, V^0_{cs}, \, \sigma^0)$ belonging to $\mathcal{M}^{us}_E$, denote by $(u(x), \, V_{cs}(x), \, \sigma  (x))$ the solution of \eqref{e:Binv:char:reduced} starting at  $(u^0, \, V^0_{cs}, \, \sigma^0)$. Then there exists a point $(u^{\infty}, \, \vec 0, \, \sigma^{\infty})$ belonging to $E$ such that 
$$
    \lim_{x \to + \infty}  \Big(   |  u(x)  - u^{\infty}|   + |V_{cs}  (x) |   + |\sigma(x) - \sigma^{\infty}  |\Big) e^{cx / 2}  = 0. 
$$
We remark that the uniformly stable manifold $\mathcal{M}^{cs}_{E}$ is contained but in general does not coincide with a given center stable manifold. Indeed, $\mathcal{M}^{us}_E$ does not contain, for example, the trajectories that converge to a point in $E$ with speed slower than $e^{- cx }$. 

The existence of the uniformly stable manifold is implied by
Hadamard-Perron theorem, which is discussed for example in  \cite{KHass}. 
The application of the uniformly stable manifold to the study of hyperbolic initial boundary values problems is derived from \cite{AnBia}.  In our case, the existence of the uniformly stable manifold is loosely speaking guaranteed by the strict hyperbolicity of $A$ (Hypothesis \ref{hyp_hyperbolic_I}), namely the fact that the eigenvalues of $A$ are uniformly separated one from each other.   
In the following, we will just recall some of the properties of the uniformly stable manifold. 

Let 
$$
    \tilde{V}^s = \Big\{    \Big(\bar{u}_k +  \sum_{i < k}  v^i_{s} \vec{\chi}_{i} (\bar{u}_k), \, \sum_{i =1}^k v^i_s  \vec{e}_i, \, \lambda_k (\bar{u}_k)   \Big), \; v^1_s \dots v^{k-1}_s \in \mathbb{R}  \Big\},
$$
where $\vec{\chi}_1$ are as before the eigenvectors of $B^{-1} ( A - E )$ and $\vec{e}_i \in \mathbb{R}^k$ are the vectors of the canonical basis 
of $\mathbb{R}^k$, $\vec{e}_1 = (1, \, 0 \dots 0)$ and so on. 
The uniformly stable manifold $M^{cs}_E$ is parameterized by $E \oplus \tilde{V}^s$ and hence has dimension $N+ k-1$. The parameterization 
$$
    \Phi_{us}: \tilde{V}^s \oplus E \to \mathbb{R}^{k+2}
$$
can be chosen in such a way that the following property is satisfied. Let 
$$
    \tilde{V}^c = \Big\{     \Big(    \bar{u}_k  +v^k \vec{r}_{k} (\bar{u}_k), \, v^k  \vec{e}_k, \, \lambda_k (\bar{u}_k)   \Big), \; v^k \in \mathbb{R}  \Big\}.
$$
If we write $\mathbb{R}^{N + k +1}$ as $\tilde{V}^c \oplus \tilde{V}^s \oplus E$ and we denote by $\pi_{us}$ the projection onto $\tilde{V}^s \oplus E$, then we can have $\pi_{us} \circ \Phi_{us}$ equal to the identity. 

By considerations analogous, but easier, to those performed to get \eqref{e:Binv:char:pcs} one deduces that a point $(u, \, V_{cs}, \, \sigma)$ belongs to $\mathcal{M}^{us}_E$ if and only if 
$$
    V_{cs} = R_s (u, \, V_s, \, \sigma) V_s,
$$
where $V_s = ( v^1_{us} \dots v^{k -1}_{us})$ are the components of $V_{cs}$ along $\vec{e}_1 \dots \vec{e}_{k-1}$. The matrix $R_s $ belongs to $ \mathbb{M}^{k \times (k -1)}$ and 
$$
     R_s \Big( \bar{u}_k, \, \vec{0}, \, \lambda_k (\bar{u}_k) \Big) = \Big( \vec{e}_1 | \dots | \vec{e}_{k -1} \Big).  
$$
We want now to compute $\Lambda_{cs} R_s V_s$. Plugging the relation $V_{sx} = \Lambda_{cs} V_s$ into \eqref{e:Binv:char:reduced} one gets 
\begin{equation}
\label{e:Binv:char:vsx}
      \Big[   R_s + D_s \Big] V_{sx} = \Big[    \Lambda_{cs} R_s -   D_u R_s R_{cs}  R_s V_s   \Big] V_s,
\end{equation}
where $D_s $ is a matrix such that $D_s ( u, \, \vec{0}, \, \sigma)$ for every $u$ and $\sigma$. Thus, 
$$
   V_{sx} = \tilde{\Lambda}_s V_s 
$$
for a suitable matrix $\tilde{\Lambda}_s$ such that $\tilde{\Lambda}_s \Big(\bar{u}_k, \, \vec{0}, \, \lambda_k (\bar{u}_k \Big) $ is a diagonal matrix with all the eigenvalues with strctly negative real part. Plugging back into \eqref{e:Binv:char:vsx} one gets 
\begin{equation}
\label{e:Binv:char:lambdars}
        \Lambda_{cs} R_s V_s = \Big\{      \Big[   R_s + D_s \Big]  \tilde{ \Lambda}_s +    D_u R_s R_{cs}  R_s V_s            \Big\}  V_s. 
\end{equation}

\item {\sl Analysis of the center component of \eqref{e:Binv:char:reduced} }

Linearizing system \eqref{e:Binv:char:reduced} around the equilibrium point $( \bar{u}_k, \, \vec{0}, \, \lambda_k (\bar{u}_k))$ one obtains \eqref{e:Binv:char:jaccs}
and hence the center space is given by
$$
    \tilde{V}^c = \Big\{ ( u, \, v_k \vec{e}_k, \, \sigma): \;  u \in \mathbb{R}^N, \; v_k, \, \sigma \in \mathbb{R} \, \Big\}.
$$
Consider the center manifold $\mathcal{M}^c$, parameterized by $\tilde{V}^c$: thanks to considerations 
similar, but easier, than those performed in the previous steps one gets that 
a point $(u, \, V_{cs}, \, \sigma)$ belongs to $\mathcal{M}^c$ if and only if 
$$
    V_{cs} = \check{r}_k (u, \, v_k, \, \sigma )  v_k,
$$
where $\check{r}_k \in \mathbb{R}^k$ and $  \check{r}_k ( \bar{u}_k, \, \vec{0}, \, \lambda_k (\bar{u}_k) ) = \vec{e}_k$.

Again, with considerations analogous, but easier than those performed at the previous step one gets 
\begin{equation}
\label{e:Binv:charcorr}
    \Lambda_{cs} \check{r}_k  v_k= \Big[    (\check{r}_k + \check{r}_{kv} v_k  ) \tilde{\phi}_k   + D_u \check{r}_k R_{cs} \check{r}_k v_k    \Big] v_k
\end{equation}
for a suitable function $\tilde{\phi}_k $ satisfying $\tilde{\phi}_k \Big( \bar{u}_k, \, \vec{0}, \, \lambda_k (\bar{u}_k ) \Big) = \vec{0}.$ Also, thanks to \eqref{e:Binv:char:Lambdacssigma},
\begin{equation}
\label{e:Binv:char:phisigma}
          \frac{\partial \tilde{\phi}_k }{\partial \sigma}     \Bigg|_{u = \bar{u}_k, \, V_{cs} = \vec{0}, \, \sigma = \lambda_k (\bar{u}_k  }   = - \bar{k} < 0.    
\end{equation}

\item {\sl Decomposition of \eqref{e:Binv:char:reduced} in center and uniformly stable component }

In this step, we will fix a trajectory $(u, \, V_{cs}, \, \sigma)$ of \eqref{e:Binv:char:reduced} and we will decompose it in a {\sl uniformly stable} and in a {\sl center} component, in the following sense. We exploits the manifolds $\mathcal{M}^{c}$ and $\mathcal{M}^{us}_E$ introduced in the previous steps and we decompose
\begin{equation}
\label{e:Binv:char:vcs}
    V_{cs} = R_s (u, \, V_s, \, \sigma) V_s + \check{r}_k (u, \, v_k, \, \sigma )v_k.
\end{equation}
Plugging this expression into the first line of \eqref{e:Binv:char:reduced} one gets 
$$
    u_x = R_{cs} R_s V_s + R_{cs} \check{r}_k v_k.
$$
From the second line of \eqref{e:Binv:char:reduced} one gets 
\begin{equation*}
\begin{split}
            V_{csx}
&             = 
            R_s V_{sx} + \Big[  D_u R_s (  R_{cs} R_s V_s + R_{cs} \check{r}_k v_k    )  + D_V R_s V_{sx }  \Big]  V_s \\
&         \quad    +  \check{r}_k v_{k x} + \Big[  D_u  \check{r}_k (  R_{cs} R_s V_s + R_{cs} \check{r}_k v_k ) + \check{r}_{kv} v_{kx} \Big] v_k \\
&          =  \Lambda_{cs}  R_s V_s + \Lambda_{cs} \check{r}_k v_k.   
\end{split}
\end{equation*}
One can prove that
$$
    D_V R_s V_{sx }   V_s = D_s (u, \, V_s, \, \sigma ) V_{sx} 
$$
for a suitable matrix $D_s$ which satisfies $D(u, \, \vec{0}, \, \sigma ) = 0$ for every $u$ and $\sigma$ (the same as in the previous step). 
Thus, exploiting \eqref{e:Binv:char:lambdars} and \eqref{e:Binv:charcorr}, we get
\begin{equation}
\label{e:Binv:char:diag}
\begin{split}
         \Big[   R_s +  D_s   \Big] V_{sx} + \Big[  \check{r}_k +   \check{r}_{kv} v_k  \Big] v_{k x }  
&         = 
            \Lambda_{cs}R_s V_s + \Lambda_{cs} \check{r}_k v_k  \\
&          \quad    - 
             D_u R_s (  R_{cs} R_s V_s + R_{cs} \check{r}_k v_k    )  V_s -  
              D_u  \check{r}_k (  R_{cs} R_s V_s + R_{cs} \check{r}_k v_k )  v_k  \\
&           = \Big\{      \Big[   R_s + D_s \Big]  \tilde{ \Lambda}_s +    D_u R_s R_{cs}  R_s V_s            \Big\}  V_s \\
&           \quad + \Big[    (\check{r}_k + \check{r}_{kv} v_k  ) \tilde \phi_k   + D_u \check{r}_k R_{cs} \check{r}_k v_k    \Big] v_k \\
&           \quad   - 
             D_u R_s (  R_{cs} R_s V_s + R_{cs} \check{r}_k v_k    )  V_s -  
              D_u  \check{r}_k (  R_{cs} R_s V_s + R_{cs} \check{r}_k v_k )  v_k \\ 
&           =     \Big[   R_s + D_s \Big]  \tilde{ \Lambda}_s V_s +    (\check{r}_k + \check{r}_{kv} v_k  ) \tilde \phi_k  v_k \\
&           \quad   - 
             \Big[ D_u R_s   R_{cs} \check{r}_k  +   
              D_u  \check{r}_k   R_{cs} R_s \Big] V_s   v_k            \phantom{\Big[}  \\             
\end{split}
\end{equation}
The functions $R_s + D_s$ and $\check{r}_k + \check{r}_k v_k$ satisfy
$$
     (R_s + D_s  )\Big( \bar{u}_k, \, \vec{0}, \, \lambda_k (\bar{u}_k) \Big) = \Big( \vec{e}_1 | \dots | \vec{e}_{k -1} \Big) \qquad 
      (\check{r}_k +  \check{r}_k v_k ) \Big( \bar{u}_k, \, \vec{0}, \, \lambda_k (\bar{u}_k) \Big) = \vec{e}_k
$$
and hence in a neighbourhood of $(   \bar{u}_k, \, \vec{0}, \, \lambda_k (\bar{u}_k)  )$ the columns of $R_s + D_s $ and $\check{r}_k+  \check{r}_k v_k  $ are all linearly independent. 
Denote by $\ell_1 \dots \ell_k$ the vectors of the dual basis and define 
\begin{equation*}
L_s (u, \, V_s, \, v_k, \, \sigma) : =
\left(
\begin{array}{ccc}
            \ell_1 \\ 
            \dots  \\
            \ell_{k -1} \\
\end{array}
\right)
\end{equation*}
Then multiplying \eqref{e:Binv:char:diag} on the left by $L_s$ one gets 
$$
    V_{sx} = \Lambda_s (u, \, v_k, \, V_s, \sigma ) V_s,  
$$
where  $\Lambda_s \in \mathbb{M}^{k-1 \times k-1}$ is given by 
$$
    \Lambda_s (u, \, v_k, \, V_s, \sigma ) =  L_s \Big[   R_s + D_s \Big]  \tilde{ \Lambda}_s - L_s  \Big[ D_u R_s   R_{cs} \check{r}_k  +   
              D_u  \check{r}_k   R_{cs} R_s \Big]   v_k . 
$$
The matrix $ \Lambda_s   (   \bar{u}_k, \, 0, \, \vec{0}, \, \lambda_k (\bar{u}_k)  )  $
is diagonal and all the eigenvalues have strictly negative real part. 

Also, multiplying \eqref{e:Binv:char:diag} on the left by $\ell_k$ one gets 
$$
     v_{kx} = {\phi}_k v_k,
$$
where the real valued function $  {\phi}_k  ( u, \, v_k, \, V_s, \sigma    ) $ is given by 
$$
    {\phi}_k  ( u, \, v_k, \, V_s, \sigma    ) = \ell_k   (\check{r}_k + \check{r}_{kv} v_k  ) \tilde \phi_k  - \ell_k 
             \Big[ D_u R_s   R_{cs} \check{r}_k  +   
              D_u  \check{r}_k   R_{cs} R_s \Big] V_s .  
$$
Also, ${\phi}_k  (   \bar{u}_k, \, 0, \, \vec{0}, \, \lambda_k (\bar{u}_k)  ) =0$. 

In conclusion, system \eqref{e:Binv:char:reduced} can be decomposed as 
\begin{equation}
\label{e:Binv:char:sys1}
\left\{
\begin{array}{llll}
            u_x = R_{cs} (u, \, v_k \check{r}_k + R_s V_s, \, \sigma) \check{r}_k v_k + R_{cs} (u, \, v_k \check{r}_k + R_s V_s, \, \sigma) R_s V_s \\
            v_{kx} = {\phi}_k (u, \, v_k, \, V_s, \, \sigma) v_k \\
                    V_{sx} = \Lambda_s (u, \, v_k, \, V_s, \, \sigma)  V_s    \\    
            \sigma_{x} =0 \\ 
\end{array}
\right.
\end{equation}
Note that, thanks to \eqref{e:Binv:char:phisigma},
\begin{equation}
\label{e:BInv:char:partialsigma}
   \frac{\partial   \phi_k  \big(   u_k     , \, v_k   , \, V_s, \, \sigma_k  \big)    }{\partial \sigma }   \bigg|_{  u=  \bar{u}_k, \, v_k = 0, \,  V_s = \vec 0, \, \sigma_k = \lambda_k (\bar{u}_k)    } =  - \bar{k}   < 0. 
\end{equation}

To study system \eqref{e:Binv:char:sys1} we will first consider the solution of \eqref{e:Binv:char:sys1} when $u_s \equiv 0$ and $V_s \equiv 0$. Fix $s_k > 0$: we will actually study the following fixed point problem, which is defined on the interval $[0, \, s_k]$: 
\begin{equation}
\label{eq:Binv:char:centerpuret} 
\left\{
\begin{array}{lll}
      u_k(\tau) = \bar{u}_k + {\displaystyle \int_0^{\tau} \check{r}_k (u_k(\xi), \, v_k(\xi), \, 0, \, 
      \sigma_k(\xi))d \xi }  \\
      v_k (\tau) = f_k [u_k, \, v_k, \, \sigma_k ] (\tau)-
      \mathrm{mon}_{[0, \, s_k]} f_k [ u_k, \, v_k, \, \sigma_k ] (\tau)   \\
      \sigma_k(\tau)= {\displaystyle \frac{1}{\bar{c}}  \frac{d}{d \tau}
      \mathrm{mon}_{[0, \, s_k]} f_k  (\tau, \, u_k, \, v_k, \, \sigma_k )}. \\
\end{array}
\right.
\end{equation}
In the previous expression, 
\begin{equation}
\label{e:Binv:char:f}
     f_k  [ u_k, \, v_k, \, \sigma_k ] (\tau) = \int_0^{\tau} \tilde{\lambda}_k [u_k, \, v_k, \, \sigma_k ] (\xi)   d \xi  
\end{equation}
where
\begin{equation}
\label{e:BInv:char:lambda}
     \tilde{\lambda}_k [u_k, \, v_k, \, \sigma_k ] (\xi)  = \phi_k  \Big(   u_k  (\xi)   , \, v_k (\xi)   , \, 0, \, \sigma_k (\xi)  \Big)   +  \bar{k} \sigma_k (  \xi  ).
\end{equation}
The constant $\bar{k}$ is defined by \eqref{e:BInv:char:partialsigma}. 

If $s_k<0$, one considers a fixed problem like \eqref{eq:Binv:char:centerpuret} but instead of the monotone concave envelope of $f_i$ one takes the monotone convex envelope:
$$
    \mathrm{monconv} f_i  (\tau)=  \sup \{ h (s): \; \mathrm{h \; is \; convex \; and \; monotone \; non \; decreasing}, \; h (y) \leq f_i (y) \; \forall \, y \in [0, \, s_k ] \} .
$$
In the following we will consider only the case $s_k >0$, the case $s_k < 0$ being entirely similar.

The link between system \eqref{eq:Binv:char:centerpuret} and system \eqref{e:Binv:char:sys1} in the case $u_s \equiv 0$, $V_s \equiv \vec{0}$ is the following. 
Let $(u_k, \, v_k, \, \sigma_k)$ solve \eqref{eq:Binv:char:centerpuret},   assume that $v_k (s_k)  < 0$ and define 
\begin{equation}
\label{e:Binv:char:alpha}
          \alpha (\tau )  =    -   \int_{\tau}^{s_k}  \frac{1}{   v_k (s) }   ds.
\end{equation} 
Let 
$$
    \underline{s}: = \max \{ \tau \in [0, \, s_k ]:   \:  v_k (\tau)  =  0\},
$$
in the following we will prove that  $\alpha (\tau) < + \infty $ if and only if $\tau > \underline{s}$. 

The function  $( u_k \circ  \alpha, \, v_k \circ \alpha, \, 0)$ solves  \eqref{e:Binv:char:sys1} in the case $u_s \equiv 0$, $V_s \equiv \vec{0}$. If $v_k (s_k) = 0$ then $(u_k(s_k), \, 0, \, 0)$ is also a trivial solution of  \eqref{e:Binv:char:sys1} in the case $u_s \equiv 0$, $V_s \equiv \vec{0}$ and the following properties are satisfied.  Note that the function $u_k \circ  \alpha$
is a steady solution of the original parabolic equation since $\sigma \equiv 0$:
$$
    A(u_k, \, u_{kx} ) u_{kx} = B (u_k) u_{k xx}.
$$ 
Moreover,
\begin{equation}
\label{e:Binv:char:limit}
         \lim_{   x \to \infty} ( u_k \circ \alpha ) (x) = u_k (  \underline{s} ).
\end{equation}
From  system \eqref{eq:Binv:char:centerpuret} we also get that  $\bar{u}_k$ is connected to $u_k (\underline{s})$ by a sequence of rarefaction waves and shocks with non negative speed.

After  finding a solution of \eqref{eq:Binv:char:centerpuret}  in the case $u_s \equiv 0$, $V_s \equiv 0$ we will also find  a solution of \eqref{eq:Binv:char:centerpuret}  in the case $u_k \equiv \bar{u}_k$, $v_k \equiv 0$ and $\sigma \equiv 0$. We will also impose 
 $$
         \lim_{   x \to \infty} u_s (x) =0
 $$
and hence we will study the fixed point problem 
\begin{equation}
\label{e:Binv:char:stable}
   \left\{
    \begin{array}{ll}
           u_{s} (x) =  \displaystyle{       -  \int_x^{+   \infty   }     \check{R}^s  \Big(  \bar{u}_k + u_s (y) , \, 0, \, V^s  (y)  \Big) V^s (y)   dy    } \\
           V_{s } (x)=   \displaystyle{     e^{    \bar{\Lambda} x         }  V^s   (0)   +    \int_0^x   e^{   \bar{\Lambda} (x - y)    } 
           \Big[  \check{\Lambda}^s   \Big(  u_s (y)  + \bar{u}_k, \, 0, \, V^s  (y) \Big) - \bar{\Lambda} \Big] V^s (y)   dy   },
    \end{array}
    \right.
\end{equation}
where 
$$
     \bar{\Lambda} =    \check{\Lambda}^s   \Big(   \bar{u}_k, \, 0, \, 0  \Big).   
$$
Note that again, being $\sigma= 0$,  $u_s$ provides a steady solution of the origin parabolic equation,
$$
     A(u_s, \, u_{sx}) u_{sx} = B(u_s) u_{sxx} 
$$

Finally, we will consider a component of perturbation, due to interaction between the purely center component, defined by  \eqref{eq:Binv:char:centerpuret} and the purely stable component, which satisfies \eqref{e:Binv:char:stable}. More precisely, we will define $(U, \, q, \, p)$ in such a way that $u= U + u_k \circ \alpha+ u_s $, $v_k + q$ and $V_s + p$ is a solution of 
\eqref{e:Binv:char:sys1}. 
\item {\bf Analysis of the purely center component}

The purely center component is the solution of system \eqref{eq:Binv:char:centerpuret}.  
\begin{lem}
\label{l:centercont}
           Fix $\delta >0$ such that $s_k \leq  \delta  < < 1  $. Then system \eqref{eq:Binv:char:centerpuret} admits a unique solution $(u_k, \, v_k, \, \sigma_k)$
           satisfying 
           \begin{enumerate}
           \item $u_k \in \mathcal{C}^0 (   [0, \,    s_k   ]    )$,  $\|  u_k   -   \bar{u}_k        \|_{   \mathcal{C}^0     }    \leq c_u \delta    $, $u_k$ is Lipschitz with constant $Lip  (u_k)  \leq 2 $.
           \item $v_k \in \mathcal{C}^0 (   [0, \,    s_k   ]    )$,  $\| v_k        \|_{   \mathcal{C}^0     }    \leq c_v   \delta    $, $v$ is Lipschitz with constant $Lip  (v_k)  \leq   \tilde{c}_v  \delta $.
           \item $\sigma_k  \in \mathcal{C}^0 (   [0, \,    s_k   ]    )$,  $\|  \sigma_k    -   \lambda_k  (   \bar{u}_k    )        \|_{   \mathcal{C}^0     }    \leq c_{\sigma}   \delta  / \eta_k  $, $\sigma_k$ 
           is Lipschitz  with constant $Lip  (\sigma_k)  \leq \tilde{c}_{\sigma} $.
           \end{enumerate}
           The constants $c_u$, $c_v$, $\tilde{c}_v$, $c_{\sigma}$, $\eta_k$ and $\tilde{c}_{\sigma}$ do not depend on $\delta$.  The function $f$ defined by \eqref{e:Binv:char:f} belongs to $\mathcal{C}^{1, \, 1}_k  (    [0, \, s_k   ])$ for a suitable constant $k$ which again does not depend on $\delta$.
           
\end{lem}
The symbol $\mathcal{C}^{1, \, 1}_k  (    [0, \, s_k   ])$ denotes the space of functions $ f \in \mathcal{C}^{1}  (    [0, \, s_k   ]) $ such that $f'$ is Lipschitz continuous and has Lipschitz constant
$Lip(f'   )  \leq k$.

\begin{proof}
The proof of the lemma relies on a fixed point argument. Only the fundamental steps are sketched here. 

Define 
\begin{equation}
\label{e:Binv:char:k:spaces}
\begin{split}
&     X_{u k} : = \big\{  u_k  \in \mathcal{C}^0 (   [0, \, s_k], \;   \mathbb{R}^{N}   ):  \; \|  u_k -  \bar{u}_k  \|_{\mathcal{C}^0}   
       \leq c_u \delta,    \;  u_k   \; \mathrm{is \;  Lipschtiz}, \;   Lip  (u_k)  \leq 2 \big\}, \\
&     X_{v k}  : = \big\{  u_k  \in \mathcal{C}^0 (   [0, \, s_k], \;   \mathbb{R}  ):  \; \|  v_k  \|_{\mathcal{C}^0} 
       \leq c_v \delta,    \;  v_k   \; \mathrm{is \;  Lipschtiz}, \;   Lip  (v_k)  \leq \tilde{c}_v   \delta  \big\}, \\       
&    X_{\sigma k }  : = \big\{  \sigma_k  \in \mathcal{C}^0 (   [0, \, s_k], \;   \mathbb{R}   ):  \; \|  \sigma_k - \lambda_k(  \bar{u}_k )   \|_{\sigma}   
       \leq c_{\sigma} \delta^2,    \;  \sigma_k   \; \mathrm{is \;  Lipschtiz}, \;   Lip  (\sigma_k)  \leq \tilde c_{\sigma}   \big\}. \\
\end{split}
\end{equation}      
In the previous expression, $ \|  \cdot    \|_{\sigma}$ is defined as follows:  
$$
     \|  \cdot    \|_{\sigma} : = \eta_k  \delta \|  \cdot    \|_{\mathcal{C}^0},
$$
for a suitable constant $\eta_k$ whose exact value does not depend on $\delta$ and will be determined in the following.  Also the constants $c_u$, $c_v$, $\tilde{c}_v$, $c_{\sigma}$ and $\tilde{c}_{\sigma}$ do not depend on $\delta$  and their exact value will be estimated in the following.

If $(u_k, \, v_k, \, \sigma_k)  \in X_{uk} \times X_{vk} \times X_{\sigma k}$ then the function $f$ defined by \eqref{e:Binv:char:f} satisfies $f \in \mathcal{C}^{1, \, 1}_k$ for a large enough constant $k$. Moreover, exploiting \eqref{e:Binv:char:lambda:piccino} and \eqref{e:BInv:char:lambda}, one gets 
\begin{equation}
\label{e:Binv:char:est:f}
\begin{split}
          \|     f  \|_{\mathcal{C}^0} 
 &         \leq \delta  \| \check{\lambda}_k \|_{\mathcal{C}^0}  \phantom{\Bigg( }  \\
&           \leq 
          \delta     \Bigg(   \Big|  \phi_k  \Big(   \bar{u}_k, \, 0, \, \lambda_k (\bar{u}_k)  \Big) \Big|     +           
          \Big\| \phi_k \Big(  u_k, \, v_k, \, \sigma_k   \Big) - \phi_k  \Big(   \bar{u}_k, \, 0, \, \lambda_k (\bar{u}_k)  \Big) \Big\|_{\mathcal{C}^0}   + 
          \Big|  \lambda_k  (   \bar{u}_k ) \Big| \\
	  &       \qquad           +  
          \| \lambda_k ( \bar u_k ) - \sigma_k \|_{\mathcal{C}^0}
          \Bigg) \\
&           \leq \delta \Bigg(    0 + \unpo (c_u \delta + c_v \delta + c_{\sigma}  \delta)  + M \delta +   \frac{c_{\sigma}}{\eta_k}  \delta             \Bigg) \\
&       \leq \frac{1}{2} c_v \delta 
\end{split}
\end{equation}
if $\delta$ is sufficiently small and $c_v$ sufficiently large.

The same computations ensure that 
\begin{equation}
\label{e:Binv:char:est:fprime}
          \|     f '   \|_{\mathcal{C}^0} \leq \frac{1}{2} \tilde{c}_v \delta
\end{equation}
for a large enough $\tilde{c}_v$. Finally, 
\begin{equation}
\label{e:Binv:char:fkappa}
\begin{split}
          |  f'(x)   - f'(y) | 
&          = |   \check{\lambda}_k \Big(  u_k (x), \, v_k  (x), \, \sigma_k (x)  \Big)  -   \check{\lambda}_k \Big(  u_k (y), \, v_k  (y), \, \sigma_k (y)  \Big)         |
          \\
&             \leq \unpo (2 + \tilde{c}_v \delta + \tilde{c}_{\sigma}  \delta ) |x- y|  \leq k |x -y|         
\end{split}
\end{equation}
for a sufficiently large constant $k$. In the previous estimate, we have exploited the following observation: since by \eqref{e:BInv:char:partialsigma}
$$
    \frac{\partial \check{\lambda }}{\partial \sigma} \Bigg|_{\Big( u = \bar{u}_k, \, v_k = 0, \, \sigma_k = \lambda_k (\bar{u}_k) \Big)} =0,
$$
then 
\begin{equation}
\label{e:Binv:char:sigma}
     \Bigg\|  \frac{\partial \check{\lambda }}{\partial \sigma}        \Bigg\|_{\mathcal{C}^0} \leq \unpo \delta. 
\end{equation}
Because of \eqref{e:Binv:char:est:f}, \eqref{e:Binv:char:est:fprime} and \eqref{e:Binv:char:sigma}, $f \in \mathcal{C}^{1 \, 1}_k$ and hence by Proposition \ref{p:es_mon} 
$\mathrm{mon}_{[0, \, s_k]} f \in  \mathcal{C}^{1 \, 1}_k$.  Moreover, 
$$
       \|   \mathrm{mon}_{[0, \, s_k]}       f  \|_{\mathcal{C}^0 }     \leq \tilde{c}_v \delta  \qquad 
     \| (  \mathrm{mon}_{      [0, \, s_k]        } f ) '   \|_{       \mathcal{C}^0} \leq \tilde{c}_v \delta.
$$
One can then conclude that the map $T$ defined by the right hand side of  \eqref{eq:Binv:char:centerpuret} maps $X_{uk} \times X_{vk} \times X_{\sigma k}$ into itself. 

Choosing $\eta_k$ large enough (but independent from $\delta$) it is possible to prove the contraction property. One has to exploit all the previous estimates and also properties \eqref{e:monfg} and \eqref{e:Binv:char:sigma}, which are needed to handle the second and the third component of $T$.

It turns out that the Lipschitz constant of $T$ (i.e., the constant in the contraction) is uniformly bounded with respect to $\delta$.

\end{proof}

In the following it will be useful to know how the solution of \eqref{eq:Binv:char:centerpuret} depends on the length $s_k$ of the interval of definition. More precisely, in order to underline the dependence from $s_k$, we will denote the
first component of the solution of   \eqref{eq:Binv:char:centerpuret}  as $u^{s_k}_k$. Let 
$$
    F^k(\bar{u}_k, \,  s_k) : = u^{s_k}_k (s_k),
$$ 
where the dependence from the initial point $\bar{u}_k$ is also made explicit.

To study how $F^k$ depends on $s_k$ we will exploit the following result:
\begin{lem}
\label{l:center:reg}
          Fix $\delta << 1$ and let $s^1_k$, $s^2_k$ such that $s^1_k \leq  s^2_k  \leq \delta$. Let $\big(   u^1_k, \, v^1_k, \, \sigma^1_k \big)$ and $\big(   u^2_s, \, v^2_k , \, \sigma^2_k                             \big)$ be  the corresponding solutions of \eqref{eq:Binv:char:centerpuret}. Then,  
          \begin{equation}
          \label{e:Binv:char:center:lip}
                    \| u^1_k - u^2_k   \|_{\mathcal{C}^0 ([0, \, s^1_k])}  + \| v^1_k - v^2_k   \|_{\mathcal{C}^0  ([0, \, s^1_k]) }  +
                    \eta_k \delta \| \sigma^1_k - \sigma^2_k \|_{\mathcal{C}^0  ([0, \, s^1_k])}    \leq L_{k}
                   \delta  |s^1_k - s^2_k | ,
          \end{equation}
          where $L_k$ is a suitable constant which does not depend on $\delta$. 
\end{lem}
\begin{proof}
In the following we will denote by $T^{s^1_k}$ and $T^{s^2_k}$ the maps defined by the right hand side of \eqref{eq:Binv:char:centerpuret}  when $s_k = s^1_k$ and $s_k = s^2_k$ respectively. Moreover, to simplify notations, we will denote by  $\big(   u^2_s, \, v^2_k , \, \sigma^2_k      \big)$  the restriction of the fixed point of $T^{s^2_k}$ to the interval $[0, \, s^1_k]$. 

Since $\big(   u^1_k, \, v^1_k, \, \sigma^1_k \big)$ is the fixed point of $T^{s^1_k}$, one then has 
$$
    \big\|   \big(   u^1_k, \, v^1_k, \, \sigma^1_k \big)         -          \big(   u^2_s, \, v^2_k , \, \sigma^2_k                             \big)               \big\|_{  X_{uk}   \times X_{v k}   \times X_{\sigma k}              }
   \leq \frac{1}{1 - K}   \big\|   \big(   u^2_k, \, v^2_k, \, \sigma^2_k \big)         -     T^{s^1_k}    \big(   u^2_s, \, v^2_k , \, \sigma^2_k                             \big)               \big\|_{  X_{uk}   \times X_{v k}   \times X_{\sigma k}              }.
$$
In the previous formula, $K$ denotes the Lipschitz constant of $T^{s^1_k}$, which turns out to be uniformly bounded with respect to $K$ as underlined in the proof of Lemma \ref{l:centercont}. Moreover, $ X_{uk} , \;  X_{v k}$   and $ X_{\sigma k}   $ denote the spaces defined by \eqref{e:Binv:char:k:spaces} when $s_k = s^1_k$.

Thus, to prove the lemma one actually reduces to prove 
\begin{equation*}
\begin{split}
       \big\|   \big(   u^2_k, \, v^2_k, \, \sigma^2_k \big)         -     T^{s^1_k}    \big(   u^2_s, \, v^2_k , \, \sigma^2_k                             \big)               \big\|_{  X_{uk}   \times X_{v k}   \times X_    {\sigma k}            }     
           =      
&                \| u^2_k - T_u^{s^1_k}  \big(   u^2_s, \, v^2_k , \, \sigma^2_k                             \big)    \|_{\mathcal{C}^0 ([0, \, s^1_k])}  \\
&         + \| v^2_k - T_v^{s^1_k}  \big(   u^2_s, \, v^2_k , \, \sigma^2_k                             \big)                  \|_{\mathcal{C}^0  ([0, \, s^1_k]) }  \\
&        +
                    \eta_k \| \sigma^2_k -     T_{\sigma}^{s^1_k} \big(   u^2_s, \, v^2_k , \, \sigma^2_k                             \big)     \|_{\mathcal{C}^0  ([0, \, s^1_k])} \\
&        \leq \unpo   \delta  |s^1_k - s^2_k |. 
\end{split}
\end{equation*}
In the previous formula, $T^{s^1_k} _u$, $T^{s^1_k} _v$ and $T^{s^1_k} _{\sigma}$ denote respectively the first, the second and the third component of $T^{s^1_k} $. 

Since $ \big(   u^2_k, \, v^2_k, \, \sigma^2_k \big) $ is the fixed point of $T^{s^2_k}$, then
\begin{equation*} 
\left\{
\begin{array}{lll}
      u^2_k(\tau) = \bar{u}_k + {\displaystyle \int_0^{\tau} \check{r}_k (u^2_k(\xi), \, v^2_k(\xi), \, 0, \, 
      \sigma^2_k(\xi))d \xi }  \\
      v^2_k (\tau) = f_k [u^2_k, \, v^2_k, \, \sigma^2_k ] (\tau)-
      \mathrm{mon}_{[0, \, s^2_k]} f_k [ u^2_k, \, v^2_k, \, \sigma^2_k ] (\tau)   \\
      \sigma^2_k(\tau)= {\displaystyle \frac{1}{\bar{c}}  \frac{d}{d \tau}
      \mathrm{mon}_{[0, \, s^2_k]} f_k  (\tau, \, u^2_k, \, v^2_k, \, \sigma^2_k )}. \\
\end{array}
\right.
\end{equation*}
Hence, to prove the lemma it is sufficient to show that 
$$
    \|    \mathrm{mon}_{[0, \, s^2_k]} f_k [ u^2_k, \, v^2_k, \, \sigma^2_k ]  -  \mathrm{mon}_{[0, \, s^1_k]} f_k [ u^2_k, \, v^2_k, \, \sigma^2_k ]     \|_{\mathcal{C}^0 [0, \, s^1_k]}
    \leq \unpo \delta (s^2_k  - s^1_k   )
$$
and that 
$$
    \|  (  \mathrm{mon}_{[0, \, s^2_k]} f_k [ u^2_k, \, v^2_k, \, \sigma^2_k ] )'  - ( \mathrm{mon}_{[0, \, s^1_k]} f_k [ u^2_k, \, v^2_k, \, \sigma^2_k ] )'    \|_{\mathcal{C}^0 [0, \, s^1_k]}
    \leq \unpo \delta (s^2_k  - s^1_k   ).
$$
This follows directly from Proposition \ref{p:monint}.
\end{proof}

In particular, the previous Lemma implies that 
\begin{equation}
 \label{e:Binv:char:center:lip2}
     |    F (\bar{u}_k, \, s^1_k ) - F  (  \bar{u}_k, \,, s^2_k )  |  \leq 3 | s^1_k -  s^2_k  |.
\end{equation}
Indeed,  assuming for example that $s^1_k \leq s^2_k$, one can write 
\begin{equation*}
\begin{split}
     |u^1_k (s_k^1)  -  u^2_k (s_k^2)  | 
&      \leq |u^1_k (s_k^1) - u^2_k (s_k^1) | + | u^2_k (s_k^1)  -  u^2_k (s_k^2)  | \\
&     \leq 
        L_k \delta (s^2_k - s^1_k) + \Big|  \int_{s^1_k}^{s^2_k}  \check{r}_k \big( u^2_k  (\tau ) , \, v^2_k   (\tau ) , \, \sigma^2_k  (\tau )  \big) d \tau  \Big| \\
&     \leq    (L_k \delta + 2) (s^2_k - s^1_k  ).
\end{split}
\end{equation*} 
Take $s^1_k = 0$ in the previous estimate: since
\begin{equation}
\label{e:Binv:klip}
\begin{split}
    |   F^k (\bar{u}_k, \, s^1_k)  -   F^k (\bar{u}_k, \, s^2_k) | = 
&     \int_0^{s^2_k}  \check{r}_k \big( u^2_k  (\tau ) , \, v^2_k   (\tau ) , \, \sigma^2_k  (\tau )  \big) d \tau \\
&      =  \int_0^{s^2_k}  \check{r}_k \big( \bar{u}_k , \, 0, \, 0 \big) d \tau \\
&   \quad      + 
       \int_{0}^{s^2_k} \big[   \check{r}_k \big( u^2_k  (\tau ) , \, v^2_k   (\tau ) , \, \sigma^2_k  (\tau )  \big)   -   \check{r}_k \big( \bar{u}_k , \, 0, \, 0 \big) \big] d \tau  \\
&   =     \check{r}_k \big( \bar{u}_k , \, 0, \, 0 \big) s^2_k + \unpo \delta s^2_k,      \phantom{\int}
\end{split}     
\end{equation}
one gets 
\begin{equation}
\label{e:Binv:kdiff}
    F^k (\bar{u}_k, \, s_k) = u^2_k (s_k) =    \bar{u}_k + \check{r}_k \big( \bar{u}_k , \, 0, \, 0 \big) s_k  + \unpo \delta^2.
\end{equation}
Thus, the function $F^k$ defined before is differentiable at $s_k = 0$ and the gradient is the eigenvector $r_k (\bar{u}_k)  =  \check{r}_k \big( \bar{u}_k , \, 0, \, 0 \big) $.

The function $F^k$ depends Lipschitz continuously  on $\bar{u}_k$ too:
\begin{lem}
\label{l:lip:ukappa}
        Fix $\delta$ such that $0 < \delta <<1$. Take $\bar{u}^1_k  \; \bar{u}^2_k  \in \mathbb{R}^N$ such that $| \bar{u}^1_k - \bar{u}^2_k |  \leq a \delta$ for a certain constant $a$. If $|s_k|  \leq \delta$, then 
        \begin{equation}
        \label{e:Binv:ukappa}
                  | F(\bar{u}^1_k, \,   s_k )  -  F( \bar{u}^2_k, \, s_k )  | \leq \tilde{ L}_k |  u^1_k - u^2_k |
        \end{equation}
        for a suitably large constant $\tilde{L}_k$.
\end{lem}
\begin{proof}
The value $F(\bar{u}^1_k, \, s_k)$ is the fixed point of the application $T$ defined by \eqref{eq:Binv:char:centerpuret}. To underline the dependence of $T$ on $\bar{u}_k$ we will write 
$$
    T:   \mathbb{R}^N  \times X_{uk}  \times X_{vk} \times X_{\sigma k}  \to X_{uk}  \times X_{vk} \times X_{\sigma k}
$$
The map $T$ depends Lipschitz continuously on $\bar{u}_k$, but one cannot apply directly the classical results (see e.g. \cite{Bre:note} ) on the dependence of the fixed point of a contraction from a parameter because of a technical difficulty. Indeed, the domain $X_{uk} \times X_{v k}  \times X_{\sigma k}$ of the contraction depends on $\bar{u}_k$.

To overcome this problem, one can proceed as follows. Define $X^1_{uk}$ and $X^2_{uk}$ the spaces defined as in \eqref{e:Binv:char:k:spaces} with $\bar{u} = \bar{u}^1_k$ and $\bar{u}_k =\bar{u}^2_k$ respectively. If $ \big(  u^1_k, \, v^1_k, \, \sigma^1_k  \big)$ is the fixed point of the application $T( \bar{u}^1_k, \, \cdot ) $ defined on $X^1_{uk}  \times X_{vk} \times X_{\sigma k}$, then 
 $$
     \|  u^1_k - \bar{u}_k \|_{\mathcal{C}^0} \leq \|  u^1_k - \bar{u}^1_k  \|_{\mathcal{C}^0}  + |  u^1_k   -  u^2_k  |  \leq 2 s_k + a \delta \leq (2 + a) \delta.
 $$
 Choosing $c_u \ge (2 +a)$, one gets that $ \big(  u^1_k, \, v^1_k, \, \sigma^1_k  \big) $ belongs to $    X^2_{uk}  \times X_{vk} \times X_{\sigma k}   $. Thus, 
 $$
     \|     \big(  u^1_k, \, v^1_k, \, \sigma^1_k  \big)  -  \big(  u^2_k, \, v^2_k, \, \sigma^2_k  \big)               \|_{         X^2_{uk}  \times X_{vk} \times X_{\sigma k}     }
     \leq  \frac{1}{1-K}    \|     \big(  u^1_k, \, v^1_k, \, \sigma^1_k  \big)  -  T \big(  \bar{u}^2_k, \, u^1_k, \, v^1_k, \, \sigma^1_k  \big)               \|_{         X^2_{uk}  \times X_{vk} \times X_{\sigma k}     }.
 $$
 In the previous expression, $ \big(  u^2_k, \, v^2_k, \, \sigma^2_k  \big)$ is the fixed point of the application $T ( \bar{u}^2_k, \, \cdot) $. One can then proceed as in \cite{Bre:note} and conclude that 
 \begin{equation}
 \label{e:Binv:char:kutot}
             \|     \big(  u^1_k, \, v^1_k, \, \sigma^1_k  \big)  -  \big(  u^2_k, \, v^2_k, \, \sigma^2_k  \big)               \|_{         X^2_{uk}  \times X_{vk} \times X_{\sigma k}     } 
            \leq   \tilde{L}_k |  \bar{u}^1_k -  \bar{u}^2_k   |  
 \end{equation}

\end{proof}
\item {\bf Analysis of the purely stable component}

The purely center component satisfies \eqref{e:Binv:char:stable}.
  
\begin{lem}
\label{l:purestable}
          Fix $\delta$ such that $|V^s (0)|  \leq \delta <   <  1$. Then system \eqref{e:Binv:char:stable} defines a contraction on the space
          $ X^s_u \times X^s_v$, where 
          $$
                X^s_u : = \Big\{   u^s    \in \mathcal{C}^0    ( [0,    \, + \infty    ), \, \mathbb{R}^N),  \;   \|   u^s    \|_{us}  \leq m_u  \delta \Big\}
                \qquad 
                 X^s_v : = \Big\{   V^s    \in \mathcal{C}^0    ( [0,    \, + \infty    ), \, \mathbb{R}^{  k -1   }),  \;   \|   V^s    \|_{vs}  \leq m_v  \delta \Big\}.    
          $$
          The constants $m_u$ and $m_v$ do not depend on $\delta$
          the norms $  \|   \;   \cdot  \;  \|_{us} $  and $  \|   \;   \cdot \;    \|_{vs} $ are defined by 
          \begin{equation}
          \label{e:Binv:char:norm}
                \|   u^s    \|_{us}  : = \eta_s   \sup \big\{  e^{c x /2  }   |  u^s  (x) |     \big\}   \qquad  
                \|   V^s    \|_{vs}   : =  \sup  \big\{   e^{   c x  / 2 }    | V^s (x)  |       \big\}, 
          \end{equation}
          where $\eta$ is a small enough constant which does not depend on $\delta$ and $c$ is defined by
          \eqref{eq_Binv_char_separation}.
\end{lem}
\begin{proof}
We know that $\bar{\Lambda}$ is has eigenvalues $\lambda_1 (  \bar{u}_k) \dots \lambda_{k -1}   (  \bar{u}_k    )$.  Relying on \eqref{eq_Binv_nonchar_separation}, one obtains 
$$
    |    e^{ \bar{\Lambda} x  }  V^s( 0 )  |      \leq e^{  - c x /2}  |   V^s (  0  )   |  \leq  e^{  - c x /2} \delta.
$$
Moreover, if $u^s \in X^s_u$ and $V^s \in X^s_v$, then 
$$
    |  \Lambda (  \bar{u}_k + u^s, \, 0, \, V^s)  - \bar{ \Lambda }     |  \leq \unpo \delta,
$$
where $\unpo$ denotes (here and in the following) some constant that depends neither on $\eta$ nor on $\delta$.  

The columns of $\check{R}^s  ( \bar{u}_k, \, 0, \, 0  ) $ are the eigenvectors $r_1 (\bar{u}_k)    \dots r_{k -1 } (\bar{u}_k)$, which are all unit vectors. Thus, if 
 $u \in X^s_u$ and $V^s \in X^s_v$, then 
$$
   |   \check{R}^s  ( \bar{u}_k  + u^s, \, 0, \, V^s  )         |   \leq \unpo.
$$ 
From the previous observations one gets that the application defined by \eqref{e:Binv:char:stable}  maps  $X^s_u  \times X^s_v$ into itself.

To prove that the application is a contraction, take $u^{1 s},   \;   u^{2 s } \in X^s_u$ and $V^{ 1 s}, \;  V^{2  s}   \in X^s_v$. Then 
\begin{equation*}
\begin{split}
&         \eta_s  \, \Big|   \displaystyle{        \int_x^{+   \infty   }     \check{R}^s  \Big(  \bar{u}_k + u^{1s}  (y) , \, 0, \, V^{1s}  (y)  \Big) V^{1s} (y)   dy    } 
                     \displaystyle{       -  \int_x^{+   \infty   }     \check{R}^s  \Big(  \bar{u}_k + u^{2s}  (y) , \, 0, \, V^{2s}  (y)  \Big) V^{2 s} (y)   dy    }
           \Big|  \\
&         \leq    \unpo \delta \|  u^{1s}   -  u^{2s}     \|_{us } +  \unpo \delta \eta_s   \|  V^{1s}   -  V^{2s}     \|_{v s }  +  \unpo \eta_s    \|  V^{1s}   -  V^{2s}     \|_{v s }.   
          \phantom{\int}\\       
\end{split}
\end{equation*}  
Choosing $\eta_s$ sufficiently small, one can suppose that $\unpo \eta_s < 1 / 4 $ in the previous expression. Moreover, one can suppose that $\delta$ is small enough to have that 
$\unpo \delta < 1 / 2$, $\unpo \delta \eta_s < 1 / 4$. 

Finally, 
\begin{equation*}
\begin{split}
&       e^{ c x / 2 }  \Big|   \displaystyle{        \int_0^x   e^{   \bar{\Lambda} (x - y)    } 
           \Big[  \check{\Lambda}^s   \Big(  u_{1s} (y)  + \bar{u}_k, \, 0, \, V^{1s}  (y) \Big) - \bar{\Lambda} \Big] V^{1s} (y)   dy   } \\
 &       \qquad           
           -   \displaystyle{        \int_0^x   e^{   \bar{\Lambda} (x - y)    } 
           \Big[  \check{\Lambda}^s   \Big(  u_{2s} (y)  + \bar{u}_k, \, 0, \, V^{2s}  (y) \Big) - \bar{\Lambda} \Big] V^{2s} (y)   dy   }  \Big|\\
&        \leq \unpo \frac{\delta}{\eta_s}  \|  u^{1s}   -  u^{2s}     \|_{u s } + \unpo \delta \|  V^{1s}   -  V^{2s}     \|_{v s } .          
\end{split}
\end{equation*}
Assuming that $\delta$ is small enough, $\unpo \delta /  \eta_s  < 1 /2$ and $\unpo \delta < 1/ 2 $. Thus the map is contraction and the constant of the contraction is less or equal to $1 /2$ uniformly for $\delta \to 0^+$.
\end{proof}

The solution of \eqref{e:Binv:char:stable} depends on the parameter $V^s (0)$. The regularity of $\Big( u_s, \, V_s \Big)$ with respect to $V^s (0)$ is discussed in the following 
lemma. 
\begin{lem}
\label{l:stable:reg}
          Fix $\delta << 1$ and let $V^1_s (0)$, $V^2_s (0)$ two initial data such that $|V^1_s (0)|$,  $|V^2_s (0)|  \leq \delta$. Let $(   u^1_s, \, V^1_s )$ and $(   u^2_s, \, V^2_s  )$   the corresponding solutions of \eqref{e:Binv:char:stable}. Then, 
          \begin{equation}
          \label{e:Binv:char:stable:lip}
                    \| u^1_s - u^2_s   \|_{us}  + \| V^{1 s }  -  V^{2 s }   \|_{v s }  \leq L_{s}
                    |V^1_s (0) - V^2_s (0) | ,
          \end{equation}
          where $L_s$ is a suitable constant which does not depend on $\delta$. Moreover, if $V^2_s (0) = \vec{0}$, then $u^2_s \equiv 0$. Also,  
          \begin{equation}
          \label{e:Binv:char:stable:diff}
                    u^1_s (0)  =  \check{R}^s (  \bar{u}_k, \, 0, \, 0) \bar{\Lambda}^{-1} V^1_s (0)  + \unpo \delta^2. 
          \end{equation}
 \end{lem}
Equation \eqref{e:Binv:char:stable:diff} guarantees, in particular, that the application  
$$
    F^s  (    \bar{u}, \, V_s (0)  ) : = u^s (0). 
$$
is differentiable with respect to $V_s (0)$ when $V_s (0) = 0$ and that  the jacobian is the matrix $ \check{R}^s (\bar{u}_k, \, 0, \, 0) \bar{\Lambda}^{-1}$, whose columns are the eigenvectors    
$r_1  (  \bar{u}_k)  / \lambda_1 (\bar{u}_k)  \dots     r_{k-1}  ( \bar{u}_k) / \lambda_{k-1} (\bar{u}_k)$.
\begin{proof}
Let $T$ the application defined by the right hand side of \eqref{e:Binv:char:stable}. To underline the dependence on the parameter $V^s (0)$, we write 
$$
    T: X_{us} \times X_{vs} \times \mathbb{R}^{k-1} \to  X_{us} \times X_{vs}. 
$$
and denote by $T_u$ and $T_v$ respectively the first and the second component of $T$.
For every $(u_s, \, V_s) \in X_{us} \times X_{vs} $,  
$$
    \Big\|  T_u  \Big(  u_s, \, V_s, \, V^1_s (0)   - T_u  \Big(  u_s, \, V_s, \, V^2_s (0)  \Big)       \Big\| _{us} = 0 
$$
and
$$ 
      \Big\|  T_v  \Big(  u_s, \, V_s, \, V^1_s (0)   - T_v  \Big(  u_s, \, V_s, \, V^2_s (0)  \Big)       \Big\| _{vs} \leq |V^1_s - V^2_s| .
$$
Hence, \eqref{e:Binv:char:stable:lip} holds with $L_s := 1 / (1-K)$, where $K$ is the Lipschitz constant of $T$ and hence it is smaller  than $1$
because $T$ is a contraction. Moreover, from the proof of Lemma \ref{l:purestable} it follows that $k$ is bounded away from $0$ uniformly with respect to $\delta$. This concludes the proof of the first part of the lemma.

To prove the second part, we observe that 
$$
    V^1_s (x)  = e^{\bar{\Lambda}  x} V^1_s (0) +  \unpo \delta^2 e^{  - c x / 2 }
$$
and that 
$$
    \big|  \check{R}^s (\bar{u}_k  + u^1_s, \, 0, \, V^1_s)    -  \check{R}^s (\bar{u}_k, \, 0, \, 0)   \big|  \leq \unpo \delta. 
$$
Hence, 
\begin{equation*}
\begin{split}
              u^1_s (0) =
&              \int_0^{  + \infty }    \check{R}^s (\bar{u}_k, \, 0, \, 0)  V^1_s (y) dy      + 
                \int_0^{+  \infty }   \Big[  \check{R}^s (\bar{u}_k  + u^1_s, \, 0, \, V^1_s) -  \check{R}^s (\bar{u}_k  + u_s, \, 0, \, V_s)  \Big]  V^1_s (y)  dy    \\  
&             =  \int_0^{  + \infty }    \check{R}^s (\bar{u}_k, \, 0, \, 0)  e^{\bar{\Lambda} y } V^1_s (0) dy      + \unpo \delta^2 \int_{+  \infty}^0  e^{   -  cy / 2  } dy + 
                 \unpo \delta^2   \int_{+  \infty}^0  e^{   -  cy / 2  } dy \\
&            =      \check{R}^s (  \bar{u}_k, \, 0, \, 0) \bar{\Lambda}^{-1} V^1_s (0)  + \unpo \delta^2.   \phantom{\int}          
\end{split}
\end{equation*}
This completes the proof of the lemma. 
\end{proof}
To study the dependence of $F^s  \big(  \bar{u}_k, \, V_s (0), \,  s_k \big)$ on $\bar{u}_k$ one can proceed as follow. Fix $\bar{u}^1_k, \; \bar{u}^2_k \in \mathbb{R}^N$ close enough  and denote by 
$\big( u^1_s, \, V^1_s \big)$ and $\big(  u^2_s, \, V^2_s \big)$ the solution of \eqref{e:Binv:char:stable}  corresponding to $\bar{u}_k = \bar{u}^1_k$ and $\bar{u}_k= \bar{u}^2_k$ respectively.  Then  classical results (see e.g. \cite{Bre:note}) on the dependence of the fixed point of a contraction on a parameter ensure 
that 
\begin{equation}
\label{e:Binv:char:sviessetot}
         \| u^1_k -   u^2_k  \|_{ us}    + \| V^{1s}  - V^{2s}   \|_{vs} \leq \tilde{L}_s |   \bar{u}^1_k   - \bar{u}^2_k |
\end{equation}
for a suitable constant $\tilde{L}^s$. In particular, 
\begin{equation}
\label{e:Binv:char:eviesse}
          |  F^s \big( \bar{u}^1_k, \, V^s (0)  \big)    -  F^s \big(  \bar{u}^2_k, \, V^s (0) \big) | \leq \tilde{L}_s |  \bar{u}^1_k - \bar{u}^2_k |. 
\end{equation}
\item {\bf The component of perturbation}

Assume that $v_k (s_k)  < 0$ and define $\alpha $ as in \eqref{e:Binv:char:alpha}:
$$
          \alpha (\tau )  =    -   \int_{\tau}^{s_k}  \frac{1}{   v_k (s)}   ds.
$$
Let 
$$
    \underline{s}: = \max \{ \tau \in [0, \, s_k ]:   \:  v_k (\tau)  =  0\},
$$
we can now prove that $\alpha (\tau) < + \infty $ if and only if $\tau > \underline{s}$. 

Indeed, for $s \ge \underline{s}$ it holds  
$$
    - v_k (s) = - v_k (s) + v_k (\underline{s})  \leq  \tilde{c}_v \delta (s - \underline{s}) 
$$
and hence 
$$
    \alpha (\tau )  =    \int_{\tau}^{s_k}    - \frac{1  }{v_k (s)  }  ds      \ge   \int_{   \tau  }^{   s_k}    \frac{1}{   \tilde{c}_v \delta (s - \underline{s}    )       } ds  = + \infty. 
$$
Let 
$$
    \beta : [0, \, + \infty [  \to     ]\underline{s}, \, s_k]
$$
the inverse function of $\alpha$, $\beta \circ \alpha (\tau) = \tau$.     
If $v_k (s_k) =0$, we set $\beta (x) = s_k  $ for every $x$.

The component of perturbation $(U, \, q, \, p)$ is the fixed point of the transformation 
\begin{equation}
\label{e:Binv:char:per}
\left\{
\begin{array}{lll}
           U (x) = \displaystyle{  \int^x_{+ \infty} \Big[    \check{R}^s \Big(  u_k \circ \beta (y) + u_s (y) + U(y), \, v_k \circ \beta (y)  + q (y), \, V_s (y) + p (y), \, 0 \Big) } \\ 
           \qquad  \qquad \qquad -
            \displaystyle{\check{R}^s \Big(  \bar{u}_k  + u_s (y) , \, 0, \, V_s (y) , \, 0 \Big) \Big] V_s (y) dy   \phantom{\int} }  \\
            \qquad  \qquad +   
           \displaystyle{ \int^x_{+ \infty}  \check{R}^s \Big(  u_k \circ \beta (y) + u_s (y) + U(y), \, v_k \circ \beta (y)  + q (y), \, V_s (y) + p (y), \, 0 \Big)  p(y) dy}  \\
            \qquad \qquad \displaystyle{    +.       \int^x_{+ \infty} \Big[    \check{r}^k \Big(  u_k \circ \beta (y) + u_s (y) + U(y), \, v_k \circ \beta (y)  + q (y), \, V_s (y) + p (y), \, 0 \Big) } \\
            \qquad \qquad \qquad  \displaystyle{-    \check{r}^k \Big(  u_k \circ \beta (y) , \, v_k \circ \beta (y), \, 0, \, 0 \Big) \Big] v_k \circ \beta (y) dy \phantom{\int} } \\
            \qquad \qquad \displaystyle{+  \int^x_{+ \infty} \check{r}^k \Big(  u_k \circ \beta (y) + u_s (y) + U(y), \, v_k \circ \beta (y)  + q (y), \, V_s (y) + p (y), \, 0 \Big)   q(y) dy} \\
            q(x) = \displaystyle{       \int_{+ \infty}^x         \Big[    \phi_k  \Big(  u_k \circ \beta (y) + u_s (y) + U(y), \, v_k \circ \beta (y)  + q (y), \, V_s (y) + p (y), \, \sigma \circ \beta  \Big)  } \\
            \qquad \qquad \qquad  \displaystyle{   - 
             \phi_k \Big(  u_k \circ \beta (y), \, v_k \circ \beta (y), \,  0 , \,  \sigma \circ \beta  \Big) \Big]  v_k \circ \beta (y)  dy     \phantom{\int}  } \\
             \qquad \qquad \displaystyle{  + \int_{+  \infty}^x   \phi_k  \Big(  u_k \circ \beta (y) + u_s (y) + U(y), \, v_k \circ \beta (y)  + q (y), \, V_s (y) + p (y), \, \sigma_k \circ \beta  \Big)  q (y) dy } \\
           p(x) =  \displaystyle{ \int_0^x e^{\bar{\Lambda }    (x - y)  }   \Big[  \Lambda \Big(  u_k \circ \beta (y) + u_s (y) + U(y), \, v_k \circ \beta (y)  + q (y), \, V_s (y) + p (y)  \Big) - \bar{\Lambda} \Big] p (y) dy   } \\
           \qquad \qquad + \displaystyle{    \int_0^x    e^{\bar{\Lambda }    (x - y) }
            \Big[     \Lambda \Big(  u_k \circ \beta (y) + u_s (y) + U(y), \, v_k \circ \beta (y)  + q (y), \, V_s (y) + p (y)  \Big)  } \\
            \qquad \qquad \qquad \qquad \qquad \displaystyle{  -    \Lambda \Big(  u_s (y) , \, 0), \, V_s (y)  \Big)  \Big]  V^s (y) dy           \phantom{\int}}. 
\end{array}
\right.
\end{equation}
In the equation for $p$ we used the same notation as in \eqref{e:Binv:char:stable}:
$
    \bar{\Lambda} = \check{\Lambda}^s (\bar{u}_k, \, 0, \, 0). 
$
\begin{lem}
\label{l:pert} 
          Let $(u_k, \, v_k,    \, \sigma_k)$ and $(u_s, \, V^s)$ be as in Lemma \ref{l:centercont} and \ref{l:purestable} respectively. Then there exists a unique solution $(U, \, q, \, p)$ of 
          \eqref{e:Binv:char:per} belonging to $X_U \times X_q \times X_p$, where the spaces $X_U$, $X_q$ and $X_p$ are defined as follows: 
          $$
              X_U : = \Big\{  U \in \mathcal{C}^0 ([0, \, +  \infty[ , \, \mathbb{R}^N ): \;  | U(x) |  \leq k_U \delta^2 e^{- cx / 2 }  \Big\},   
          $$
          $$
              X_p : = \Big\{  p \in \mathcal{C}^0 ([0, \, +  \infty[ , \, \mathbb{R}^{k-1} ): \;  | p(x)| \leq k_p \delta^2  e^{- cx / 2 }                                    \Big\} ,
          $$
           and     
         $$
              X_q : = \Big\{  q \in \mathcal{C}^0 ([0, \, +  \infty[ , \, \mathbb{R} ): \;  |   q  (x)  |_{q} \leq k_q \delta^2  e^{- cx / 2 }.                                       \Big\}
         $$
         In the previous expressions, the constants $k_U, \; k_q$ and $k_p$ do not depend on $\delta$.                  
         
         The space $X_U$ is endowed with the norm 
         $$
            \| U  \|_{U}= \eta_U  \sup_x \Big\{  e^{  c x  /  4 }   |U ( x ) |     \Big\}
         $$
         where $\eta_U$ is a suitable constant which does not depend on $\delta$. 
         The spaces $X_p$ and $X_q$ are endowed respectively with the norms
         $$
            \| p  \|_{p}=   \sup_x \Big\{  e^{  c x  /  4 }  |p ( x ) |      \Big\}
         $$
         and 
         $$
            \| q  \|_{q}=   \sup_x \Big\{  e^{  c x  /  4 }  |q ( x ) |      \Big\}.
         $$
         \end{lem}
\begin{proof}
It is enough to prove that the right hand side of \eqref{e:Binv:char:per} defines a contraction on  $X_U \times X_q \times X_p$.  As in the proof of Lemma \ref{l:purestable}, one observes that $\bar{\Lambda}$ is a diagonal matrix with all the eigenvalues less or equal to $-c$. Thus, 
$$
   | e^{  \bar{\Lambda} x} \vec{\xi} | \leq e^{ - c x} |  \vec{\xi} | \quad \forall, \; \vec{\xi} \in \mathbb{R}^{  k-1 }. 
$$
Moreover, from Lemma \ref{l:purestable} one gets that for every $x$, 
$$
    |V_s (x)| \leq m_v \delta e^{ - c x / 2} 
    \qquad |u_s (x)|  \leq m_u \delta e^{ - c x / 2}.  
$$
Another useful observation is the following:
\begin{equation}
\label{e:Binv:char:changevar}
          \int_{  +  \infty} ^0  v_k \circ \beta (y) dy \leq s_k \leq \delta. 
\end{equation} 
Indeed, if $v_k (s_k) < 0$, then one can exploit the change of variable $\tau = \beta (y)$. Recalling that $\beta$ is the inverse function of $\alpha$, which is defined by 
\eqref{e:Binv:char:alpha},  one has 
$$
          \int_{  +  \infty} ^0  v_k \circ \beta (y) dy \leq s_k  = \int_{\underline{s}}^{  s_k }  v_k (\tau) \alpha' (\tau) d \tau =
          \int_{\underline{s}}^{  s_k }  d \tau = (s_k - \underline{s}) \leq s_k.	
$$
If $v_k (s_k) =0$, then $v_k \circ \beta \equiv 0$ and hence \eqref{e:Binv:char:changevar} is trivial.  

One can also assume 
\begin{equation}
\label{e:Binv:char:phik}
         |\phi_k   \big(  \bar{u}_k , 0, \, 0, \, \lambda_k ( \bar{u}_k  \big)| \leq K \delta
\end{equation}
and hence 
$$
     |\phi_k   \big(    u_k \circ \beta (y) + u_s (y) + U(y), \, v_k \circ \beta (y)  + q (y), \, V_s (y) + p (y), \, \sigma_k \circ \beta     \big)| \leq    \unpo  \delta.
$$
Moreover, it holds
$$
         \Big|     \Lambda \Big(  u_k \circ \beta (y) + u_s (y) + U(y), \, v_k \circ \beta (y)  + q (y), \, V_s (y) + p (y)  \Big)
           -    \bar{\Lambda}   \Big| \leq \unpo \delta
$$
Since $
    \bar{\Lambda} = \check{\Lambda}^s (\bar{u}_k, \, 0, \, 0), 
$
this follows from the regularity of $\Lambda$ and from the estimates
$$
    |u_k - \bar{u}_k  |, \; |u_s|, \; |v_k|, \; |q|, \;  | V_s|, \;  |p| \leq \unpo \delta.  
$$

Exploiting the previous observations,  by  direct check one gets that the right hand side of \eqref{e:Binv:char:per} belongs to $X_U \times X_q \times X_p$.  

To prove that the map is actually a contraction, one has to exploit again the previous remarks, combined with the following: for every $x$,  
$$
     | U^1 - U^2  | ( x ) \leq \frac{1}{ \eta_U  }  \|  U^1 -   U^2   \|_U e^{- cx / 4 }, \quad 
      | q^1 - q^2  | ( x ) \leq   \|  q^1 -   q^2   \|_q e^{- cx / 4 } 
$$
and 
$$      
           | p^1 - p^2  | ( x ) \leq   \|  p^1 -   p^2   \|_p e^{- cx / 4 } . 
$$
One can then check that the map is a contraction, provided that $\eta_U$ is sufficiently small.
\end{proof}         

Since $(u_k, \, v_k, \, \sigma_k)$ and $(u_s, \, V_s)$ depend on the parameters  $s_k$ and $V_s (0)$ respectively, then also $(U, \, q, \, p)$ does. Let 
$$
    F^p \big( \bar{u}_k, \, V_s (0), \, s_k  \big) : =  U(0).
$$
The following lemma concerns the regularity of $F^p$ with respect to $\big(V_s (0), \, s_k  \big)$. 
\begin{lem}
\label{l:pert:reg}
          The function $F^p$ is differentiable with respect to  $\big(V_s (0), \, s_k  \big)$ at $\big(V_s (0), \, s_k  \big)= \vec{0}$ and the jacobian 
          is the null matrix. 
          
          Also, let $ \big(V^1_s (0), \, s^1_k  \big)$ and  $\big(V^2_s (0), \, s^2_k  \big) $ such that   $\big| \big( V^1_s (0), \, s^1_k  \big) \big|$, 
          $\big| \big( V^2_s (0), \, s^2_k  \big) \big| \leq \delta$. 
          Then, 
          \begin{equation}
          \label{e:Binv:char:per:lip}
              \big|  F^p   \big(\bar u^1_k, \, V^1_s (0), \, s^1_k  \big) -   F^p \big(\bar u^2_k \, V^2_s (0), \, s^2_k  \big) \big| \leq L_p \Big( |s^1_k -s^2_k |    +   |  V^1_s(0)  - V^2_s (0)   |\Big).
          \end{equation}
\end{lem}
\begin{proof}
By construction, 
\begin{equation}
\label{e:Binv:char:per:diff}
    |U(0)| \leq \unpo \delta^2.
\end{equation}
This implies the differentiability at $\big(V_s (0), \, s_k  \big)  =      \vec{0}$.

To prove the second part of the lemma, we will focus only on the dependence from $V_s (0)$ and $s_k$. The proof of the Lipschitz continuous dependence from $\bar u_k$ is completely analogous. 

We can observe that  $(U, \, q, \, p)$  is the fixed point of a contraction $T$, defined by the right hand side of \eqref{e:Binv:char:per}. To underline the dependence from $\big(V_s (0), \, s_k  \big)$ we write 
$$
    T: X_u \times X_q \times X_p \times \mathbb{R}^{k - 1} \times \mathbb{R}   \to  X_u \times X_q \times X_p 
$$
and denote by $T_U$, $T_q$ and $T_p$ respectively the first, the second and the third component of $T$.
Since the Lispchitz constant  of $T$ is not only smaller than $1$ but also bounded away from $0$ uniformly with respect to $\delta$, then to prove  \eqref{e:Binv:char:per:lip} it is enough to prove that,
for every $(U, \, q, \, p)$, 
\begin{equation}
\label{e:Binv:char:per:lip2}
\begin{split}
&             \big\|    T_U \big(  V^1_s (0), \, s^1_k,   \, U,   \, p, \, q            \big)      -        T_U \big(  V^2_s (0), \, s^2_k,   \, U,   \, p, \, q            \big)               \big\|_U   \leq 
                  \unpo       \Big( |s^1_k -s^2_k |    +   |  V^1_s(0)  - V^2_s (0)   |\Big) \\
&                   \big\|    T_q \big(  V^1_s (0), \, s^1_k,   \, U,   \, p, \, q            \big)      -        T_q \big(  V^2_s (0), \, s^2_k,   \, U,   \, p, \, q            \big)               \big\|_q   \leq 
                  \unpo       \Big( |s^1_k -s^2_k |    +   |  V^1_s(0)  - V^2_s (0)   |\Big)  \\
&                  \big\|    T_p  \big(  V^1_s (0), \, s^1_k,   \, U,   \, p, \, q            \big)      -        T_p   \big(  V^2_s (0), \, s^2_k,   \, U,   \, p, \, q            \big)               \big\|_p \leq 
                  \unpo       \Big( |s^1_k -s^2_k |    +   |  V^1_s(0)  - V^2_s (0)   |\Big).                   
\end{split}
\end{equation}
Let 
$$
    \alpha_1 (\tau) = - \int_{\tau}^{s^1_k} \frac{1}{v^1_k  (s)} ds    \qquad 
    \alpha_2 (\tau) = - \int_{\tau}^{s^2_k} \frac{1}{v^2_k  (s)} ds 
$$
and $\beta^1$ and $\beta^2$ the inverse functions. Just to fix the ideas, let us assume that $s^1_k < s^2_k$, then 
$\beta_1 (0) = s^1_k < \beta_k^2 (0)  = s_k^2$. Exploiting Lemma \ref{l:center:reg} one gets that, as far as $\beta_1 (x) \leq \beta_2 (x)$, it holds  
\begin{equation*}
\begin{split}
           \frac{d  (\beta_2 - \beta_1)}{d x}
&         = v^2_k   \circ \beta^2 -  v^1_k \circ \beta^1 = v^2_k   \circ \beta^2  - v^2_k   \circ \beta^1 + v^2_k   \circ \beta^1 - v^1_k   \circ \beta^1   \\
&        \leq \|   v^2_k   - v^1_k   \|_{\mathcal{C}^0} + \tilde{c}_v \delta (\beta^2 - \beta^1 )    \phantom{\int}\\
&        \leq L_k \delta (s^2_k - s^1_k ) + \tilde{c}_v \delta (\beta^2 - \beta^1 )   \phantom{\int}  \\
\end{split}
\end{equation*}
and hence 
\begin{equation}
\label{e:Binv:char:exp:delta}
          |   \beta^1 (x) - \beta^2 (x) |  \leq   \Big[  \frac{L_k}{  \tilde{c}_v}   (s^2_k - s^1_k ) +   (s^2_k - s^1_k )   \Big]    e^{\displaystyle{\tilde{c}_v \delta x }}  -  
         \frac{L_k}{  \tilde{c}_v}   (s^2_k - s^1_k ) \leq \unpo (s^2_k - s^1_k )    e^{\displaystyle{\tilde{c}_v \delta x }}
\end{equation}
If $\beta^1 > \beta^2$, then 
$$
     \beta^1 (x) - \beta^2 (x) \leq   \frac{L_k}{ \tilde{c}_v }  \Big[     e^{\displaystyle{\tilde{c}_v \delta (x - \bar{x} ) }} - 1    \Big]   (s^2_k - s^1_k ),   
$$
where 
$$
     \bar{x} = \max \{   y <  x:    \;     \beta^1 (y)  = \beta^2 (y)   \}.
$$
Thus, estimate \eqref{e:Binv:char:exp:delta} still holds.

Exploiting \eqref{e:Binv:char:exp:delta}, one gets 
\begin{equation}
\label{e:Binv:char:exp:delta:u}
\begin{split}
           | u^2_k   \circ \beta^2 (x)-  u^1_k \circ \beta^1 (x)  |  
&        \leq     | u^2_k   \circ \beta^2  - u^2_k   \circ \beta^1 |  +   |  u^2_k   \circ \beta^1 (x) - u^1_k   \circ \beta^1 (x) |   \\
&        \leq \|   u^2_k   - u^1_k   \|_{\mathcal{C}^0} + 2 |  \beta^2 - \beta^1 |     \leq L_k \delta (s^2_k - s^1_k ) +   \unpo (s^2_k - s^1_k )    e^{\displaystyle{\tilde{c}_v \delta x }} 
          \phantom{\int}  \\
&       \leq  \unpo (s^2_k - s^1_k )    e^{\displaystyle{\tilde{c}_v \delta x }}      \phantom{\int}  \\
\end{split}
\end{equation}
and, by analogous considerations, 
\begin{equation}
\label{e:Binv:char:exp:delta:sigmav}
                 | v^2_k   \circ \beta^2 (x)-  v^1_k \circ \beta^1 (x)  |,   \;         | \sigma^2_k   \circ \beta^2 (x)-      \sigma^1_k \circ \beta^1 (x)  |    
                 \leq \unpo (s^2_k - s^1_k )    e^{\displaystyle{\tilde{c}_v \delta x }}  . 
\end{equation}

To prove \eqref{e:Binv:char:per:lip2}, the most complicated terms to handle are in the form  
\begin{equation}
\label{e:Binv:char:badterm}
\begin{split}
&       \Big|   \int^x_{+ \infty} \Big[    \check{r}^k \Big(  u^1_k \circ \beta^1 (y) + u^1_s (y) + U(y), \, v^1_k \circ \beta^1 (y)  + q (y), \, V^1_s (y) + p (y), \, 0 \Big) \\
&         \qquad \qquad             -    \check{r}^k \Big(  u^1_k \circ \beta^1 (y) , \, v^1_k \circ \beta^1 (y), \, 0, \, 0 \Big) \Big] v^1_k \circ \beta^1 (y) dy  \\
&      \qquad  -     \int^x_{+ \infty} \Big[    \check{r}^k \Big(  u^2_k \circ \beta^2 (y) + u^2_s (y) + U(y), \, v^2_k \circ \beta^2 (y)  + q (y), \, V^2_s (y) + p (y), \, 0 \Big) \\
&      \qquad \qquad     -    \check{r}^k \Big(  u^2_k \circ \beta^2 (y) , \, v^2_k \circ \beta^2 (y), \, 0, \, 0 \Big) \Big] v^2_k \circ \beta^2 (y) dy   \Big|  \\     
&      \leq      \Big|   \int^x_{+ \infty} \Big[    \check{r}^k \Big(  u^1_k \circ \beta^1 (y) + u^1_s (y) + U(y), \, v^1_k \circ \beta^1 (y)  + q (y), \, V^1_s (y) + p (y), \, 0 \Big) \\
&       \qquad \qquad             -    \check{r}^k \Big(  u^1_k \circ \beta^1 (y) , \, v^1_k \circ \beta^1 (y), \, 0, \, 0 \Big) \Big]  \Big(  v^1_k \circ \beta^1 (y)  - v^2_k \circ \beta^2 (y) \Big) dy  \Big| \\
&       +     \Big|   \int^x_{+ \infty}    \Bigg\{   \Big[    \check{r}^k \Big(  u^1_k \circ \beta^1 (y) + u^1_s (y) + U(y), \, v^1_k \circ \beta^1 (y)  + q (y), \, V^1_s (y) + p (y), \, 0 \Big) \\
&         \qquad \qquad             -    \check{r}^k \Big(  u^1_k \circ \beta^1 (y) , \, v^1_k \circ \beta^1 (y), \, 0, \, 0 \Big)  \Big]   \\
&       \quad \qquad -   \Big[    \check{r}^k \Big(  u^2_k \circ \beta^2 (y) + u^2_s (y) + U(y), \, v^2_k \circ \beta^2 (y)  + q (y), \, V^2_s (y) + p (y), \, 0 \Big) \\
&      \qquad \qquad     -    \check{r}^k \Big(  u^2_k \circ \beta^2 (y) , \, v^2_k \circ \beta^2 (y), \, 0, \, 0 \Big) \Big]   \Bigg\}  v^2_k \circ \beta^2 (y) dy   \Big| \\
\end{split}
\end{equation}
Exploiting \eqref{e:Binv:char:exp:delta:sigmav} the first term in the previous sum can be estimated as follows: 
\begin{equation*}
\begin{split}
&        \Big|   \int^x_{+ \infty} \Big[    \check{r}^k \Big(  u^1_k \circ \beta^1 (y) + u^1_s (y) + U(y), \, v^1_k \circ \beta^1 (y)  + q (y), \, V^1_s (y) + p (y), \, 0 \Big) \\
&       \qquad \qquad             -    \check{r}^k \Big(  u^1_k \circ \beta^1 (y) , \, v^1_k \circ \beta^1 (y), \, 0, \, 0 \Big) \Big]  \Big(  v^1_k \circ \beta^1 (y)  - v^2_k \circ \beta^2 (y) \Big) dy  \Big| \\
&       \leq \int^x_{   +  \infty}     e^{\displaystyle{  - c y / 2  }} \Big[ m_u \delta + k_u \delta^2 + k_q \delta^2 + m_v \delta + k+p \delta^2  \Big] \unpo (s^2_k - s^1_k )    e^{\displaystyle{\tilde{c}_v \delta y }}  dy \\
&      \leq \unpo   (s^2_k - s^1_k )    \delta  e^{\displaystyle{  - c x / 4  }} ,
\end{split}
\end{equation*}
where we have supposed that $\delta$ is small enough to obtain $\tilde{c}_v \delta \leq c / 4$. 

To handle with the second term in \eqref{e:Binv:char:badterm} one can observe that 
\begin{equation*}
\begin{split}
&        \Big[    \check{r}^k \Big(  u^1_k \circ \beta^1 (y) + u^1_s (y) + U(y), \, v^1_k \circ \beta^1 (y)  + q (y), \, V^1_s (y) + p (y), \, 0 \Big) \\
&         \qquad             -    \check{r}^k \Big(  u^1_k \circ \beta^1 (y) , \, v^1_k \circ \beta^1 (y), \, 0, \, 0 \Big)  \Big]  \\
&       =  \int_0^1 D_u \check{r}^k \Big(  u^1_k \circ \beta^1 (y) + t u^1_s (y) +  t U(y), \, v^1_k \circ \beta^1 (y)  + t q (y), \,  t  V^1_s (y) + t p (y), \, 0  \Big)  \Big[ u^1_s (y) + U (y)  \Big] dt \\
&      \qquad + \int_0^1 D_v   \check{r}^k \Big(  u^1_k \circ \beta^1 (y) + t u^1_s (y) + t U(y), \, v^1_k \circ \beta^1 (y)  + t q (y), \, t V^1_s (y) + t p (y), \, 0  \Big)  \Big[ q (y)  \Big] dt \\          
&       \qquad +  \int_0^1 D_V \check{r}^k \Big(  u^1_k \circ \beta^1 (y) + t u^1_s (y) +  t U(y), \, v^1_k \circ \beta^1 (y)  + t q (y), \, t V^1_s (y) + t p (y), \, 0  \Big)  
        \Big[ V^1_s (y) + p (y)  \Big] dt .\\
\end{split}
\end{equation*}
Moreover, 
\begin{equation*}
\begin{split}
&       \Bigg|   \int_{+ \infty }^x \Bigg\{   \int_0^1 D_u \check{r}^k \Big(  u^1_k \circ \beta^1 (y) + t u^1_s (y) +  t U(y), \, v^1_k \circ \beta^1 (y)  + t q (y), \,  t  V^1_s (y) + t p (y), \, 0  \Big)  \Big[ u^1_s (y)        + U (y)  \Big] dt \\
&       \quad - D_u \check{r}^k \Big(  u^2_k \circ \beta^2 (y) + t u^2_s (y) +  t U(y), \, v^2_k \circ \beta^2 (y)  \\
&       \qquad \qquad + t q (y), \,  t  V^2_s (y) + t p (y), \, 0  \Big)  \Big[ u^2_s (y) + U (y)  \Big] dt 
                  \Bigg\} v_k^2 \circ \beta^2 (y) dy   \Bigg| \\
 &      \leq   \Bigg|  \int_{+ \infty }^x \Bigg\{   \int_0^1 D_u \check{r}^k \Big(  u^1_k \circ \beta^1 (y) + t u^1_s (y) +  t U(y), \, v^1_k \circ \beta^1 (y)  \\
 &       \qquad \qquad + t q (y), \,  t  V^1_s (y) + t p (y), \, 0  \Big) 
         \Big[   u^1_s (y) - u_s^2 ( y) \Big]  dt \Bigg\}   v_k^2 \circ \beta^2 (y) dy   \Bigg| \\
 &      +  \Bigg|    \int_{+ \infty }^x \Bigg\{   \int_0^1 \Big[ D_u \check{r}^k \Big(  u^1_k \circ \beta^1 (y) + t u^1_s (y) +  t U(y), \, v^1_k \circ \beta^1 (y)  + t q (y), \,  t  V^1_s (y) + t p (y), \, 0  \Big) \\
  &     \quad    - D_u \check{r}^k \Big(  u^2_k \circ \beta^2 (y) + t u^2_s (y) +  t U(y), \, v^2_k \circ \beta^2 (y)  \\
  &     \qquad \qquad + t q (y), \,  t  V^2_s (y) + t p (y), \, 0  \Big)  \Big] [u^2_s (y) + U(y)]     dt \Bigg\}   
           v_k^2 \circ \beta^2 (y) dy   \Bigg| \\                          
&        \leq        \int_{+ \infty }^x \Bigg\{   \int_0^1   \unpo L_s | V^1_s (0)  - V^2_s (0) |  e^{\displaystyle{  - c y / 2   }} dt  \Bigg\} c_v \delta dy \\
&        \quad        +    \int_{+ \infty }^x \Bigg\{   \int_0^1 \Big[  \unpo (s^2_k - s^1_k)   e^{\displaystyle{\tilde{c}_v \delta y }}  +  
          L_s | V^1_s (0)  - V^2_s (0) |  e^{\displaystyle{  - c y / 2   }}  \Big] \delta   e^{\displaystyle{- c y / 2 }} dt \Bigg\} \delta \\
&      \leq   \unpo \delta \Big(  | V^1_s (0)  - V^2_s (0) |    +   (s^2_k - s^1_k)    \Big)   e^{\displaystyle{- c x /4 }} .   \phantom{\Bigg\{}    
\end{split}
\end{equation*}
By analogous considerations, one can conclude that the second term in  \eqref{e:Binv:char:badterm}  can be bounded by $ \unpo \delta \Big(  | V^1_s (0)  - V^2_s (0) |    +   (s^2_k - s^1_k)    \Big)   e^{\displaystyle{- c x /4 }} $. 

One can also perform similar estimates on all the other terms that appear in
 \begin{equation*}
 \begin{split}     
      \big\|    T_U \big(  V^1_s (0), \, s^1_k,   \, U,   \, p, \, q            \big)      -     &   T_U \big(  V^2_s (0), \, s^2_k,   \, U,   \, p, \, q            \big)               \big\|_U, \\
&                            \big\|    T_q \big(  V^1_s (0), \, s^1_k,   \, U,   \, p, \, q            \big)      -        T_q \big(  V^2_s (0), \, s^2_k,   \, U,   \, p, \, q            \big)               \big\|_q  \\    
\end{split}
\end{equation*}
and 
$$
                 \big\|    T_p  
                    \big(  V^1_s (0), \, s^1_k,   \, U,   \, p, \, q            \big)      -        T_p   \big(  V^2_s (0), \, s^2_k,   \, U,   \, p, \, q            \big)               \big\|_p 
$$                    
and conclude that \eqref{e:Binv:char:per:lip2} holds. 
\end{proof}
\end{enumerate}

Set $V_s (0) = (s_1 \dots s_{k -1})$ and define 
\begin{equation}
\label{e:Binv:char:F}
     F (\bar{u}_k, \, s_1, \dots s_k ) : =   F^k (  \bar{u}_k, \, s_k )  +  F^s ( \bar{u}_k, \, s_1, \dots s_{k - 1} ) + F^p ( \bar{u}_k, \, s_1, \dots s_{k} ).
\end{equation}     
Combining \eqref{e:Binv:char:center:lip}, \eqref{e:Binv:char:stable:lip} and \eqref{e:Binv:char:per:lip} one gets that $F$ is Lipschitz continuous with respect to $(s_1 \dots s_k)$:
\begin{equation}
\label{e:Binv:char:totlips}
         | F(\bar{u}_k, \, s_1, \dots s_k) - F ( \bar{u}_k, \, \tilde{s}_1, \dots \tilde{s}_k )   | \leq L \big(  |s_1 - \tilde{s}_1| + \dots | s_k - \tilde{s}_k  |  \big)
\end{equation}
for a suitable constant $L$.

Moreover, $F$ is differentiable at $(s_1 \dots s_k) = \vec{0}$ because of \eqref{e:Binv:kdiff} , \eqref{e:Binv:char:stable:diff} and \eqref{e:Binv:char:per:diff}. The colums of the jacobian matrix are the eigenvectors $r_1 (\bar{u}_k) \dots r_k (\bar{u}_k  )$.  

Finally, $F$ is also Lipschitz continuous with respect to $\bar{u}_k$, namely
\begin{equation}
\label{e:Binv:char:totlipu}
          | F ( \bar{u}^1_k,  \,  s_1 \dots s_k ) -  F ( \bar{u}^2_k, \, s_1 \dots s_k ) | \leq \tilde{L} | \bar{u}^1_k - \bar{u}^2_k  |.
\end{equation}
Define 
\begin{equation}
\label{e:Binv:char:phi}
         \phi ( \bar{u}_0, \, s_1 \dots s_N) : = F \Big( T^{k +1}_{s_{  k + 1}  } \circ \dots \circ T^N_{s_N} (\bar{u}_0) , \, s_1 \dots s_k      \Big). 
\end{equation}
We recall that the composite map $ T^{k +1}_{s_{  k + 1}  } \circ \dots \circ T^N_{s_N} (\bar{u}_0)  $ is Lipschitz continuous and differentiable at $(s_{k+1} \dots s_N) = (0 \dots 0)$, with 
the columns of the jacobian given by the eigenvectors $r_{k +1}  (\bar{u}_0) \dots r_N (\bar{u}_0)$.  This was proved in \cite{Bia:riemann}. 

Combining this with \eqref{e:Binv:char:totlipu}  and \eqref{e:Binv:char:phi} one gets that $\phi$ is Lipschitz continuous. Moreover, it is differentiable at $(s_1 \dots s_N)$ and the columns of the jacobian are the eigenvectors $r_1 (\bar{u}_0) \dots r_N (\bar{u}_0)$. Thus, thanks to \eqref{eq_Binv_direct_sum}, the jacobian is invertible and hence, exploiting the extension of the implicit function theorem discussed in \cite{Cl} (page 253), the map $\phi ( \bar{u}_0, \cdot) $ is invertible in a neighbourhood of $(s_1 \dots s_N) = ( 0 \dots 0 )$. 

Thus, if $\bar{u}_b$ is fixed and is sufficiently close to $\bar{u}_0$, then the values of $s_1 \dots s_N$ are uniquely determined by the equation 
\begin{equation}
\label{e:Binv:char:s1sN}
    \bar{u}_b = \phi (\bar{u}_0, \, s_1 \dots s_N).
\end{equation}
We will take $\bar{u}_b$ equal to the boundary datum imposed on the parabolic approximation \eqref{eq_Binv_the_system}. Once $(s_1 \dots s_N)$ are determined one can determine the values of the hyperbolic limit  $u (t, \, x)$ a.e. $(t, \, x)$. This can be done glueing together pieces like \eqref{e:Binv:s+:limit}. Here we will just make a remark concerning the trace $\bar{u}$ of the hyperbolic limit on the axis $x =0$. Let $s_k$ determined by equation \eqref{e:Binv:char:s1sN} and let $(u_k, \, v_k, \, \sigma_k)$ solve the fixed  point problem \eqref{eq:Binv:char:centerpuret}. Define
$$
    s: = \min \{ s \in \,  [0, \, s_k] :  \; \sigma_k (s) = 0    \}, 
$$
then the trace $\bar{u}$ of the hyperbolic limit on $x = 0$ is  
\begin{equation}
\label{e:Binv:char:trace}
    \bar{u} = u_k (\bar{s})
\end{equation}

The following theorem collects the results discussed in this section. 
\begin{teo}
\label{t:Binv:char}
             Let Hypotheses \ref{hyp_invertible}, \ref{hyp_hyperbolic_I}, \ref{hyp_convergence}, \ref{hyp_small_bv}, \ref{hyp_stability} , \ref{hyp_finite_pro} and \ref{hyp_char} hold. Then 
             there  exists $\delta >0$ small enough such that the following holds. If $|  \bar{u}_0 - \bar{u}_b| << \delta$, then the limit of the parabolic approximation 
             \eqref{eq_Binv_the_system} 
            satisfies 
            $$
               \bar{u}_b = \phi (\bar{u}_0, \, s_1 \dots s_N)
        $$
         for a suitable vector $(s_1 \dots s_N)$. The map $\phi$ is given by 
         $$
               \phi ( \bar{u}_0, \, s_1 \dots s_N) : = F \Big( T^{k +1}_{s_{  k + 1}  } \circ \dots \circ T^N_{s_N} (\bar{u}_0) , \, s_1 \dots s_k      \Big),
         $$
         where $T^{k+1}_{s_{k+1}} \dots T^N_{s_N}$ are the curves of admissible states defined in \cite{Bia:riemann}, while the function $F$ is defined by \eqref{e:Binv:char:F}. 
         The map $\phi$ is locally invertible: given $\bar{u}_0$ and $\bar{u}_b$, one can determine uniquely 
         $     (s_1 \dots s_N)$. Once $     (s_1 \dots s_N)$ are known, the value $u(t, \, x)$ assumed by the limit is determined for a.e. $(t, \, x)$. In particular, the trace $\bar{u}$ of the hyperbolic limit in the axis $x =0$ is given by \eqref{e:Binv:char:trace}.
\end{teo}

\section{The characterization of the hyperbolic limit in the case of a singular viscosity matrix}
\label{sec_B_non_invertible}

The aim of this section is to describe the limit of the parabolic approximation 
\eqref{eq_introduction_parabolic_approx}  when the viscosity matrix $\tilde{B}$ is not invertible. The precise hypotheses assumed in this case are introduced in Sections 
\ref{subsec_hypotheses_introduction} and
\ref{subsec_hypotheses_non_inv}. In particular, they guarantee
that it is sufficient to study the family of systems
\begin{equation}
\label{eq_Bsing_parabolic_approx}
      \left\{
      \begin{array}{ll}
             E(\ue) \ue_t + A(\ue, \, \ee \ue_x) \ue_x = \ee
             B(\ue) \ue_{xx} \qquad \ue \in \mathbb{R}^N \\
             \B(\ue (t, \, 0)) \equiv \bar g \qquad
             \ue (0, \, x) \equiv \bar{u}_0, \\
      \end{array}
      \right.
\end{equation}
with $E$, $A$ and $B$ satisfying Hypothesis \ref{hyp_Bsing}, $\B$
given by Definition \ref{def_bc}.

The exposition is organized as follows. 
In Section
 \ref{sub_Bsing_preliminary}  we introduce some preliminary results, while
 in Section \ref{sub_Bsing_Riemann_solvernonchar} we give a characterization of the limit of \eqref{eq_Bsing_parabolic_approx}. More precisely,
 in Section \ref{subsub_Bsing_travelling_waves} we briefly recall the characterization of the limit in the case of the Cauchy problem:
\begin{equation*}
      \left\{
      \begin{array}{ll}
             E(\ue) \ue_t + A(\ue, \, \ee \ue_x) \ue_x = \ee
             B(\ue) \ue_{xx} \\
             \ue (0, \, x) =
             \left\{
             \begin{array}{ll}
                       \bar{u}_0 &  x \ge 0 \\
                       u^- & x < 0  \\
              \end{array}
              \right.
              \\
      \end{array}
      \right.
\end{equation*} 
 We refer to \cite{Bia:riemann} for the complete analysis. In Section 
 \ref{subsub_Bsing_dimension} we introduce a lemma on the
 dimension of the stable manifold of the equation 
 $$
    A(u, \,  u_x) u_x = 
             B(u) u_{xx},
 $$
 which is the equation satisfied by the steady solutions of \eqref{eq_Bsing_parabolic_approx}. 
 Such a lemma
 gives an answer to question which was left open in \cite{Rousset:char}. In Sections \ref{subsub_Bsing_nonchar}  we give a characterization of the limit in the non characteristic case, i.e. when none of the eigenvalues of $E^{-1}A$ can attain the value zero. The boundary characteristic case occurs when an eigenvalue of $E^{-1}A$ can attain the value $0$ and it is discussed 
 in Section \ref{subsub_Bsing_char}. Finally, in Section \ref{sub_Bsing_boundary_datum} we discuss a transversality lemma, a technical result which is needed to complete the analysis in Sections  \ref{subsub_Bsing_nonchar} and  \ref{subsub_Bsing_char}.

\subsection{Preliminary results}
\label{sub_Bsing_preliminary} This section introduces the
preliminary results needed in Section
\ref{sub_Bsing_Riemann_solvernonchar}. The exposition is organized
as follows: Section \ref{subsub_Bsing_prel_lemmata} gives some
results about the structure of the matrix $A(u, \, u_x)$, all
implied by Hypothesis \ref{hyp_Bsing} and in particular by
Kawashima condition. In Section \ref{subsub_Bsing_prel_normal} we rewrite equations 
\begin{equation}
\label{e:Bsing:pre:tw}   
   E(u) u_t + A(u, \,  u_x) u_x = 
             B(u) u_{xx}
\end{equation}
and 
\begin{equation}
\label{e:Bsing:pre:bl}
   A(u, \,  u_x) u_x = 
             B(u) u_{xx}
\end{equation}
in a form easier to handle. 

Some notations have to be introduced: given a real parameter
$\sigma$, let
\begin{equation}
\label{eq_Bsing_prel_Esigma}
      A(u, \, u_x, \, \sigma) := A(u, \, u_x) - \sigma E(u).
\end{equation}
Consequently, both systems \eqref{e:Bsing:pre:tw} and
\eqref{e:Bsing:pre:bl} may be written in the form
\begin{equation}
\label{eq_Bsing_prel_system_gen}
      A(u, \, u_x, \, \sigma) u_x = B(u) u_{xx}.
\end{equation}
Finally, let
\begin{equation}
\label{eq_Bsing_prel_A_block} \left(
\begin{array}{ll}
      A_{11}(u, \, \sigma) & A_{12}(u, \, \sigma) \\
      A_{21}(u, \, \sigma) & A_{22}(u, \, u_x, \, \sigma)
\end{array}
\right)= \left(
\begin{array}{ll}
      A_{11}(u) - \sigma E_{11}(u) & A_{12}(u) - \sigma E_{12}(u) \\
      A_{21}(u) - \sigma E_{21}(u) & A_{22}(u, \, u_x) -
      \sigma E_{22}(u)
\end{array}
\right)
\end{equation}
and
\begin{equation}
\label{eq_Bsing_prel_wz}
      \Xi = (w, \, z)^T, \quad w
      \in \mathbb{R}^{N-r}, \; z \in \mathbb{R}^{r}
\end{equation}
be the block decompositions of $A(u, \, u_x, \, \sigma)$
 and of $\Xi \in \mathbb{R}^N$ corresponding to \eqref{eq_hypotheses_Bsing_B}.
\subsubsection{Some results about the structure of $A(u, \, u_x, \, \sigma)$ }
\label{subsub_Bsing_prel_lemmata} Even if it is not explicitly
specified in the statement, all the results of this section
presuppose that Hypothesis \ref{hyp_Bsing} holds.

The first observation is immediate, but useful:
\begin{lem}
\label{lem_Bsing_symmetry}
      Let $A_{i j}(u, \, u_x, \, \sigma)$ be the blocks of $A(u, \, u_x, \, \sigma)$
      defined by
      \eqref{eq_Bsing_prel_A_block}. Then
      \begin{equation*}
             A_{12} (u, \, \sigma) = A_{21}^T(u, \, \sigma)
              \qquad A_{11}^T (u, \, \sigma ) = A_{11}(u, \, \sigma)
      \end{equation*}
\end{lem}

Some further notations are required: let $P_0(u, \, \sigma)$
denote the projection of $\mathbb{R}^{N-r}$ on $\text{ker} \,
A_{11}(u, \, \sigma)$. Thanks to the third condition in Hypothesis
\ref{hyp_Bsing}, the dimension of this subspace is constant in
$u$: we will denote it by $q$. The operator $P_0$ can be
identified with a matrix $P_0 (u, \, \sigma ) \in \mathbb{M}^{q
\times (N-r)}$. Similarly, the projection on the subspace
orthogonal to $\text{ker} A_{11}(u, \, \sigma)$ is identified with
a matrix $P_{\bot} (u, \, \sigma ) \in \mathbb{M}^{(N-r-q) \times
(N-r)}$. The decomposition \eqref{eq_Bsing_prel_wz} is hence
refined setting
\begin{equation}
      \bar{w} = P_0  w
      \qquad
      \tilde{w} = P_{\bot} w
      \qquad
      \bar{w} \in \mathbb{R}^{q}, \; \; \tilde{w} \in
      \mathbb{R}^{N-r-q}.
\end{equation}

Using the previous notations and recalling Lemma
\ref{lem_Bsing_symmetry}, the product $A (u, \, u_x, \, \sigma)
\vec{\Xi}$ writes
\begin{equation}
\label{eq_Bsing_prel_A_blockII}
      \left(
             \begin{array}{cc}
                   A_{11}(u, \, \sigma)  & A_{21}(u, \, \sigma)^T \\
                   A_{21}(u, \, \sigma) & A_{22}(u, \, u_x, \, \sigma) \\
             \end{array}
      \right)
       \left(
             \begin{array}{cc}
                   w \\
                   z \\
             \end{array}
      \right) =
      \left(
             \begin{array}{ccc}
      0 & 0                                &  \big( A^{I}_{21} \big) ^T  (u, \, \sigma)\\
      0 & \tilde{A}_{11}(u, \, \sigma)     & \big( A^{II}_{21}\big)^T (u, \, \sigma)\\
      A^{I}_{21}(u) &  A^{II}_{21}(u, \, \sigma) & A_{22}(u, \, \sigma, \, u_x) \\
             \end{array}
      \right)
      \left(
             \begin{array}{ccc}
                   \bar{w} \\
                   \tilde{w}\\
                    z \\
             \end{array}
      \right),
\end{equation}
where
\begin{equation*}
       A^{I}_{21}  (u, \, \sigma) = A_{21}  (u, \, \sigma)
       P_0^T
      (u, \, \sigma)  \in \mathbb{M}^{r \times q}
      \qquad \qquad
       A^{II}_{21}  (u, \, \sigma) =  A_{21}  (u, \, \sigma)
       P_{\bot}^T
      (u, \, \sigma)  \in \mathbb{M}^{r \times (N-r-q)}
\end{equation*}
and
\begin{equation*}
       \tilde{A}_{11}  (u, \, \sigma) =    P_{\bot}
      (u, \, \sigma)  A_{11}  (u, \, \sigma)
      P_{\bot}^T(u, \, \sigma) \in \mathbb{M}^{(N-r-q) \times (N-r-q)}
\end{equation*}
It is now possible to state the following result, a consequence of
Kawashima condition:
\begin{lem}
\label{lem_Bsing_kawashima_implies}
      The matrix $A^{I}_{21}  (u, \, \sigma) \in \mathbb{M}^{r \times q}$
      has rank $q$. Thus, in
      particular, $q \leq r$.
\end{lem}
\begin{proof}
First of all, one observes that
\begin{equation*}
      E^{-1} A(u, \, 0) \Xi = \sigma \Xi
      \iff A (u, \, 0, \, \sigma ) \Xi =0
\end{equation*}
and hence Kawashima condition writes
\begin{equation}
\label{eq_Bsing_prel_kaw_equival}
      \text{ker} \big( A (u, \, 0, \, \sigma ) \big) \cap
      \text{ker} \big( B(u) \big) = \big\{ 0 \big\}
      \qquad \forall \, u, \, \sigma.
\end{equation}
We claim that this implies
\begin{equation}
\label{eq_Bsing_kawashima_implies}
       {A}^I_{21} \bar{w} = 0 \Rightarrow \bar{w} =0.
\end{equation}
Indeed, suppose by contradiction that there exists $\bar{w} \neq
0$ such that $A^I_{21} \bar{w} = 0$. Then
\begin{equation*}
      \left(
      \begin{array}{ccc}
             0        & 0              & (A_{21}^I)^T \\
             0        & \tilde{A}_{11} & ({A}^II_{21})^T \\
             A^I_{21} & A^II_{21}      & A_{22} \\
      \end{array}
      \right)
      \left(
      \begin{array}{lll}
             \bar{w} \\
             0 \\
             0
      \end{array}
      \right)  =
       \left(
      \begin{array}{lll}
             0  \\
             0  \\
             0  \\
      \end{array}
      \right) =
      \left(
      \begin{array}{ccc}
             0 & 0 & 0 \\
             0 & 0 &  0 \\
             0 & 0 & b
      \end{array}
      \right)
      \left(
      \begin{array}{lll}
             \bar{w} \\
             0 \\
             0
      \end{array}
      \right),
\end{equation*}
which contradicts \eqref{eq_Bsing_prel_kaw_equival}. Hence
\eqref{eq_Bsing_kawashima_implies} holds.

From \eqref{eq_Bsing_kawashima_implies} one then deduces that
$A^{I}_{21}$ has $q$ independent columns and hence it admits a
submatrix of rank $q$.
\end{proof}
In the following we will isolate an invertible submatrix of
$A^{I}_{21}(u, \, \sigma)$ and we will denoted it by $a_{11}(u, \,
\sigma)$. More precisely, $a_{11}$ is obtained from $A^{I}_{21}(u,
\, \sigma)$ selecting $q$ independent rows. By continuity, one can
suppose that the rows selected are the same for every $u$.

With considerations analogous to the ones that lead to
\eqref{eq_Bsing_prel_A_blockII}, one can write
\begin{equation*}
      z = \bar{P}^T (\sigma) \bar{z}+ \tilde{P}^T(\sigma) \tilde{z},
\end{equation*}
where $\bar{P} \in \mathbb{M}^{q \times r}$ is the projection such
that $a_{11}= A^I_{21} \bar{P}^T$ and $\tilde{P} \in
\mathbb{M}^{(r- q ) \times r}$ is the projection on the subspace
orthogonal to the image of $\mathbb{R}^r$ through $\bar{P}$.

Thanks to the previous considerations, one obtains that $A(u, \,
u_x, \, \sigma ) \Xi$ writes
\begin{equation}
\label{eq_Bsing_Adivisa} \left(
\begin{array}{cccc}
       0 & 0 & a_{11}^T (u, \, \sigma)& a_{21}^T(u, \, \sigma)  \\
       0 & \tilde{A}_{11}(u, \, \sigma)
      & a_{12}^T(u, \, \sigma) & a_{22}^T(u, \, \sigma) \\
      a_{11}(u, \, \sigma) & a_{12}(u, \, \sigma) &
      \alpha_{11}(u, \, u_x, \, \sigma) &
      \alpha_{21}^T (u, \, u_x, \, \sigma)\\
      a_{21}(u) & a_{22}(u, \, \sigma)
      & \alpha_{21}(u, \, u_x, \, \sigma) &
      \alpha_{22}(u, \, u_x, \, \sigma)
\end{array}
\right) \left(
\begin{array}{cccc}
      \bar{w} \\
      \tilde{w} \\
      \bar{z} \\
      \tilde{z},
\end{array}
\right)
\end{equation}
with
\begin{equation*}
\begin{split}
&       \qquad \qquad \qquad
       a_{12} \in \mathbb{M}^{q \times (N-r-q)}, \; \;
       a_{22} \in \mathbb{M}^{(r-q) \times (N-r-q)} \\
&      \qquad \qquad \alpha_{11} \in \mathbb{M}^{q \times q}, \; \;
       \alpha_{21} \in \mathbb{M}^{(r- q) \times q}, \; \;
       \alpha_{22} \in \mathbb{M}^{(r-q) \times (r-q)}  \\
&      \qquad \qquad \qquad \qquad \qquad \qquad
       \bar{z} \in \mathbb{R}^{q} \; \; \tilde{z} \in
       \mathbb{R}^{r-q} \\
\end{split}
\end{equation*}
The corresponding decomposition of $B(u)$ is
\begin{equation}
\label{eq_Bsing_Bdivisa} \left(
\begin{array}{cccc}
      0 & 0 & 0 & 0  \\
      0 & 0 & 0 & 0 \\
      0 & 0 & b_{11}(u) & b_{12}(u) \\
      0 & 0 & b_{21}(u) & b_{22}(u) \\
\end{array}
\right)
\end{equation}
with
\begin{equation*}
      b_{11} \in \mathbb{M}^{q \times q}, \; \;
      b_{12} \in \mathbb{M}^{q \times (r-q)}, \; \;
      b_{21} \in \mathbb{M}^{(r- q) \times q}, \; \;
      b_{22} \in \mathbb{M}^{(r-q) \times (r-q)}.
\end{equation*}
It is worth underling explicitly that formulations
\eqref{eq_Bsing_Adivisa} and \eqref{eq_Bsing_Bdivisa} include in
particular the limit cases $q=0$ and $q=N-r$: indeed, it is enough
to assume $w= \tilde{w}$, $z = \tilde{z}$ and $w= \bar{w}$, $z =
\bar{z}$ respectively.

\begin{lem}
\label{lem_Bsing_inv}
      The matrix
      \begin{equation*}
            a_{21} a_{11}^{-1} \Big(  b_{11}(a_{11}^T)^{-1}
            a_{21}^T  - b_{12} \Big) - b_{21} (a_{11}^T)^{-1} a_{21}^T
            + b_{22}
      \end{equation*}
       is invertible.
\end{lem}
\begin{proof}
By contradiction, assume that there exists an $(r-q)$-dimensional
vector $\vec{\xi} \neq 0$ such that
\begin{equation*}
      \Big( a_{21} a_{11}^{-1} \Big(  b_{11}(a_{11}^T)^{-1}
      a_{21}^T  - b_{12} \Big) - b_{21} (a_{11}^T)^{-1} a_{21}^T
      + b_{22}
              \Big)\vec{ \xi} = 0.
\end{equation*}
It follows that
\begin{equation}
\label{eq_Bsing_trasposta2}
      \left(
             a_{21} a_{11}^{-1} , \,    - I_{r-q} 
      \right)
      \left(
      \begin{array}{cc}
              b_{11} &  b_{12} \\
            b_{21} &  b_{22}
      \end{array}
      \right)
      \left(
      \begin{array}{cc}
             \big( a_{11}^{T} \big)^{-1}a_{21}^T  \\
           -  I_{r-q}
      \end{array}
      \right)
       \vec \xi  = 0 
\end{equation}
Define 
\begin{equation*}
      D =
      \left(
      \begin{array}{cc}
             ( a_{11}^{-1})^T a_{21}^T  \\
             I_{r-q}
      \end{array}
      \right),
\end{equation*}
then \eqref{eq_Bsing_trasposta2} implies 
\begin{equation*}
       0 = \langle D \vec \xi , \,  b D \vec \xi \rangle .
\end{equation*}
Since $D\, \xi \neq \vec 0$, this is a contradiction because of the second condition in Hypothesis \ref{hyp_Bsing}.
\end{proof}
\subsubsection{The explicit form of system \eqref{eq_Bsing_prel_system_gen}}
\label{subsub_Bsing_prel_normal} The following result guarantees
that system \eqref{eq_Bsing_prel_system_gen}, and hence systems
\eqref{e:Bsing:pre:tw} and
\eqref{e:Bsing:pre:bl} can be reduced to an
explicit form like \eqref{eq_Bsing_explicit_form}
\begin{lem}
\label{lem_Bsing_explicit_form}
       There exists a sufficiently small constant $\delta >0$ such that if $|u - \bar{u}_0|< \delta$ and $|u_x| < \delta$, then the following holds. System \eqref{eq_Bsing_prel_system_gen} can be rewritten as 
      \begin{equation}
      \label{eq_Bsing_explicit_form}
      \left\{
      \begin{array}{ll}
            u_x = \Big( \bar{w}(u, \, \tilde{z}, \, \sigma), \;
                  \tilde{w}(u, \, \tilde{z}, \, \sigma), \;
                  \bar{z}(u, \, \tilde{z}, \, \sigma), \; \tilde{z}
                  \Big)^T \\
            \tilde{z}_x = f(u, \, \tilde{z}, \, \sigma)
      \end{array}
      \right.
      \end{equation}
      for suitable functions
      \begin{equation*}
      \begin{split}
      &      \bar{w}:
             \mathbb{R}^{N} \times \mathbb{R}^{r-q} \times \mathbb{R} \to
             \mathbb{R}^{q} \qquad
             \tilde{w}: \mathbb{R}^{N} \times \mathbb{R}^{r-q} \times \mathbb{R}
             \to
             \mathbb{R}^{N-r-q} \qquad
             \bar{z}: \mathbb{R}^{N} \times \mathbb{R}^{r-q} \times \mathbb{R}
              \to
             \mathbb{R}^{q} \\
      &     \qquad \qquad \qquad \qquad
            f: \mathbb{R}^{N} \times \mathbb{R}^{r-q} \times \mathbb{R} \to
            \mathbb{R}^{r-q}. \\
      \end{split}
      \end{equation*}
      For every $\sigma$ and for every $u$ the following holds. 
      \begin{equation}
      \label{e:Bsing:prel:zerozero}
            \bar w (u, \, \vec 0, \, \sigma ) = \vec 0 \qquad \tilde w (u, \, \vec 0, \, \sigma ) = \vec 0
	    \qquad \bar z  (u, \, \vec 0, \, \sigma ) = \vec 0 \qquad f  (u, \, \vec 0, \, \sigma ) = \vec 0  
      \end{equation}
     Also, these functions satisfy
      \begin{equation}
      \label{eq_Bsing_barz}
      \begin{split}
      &      D_u \bar{z} \big|_{u=\bar{u}_0, \, \tilde{z}=0, \, \sigma} = 0
            \qquad
            D_{\sigma} \bar{z} \big|_{u=\bar{u}_0, \, \tilde{z}=0, \, \sigma} = 0
            \qquad
            D_{\tilde{z}} \bar{z} \big|_{u=\bar{u}_0, \, \tilde{z}=0} =
            - (a_{11}^T)^{-1} a_{21}^T  \phantom{\bigg(}\\
       \end{split}
      \end{equation}
      and
      \begin{equation}
      \label{eq_Bsing_tildew}
      \begin{split}
      &     \qquad \qquad
            D_u \tilde{w} \big|_{u=\bar{u}_0, \, \tilde{z}=0, \, \sigma} = 0
            \qquad
            D_{\sigma} \tilde{w} \big|_{u=\bar{u}_0, \, \tilde{z}=0, \, \sigma} = 0
            \\
      &     D_{\tilde{z}} \tilde{w} \big|_{u=\bar{u}_0, \, \tilde{z}=0, \, \sigma}
            = - \tilde{A}_{11}^{-1} \big( a_{12}^T D_{\tilde{z}}
            \bar{z} +  a_{22}^T \big) =
            - \tilde{A}_{11}^{-1} \big( - a_{12}^T  (a_{11}^T)^{-1} a_{21}^T
            +  a_{22}^T \big) . \\
      \end{split}
      \end{equation}
      Moreover,
      \begin{equation}
      \label{eq_Bsing_barw}
      \begin{split}
      &     \qquad \qquad \qquad \qquad \qquad \qquad
            D_u \bar{w} \big|_{u=\bar{u}_0, \, \tilde{z}=0, \, \sigma} = 0
            \qquad
            D_{\sigma} \bar{w} \big|_{u=\bar{u}_0, \, \tilde{z}=0, \, \sigma} = 0
            \\
      &     D_{\tilde{z}} \bar{w} \big|_{u=\bar{u}_0, \, \tilde{z}=0} =
            - a_{11}^{-1}
            \Big(
                 a_{12} D_{\tilde{z}} \tilde{w} +
                 \alpha_{11} D_{\tilde{z}} \bar{z} +
                 \alpha_{12}^T
            \Big) +
            a_{11}^{-1} b_{11} D_{\tilde{z}}
            \bar{z} D_{\tilde{z}} f +
            a_{11}^{-1} b_{12} D_{\tilde{z}} f  \\
      &     \qquad \qquad \qquad = a_{11}^{-1} \Big( a_{12}
            \tilde{A}_{11}^{-1} \big(
                            a_{22}^T - a_{12}^T
                            (a_{11}^T)^{-1}
                            a_{21}^T
            \big)+ \alpha_{11}(a_{11}^T)^{-1} a_{21}^T -
            \alpha_{12} \Big)
            + a_{11}^{-1} \bigg( b_{12} \\
      &     \qquad \qquad \qquad \;  - b_{11}(a_{11}^T)^{-1}
            a_{21}^T  \bigg) \bigg( a_{21} a_{11}^{-1} \Big(  b_{11}(a_{11}^T)^{-1}
            a_{21}^T  - b_{12} \Big) - b_{21} (a_{11}^T)^{-1}
            a_{21}^T \\
       &     \qquad \qquad  \qquad \;
            + b_{22} \bigg)^{-1} \bigg( a_{21} a_{11}^{-1} \Big( a_{12}
            \tilde{A}_{11}^{-1} \big(
                            a_{22}^T -
            a_{12}^T
                            (a_{11}^T)^{-1}
                             a_{21}^T
           \big)
           +  \alpha_{11}(a_{11}^T)^{-1} a_{21}^T -
           \alpha_{12} \Big) \\
      &    \qquad \qquad \qquad \;
           + a_{22} \tilde{A}_{11}^{-1}
           \big( a_{12}^T (a_{11}^T)^{-1} a_{21}^T
           - a_{22}^T - \alpha_{21} (a_{11}^T)^{-1} a_{21}^T
           + \alpha_{22 } \big) \bigg) \\
      \end{split}
      \end{equation}
      and
      \begin{equation}
      \label{eq_Bsing_f}
      \begin{split}
      &      D_u f \big|_{u=\bar{u}_0, \, \tilde{z}=0, \, \sigma} = 0
            \qquad
            D_{\sigma} f \big|_{u=\bar{u}_0, \, \tilde{z}=0, \, \sigma} =
            0 \\
      &     \qquad
            D_{\tilde{z}} f \big|_{u=\bar{u}_0, \, \tilde{z}=0} =
            \bigg( b_{21} D_{\tilde{z}} \bar{z}  + b_{22} - a_{21}
            \Big( a_{11}^{-1} b_{11} D_{\tilde{z}} \bar{z}
             +
            a_{11}^{-1} b_{12} \Big) \bigg)^{-1}
            \bigg(
                 - a_{21} a_{11}^{-1}
            \Big(
                 a_{12} D_{\tilde{z}} \tilde{w} +
                 \alpha_{11} D_{\tilde{z}} \bar{z}  \\
      &     \qquad \qquad \qquad \qquad \qquad +
                 \alpha_{12}^T
            \Big) +
            a_{22} D_{\tilde{z}} \tilde{w}  +
            \alpha_{21} D_{\tilde{z}} \bar{w} +
            \alpha_{22} \bigg) \\
      &     \;  \qquad \qquad \qquad \qquad
            = \bigg( a_{21} a_{11}^{-1} \Big(  b_{11}(a_{11}^T)^{-1}
            a_{21}^T  - b_{12} \Big) - b_{21} (a_{11}^T)^{-1} a_{21}^T
            + b_{22} \bigg)^{-1} \bigg( a_{21} a_{11}^{-1} \Big( a_{12}
            A_{11}^{-1} \big(
                            a_{22}^T \\
      &     \qquad \qquad  \qquad \qquad \qquad - a_{12}^T
                            (a_{11}^T)^{-1}
                             a_{21}^T
           \big)
           +  \alpha_{11}(a_{11}^T)^{-1} a_{21}^T -
           \alpha_{21}^T \Big)
           + a_{22} \tilde{A}_{11}^{-1}
           \Big( a_{12}^T (a_{11}^T)^{-1} a_{21}^T \\
      &     \qquad \qquad  \qquad \qquad \qquad
           - a_{22}^T \Big) - \alpha_{21} (a_{11}^T)^{-1} a_{21}^T
           + \alpha_{22 }  \bigg).  \\
      \end{split}
      \end{equation}
\end{lem}
\begin{proof}
Using notation \eqref{eq_Bsing_Adivisa} and
\eqref{eq_Bsing_Bdivisa}, equation rewrites as
\begin{equation*}
\begin{split}
&  \left(
\begin{array}{cccc}
      0 & 0 & a_{11}^T (u, \, \sigma )& a_{21}^T(u, \, \sigma)  \\
      0 & \tilde{A}_{11}(u, \, \sigma) & a_{12}^T(u, \, \sigma) & a_{22}^T(u\, \sigma) \\
      a_{11}(u, \, \sigma) & a_{12}(u, \, \sigma) &
      \alpha_{11}(u, \, \bar{w},  \, \tilde{w}, \, \bar{z}, \, \tilde{z}, \, \sigma) &
      \alpha_{21}^T (u, \, \bar{w},  \, \tilde{w}, \, \bar{z}, \, \tilde{z}, \, \sigma)\\
      a_{21}(u, \, \sigma) & a_{22}(u, \, \sigma)
      & \alpha_{21}(u, \, \bar{w},  \, \tilde{w}, \, \bar{z}, \, \tilde{z}, \, \sigma) &
      \alpha_{22}(u, \, \bar{w},  \, \tilde{w}, \, \bar{z}, \, \tilde{z}\, \sigma) \\
\end{array}
\right) \left(
\begin{array}{cccc}
      \bar{w}  \\
      \tilde{w} \\
      \bar{z} \\
      \tilde{z}
\end{array}
\right) \\
& =  \left(
\begin{array}{cccc}
      0 & 0 & 0 & 0  \\
      0 & 0 & 0 & 0 \\
      0 & 0 & b_{11}(u) & b_{12}(u) \\
      0 & 0 & b_{21}(u) & b_{22}(u) \\
\end{array}
\right) \left(
\begin{array}{cccc}
      \bar{w}_x  \\
      \tilde{w}_x \\
      \bar{z}_x \\
      \tilde{z}_x
\end{array}
\right) \\
\end{split}
\end{equation*}
Hence from the first line it follows
\begin{equation*}
      a_{11}^T(u, \, \sigma) \bar{z}+
      a_{21}^T(u, \, \sigma) \tilde{z} =0
\end{equation*}
and therefore
\begin{equation*}
      \bar{z}(u, \, \tilde{z}) = - \big( a_{11}^T(u, \, \sigma) \big)^{-1}
      a_{21}^T(u, \, \sigma) \tilde{z} .
\end{equation*}
The second line then reads
\begin{equation*}
      \tilde{A}_{11}(u, \, \sigma)
      \tilde{w}- a_{12}^T (u)\big( a_{11}^T(u, \, \sigma) \big)^{-1}
      a_{21}^T(u, \, \sigma) \tilde{z}
      + a_{22}^T (u, \, \sigma) \tilde{z} =0
\end{equation*}
and hence thanks to the invertibility of $\tilde{A}_{11}(u)$
\begin{equation*}
      \tilde{w}(u, \, \tilde{z}) =
      \tilde{A}^{-1}_{11}(u, \, \sigma)
      \Big( a_{12}^T (u)\big( a_{11}^T(u, \, \sigma) \big)^{-1}
      a_{21}^T(u, \, \sigma)
      + a_{22}^T (u, \, \sigma) \Big) \tilde{z}.
\end{equation*}
Moreover, one has
\begin{equation*}
      \big( \bar{z}(u, \, \tilde{z}) \big)_x =
      - \big( a_{11}^T(u, \, \sigma) \big)^{-1}
      a_{21}^T(u, \, \sigma) \tilde{z}_x -
      D_u \Big( \big( a_{11}^T(u, \, \sigma) \big)^{-1}
      a_{21}^T(u, \, \sigma) \Big)
      u_x =
      \bar{z}_x (u, \, \bar{w}, \, \tilde{z}, \, \tilde{z}_x)
\end{equation*}
and hence the third line reads
\begin{equation*}
      a_{11}(u) \bar{w} + a_{12}(u) \tilde{w}(u, \, \tilde{z}) +
      \alpha_{11} (u, \, \bar{w}, \, \tilde{z}) \bar{z}(u, \,
      \tilde{z})+ \alpha_{21}^T(u, \, \bar{w}, \, \tilde{z})
      \tilde{z} = b_{11}(u) \bar{z}_x
      (u, \, \bar{w}, \, \tilde{z}, \, \tilde{z}_x) +
      b_{12}(u) \tilde{z}_x
\end{equation*}
and therefore, thanks to the invertibility of $a_{11}$, in a small
enough neighborhood of \\ ${(u=\bar{u}_0, \, \bar{w}=0, \,
\tilde{z}=0, \, \tilde{z}_x=0)}$ it is implicitly defined a map
$\bar{w}= \bar{w}(u, \, \tilde{z}, \, \tilde{z}_x)$ such that
\begin{equation*}
\begin{split}
&     D_u \bar{w} \big|_{u=\bar{u}_0, \, \tilde{z}=0, \,
      \tilde{z}_x=0}=0 \qquad \qquad
      D_{\tilde{z}_x} \bar{w} \big|_{u=\bar{u}_0, \, \tilde{z}=0, \,
      \tilde{z}_x=0}= a_{11}^{-1} \Big( b_{12} - b_{11}(a_{11}^T)^{-1}
      a_{21}^T  \Big) \\
&     D_{\tilde{z}} \bar{w} \big|_{u=\bar{u}_0, \, \tilde{z}=0, \,
      \tilde{z}_x=0}= a_{11}^{-1} \Big( a_{12}
      A_{11}^{-1} \big(
                      a_{22}^T - a_{12}^T
                      (a_{11}^T)^{-1}
                      a_{21}^T
      \big)+ \alpha_{11}(a_{11}^T)^{-1} a_{21}^T -
      \alpha_{21}^T \Big). \\
\end{split}
\end{equation*}
Then the fourth line reads
\begin{equation*}
      a_{21}(u) \bar{w}(u, \, \tilde{z}, \, \tilde{z}_x )+
      a_{22}(u) \tilde{w}(u, \, \tilde{z})+
      \alpha_{21} ( u, \, \tilde{z}, \, \tilde{z}_x)
      \bar{z}_x ( u, \, \tilde{z}, \, \tilde{z}_x) +
      \alpha_{22} ( u, \,  \tilde{z}, \, \tilde{z}_x )
      \tilde{z} =
      b_{21} (u) \bar{z}_x
      ( u, \,  \tilde{z}, \, \tilde{z}_x ) + b_{22}
      (u) \tilde{z}_x
\end{equation*}
Hence thanks to Lemma \ref{lem_Bsing_inv}, it is implicitly
defined a map $\tilde{z}_x = f(u, \, \tilde{z})$, which satisfies
the hypotheses described in the statement of the lemma.
\end{proof}
In order to simplify the notations, we set
\begin{equation}
\label{eq_Bsing_prel_underline}
\begin{split}
&     \underline{a} (u, \, \tilde{z}, \, \sigma): =
      a_{21} a_{11}^{-1} \Big( a_{12}
            \tilde{A}_{11}^{-1} \big(
                            a_{22}^T -
                            a_{12}^T
                            (a_{11}^T)^{-1}
                             a_{21}^T
           \big)
           +  \alpha_{11}(a_{11}^T)^{-1} a_{21}^T -
           \alpha_{21}^T \Big)
           + a_{22} \tilde{A}_{11}^{-1}
           \big( a_{12}^T (a_{11}^T)^{-1} a_{21}^T \\
     &     \qquad \qquad
           - a_{22}^T - \alpha_{21} (a_{11}^T)^{-1} a_{21}^T
           + \alpha_{22 } \big)  \in \mathbb{M}^{(r-q) \times (r-q)} \\
     &     \underline{b}(u,  \, \sigma):=
           a_{21} a_{11}^{-1} \Big(  b_{11}(a_{11}^T)^{-1}
            a_{21}^T  - b_{12} \Big) - b_{21} (a_{11}^T)^{-1} a_{21}^T
            + b_{22} \in \mathbb{M}^{(r-q) \times (r-q)}
\end{split}
\end{equation}
With this notations,
\begin{equation*}
       D_{\tilde{z}} f \Big|_{u= \bar{u}_0, \, \tilde{z}=0, \, \sigma}=
       \underline{b}^{-1}(\bar{u}_0, \, 0, \, \sigma )
       \underline{a}(\bar{u}_0, \, 0, \, \sigma ).
\end{equation*}
Also, if we consider the jacobian
$$
  D_{\tilde z} \Big( \bar w(u, \, \tilde z, \, \sigma) , \, \tilde w (u, \, \tilde z, \, \sigma), \, 
  \bar z (u, \, \tilde z, \, \sigma), \, \tilde z (u, \, \tilde z, \, \sigma) \Big)
$$
and we compute it at the point $\Big( u= \bar{u}_0, \, \tilde{z}=0, \, \sigma)$, we get  
\begin{equation}
\label{e:Bsing:dzu}
       \left(
       \begin{array}{cccc}
                a_{11}^{-1} a_{12} \tilde{A}_{11}^{-1} \Big( a_{22}^T - a_{12}^T
                   (a_{11}^T)^{-1}   a_{21}^T  \Big)\vec{\xi} +
                   a_{11}^{-1} \alpha_{11}(a_{11}^T)^{-1} a_{21}^T \vec{\xi}-
                   a_{11}^{-1} \alpha_{21}^T \vec{\xi} +
                   \Big( a_{11}^{-1} b_{12} - a_{11}^{-1}
                   b_{11} (a_{11}^{-1})^T a_{21}^T \Big) \underline{b}^{-1} \underline{a} \\
	  \tilde{A}_{11}^{-1} \Big( a_{12}^T (a_{11}^T)^{-1} a_{21}^T
                   - a_{22}^T \Big) \\	   
       - (a_{11}^T)^{-1} a_{21}^T \\
       I_{r - q}
       \end{array}
       \right)
\end{equation}
Finally, the proof of Lemma \ref{lem_Bsing_explicit_form} implies the
following result:
\begin{lem}
\label{lem_Bsing_exchange_eigenvalues}
       Given $\vec{\xi} \in \mathbb{R}^{r-q}$, the condition
       $\underline{a}(u, \, 0, \, \sigma ) \vec{\xi} =
       \lambda \underline{b}(u, \sigma ) \vec{\xi} $
       is equivalent to ${A (u, \, 0, \, \sigma) \; \Xi =
       \lambda B(u) \; \Xi}$,
       where $\Xi \in \mathbb{R}^N$
       is given by
       \begin{equation}
       \label{eq_Bsing_prel_exchange}
       \Xi =
       \left(
       \begin{array}{cccc}
              a_{11}^{-1} a_{12} \tilde{A}_{11}^{-1} \Big( a_{22}^T - a_{12}^T
                   (a_{11}^T)^{-1}   a_{21}^T  \Big)\vec{\xi} +
                   a_{11}^{-1} \alpha_{11}(a_{11}^T)^{-1} a_{21}^T \vec{\xi}-
                   a_{11}^{-1} \alpha_{21}^T \vec{\xi} +
                   \Big( a_{11}^{-1} b_{12} - a_{11}^{-1}
                   b_{11} (a_{11}^{-1})^T a_{21}^T \Big) \underline{b}^{-1} \underline{a} \vec{\xi} \\
              \tilde{A}_{11}^{-1} \Big( a_{12}^T (a_{11}^T)^{-1} a_{21}^T
                   - a_{22}^T \Big)\vec{\xi} \\
              - (a_{11}^T)^{-1} a_{21}^T \vec{\xi}\\
              \vec{\xi}.
               \end{array}
               \right)
       \end{equation}
       The case $\lambda =0$ and hence $\underline{a}\vec{\xi} =0$
       is, in particular, included in the previous formulation.
\end{lem}
\begin{rem}
For simplicity, in the following we assume that
$\underline{b}^{-1} \underline{a}$ is diagonalizable. Hence, in
particular, one can find $(r-q)$ independent vectors $\Xi_1, \dots
\Xi_{r-q}$ such that
\begin{equation*}
      (A- \lambda_i B) \Xi_i =0.
\end{equation*}
To handle the general case, one should work with generalized eigenspaces.
\end{rem}

\subsection{The hyperbolic limit in the case
of non invertible viscosity matrix}
\label{sub_Bsing_Riemann_solvernonchar} 

\subsubsection{The hyperbolic limit in the case of the Cauchy problem}
\label{subsub_Bsing_travelling_waves} 
In this section for completeness we review the construction of the Riemann solver for the Cauchy problem in the case the viscosity matrix  $B$ is not invertible. We refer to 
\cite{Bia:riemann} for the complete analysis. Also, some of the steps in the construction are the same as in the case 
of an invertible viscosity matrix, which is discussed in Section \ref{subsub_Binv_speed+}. Thus, here we will focus only on the points where the singularity of the viscosity matrix plays an important role.

The goal is the characterization of the limit of 
\begin{equation}
\label{e:Bsing:s+:app}
\left\{
\begin{array}{ll}
          E( \ue) \ue_t + A (\ue, \, \ue_x ) \ue_x = \ee B(  \ue ) \ue_{xx} \\
          \ue (0, \,,  x) = 
          \left\{
          \begin{array}{ll}
                     u^-    &  x \leq 0\\
                     \bar{u}_0     & x >0  \\
          \end{array}
          \right.
\end{array}
\right.
\end{equation}
under Hypotheses \ref{hyp_Bsing}, \ref{hyp_hyperbolic_I}, \ref{hyp_small_bv} and \ref{hyp_stability}.

The construction is as follows. Thanks to Lemma \ref{lem_Bsing_explicit_form}, the equation of travelling waves 
\begin{equation}
\label{e:Bsing:cauchy:tw}
    B (U ) U'' = \Big( A(U, \, U' )  - \sigma E  ( U )     \Big) U',
\end{equation}    
is equivalent to the system 
\begin{equation}
\label{eq_Binv_system_travdue} \left\{
\begin{array}{lll}
      u' = \Big( \bar{w} (u, \, \tilde{z}, \, \sigma), \; \tilde{w} (u, \, \tilde{z}, \, \sigma), \; \bar{z} (u, \, \tilde{z}, \, \sigma), \; \tilde{z}  \Big) \\
      \tilde{z}' = f (u, \, \tilde{z}, \, \sigma) \\
      \sigma' = 0.
\end{array}
\right.
\end{equation}
When $\tilde{z} =\vec{0} $ then for every $u$ and $\sigma$ $f (u, \, \vec{0}, \, \sigma) = \vec{0}$, 
$\bar{w} (u, \, \vec{0}, \, \sigma) = \vec{0}$, $ \tilde{w} (u, \, \vec{0}, \, \sigma) = \vec{0}$, $ \bar{z} (u, \, \vec{0}, \, \sigma)=\vec{0}$. Thus, $(\bar{u}_0, \, \vec{0}, \, \lambda_i (\bar{u}_0))$ is an equilibrium for 
\eqref{eq_Binv_system_trav}, where $\lambda_i (\bar{u}_0)$ is an eigenvalue of $E^{-1}(\bar{u}_0) A(\bar{u}_0, \, 0)$. Linearizing around such an equilibrium point one obtains  
\begin{equation}
\label{e:Bsing:tv:jac}
\left(
\begin{array}{ccc}
       0 & \underline{c} & 0 \\
       0 & \underline{b}^{-1} \bar{a}(\bar u_0, \,
           0, \, \sigma_i) & 0 \\
       0 & 0 & 0 \\
\end{array}
\right),
\end{equation}
where the exact expression of the matrix $\underline{c} \in
\mathbb{R}^{N \times (r-q)}$ is not important here. Let $\Xi_i (\bar{u}_0)$ denote the eigenvector
corresponding to $\lambda_i (\bar{u}_0)$, then $\Xi_i (\bar{u}_0$ is in the form given by Lemma \ref{lem_Bsing_exchange_eigenvalues} for some vector $\vec{\xi} (\bar{u}_0) \in \mathbb{R}^{r -q}$. 
Then the eigenspace of \eqref{e:Bsing:tv:jac} associated to the eigenvalue $0$ is 
$$
  V^c : = \Big\{ \Big(u, \; v_i \vec{\xi}_i (\bar{u}_0), \, \sigma \Big): \; u \in \mathbb{R}^N, \; v_i, \, \sigma \in \mathbb{R} \Big\}.
$$
One can then proceed as in the case of an invertible viscosity matrix (see Section \ref{subsub_Binv_speed+}). Fix a center manifold $\mathcal{M}^c$, then one can verify that a point$(u, \, \tilde{z}, \, \sigma)$ belongs to $\mathcal{M}^c$ if and only if 
$$
   \tilde{z} = v_i \tilde{\xi}_i (u, \, v_i, \, \sigma),
$$ 
where $\tilde{\xi}_i \in \mathbb{R}^{r -q}$ is a suitable vector valued function such that  
$$
  \tilde{\xi}_i \Big( \bar{u}_0, \, 0, \, \lambda_i (\bar{u}_0 )  \Big) = \vec{\xi}_i (\bar{u}_0). 
$$
Thus, on the center manifold 
$$
  u' = \Big( \bar{w} (u, \, \tilde{z}, \, \sigma_i), \; \tilde{w} (u, \, \tilde{z}, \, \sigma_i), \; \bar{z} (u, \, \tilde{z}, \, \sigma_i), \; \tilde{z}  \Big) = \tilde{\Xi}_i (u, \, v_i, \, \sigma_i) v_i
$$
for some vector valued function $\tilde{\Xi}_i \in \mathbb{R}^N$ such that 
$$
  \tilde{\Xi}_i \Big( \bar{u}_0, \, 0, \, \lambda_i (\bar{u}_0 )  \Big) = {\Xi}_i (\bar{u}_0).
$$
Plugging the relation $u' = \tilde{\Xi}_i v_i$ into the equation 
$$
  B(u) u'' = [ A(u, \, u') - \sigma E(u) ] u'
$$
one obtains 
$$
  [B \tilde{\Xi}_i + \tilde{\Xi}_{iv}   v_i ] v_{ix} = \Big( [ A - \sigma E ] \tilde{\Xi}_i - D_u \tilde{\Xi}_i 
  \tilde{\Xi}_i v_i \Big) v_i.
$$
Taking both members dot product $\tilde{\Xi}_i$ one gets
$$
  c_i v_{ix} = a_i v_i
$$
for suitable functions $c_i(u, \, v_i, \, \sigma)$ and $a_i(u, \, v_i, \, \sigma)$. Thanks to Lemma \ref{lem_Bsing_exchange_eigenvalues}, $\Xi_i(\bar{u}_0)$ is not in the kernel of $B$ and hence 
$$
  c_i \Big( \bar{u}_0, \, 0, \, \lambda_i (\bar{u}_0) \Big) = \langle \Xi_i(\bar{u}_0), \, B (\bar{u}_0) \Xi_i (\bar{u}_0) \rangle > 0 
$$
Thus, $c_i>0$ in a neighborhood and hence one can introduce $\phi_i:= a_i / c_i$. Also, 
\begin{equation}
\label{e:Bsing:tv:phi}
       \frac{\partial \phi_i}{\partial \sigma } \Bigg|_{ ( \bar{u}_0, \, 0, \, \lambda_i (\bar{u}_0))  } = - 
       \langle \Xi_i(\bar{u}_0), \, E \Xi_i(\bar{u}_0) \rangle < 0.
\end{equation}

Thus, system the solutions of \eqref{eq_Binv_system_trav} laying on $\mathcal{M}^c$ satisfy  
\begin{equation}
\label{eq_Bsing_dynamic_center_man} \left\{
\begin{array}{lll}
        u_i' = v_i \tilde{\Xi}_i (u, \, v_i, \, \sigma_i) \\
        v_i' = \phi_i(u, \, v_i, \, \sigma_i) \\
        \sigma_{i}' =0 \\
\end{array}
\right.
\end{equation}
with $\phi$ satisfying \eqref{e:Bsing:tv:phi}. One can therefore apply the results in \cite{Bia:riemann}. Proceeding as in Section \ref{subsub_Binv_speed+} one eventually defines the $i$-th curve of admissible states 
$T^i_{s_i} \bar{u}_0$, which satisfies 
$$
  \frac{\partial T^i \bar{u}_0}{\partial s_i} \Bigg|_{s_i =0} = \Xi_i (\bar{u}_0)
$$
and all the other properties listed in Section \ref{subsub_Binv_speed+}.  

Consider the composite function  
$$
    \psi (\bar{u}_0, \, s_1 \dots s_N) = T^1_{s_1} \circ \dots T^N_{s_N} \bar{u}_0
$$  
With the previous expression we mean that the starting point for $T^{N-1}_{s}$ is  $T^N_{s_N} \bar{u}_0$. It turns out that 
the map $\phi ( \bar{u}_0, \cdot) $ is invertible in a neighbourhood of $(s_1 \dots s_N) = ( 0 \dots 0 )$. 
In other words, if $u^-$ is fixed and is sufficiently close to $\bar{u}_0$, then the values of $s_1 \dots s_N$ are uniquely determined by the equation 
\begin{equation}
\label{e:Binv:char:s1sNdue}
    u^- = \psi (\bar{u}_0, \, s_1 \dots s_N).
\end{equation}
Taking the same $u^-$ as in \eqref{e:Bsing:s+:app}, one obtains the parameters $(s_1 \dots s_N)$ which can be used to reconstruct the hyperbolic limit $u$ of \eqref{e:Bsing:s+:app}. Indeed, once $(s_1 \dots s_N)$  are known then $u$ can be obtained in the way described in Section \ref{subsub_Binv_speed+}. 

\begin{rem}
\label{rem_twave}
          To construct the curves $T^i_{s_i} \bar{u}_0$ one proceeds as follows. Suppose that $s_i >0$, then one considers the fixed point problem
          \begin{equation}
          \label{e:Bsing:tw:rem}
 \left\{
\begin{array}{lll}
      u(\tau) = \bar{u}_k + {\displaystyle \int_0^{\tau} \tilde{r}_k(u(\xi), \, v_k(\xi), \,
      \sigma_k(\xi))d \xi }  \\
      v_k (\tau) = f_k (\tau, \, u, \, v_k, \, \sigma_k ) -
      \mathrm{conc} f_k (\tau, \, u, \, v_k, \, \sigma_k )\\
      \sigma_k(\tau)=  {  \frac{1}{c_E  (\bar{u}_0) }  \displaystyle \frac{d}{d \tau}
      \mathrm{conc} f_k (\tau, \, u, \, v_k, \, \sigma_k )}. \\
\end{array}
\right.
\end{equation}
          Denote by $(u_i, \, v_i, \, \sigma_i)$ the solution of \eqref{e:Bsing:tw:rem}, whose existence is ensured by results in \cite{Bia:riemann}. Then 
          we define 
          $$
              T^i_{s_i} \bar{u}_0 = u_i (s_i).
          $$
          The value $u_i (s_i)$ is connected by $\bar{u}_0$ by a sequence of rarefaction and shocks. This is actually true only if system \eqref{e:Bsing:tv:phi}
          is equivalent to \eqref{e:Bsing:cauchy:tw}, i.e. if Lemma \ref{lem_Bsing_explicit_form} holds true. The problem here is that in \eqref{e:Bsing:tv:phi} $\sigma$ is a parameter, while in \eqref{e:Bsing:tw:rem} can vary. Actually, it turns out that we can still apply Lemma \ref{lem_Bsing_explicit_form} provided that the dimension of the kernel of   
          $A_{11}(u) - \sigma E_{11}(u)$ does not vary in a
      small enough neighborhood of $\sigma= \lambda_i(\bar u_0)$, i.e. if the kernel of  $A_{11}(u) - \lambda_i(\bar u_0) E_{11}(u)$
         is trivial. This is certainly a restrictive condition. In \cite{BiaSpi:per} the authors will address the problem of finding less restrictive conditions. 
  \end{rem}

\subsubsection{A lemma on the dimension of the stable manifold}
\label{subsub_Bsing_dimension} The aim of this section is to
introduce Lemma \ref{lem_Bsing_crucial}, which determines the
dimension of the stable manifold of the system
\begin{equation}
\label{eq_Bsing_riemann_boundary_leyers}
       B(U) U''= A(U, \, U') U',
\end{equation}
satisfied by the steady solutions of 
$$
  E(u) u_t + A(u, \, u_x) u_x = u_{xx}.
$$
In this section we thus give an
answer to a question introduced in \cite{Rousset:char}: more
precisely,  it is shown that the condition called there Hypothesis
(H5) is actually a consequence of Kawashima condition.

Thanks to Lemma \eqref{lem_Bsing_explicit_form}, system
\eqref{eq_Bsing_riemann_boundary_leyers} is equivalent to
\begin{equation}
\label{eq_Bsing_reduced}
 \left\{
      \begin{array}{ll}
            u_x = \Big( \bar{w}(u, \, \tilde{z}, \, 0), \;
                  \tilde{w}(u, \, \tilde{z}, \, 0), \;
                  \bar{z}(u, \, \tilde{z}, \, 0), \; \tilde{z}
                  \Big)^T \\
            \tilde{z}_x = f(u, \, \tilde{z}, \, 0). \phantom{\Big(}
      \end{array}
      \right.
\end{equation}
Hence, linearizing around the equilibrium point $(\bar{u}, \, 0)$,
one finds that the Jacobian is represented by the matrix
\begin{equation*}
\left(
\begin{array}{cc}
      0 & \underline{m} \\
      0 & \underline{b}^{-1} \underline{a} (\bar{u}, \, 0, \, 0)
\end{array}
\right)
\end{equation*}
where $\underline{a}$ and $\underline{b}$ are defined by
\eqref{eq_Bsing_prel_underline} and $\underline{m}$ is a $N \times
(r-q)$ dimensional matrix whose exact expression is not important
here. Hence the stable manifold of system \eqref{eq_Bsing_reduced}
is tangent in $(\bar{u}, \, 0)$ to
\begin{equation}
\label{eq_Bsing_stable_space}
      \tilde{V}^s = \Big\{\Big( \bar{u} + \sum_{i=1 }^{d}
      v_i \underline m \vec{\theta}_i (\bar{u}), \, \sum_{i=1 }^{d}
      \vec{\theta}_i(\bar{u})\Big): \; v_i \in \mathbb{R} \Big\}
      \subseteq \mathbb{R}^{N + r-q},
\end{equation}
In the previous expression $\vec{\theta}_1, \dots, \vec{\theta}_d   \in
\mathbb{R}^{r-q}$ are the eigenvectors of
$\underline{b}^{-1} \underline{a}$ associated to 
eigenvalues with strictly negative real part. What we want to determine in 
Lemma \eqref{lem_Bsing_explicit_form} is the number $d$, i.e. the number of eigenvalues of 
$\underline{b}^{-1} \underline{a}$ with strictly negative real part.  
 
Before stating the lemma, we recall some notations:
$k-1$ denotes the number of eigenvalues of
$E^{-1}(u)A(u, \, 0)$, that satisfy
\begin{equation*}
\label{eq_Bsing_separation_speed}
       \lambda_1(u) < \lambda_2(u) < \dots
       < \lambda_{k-1}(u) < -c <0
\end{equation*}
for some constant $c >0$. The number $(k-1)$ does not depend on
$u$ and $u_x$ because of the strict hyperbolicity (Hypothesis
\ref{hyp_hyperbolic_I}). Moreover, let
\begin{equation}
\label{eq_Bsing_block}
      A({u}, \, u_x)=
      \left(
      \begin{array}{cc}
             A_{11}({u}) & A_{21}^T(u) \\
             A_{21}({u}) & A_{22}(u, \, u_x) \\
      \end{array}
      \right)
      \qquad \qquad
      E(u)=
      \left(
      \begin{array}{cc}
             E_{11}({u}) & E_{21}^T(u) \\
             E_{21}({u}) & E_{22}(u) \\
      \end{array}
      \right)
\end{equation}
be the block decompositions of $A$ and $E$ corresponding to the block decomposition 
$$
   B(u)=
      \left(
      \begin{array}{cc}
             0 & 0 \\
             0 & b(u) \\
      \end{array}
      \right)
$$
Also, $n_{11}$ denotes the
number of negative eigenvalues of $A_{11}$. Since $E$ is positive
definite, also $E_{11}$ is and hence by Lemma
\ref{lem_Binv_dimension} $n_{11}$ is equal to the number of
strictly negative eigenvalues of $E_{11}^{-1}(u)A_{11}(u)$. Hence,
by the third assumption in Hypothesis \ref{hyp_Bsing}, $n_{11}$
does not depend on $u$. Finally, $q$ denotes the
dimension of the kernel of $A_{11}(u)$, which as well does not
depend on $u$.

Before introducing Lemma \ref{lem_Bsing_crucial}, it is necessary
to state a preliminary result, which is exploited in the
proof of the lemma:
\begin{lem}
\label{lem_Bsing_rank}
      Let Hypothesis \ref{hyp_Bsing} hold. Then there exists $\delta >0$ such that if $u$ belongs to a neighborhood of $\bar u_0$ of size $\delta$ then the following holds. The polynomial 
      \begin{equation}
            \det \Big( A(\bar{u}, \, 0) - \mu B(\bar{u})\Big)
      \end{equation}
      has degree $(r-q)$ with respect to $\mu$.
\end{lem}
The proof follows immediately from Lemma
\ref{lem_Bsing_exchange_eigenvalues}, which guarantees that
\begin{equation*}
      \det \Big( A(\bar{u}, \, 0) - \mu B(\bar{u})\Big) =0 \,
      \Longleftrightarrow \,
      \det \Big( \underline{b}^{-1} (\bar{u}, \, 0)\underline{a}(\bar{u}, \, 0, \, 0) -
      \mu I_{r-q} \Big) =0,
\end{equation*}
where the matrices $\underline{a}, \; \underline{b} \in
\mathbb{M}^{(r-q) \times (r-q)}$ are defined by
\eqref{eq_Bsing_prel_underline} and $\underline{b}$ is invertible
by Lemma \ref{lem_Bsing_inv}. The symbol $I_{r-q}$ denotes the
identity matrix of dimension $(r-q)$ .

\begin{lem}
\label{lem_Bsing_crucial}
       Let Hypotheses \ref{hyp_Bsing} and \ref{hyp_hyperbolic_I}
       hold. Then the dimension $d$ of the stable space
       \eqref{eq_Bsing_stable_space} of system
       \eqref{eq_Bsing_riemann_boundary_leyers}
       is $d=k-1 -n_{11}- q$.
\end{lem}
\begin{figure}
\begin{center}
\caption{the behavior of the roots of
\eqref{eq_Bsing_prel_determinant_perturbation} when the
perturbation $\sigma \to 0^+$} \psfrag{a}{$N-k$ roots: $Re
\mu_i(\sigma)
>0$} \psfrag{b}{$k-1$ roots: $Re \mu_i(\sigma) <0$ }
\psfrag{d}{$n_{11}+q$ eigenvalues:$Re \mu_i <0$, $| \mu(\sigma)|
\to + \infty $} \psfrag{c}{$ r -q $ eigenvalues: $\mu_i(\sigma)
\to \mu(0)$} \psfrag{e}{$\mu_k (\sigma ) \equiv 0$}
\label{fig_eigenvalues_II}
\includegraphics[scale=0.4]{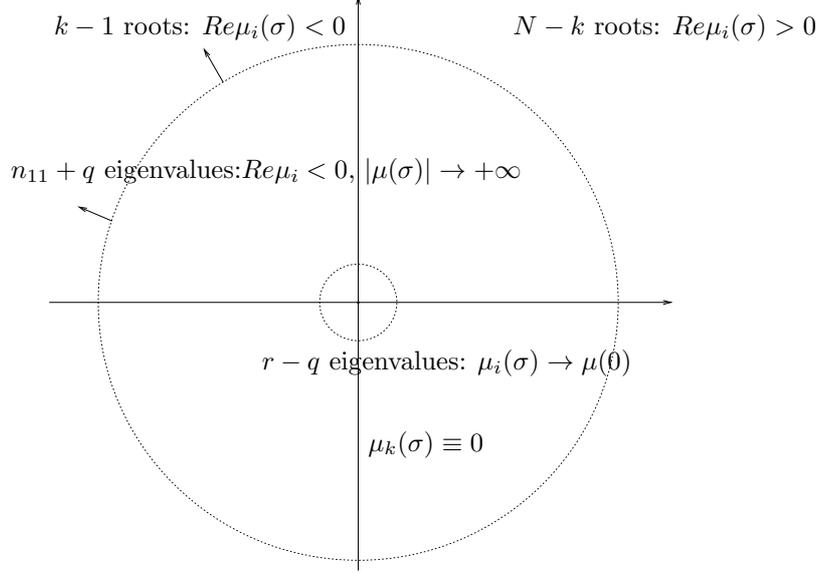}
\end{center}
\end{figure}

\begin{proof}
Since the state $\bar{u}$ is fixed, in the following for
simplicity the matrices $A(\bar{u}, \, 0)$ and $B(\bar{u})$ will
be denoted by $A$ and $B$ respectively. To prove the lemma one has
to show that the equation
\begin{equation}
\label{eq_Bsing_prel_determinant}
       \det \, \Big( A - \mu_{\underline{i}} B \Big) =0
\end{equation}
admits exactly $(k-1 - n_{11}-q)$ roots with negative real part.

The proof is organized into six steps:
\begin{enumerate}
\item First, it is introduced a perturbation technique. More
precisely, it is considered the problem
\begin{equation}
\label{eq_Bsing_prel_determinant_perturbation}
       \det \, \Big( A - \mu_i(\sigma) (B + \sigma I) \Big) =0.
       \qquad \sigma >0,
\end{equation}
where $\sigma$ is a positive parameter which is allowed to go to
$0$. The symbol $I$ denotes the $N$-dimensional identity matrix.
Just to fix the ideas and to consider the most general case, it
will be supposed that the matrix $A$ is singular:
\begin{equation*}
      \lambda_1 < \dots< \lambda_{k-1} < 0 = \lambda_k < \lambda_{k+1} < \dots <
      \lambda_N
\end{equation*}
Once the value of the parameter $\sigma$ is fixed, the solutions
$\mu_i (\sigma)$ of \eqref{eq_Bsing_prel_determinant_perturbation}
are the roots of a polynomial of degree $N$. Moreover, the matrix
$(B + \sigma I)$ is positive definite and hence Lemma
\ref{lem_Binv_dimension} guarantees that $(k-1)$ roots have
strictly negative real part, $(N -k)$ roots have strictly positive
real part and one root is zero.

If one let the parameter $\sigma$ vary,
\eqref{eq_Bsing_prel_determinant_perturbation} defines $N$
algebraic functions $\mu_1(\sigma), \dots \mu_N(\sigma)$: since
the coefficients of the polynomial
\eqref{eq_Bsing_prel_determinant_perturbation} are analytic with
respect to $\sigma$, it is known (see for example
\cite{Knopp:book}) that the functions $\mu_i$ are analytic
everywhere in the complex plane except for some so called
exceptional points that are certainly of finite number in every
compact subset. Moreover, for every fixed function $\mu_i$ only
two behaviors are possible in an exceptional point $\sigma =
\bar{\sigma}$:
\begin{itemize}
\item $\mu_i(\sigma)$ is continuous in $\bar{\sigma}$. \item
$\lim_{\sigma \to \bar{\sigma}} |\mu_i (\sigma)|= + \infty$.
\end{itemize}
In the case of the algebraic functions defined by
\eqref{eq_Bsing_prel_determinant_perturbation}, the point $\sigma
=0$ is an exceptional point since in $\sigma=0$ the degree of the
polynomial drops from $N$ to $r-q$: we will denote by
$\mu_{\underline{i}}(0)$, $\underline{i}=1, \dots (r-q)$ the
$(r-q)$ roots of \eqref{eq_Bsing_prel_determinant}.

Thanks to the known results recalled before, however, one deduces
that there are
\begin{itemize}
\item $r-q$ functions $\mu_{\underline{i}}(\sigma)$ such
      that
      \begin{equation*}
             \lim_{\sigma \to 0} \mu_{\underline{i}}(\sigma) =
             \mu_{\underline{i}}(0) \qquad i=1, \dots r.
      \end{equation*}
\item $N-r+q$ functions $\mu_{\underline{i}}(\sigma)$ such that
      \begin{equation*}
             \lim_{\sigma \to 0} |\mu_{\underline{i}}(\sigma)| =
             + \infty \qquad i=r+1, \dots N.
      \end{equation*}
\end{itemize}
The situation is summarized in Figure \ref{fig_eigenvalues_II}.

Moreover, restricting to a small enough neighborhood, it is not
restrictive to suppose that $\sigma=0$ is the only singular point.
\item As a second step we prove that to obtain the lemma it is
sufficient to show that the number of functions
$\mu_{\underline{i}}(\sigma)$ continuous in $\sigma =0$ and such
that $Re \mu_{\underline{i}}(\sigma) < 0$ when $\sigma \in
\mathbb{R}^+$ is equal to $k-1 - n_{11} - q$. Since the number of
roots of \eqref{eq_Bsing_prel_determinant_perturbation} with
negative is equal to $k-1$, this is equivalent to show that there
are exactly $n_{11}+q$ functions $\mu_{\underline{i}}(\sigma)$
such that $Re \mu_{\underline{i}}(\sigma) < 0$ and $
|\mu_{\underline{i}}(\sigma)| \uparrow + \infty$ when $\sigma \to
0^+$.

The only case that has to be excluded is the possibility that in
the limit a function $\mu_{\underline{i}}(\sigma)$ with $Re
\mu_{\underline{i}}(\sigma) < 0$ has zero real part.

First, one observes that equation
\eqref{eq_Bsing_prel_determinant} does not admit purely imaginary
roots: by contradiction, assume that there exists $\mu \in
\mathbb{R}$, $\mu \neq 0$ and $v = v_r + i v_i$, $v_r, \, v_i \in
\mathbb{R}^N$ such that
\begin{equation}
\label{eq_Bsing_prel_pflem_zeroim}
      A (v_r + i v_i) =  i \mu B (v_r  + i v_i).
\end{equation}
Just to fix the ideas, one can assume $\mu >0$: as in the proof of
Lemma \ref{lem_Binv_dimension}, it turns out that
\begin{equation}
       0 \ge - \mu \langle B v_i, \, v_i \rangle =
       \langle A v_r, \, v_i  \rangle = \langle A v_i, \, v_r
       \rangle= \mu \langle B v_r, \, v_r  \rangle \ge 0
\end{equation}
because $B$ is positive semidefinite. Hence the only possibility
is that $v_i$ and $v_r$ both belong to $\text{ker} \, B$, but this
contradicts \eqref{eq_Bsing_prel_pflem_zeroim} and hence there are
no purely imaginary roots of \eqref{eq_Bsing_prel_determinant}.

Hence one is left to exclude the possibility that a function
$\mu_{\underline{i}}(\sigma)$ continuous in $\sigma =0$ and such
that $Re \mu_{\underline{i}}(\sigma) < 0$ converges to $0$. From
the strict hyperbolicity of the matrix $E^{-1}A$ (Hypothesis
\ref{hyp_hyperbolic_I}) and from Lemma \ref{lem_Binv_dimension} it
follows that $0$ is an eigenvalue of $A$ with multiplicity one and
hence $\mu_{\underline{k}}=0$ is a root of
\eqref{eq_Bsing_prel_determinant} with multiplicity one. Moreover,
again from Lemma \ref{lem_Binv_dimension}, one deduces that the
algebraic function $\mu_{\underline{k}}(\sigma)$ is identically
equal to $0$. Since the multiplicity of $0$ as a root of
\eqref{eq_Bsing_prel_determinant} is one, it cannot happen that
other functions $\mu_{\underline{i}}(\sigma)$, $\underline{i}\neq
\underline{k}$ converge to zero when $\sigma \to 0^+$.
\item As a third step, we perform a change of variable. More
precisely, the number of solutions $\mu_{\underline{i}}(\sigma)$
of \eqref{eq_Bsing_prel_determinant_perturbation} such that $Re
\mu_{\underline{i}}(\sigma) < 0$ and $
|\mu_{\underline{i}}(\sigma)| \uparrow + \infty$ when $\sigma \to
0^+$ is equal to the number of roots $x_{\underline{i}} (\sigma)$
of
\begin{equation}
\label{eq_Bsing_prel_pflem_over}
      \det \Big( B + \sigma I - x_{\underline{i}}(\sigma) A \Big) =0
      \qquad \underline{i}=1, \dots N-r + q
\end{equation}
such that $Re x_i (\sigma) < 0$ and $x_i (\sigma) \to 0$ when
$\sigma \to 0^+$. This is clear if one defines
${x_{\underline{i}}(\sigma): = 1 / \mu_{\underline{i}}(\sigma)}$.
\item In the fourth part of the proof, we study the inverse
problem to \eqref{eq_Bsing_prel_pflem_over}, i.e. the eigenvalue
problem
\begin{equation}
\label{eq_Bsing_prel_pflem_inverse}
      \det \Big(B  - x A + \sigma_{\underline{j}}(x) I \Big) =0
      \qquad \underline{j}=1, \dots N.
\end{equation}
The meaning of the notation is that in this case one fixes $x$ and
then finds $\sigma $.

More precisely, it will be studied the behavior of the $N-r$
eigenvalues $- \sigma_{\underline{j}}(x)$ such that
$\sigma_{\underline{j}}(x) \to 0$ when $x \to 0$. From the
analysis in \cite{BiaHanNat} it follows that
\begin{equation*}
      F(x) = - A_{11} x - A_{21}^T b^{-1} A_{21} x^2 + o(x^2)
      \qquad x \to 0
\end{equation*}
where $F(x)$ denotes the projection of $B - x A$ on the
generalized subspaces converging to the kernel of $B$. The blocks
$A_{11}$ and $A_{21}$ are defined by \eqref{eq_Bsing_block}. Let
$F_0(x)$ be the projection of $F(x)$ on the kernel of $A_{11}$ and
$F_{\bot}(x)$ be the projection on the subspace orthogonal to
$\text{ker} A_{11}$.

Thus again the analysis in \cite{BiaHanNat} guarantees that
\begin{equation*}
      F_0 (x)= - (A^I_{21})^T b^{-1} A^I_{21} x^2 + o(x^2)
      \qquad x \to 0
\end{equation*}
and
\begin{equation*}
      F_{\bot}(x) = - \tilde{A}_{11} x + o(x)
      \qquad x \to 0.
\end{equation*}
In the previous expression it has been used the notations
introduced in \eqref{eq_Bsing_prel_A_blockII}:
\begin{equation*}
       \left(
             \begin{array}{cc}
                   A_{11}(u, \, \sigma)  & A_{21}(u\, \sigma)^T \\
                   A_{21}(u\, \sigma) & A_{22}(u, \, u_x, \, \sigma) \\
             \end{array}
      \right)=
      \left(
             \begin{array}{ccc}
      0 & 0                                &  \big( A^{I}_{21} \big) ^T  (u, \, \sigma)\\
      0 & \tilde{A}_{11}(u, \, \sigma)     & \big( A^{II}_{21}\big)^T (u, \, \sigma)\\
      A^{I}_{21}(u) &  A^{II}_{21}(u, \, \sigma) & A_{22}(u, \, \sigma, \, u_x) \\
             \end{array}
      \right)
\end{equation*}
From Kawashima condition it follows that the matrix $(A^I_{21})^T
b^{-1} A^I_{21}$ is positive definite: one can proceed in the same
way as in the proof of Lemma \eqref{lem_Bsing_kawashima_implies}.

Hence projecting on suitable subspaces (we refer again to
\cite{BiaHanNat} for the precise computations) it follows that the
following expansions hold:
\begin{equation}
\label{eq_Bsing_prel_pflem_exp_nonzero}
      \sigma_{\underline{j}} (x) = \mu_{\underline{j}} x^2 + o(x)
      \qquad x \to 0, \qquad j= 1, \dots q
\end{equation}
and
\begin{equation}
\label{eq_Bsing_prel_pflem_exp_zero}
      \sigma_{\underline{j}} (x) = \lambda_{\underline{j}} x + o(x)
      \qquad x \to 0 \qquad j= q+1, \dots N-r.
\end{equation}
In the previous expression, $\mu_{\underline{j}} >0$ is an
eigenvalue of $(A^I_{21})^T b^{-1} A^I_{21}$ and
$\lambda_{\underline{j}} \neq 0$ is an eigenvalue of
$\tilde{A}_{11}$.
\item As a fifth step, we analyze equation
\eqref{eq_Bsing_prel_pflem_exp_zero}: for every fixed
$\underline{j}=q+1, \dots N-r$, the expansion guarantees that the
function $\sigma_{\underline{j}} (x)$ is invertible in a
neighborhood of $x=0$. The inverse function
$x_{\underline{j}}(\sigma)$ satisfies
\eqref{eq_Bsing_prel_pflem_over} and because of
\eqref{eq_Bsing_prel_pflem_exp_zero} the condition that $Re
x_{\underline{j}}(\sigma) < 0$ when $\sigma \in \mathbb{R}^+$ is
satisfied if and only if $\lambda_{\underline{j}}$ is an
eigenvalue of $\tilde{A}_{11}$ with negative real part. This
implies that there are exactly $n_{11}$ functions
$x_{\underline{j}}(\sigma)$ satisfying
\eqref{eq_Bsing_prel_pflem_exp_zero} and such that $Re
x_{\underline{j}}(\sigma) < 0$ when $\sigma \in \mathbb{R}^+$.

Moreover, there are $N-r-q$ functions $x_{\underline{j}}(\sigma)$
satisfying \eqref{eq_Bsing_prel_pflem_over} and
\eqref{eq_Bsing_prel_pflem_exp_zero}.
\item  Since there are exactly $2q$
functions $x_{\underline{j}}(\sigma)$ satisfying
\eqref{eq_Bsing_prel_pflem_over} and
\eqref{eq_Bsing_prel_pflem_exp_nonzero}, one deduces that there
are $q$ functions with strictly positive real part and $q$
functions with strictly negative real part. The possibility of a
root with zero real part is excluded on the basis of the same
considerations as in step 2.

Hence from steps 5 and 6 one obtains that there are exactly
$n_{11}+q$ solutions $x_{\underline{i}}(\sigma)$ of
\eqref{eq_Bsing_prel_pflem_over} such that $Re
x_{\underline{i}}(\sigma) < 0$ and $x_{\underline{i}}(\sigma) \to
0$ when $\sigma \to 0^+$. Thanks to the considerations in steps 1,
2 and 3 this concludes the proof of the lemma.
\end{enumerate}
\end{proof}
To introduce Lemma \ref{lem_Bsing_riemann_intersection}, it is useful to
introduce some further notations: let
\begin{equation}
\label{eq_Bsing_riemann_stable}
       V^s(u) : = \text{span}  \langle \Theta_{n_{11}-q+1} , \dots \Theta_{k-1} \rangle
\end{equation}
be the space generated by the vectors that satisfy
\begin{equation}
\label{eq_Bsing_riemann_Theta}
      \Big( A(\bar{u}, \, 0) - \mu_i (\bar{u}, \, 0) B(\bar{u})
      \Big) \, \Theta_i (\bar{u}, \, 0) =0 \qquad
      \mu_i (\bar{u}, \, 0) < 0.
\end{equation}
Let
\begin{equation}
\label{eq_Bsing_riemann_unstable}
       V^u (u): = \text{span}  \langle \Xi_{k+1} , \dots \Xi_{N} \rangle
\end{equation}
the subspace generated by the $\Xi_i$ satisfying
\begin{equation}
\label{eq_Bsing_riemann_Xi}
      \Big( A(u, \, 0) - \lambda_i E(u) \Big) \Xi_i =0 \qquad
      \lambda_i >0.
\end{equation}
When $A(u, \, 0)$ is invertible the dimension of $V^u$ is equal to
$(N-n)$, where $n$ is the number of negative eigenvalues of $A$
defined as in \eqref{eq_Binv_nonchar_separation}. On the other
side, when $A(u, \, 0)$ is singular the dimension of $V^u$ is
equal to $(N-k)$, where $(k-1)$ is the number of strictly negative
eigenvalues of $A$ as in \eqref{eq_Binv_char_separation}. In this
case, the following subspace is non trivial:
\begin{equation}
\label{eq_Bsing_riemann_center}
       V^c (u): = \text{span}  \langle \, \Theta_{k} \rangle,
       \; \; \mathrm{where} \;
       A(u, \, 0) \, \Theta_k =0.
\end{equation}
The proof of the following result is analogous to that of Lemma
7.1 in \cite{BenSerreZum:Evans}:
\begin{lem}
\label{lem_Bsing_riemann_intersection}
      Let Hypotheses \ref{hyp_Bsing} and \ref{hyp_hyperbolic_I}
      hold. Then
      \begin{equation*}
            V^s \cap V^u = \big\{ 0 \big\} \qquad
            V^s \cap V^c = \big\{ 0 \big\} \qquad
            V^c \cap V^u = \big\{ 0 \big\},
      \end{equation*}
      where $V^s$, $V^c$ and $V^u$ are defined by
      \eqref{eq_Bsing_riemann_stable}, \eqref{eq_Bsing_riemann_center}
      and \eqref{eq_Bsing_riemann_unstable} respectively.
\end{lem}
\subsubsection{The hyperbolic limit in the non characteristic case}
\label{subsub_Bsing_nonchar} 
In this section we will provide a characterization of the limit of the parabolic approximation 
\eqref{eq_Bsing_parabolic_approx}  when the boundary  non characteristic, i.e. when none of the eigenvalues of $E^{-1} (u) A( u, \, u_x)$ can attain the value $0$
(Hypothesis \ref{hyp_nonchar}). As in Section \ref{subsubsec_Binv_nonchar}, $n$ will denote the number of eigenvalues of $E^{-1} (u) A( u, \, u_x)$  with strictly negative
negative real part and $N-p$ the number of eigenvalues with strictly positive real part. Also, we recall that $n_{11}$ is the number of strictly negative eigenvalues of $A_{11}(u)$, while 
the eigenvalue $0$ has multiplicity $q$. 

To give a characterization of the limit of  \eqref{eq_Bsing_parabolic_approx} we will follow the same steps Section \ref{subsubsec_Binv_nonchar}. Thus, in the exposition we will focus only on the points at which the singularity of the viscosity matrix $B$ plays an important role. 

The characterization of the limit works as follows. We will construct a map 
$ \phi( \bar{u}_0, s_{n_{11} + q}  \dots s_N)$ which describes, as $(s_{n_{11} + q}  \dots s_N)$ vary, states that can be connected to $\bar{u}_0$.
We will compose $\phi$ with the function $\beta$, which is used to assign the boundary condition in  \eqref{eq_Bsing_parabolic_approx}. We will the show that the composite map is locally invertible. This means that,  given $\bar{u}_0$ and $\bar{g}$ such that $|\beta (\bar{u}_0 -  \bar{g}|$  is sufficiently small, the  values of $(s_{n_{11} + q}  \dots s_N)$ are uniquely determined by the equation
$$
    \bar{g} = \phi ( \bar{u}_0, \, s_1 \dots s_N).
$$
Once $(s_{n_{11} + q} \dots s_N)$ are known the limit of \eqref{eq_Binv_the_system}  is completely characterized. The construction of the map $\phi$ is divided in some steps:  
\begin{enumerate}
 \item {\sl Waves with positive speed }

Consider the Cauchy datum $\bar{u}_0$, fix $(N-n)$ parameters $(s_{n+1} \dots s_N)$
and consider the value
$$
    \bar{u} = T^{n}_{s_{n}} \circ \dots T^N_{s_N} \bar{u}_0.
$$
The curves $T^{n}_{s_{n}} \dots T^N_{s_N}$ are, as in Section \ref{subsub_Bsing_travelling_waves}, the {\sl curves of admissible states} introduced in \cite{Bia:riemann}. 
The state $\bar{u}_0$ is then connected to $\bar{u}$ by a sequence of rarefaction and travelling waves with positive speed.

\item {\sl Boundary layers }

We have now to characterize the set of values $u$ such that the following problem admits a solution:
\begin{equation}
\label{e:Bsing:nonchar:aubu}
\left\{
\begin{array}{lllll}
      A(U, U_x) U_x = B(U) U_{xx} \\
      U(0) = u  \\
      \lim_{x \to + \infty} U(x) = \bar{u}.
\end{array}
\right.
\end{equation}
Because of Lemma \ref{lem_Bsing_explicit_form}, one has to study 
\begin{equation}
\label{eq_Bsing_stableman} \left\{
\begin{array}{ll}
       U_x = \Big( \bar{w} (U, \tilde{z}, \, 0), \, \tilde{w} (U, \tilde{z}, \, 0), \, \bar{z} (U, \tilde{z}, \, 0) , \, \tilde{z} \Big)  \\
       \tilde{z}_x = f (U, \tilde{z}, \, 0)  \\
\end{array}
\right.
\end{equation}
Consider the equilibrium point $(\bar{u}, \, 0)$, linearize at that point and denote by $V^s$ the stable space, i.e. the eigenspace associated to the eigenvalues with strictly negative real part.  Thanks to Lemma
\ref{lem_Bsing_crucial}, the dimension of $V^s$ is equal  to $n - n_{11} - q$.
$$
    V^s = \Big\{    \big( \bar{u}  + \sum_{i = 1}^n  \frac{x_i}{  \mu_i  (\bar{u}) } \vec{\xi}_i   (\bar{u}), \,   \sum_{i = 1}^n  x_i \vec{\xi}_i  (\bar{u})  \big), \; x_1 \dots x_n \in \mathbb{R}  \Big\},
$$
 where $\mu_1  (\bar{u}) \dots \mu_n   (\bar{u})$ are the eigenvalues of $\underline{b}^{-1} (\bar{u}, \, \vec{0}, \, 0)   \underline{a}(  (\bar{u}), \, \vec{0}, \, 0) $ with negative real part  and $\vec{\xi}_1  (\bar{u})  \dots \vec{\xi}_n$ are the  
 corresponding eigenvectors. 
 
 Denote by $\mathcal{M}^s$ the stable manifold, which is parameterized by $V^s$. Also, denote by $\phi_s$ a parameterization of $\mathcal{M}^s$:
 $$
     \phi_s :   V^s \to \mathbb{R}^N. 
 $$
 Let $\pi_u$ be the projection 
\begin{equation*}
\begin{split}
     \pi_u: \,
&    \erreN \times \mathbb{R}^{r -q}
      \to \erreN \\
&    (u, \, \tilde{z}) \mapsto u 
\end{split}
\end{equation*}
One can then proceed as in the proof of Section \ref{subsubsec_Binv_nonchar} and conclude that system  \eqref{e:Bsing:nonchar:aubu} admits a solution if $u \in \pi_u \big( \phi_s  (s_{n_{11} + q +1 } \dots s_n)  \big)$ for some $s_{n_{11} + q + 1} \dots s_n$. 
Also, thanks to \eqref{e:Bsing:dzu} the columns of the jacobian of $\pi_u \circ \phi_s$ computed at  $s_{n_{11} +q +1} =0 \dots s_n = 0$ are ${\Xi}_{n_{11} +q +1} \dots {\Xi}_n$. 

Note that the map $\pi_u \circ \phi_s $ actually depends also on the point $\bar{u}$ and it does in a Lipschitz continuos way:
$$
    |  \pi_u \circ \phi_s  ( \bar{u}_1, \, s_{n_{11} +q +1}  \dots s_n ) - 
    \pi_u \circ \phi_s ( \bar{u}_2, \, s_{n_{11} +q +1}  \dots s_n ) | \leq L | \bar{u}_1 - \bar{u}_2  |.
$$

\item {\sl Conclusion}

Define the map $\phi$ as follows:
\begin{equation}
\label{e:Bsing:nonchra:solver}
    \phi ( \bar{u}_0, \, s_{n_{11} +q +1}  \dots s_ N) = \pi_u \circ \phi_s \Big(   T^{N-k}_{s_{N-k}} \circ \dots T^N_{s_N} \bar{u}_0, \, s_{n_{11} +q +1}  \dots s_n  \Big)
\end{equation}
From the previous steps it follows that $\phi$ is Lipschitz continuos and that it is differentiable at $s_{n_{11}+q +1} =0 \dots s_N =0$. Also, the columns of the jacobian   
are $\Xi_{n_{11}+q +1}(\bar{u}_0) \dots \Xi_n  (\bar{u}_0)  , \, \Theta_{n + 1} (\bar{u}_0)  \dots \Theta_N $, where 
$$
  \Big( A (\bar{u}_0, \, 0) - \lambda_i (\bar{u}_0 E(\bar{u}_0) \Big) \Theta_i (\bar{u}_0) =0
$$
for $\lambda_i(\bar{u}_0) > 0$ and 
$$
  \Big( A (\bar{u}_0, \, 0) - \mu_i (\bar{u}_0 E(\bar{u}_0) \Big) \Xi_i (\bar{u}_0) = B(\bar{u}_0) \Xi_i (\bar{u}_0)
$$
for $\mu_i(\bar{u}_0)$ with strictly negative real part.

In the case of an invertible viscosity matrix (Section \ref{subsubsec_Binv_nonchar}) the definition of the map $\phi$ is the final step in the construction. Here, instead, one has to take into account the function $\B$, which is used to assign the boundary condition and is defined in Section \ref{subsub_hyp_datum}. Consider 
$$
   \B \circ \phi ( \bar{u}_0, \, s_{n_{11} +q +1}  \dots s_N ).
$$
Thanks to the regularity of $\B$ and to the previous remarks, $\B \circ \phi$ is Lipschitz continous and differentiable at $s_{n_{11} +q +1}= 0  \dots s_N=0$. Denote by 
$$
   \mathcal{V}(\bar{u}_0) = \mathrm{span} \langle \Xi_1(\bar{u}_0) \dots \Xi_n  (\bar{u}_0)  , \, \Theta_{n + 1} (\bar{u}_0)  \dots \Theta_N (\bar u_0 \rangle.  
$$
Lemma \ref{lem_trasversal}, which is introduced in Section \ref{sub_Bsing_boundary_datum}, ensures that 
for every $\vec{V} \in \mathcal{V}(\bar{u}_0)$
$$
  D \B (\bar{u}_0) \vec{V} = 0 \implies \vec{V}= \vec{0} 
$$
Thus, the jacobian of $\B \circ \phi$ at $s_{n_{11}+q +1} =0 \dots s_N =0$ is an invertible matrix. 
Thanks to 
the extension of the implicit function theorem discussed in \cite{Cl} (page 253) one can conclude that the map $\B \circ \phi ( \bar{u}_0, \cdot) $ is invertible in a neighbourhood of $(s_{n_{11}+q +1} \dots s_N) = ( 0 \dots 0 )$. 
In particular, if one takes $\bar{u}_b$ as in \eqref{eq_Binv_the_system} and assumes that $| \B (\bar{u}_0) - \bar{g} |$ is sufficiently small, then the values of $s_{n_{11}+q +1} \dots s_N$ are uniquely determined by the equation 
\begin{equation}
\label{eq_Bsing_noncharsolver}
      \bar{g} =   \B \circ \phi ( \bar{u}_0, \, s_{n_{11}+q +1} \dots s_n )
\end{equation}
Once the values of $s_{n_{11}+q +1} \dots s_N$ are known, then the limit $u (t, \, x )$ can be reconstructed. In particular, the trace of $u$ on the axis $x =0$ is given by
\begin{equation}
\label{e:Bsing:nonchar:trace}
    \bar{u} : =  T^{n+1}_{s_{n+1}} \circ \dots T^N_{s_n} \bar{u}_0.
\end{equation}
Also, the self similar function $u$ can be obtained gluing together pieces like \eqref{e:Binv:s+:limit} . 
\end{enumerate}

Here is a summary of the results obtained in this section:
\begin{teo}
\label{pro_Bsing_nonchar}
           Let Hypotheses \ref{hyp_Bsing}, \ref{hyp_hyperbolic_I}, \ref{hyp_convergence}, \ref{hyp_small_bv}, \ref{hyp_stability} , \ref{hyp_finite_pro} and \ref{hyp_nonchar} hold. Then 
             there  exists $\delta >0$ small enough such that the following holds. If $|  \B(\bar{u}_0) - \bar{g}| < \delta$, then the limit of the parabolic approximation 
             \eqref{eq_Binv_the_system} 
            satisfies 
           $$
                \bar{g} = \B \phi (\bar{u}_0, \,  s_{n_{11}+q +1} \dots s_N)
           $$
         for a suitable vector $( s_{n_{11}+q +1} \dots s_N)$. The map $\phi$ is defined by \eqref{eq_Bsing_parabolic_approx}. Given $\bar{u}_0$ and $\bar{g}$, one can invert $\B \circ \phi$ and determine uniquely 
         $     ( s_{n_{11}+q +1} \dots s_N)$. Once  $     ( s_{n_{11}+q +1} \dots s_N)$ are known the value $u(t, \, x)$ assumed by the limit function is determined a.e. $(t, \, x)$. In particular, the trace $\bar{u}$ of the hyperbolic limit in the axis $x =0$ is given by \eqref{e:Bsing:nonchar:trace}.     
\end{teo}

\subsubsection{The hyperbolic limit in the boundary characteristic case}
\label{subsub_Bsing_char} This section deals with the
limit of the parabolic approximation \eqref{eq_Bsing_parabolic_approx} in the boundary 
characteristic case, i.e. when one of the eigenvalues of $E^{-1} A$ can attain the value $0$. More precisely, we will assume Hypothesis \ref{hyp_char}, which is introduced at the beginning of Section \ref{subsubsec_Binv_char}. 
As in Section \ref{subsubsec_Binv_char}, $k-1$ is the number of eigenvalues of $E^{-1} (\bar{u}_0) A(\bar{u}_0, \, 0)$
that are less or equal then $-c$, where $c$ is the separation speed introduced in Hypothesis \ref{hyp_hyperbolic_I}. Also, $n_{11}$ is the number of strictly negative eigenvalues of $A_{11} (\bar{u}_0)$, while $q$ is the dimension of the kernel of $A_{11} (\bar{u}_0)$.

To give a characterization of the limit of  \eqref{eq_Bsing_parabolic_approx} one can follow the same steps Section \ref{subsubsec_Binv_char}. Also, in Section \ref{subsub_Bsing_nonchar} we explain how to tackle the difficulties due to the fact that the viscosity matrix is not invertible. Thus, in the following we give only a quick overview of the key points in the characterization.
 
The characterization works as follows. We construct a map 
$ \phi( \bar{u}_0, s_{n_{11} + q}  \dots s_N)$ such that the following holds. As $(s_{n_{11} + q}  \dots s_N)$ vary, $ \phi( \bar{u}_0, s_{n_{11} + q}  \dots s_N)$ 
 describes states that can be connected to $\bar{u}_0$.
We compose $\phi$ with the function $\beta$, which is used to assign the boundary condition in  \eqref{eq_Bsing_parabolic_approx}. We then show that the composite map is locally invertible. This means that,  given $\bar{u}_0$ and $\bar{g}$ such that $|\beta (\bar{u}_0 -  \bar{g}|$  is sufficiently small, then the  values of $(s_{n_{11} + q}  \dots s_N)$ are uniquely determined by the equation
$$
    \bar{g} = \B \circ \phi ( \bar{u}_0, \, s_1 \dots s_N).
$$
Once $(s_{n_{11} + q} \dots s_N)$ are known the limit of \eqref{eq_Binv_the_system} is completely characterized and it is obtained gluing together pieces like \ref{e:Binv:s+:limit}. The construction of the map $\phi$ is divided in some steps:  
\begin{enumerate}
 \item {\sl Waves with positive speed }

Given the Cauchy datum $\bar{u}_0$, we fix $(N-k)$ parameters $(s_{k+1} \dots s_N$
and consider the value
$$
    \bar{u}_k = T^{k+1}_{s_{k+1}} \circ \dots T^N_{s_N} \bar{u}_0.
$$
The curves $T^{k+1}_{s_{k+1}} \dots T^N_{s_N}$ are, as in Section \ref{subsub_Bsing_travelling_waves}, the {\sl curves of admissible states} introduced in \cite{Bia:riemann}. 
The state $\bar{u}_0$ is then connected to $\bar{u}$ by a sequence of rarefaction and travelling waves with positive speed.

\item {\sl Analysis of the center stable manifold}

Consider the equation satisfied by travelling waves:
$$
  B(U) U'' = \big( A -\sigma E \big) U'.
$$
Thanks to Lemma \ref{lem_Bsing_explicit_form}, this is equivalent to system 
\begin{equation}
\label{eq_Bsing_centerstableman} \left\{
\begin{array}{lll}
       U_x = \Big( \bar{w} (U, \tilde{z}, \, 0), \, \tilde{w} (U, \tilde{z}, \, 0), \, \bar{z} (U, \tilde{z}, \, 0) , \, \tilde{z} \Big)  \\
       \tilde{z}_x = f (U, \tilde{z}, \, 0)  \\
       \sigma_x =0 \\
\end{array}
\right.
\end{equation}
The point $( \bar{u}_k, \, \vec{0}, \, \lambda_k )$ is an equilibrium. One can than define an invariant {\sl center stable} $\mathcal{M}^{cs}$ manifold with the same properties listed in Section \ref{subsubsec_Binv_char}. Thanks to Lemma \ref{lem_Bsing_crucial}, the dimension of every center stable manifold is $k-n_{11} -q$. Also, one can proceed again as in Section \ref{subsubsec_Binv_char} and find the equation satisfied by the solutions of \ref{subsub_Bsing_travelling_waves} laying on $\mathcal{M}^{cs}$. Moreover, every solution laying on $\mathcal{M}^{cs}$ can be decomposed in a purely center component, a purely stable component and a component of perturbation, in the same way described in Section \ref{subsubsec_Binv_char}. 
Eventually, one is able to define a map $F(\bar{u}_k, \, s_{n_{11} +q}, \dots s_K)$ which is Lischitz continous with respect to both $\bar{u}_k$ and $s_{n_{11} +q}, \dots s_K$. Also, it is 
differentiable at $s_{n_{11} +q}=0, \dots s_K=0$ and the columns of the jacobian matrix are the vectors $\Xi_{n_{11} +  q} (\bar{u}_k) \dots \Xi_{k} (\bar{u}_k)$. For every $i = n_{11} + q \dots K$ it holds 
$$
    \Big[ A (\bar{u}_k, \, 0) - \mu_i E(\bar{u}_k \Big] \Xi_i = B (\bar{u}_k ) \Xi_i
$$
with the real part of $\mu_i$ less or equal to zero. 

\item {\sl Conclusion} \\

Define the map $\phi$ as 
\begin{equation}
\label{e:Bsing:chra:solver}
    \phi ( \bar{u}_0, \, s_{n_{11} +q +1}  \dots s_N ) = \pi_u \circ \phi_s \Big(   T^{N-k}_{s_{N-k}} \circ \dots T^N_{s_N} \bar{u}_0, \, s_{n_{11} +q +1}  \dots s_n  \Big)
\end{equation}
From the previous steps it follows that $\phi$ is Lipschitz continuos and that it is differentiable at $s_{n_{11}+q +1} =0 \dots s_N =0$. Also, the columns of the jacobian   
are $\Xi_1(\bar{u}_0) \dots \Xi_k  (\bar{u}_0)  , \, \Theta_{k + 1} (\bar{u}_0)  \dots \Theta_N $, where 
$$
  \Big( A (\bar{u}_0, \, 0) - \lambda_i (\bar{u}_0 E(\bar{u}_0) \Big) \Theta_i (\bar{u}_0) =0
$$
with $\lambda_i >0$ for every $i$. The vectors $\Xi$ are as before.

To take into account the function $\beta$, which is used to assign the boundary condition and is defined in Section \ref{subsub_hyp_datum}., one considers  
$$
   \B \circ \phi ( \bar{u}_0, \, s_{n_{11} +q +1}  \dots s_N ).
$$
Thanks to the regularity of $\B$ and to the previous remarks, $\B \circ \phi$ is Lipschitz continous and differentiable at $s_{n_{11} +q +1}= 0  \dots s_N =0$. Denote by 
$$
   \mathcal{V}(\bar{u}_0) = \mathrm{span} \langle \Xi_1(\bar{u}_0) \dots \Xi_k (\bar{u}_0)  , \, \Theta_{k + 1} (\bar{u}_0)  \dots \Theta_N  \rangle.  
$$
Lemma \ref{lem_trasversal}, which is introduced in Section \ref{sub_Bsing_boundary_datum}, ensures that 
for every $\vec{V} \in \mathcal{V}(\bar{u}_0)$
$$
  D \B (\bar{u}_0) \vec{V} = 0 \implies \vec{V}= \vec{0} 
$$
Thus, the jacobian of $\B \circ \phi$ at $s_{n_{11}+q +1} =0 \dots s_N =0$ is an invertible matrix. 
Thanks to 
the extension of the implicit function theorem discussed in \cite{Cl} (page 253) one can conclude that the map $\B \circ \phi ( \bar{u}_0, \cdot) $ is invertible in a neighbourhood of $(s_{n_{11}+q +1} \dots s_N) = ( 0 \dots 0 )$. 
In particular, if one takes $\bar{u}_b$ as in \eqref{eq_Binv_the_system} and assumes that $| \B (\bar{u}_0) - \bar{g} |$ is sufficiently small, then the values of $s_{n_{11}+q +1} \dots s_N$ are uniquely determined by the equation 
\begin{equation}
\label{eq_Bsing_charsolver}
      \bar{g} =   \B \circ \phi ( \bar{u}_0, \, s_{n_{11}+q +1} \dots s_n )
\end{equation}
Once the values of $s_{n_{11}+q +1} \dots s_N$ are known, then the limit $u (t, \, x )$ can be reconstructed gluing together pieces like \ref{e:Binv:s+:limit}.  In particular, the value of the trace $\bar{u}$ on the axis $x =0$ can be determined in the same way described in Section \ref{subsubsec_Binv_char}
\end{enumerate}

Here is a summary of the results obtained in this section:
\begin{teo}
\label{pro_Bsing_char}
           Let Hypotheses \ref{hyp_Bsing}, \ref{hyp_hyperbolic_I}, \ref{hyp_convergence}, \ref{hyp_small_bv}, \ref{hyp_stability} , \ref{hyp_finite_pro} and \ref{hyp_char} hold. Then 
             there  exists $\delta >0$ small enough such that the following holds. If $|  \B(\bar{u}_0) - \bar{g}| < \delta$, then the limit of the parabolic approximation 
             \eqref{eq_Binv_the_system} 
            satisfies 
           $$
                \bar{g} = \B \circ \phi (\bar{u}_0, \,  s_{n_{11}+q +1} \dots s_N)
           $$
         for a suitable vector $( s_{n_{11}+q +1} \dots s_N)$. The map $\phi$ is defined by \eqref{e:Bsing:chra:solver}. Given $\bar{u}_0$ and $\bar{g}$, one can invert $\B \circ \phi$ and determine uniquely 
         $     ( s_{n_{11}+q +1} \dots s_N)$. Once $     ( s_{n_{11}+q +1} \dots s_N)$ are known, then the value $u(t, \, x)$ assumed by the limit function is determined a.e. $(t, \, x)$.  
 \end{teo}

\subsection{A transversality lemma}
\label{sub_Bsing_boundary_datum} 
In the first part of this section we state and prove Lemma 
\ref{lem_trasversal}. It is a technical result and it guarantees that the map $\B \circ \phi$ that appears in both 
Theorems \ref{pro_Bsing_nonchar} and 
\ref{pro_Bsing_char} is indeed locally invertible. In the second part of the section we introduce a new definition for the map $\B$ which is used to assign the boundary condition in 
\begin{equation}
\label{eq_Bsing_bd_system}
\left\{
\begin{array}{lll}
      E(u) u_t + A(u, \, u_x) u_x = B(u) u_{xx} \\
      \B(u(t, \, 0)) = \bar g \qquad
      u(0, \, x)= \bar u_0. \\
\end{array}
\right.
\end{equation}
This definition is an extension of Definiton \ref{def_bc} and it guarantees both the local invertibily of the map $\B \circ \phi$ and the well posedness of the initial boundary value problem \eqref{eq_Bsing_bd_system}. 

Before stating Lemma \ref{lem_trasversal}, we recall some notations. We denote by $A_{11}$, $E_{11}$ the blocks of $A$, $E$ defined
by \eqref{eq_Bsing_block}. Let $\vec{\zeta}_i(u) \in
\mathbb{R}^{N-r}$ be an eigenvector of $E^{-1}_{11}(u) A_{11}(u)$
associated to an eigenvalue with non positive real part. Also, let
$Z_i (u)\in \mathbb{R}^N$ be defined by
 \begin{equation}
 \label{eq_bsing_bd_zeta}
       Z_i:=
       \left(
       \begin{array}{cc}
              \vec{\zeta}_i \\
              0
       \end{array}
       \right)
\end{equation}
and
\begin{equation*}
             \mathcal{Z}(u):= \text{span}  \langle Z_1(u),
             \dots, Z_{n_{11}+q}(u)
             \rangle.
\end{equation*}

We define the complementary subspace
\begin{equation*}
      \mathcal{W}(u):= \text{span}  \Big\langle
      \left(
      \begin{array}{cc}
             0 \\
             \vec{e}_1 \\
      \end{array}
      \right),
      \dots,
      \left(
      \begin{array}{cc}
             0 \\
             \vec{e}_r \\
      \end{array}
      \right),
      \left(
      \begin{array}{cc}
             \vec{w}_1 \\
             0 \\
      \end{array}
      \right),
      \dots,
      \left(
      \begin{array}{cc}
             \vec{w}_{N-r-n_{11}-q} \\
             0 \\
      \end{array}
      \right)
      \Big\rangle,
\end{equation*}
where $\vec{e}_i \in \mathbb{R}^r$ are the vectors of a basis in
$\mathbb{R}^r$ and $\vec{w}_j \in \mathbb{R}^{N-r}$ are the
eigenvectors of $E^{-1}_{11} A_{11}$ associated to eigenvalue with strictly
positive real part. By Definition \ref{def_bc}, the
function $\B(u)$ is the component $u_1$ of the vector $u$ in the
decomposition
\begin{equation}
\label{E:boundecom1}
      u = u_1 + u_2, \quad u_1 \in \mathcal{W}({u}), u_2 \in \mathcal{V}(u).
\end{equation}

Also, we define the space $\mathcal V (u) $ as 
\begin{equation}
\label{eq_bsing_bd_defV}
            \mathcal{V}(u): = \text{span}  \langle \Theta_{n_{11}+q+1}(u), \dots,
            \Theta_{k-1}(u), \Theta_k(u), \,
            \Xi_{k+1}(u), \dots \Xi_{N}(u)
            \rangle,
\end{equation}
when the matrix $A(u, \, 0)$ can be singular (boundary characteristic case). The vectors $\Theta_i$ satisfy
$$
  A(u, \, 0) \Theta_i(u) = \mu_i B \Theta_i \qquad Re (\mu_i ) < 0,
$$
while the vectors $\Xi_i$ satisfy
$$
   \Big( A (u, \, 0) - \lambda_i (u)E (u) \Big) \Xi_i = \vec 0 \qquad \lambda_i (u) \ge 0
$$
In the case of a non characteristic boundary (i.e when $A(u, \, 0)$ is always invertible) the subspace $\mathcal{V}(u)$ is defined as
\begin{equation*}
            \mathcal{V}(u): = \text{span}  \langle \Theta_{n_{11}+q+1}(u, \, 0), \dots,
            \Theta_{n}(u, \, 0), \,
            \Xi_{n+1}, \dots \Xi_{N}(u, \, 0)
            \rangle,
\end{equation*}
With $n$ we denote the number of
negative eigenvalues of $A$, according to
\eqref{eq_Binv_nonchar_separation}. According to the results in Sections \ref{subsub_Bsing_nonchar} and \ref{subsub_Bsing_char}, $\mathcal{V}(\bar u_0)$ is the space generated
by the columns of the jacobian of the map $\phi$ at $(s_{n_{11} + q+ 1}
\dots s_N )= \vec{0}$.

The following lemma guarantees that, even if $\mathcal{W}(u)$ and
$\mathcal{V}(u)$ do not coincide in general, nevertheless they are
transversal to the same subspace $\mathcal{Z}(u)$. This result is already known, but we repeat the proof for completeness. 
\begin{lem}
\label{lem_trasversal}
      The following holds: 
      \begin{equation*}
            \mathcal{V}(\bar u_0) \oplus \mathcal{Z}(\bar u_0) = \mathbb{R}^N
            \qquad
            \mathcal{W}(\bar u_0) \oplus \mathcal{Z}(\bar u_0) = \mathbb{R}^N
      \end{equation*}
\end{lem}
\begin{proof}
     The second part of the statement is trivial.

     In the proof of the first part for simplicity we consider only the boundary characteristic case 
     ($A(u, \, 0)$ can be singular), being the other case absolutely analogous. Also, we write $A$ instead of $A (\bar(u_0, \, )$.
     The proof is organized into four steps:
     \begin{enumerate}
     \item First, we prove that if $V$ is a non zero
     vector in $\text{span}  \langle \Xi_{k+1},
           \dots , \Xi_{N} \rangle$, then
           \begin{equation}
           \label{eq_Bsing_bd_V}
                 \langle V, \, A V \rangle >0.
           \end{equation}
           The steps are analogous to those in the proof of Lemma 7.1 in
     \cite{BenSerreZum:Evans}, but for
     completeness we repeat them.

     If $V \in \text{span}  \langle \Xi_{k+1}, \dots , \Xi_{N} \rangle$,
     then  $V$ belongs
    to the unstable manifold of the linear system
     \begin{equation*}
           u_x = A u
     \end{equation*}
     and hence there exists a solution of
     \begin{equation*}
     \left\{
     \begin{array}{ll}
           u_x = A u \\
           u(0) = V \qquad \lim_{x \to - \infty} u(x ) =0.
     \end{array}
     \right.
     \end{equation*}
     Thanks to the symmetry of $A$, such a solution satisfies
     \begin{equation*}
           \frac{d}{dx} \langle u, \, A u \rangle =
           2 \langle u_x, \, A u \rangle = 2 |Au|^2 > 0
     \end{equation*}
     and hence to conclude it is enough to observe that
     \begin{equation*}
            \lim_{x \to - \infty} \langle u(x), \, A u(x) \rangle
            =0.
     \end{equation*}
\item if
\begin{equation*}
            Z = \sum_{i=1}^{n_{11} + q} x_i Z_i,
\end{equation*}
then
 \begin{equation}
 \label{eq_Bsing_bd_>01}
           \langle Z, \, A Z \rangle \leq 0.
\end{equation}
Indeed, the matrices $A_{11}$, $E_{11}$ are symmetric and $E_{11}$
is positive definite, hence they admit an orthogonal basis of
eigenvectors $\zeta_1 \dots \zeta_{N-r}$ such that $\langle
\zeta_j, \, A_{11} \zeta_i \rangle= \eta_i \delta_{i \, j} \,
\langle \zeta_i, \, E_{11} \zeta_i \rangle$. In particular,
\begin{equation*}
           \langle Z, \, A Z \rangle = \sum_{i=1}^{n_{11} + q}
           \eta_i x_i^2 \langle \zeta_i, E_{11} \zeta_i \rangle \leq 0.
\end{equation*}
\item if
     \begin{equation*}
           \Theta  = \sum_{j= n_{11}+q+1}^{k} y_j \Theta_j,
     \end{equation*}
     then
     \begin{equation}
     \label{eq_Bsing_bd_Theta_leq}
           \langle A \Theta, \, \Theta \rangle \leq 0.
     \end{equation}
     Indeed, by Lemma \ref{lem_Bsing_explicit_form} one deduces that the following
     system admits a solution:
     \begin{equation*}
     \left\{
     \begin{array}{ll}
           B u_x = A u \\
           u(0) = \Theta \qquad \lim_{x \to + \infty} u(x ) =0.
     \end{array}
     \right.
     \end{equation*}
     Hence, by considerations analogous to those performed in the first step, one
     concludes that \eqref{eq_Bsing_bd_Theta_leq} holds true.
\item it holds
      \begin{equation}
      \label{eq_Bsing_bd_third_step}
            \text{span}  \langle Z_1, \dots , Z_{n_{11}+q},
            \Theta_{n_{11}+q+1}, \dots , \Theta_{k}
            \rangle \cap \text{span}  \langle \Xi_{k+1},
            \dots , \Xi_{N} \rangle = \{ 0 \}.
      \end{equation}
     To prove it, it is enough to show that if
     \begin{equation}
     \label{e:Bsing:lt:notin}
           Z \in \text{span}  \langle Z_1, \dots , Z_{n_{11}+q} \rangle,
           \qquad
           \Theta \in \text{span}  \langle \Theta_{n_{11}+q+1}, \dots , \Theta_{k}
           \rangle, \qquad \Big(Z + \Theta \Big) \neq 0,
     \end{equation}
     then
     \begin{equation*}
           \Big( Z + \Theta \Big) \notin \text{span}  \langle \Xi_{k+1},
            \dots , \Xi_{N} \rangle.
     \end{equation*}
     If 
     $$
         \Theta = \sum_{i = {n_{11}} + q +1}^{N} x_i \Theta_i, 
     $$
     then 
     \begin{equation*}
           \langle A Z , \, \Theta \rangle = \langle Z, \, A \Theta \rangle =
           = \sum_{i = { n_{11} }+ q +1}^{N} x_i \langle Z, \, A \Theta_i \rangle =
	    \sum_{i = {n_{11} } + q +1}^{N} \mu_i x_i \langle Z, \, B \Theta_i \rangle =0
     \end{equation*}
     and
     \begin{equation*}
           \langle Z + \Theta, \, AZ + A \Theta \rangle = \langle
           Z, \, AZ  \rangle + 2 \langle Z, \, A\Theta \rangle +
           \langle \Theta, \, A \Theta \rangle,
     \end{equation*}
     from \eqref{eq_Bsing_bd_>01}, \eqref{eq_Bsing_bd_V} and
     \eqref{eq_Bsing_bd_Theta_leq} it follows that
     \eqref{e:Bsing:lt:notin} holds true.
    \item  to conclude, it is enough to show that
    for every $i=1, \dots, (n_{11}+q)$
     \begin{equation}
     \label{eq_Bsing_bd_Theta}
           Z_i \notin \, \text{span}
           \langle \Theta_{n_{11}+q+1}(u, \, 0), \dots,
           \Theta_k(u, \, 0) \rangle.
     \end{equation}
     Suppose by contradiction
     that
     \begin{equation}
     \label{eq_bsing_bd_lin_comb}
            Z_i = \sum_{j= n_{11} + q+ 1}^{k} c_j \Theta_j
     \end{equation}
     for suitable numbers $c_{n_{11}+q+1}, \dots c_k$. Because
     of Lemma \ref{lem_Bsing_exchange_eigenvalues}, for
     every $j$, $\Theta_j$ has the structure described by
     \eqref{eq_Bsing_prel_exchange}, where $\vec{\xi}_j$ is an
     eigenvector of $\underline{b}^{-1} \underline{a}$. The
     matrices $\underline{b}$ and $\underline{a}$ have dimension
     $(r-q)$ and are defined by \eqref{eq_Bsing_prel_underline}.
     Considering the last $(r-q)$ lines of the equality
     \eqref{eq_bsing_bd_lin_comb} one obtains
     \begin{equation*}
           0 = \sum_j c_j \vec{\xi}_j,
     \end{equation*}
     which implies $c_j =0$ for every $j$, because the $\vec{\xi}_j$
     are all independent. Hence \eqref{eq_bsing_bd_lin_comb}
     cannot hold and \eqref{eq_Bsing_bd_Theta} is proved.
\end{enumerate}
\end{proof}

From Lemma \ref{lem_trasversal} we deduce that the map $\B \circ \phi$ which appears in Theorems 
\ref{pro_Bsing_nonchar} and \ref{pro_Bsing_char} is locally invertible. We proceed as follows. 

By construction the kernel of the jacobian $D \, \B$ is $\mathcal Z ( \bar u_0) $. Thus, 
$$
   D \, \B (\bar u_0 ) V = \vec{0} \; \implies \; V \in \mathcal Z (\bar u_0). 
$$
The columns of the matrix 
$$
  D \, \B ( \bar u_0)  \, D \, \phi (\bar u_0) 
$$
are $D \, \B ( \bar u_0) \Theta_{n_{11} + q+1} \dots D \, \B ( \bar u_0) \Theta_{k-1}, D \, \B ( \bar u_0) \Xi_k \dots D \, \B ( \bar u_0) \Xi_N $. To prove that the columns are all indipendent it is enough to show that 
\begin{equation}
\label{e:Bsing:ts:indip}
  \sum_{i = { n_{11} }+ q +1}^{k-1} x_i  D \, \B ( \bar u_0)  \Theta_i + 
  \sum_{i = k }^{N} x_i  D \, \B ( \bar u_0)  \Xi_i = \vec{0} \; \implies \; x_{n_{11}+q+1} = \dots x_N =0.
\end{equation}
Since $ \Theta_{n_{11} + q+1} \dots \Theta_{k-1}, \, \Xi_k \dots \Xi_N $ are all indipendent, to prove \eqref{e:Bsing:ts:indip} is is enough to show that if $V \in \mathcal V (\bar u_0)$ and 
$$
   D \, \B ( \bar u_0) V = \vec 0 
$$
then $V = \vec 0$. This a consequence of Lemma \ref{lem_trasversal}, which states that 
$$
  \mathcal V ( \bar u_0 ) \oplus \mathcal Z (\bar u_0 ) = \mathbb R^N
$$
\vspace{1cm}

We now introduce a generalization of Definition \ref{def_bc}.  
\begin{say}
\label{def_gen}
       The function 
       $$
         \B: \mathbb R^N \to \mathbb R^{N - n_{11} -q}
	$$
       is  smooth and satifies 
        \begin{equation}
       \label{eq_bd_trasv}
              \mathrm{kernel} \Big( D \B (\bar u_0)  \Big) \oplus \mathcal W (\bar u_0) = \mathbb R^N 
       \end{equation} 
       and 
       \begin{equation}
       \label{e:Bsing_trs:rie}
               \mathrm{kernel} \Big( D \B (\bar u_0)  \Big) \oplus \mathcal V (\bar u_0) = \mathbb R^N 	
       \end{equation}
\end{say}
Thanks to \eqref{eq_bd_trasv}, the initial boundary value problem 
\begin{equation*}
\left\{
\begin{array}{lll}
      E(u) u_t + A(u, \, u_x) u_x = B(u) u_{xx} \\
      \B(u(t, \, 0)) = \bar{g} \qquad
      u(0, \, x)= \bar{u}_0 \\
\end{array}
\right.
\end{equation*}
is well posed. This is a consequence of the same considerations as in Section \ref{subsub_hyp_datum}.

Also, because of \eqref{e:Bsing_trs:rie} the matrix 
$$
   D \, \B (\bar u_0) \, D \, \phi (\bar u_0) 
$$
is invertible. Thus, the function $\B \circ \phi$ which appears in Theorems \ref{pro_Bsing_nonchar} and \ref{pro_Bsing_char} is locally invertible. In other words, the analysis in Section \ref{sec_B_non_invertible} is still valid if we use Definition \ref{def_gen} instead of Definition \ref{def_bc}.

\bibliography{biblio}
\end{document}